\documentclass[fleqn,10pt]{article}
\usepackage{amsfonts}
\usepackage{amsmath}
\usepackage{amscd}
\usepackage{amsthm}
\usepackage{amssymb}
\usepackage{a4wide}
\usepackage[english]{babel}
\usepackage{graphicx}
\newtheorem{theorem}{Theorem}
\newtheorem{lemma}[theorem]{Lemma}
\newtheorem{definition}[theorem]{Definition}
\newtheorem{corollary}[theorem]{Corollary}

\newtheorem{Remark}[theorem]{Remark}
\newtheorem{Result}[theorem]{Result}
\newenvironment{remark}{\begin{Remark}\rm}{\end{Remark}}

\newcommand{\ve}{\varepsilon}

\newcommand{\diev}{\textrm{div}}

\newcommand{\ul}{\textbf}

\newcommand{\R}{\mathbb{R}}

\newcommand{\LLL}{\mbox{\boldmath$\lambda$}}
\newcommand{\www}{\mbox{\boldmath$\omega$}}

\newcommand{\PPPP}{\mbox{\boldmath$\psi$}}

\newcommand{\desda}{\Leftrightarrow}
\DeclareMathAlphabet\gothic{U}{euf}{m}{n}
\begin{document}
\title{\textbf{Left-invariant Stochastic Evolution Equations on $SE(2)$ and its Applications to Contour Enhancement and Contour Completion via Invertible Orientation Scores.}}
\author{\mbox{Remco Duits and Erik Franken} \\ \\
\centerline{{\small Eindhoven University of Technology}} \\
\centerline{{\small Department of Mathematics and Computer Science, CASA applied analysis, }} \\
\centerline{{\small Department of Biomedical Engineering, BMIA Biomedical image analysis.}}
\\
\centerline{{\small e-mail: R.Duits@tue.nl, E.M.Franken@tue.nl }}
\\
\centerline{{\small available on the web at:
http://www.win.tue.nl/casa/research/casareports/2007.html }}}
\maketitle
\begin{abstract}
We provide the explicit solutions of linear, left-invariant, (convection)-diffusion equations and the corresponding resolvent equations on the 2D-Euclidean motion group $SE(2)=\R^2 \rtimes \mathbb{T}$. These diffusion equations are forward Kolmogorov equations for well-known stochastic processes for contour enhancement and contour completion. The solutions are given by group-convolution with the corresponding Green's functions which we derive in explicit form. We have solved the Kolmogorov equations for stochastic processes on contour \emph{completion}, in earlier work \cite{DuitsR2006AMS}. Here we mainly focus on the Forward Kolmogorov equations for contour \emph{enhancement} processes which, in contrast to the Kolmogorov equations for contour completion, do not include convection. The Green's functions of these left-invariant partial differential equations coincide with the heat-kernels on $SE(2)$. Nevertheless, our exact formulae do not seem to appear in literature. Furthermore, by approximating the left-invariant basis of the generators on $SE(2)$ by left-invariant generators of a Heisenberg-type group, we derive approximations of the Green's functions.

The Green's functions are used in so-called completion distributions on $SE(2)$ which are the product of a forward resolvent evolved from a source distribution on $SE(2)$ and a backward resolvent evolution evolved from a sink distribution on $SE(2)$. Such completion distributions on $SE(2)$ represent the probability density that a random walker from a forward proces collides with a random walker from a backward process. On the one hand, the modes of Mumford's direction process (for contour completion) coincides with elastica curves minimizing $\int \kappa^{2} + \epsilon {\rm ds}$,
and they are closely related to zero-crossings of two left-invariant derivatives of the completion distribution.
On the other hand, the completion measure for the contour enhancement proposed by Citti and Sarti, \cite{Citti} concentrates on the geodesics
minimizing $\int \sqrt{\kappa^{2} + \epsilon} {\rm ds}$ if the expected life time $1/\alpha$ of a random walker in $SE(2)$ tends to zero.

This motivates a comparison between the geodesics and elastica. For reasonable parameter settings they turn out to be quite similar. However, we apply the results by Bryant and Griffiths\cite{Bryant} on Marsden-Weinstein reduction
on Euler-Lagrange equations associated to the elastica functional, to the case of the geodesic functional. This yields rather simple practical analytic solutions for the geodesics,
which in contrast to the formula for the elastica, do not involve special functions.

The theory is directly
motivated by several medical image analysis applications where
enhancement of elongated structures, such as catheters and bloodvessels, in noisy medical image data is required. Within this article we show how the left-invariant evolution processes can be used for automated contour enhancement/completion using a so-called orientation score, which is obtained from a grey-value image by means of a special type of unitary wavelet transformation. Here the (invertible) orientation score serves as both the source and sink-distribution in the completion distribution.

Furthermore, we also consider \emph{non-linear adaptive evolution equations} on orientation scores. These non-linear evolution equations are practical improvements of the celebrated standard ``coherence enhancing diffusion''-schemes on images as \emph{they can cope with crossing contours}.
Here we employ differential geometry on $SE(2)$ to include curvature in our non-linear diffusion scheme on orientation scores. Finally, we use the same differential geometry for a morphology theory on orientation scores yielding automated erosion towards geodesics/elastica.
\end{abstract}

\tableofcontents

\section{Invertible Orientation Scores \label{ch:OS}}

In many image analysis applications an object~$U_f \in \mathbb{L}_{2} (SE(2))$ defined on the 2D-Euclidean motion group $SE(2)=\R^2 \rtimes \mathbb {T}$ is constructed from a 2D-grey-value image $f \in \mathbb{L}_{2} (\R^2)$. Such an object provides an overview of all local orientations in an image. This is important for image analysis and perceptual organization,  \cite{Kali99a}, \cite{ForssenThesis}, \cite{Medi00}, \cite{Fels04a},  \cite{Duits2004PRIA}, \cite{Williams}, 
\cite{August2003} and is inspired by our own visual system, in which receptive fields exist that are tuned to various locations and orientations, \cite{Tso90}, \cite{Bosking}. In addition to the approach given in the introduction other schemes to construct
$U_f:\R^2 \rtimes \mathbb{T} \to \mathbb{C}$ from an image $f:\R^2 \to \R$ exist, but only few methods put emphasis on the stability of the inverse transformation $U_f \mapsto f$.

In this section we provide an example on how to obtain such an object~$U_{f}$ from an image~$f$. This leads to the concept of invertible orientation scores, which we developed in previous work, \cite{DuitsRThesis}, \cite{Duits2004PRIA}, \cite{Duits2005IJCV}, and which we briefly explain here.

An orientation score $U_f:\R^2 \rtimes \mathbb{T} \to \mathbb{C}$ of an image $f:\R^2 \to \R$ is obtained by means of an anisotropic convolution kernel~$\check{\psi}:\R^2 \to \mathbb{C}$ via
\[
U_{f}(g)=\int_{\R^2} \overline{\psi(R_{\theta}^{-1}(\ul{y}-\ul{x}))}
f(\ul{y})\; {\rm d}\ul{y},\ \ \  g=(\ul{x},e^{i\theta}) \in G=\R^2 \rtimes \mathbb{T}, R_{\theta} \in \textrm{SO}(2),
\]
where $\psi(-\ul{x})=\check{\psi}(\ul{x})$. Assume $\psi \in \mathbb{L}_{2}(\R^2) \cap \mathbb{L}_{1}(\R^2)$, then the transform $\mathcal{W}_{\psi}$ which maps image $f \in \mathbb{L}_{2}(\R^2)$ onto its orientation score $U_{f} \in \mathbb{L}_{2}(\R^2 \rtimes \mathbb{T})$ can be re-written as
\[
U_{f}(g)=(\mathcal{W}_{\psi}f)(g)=(\mathcal{U}_{g}\psi,f)_{\mathbb{L}_{2}(\R^2)},
\]
where $g \mapsto \mathcal{U}_{g}$ is a unitary (group-)representation of the Euclidean motion group $SE(2)=\R^2 \rtimes \mathbb{T}$ into $\mathbb{L}_{2}(\R^2)$ given by $\mathcal{U}_{g}f(\ul{y})= f(R_{\theta}^{-1}(\ul{y}-\ul{x}))$ for all $g=(\ul{x},e^{i\theta}) \in SE(2)$ and all $f \in \mathbb{L}_{2}(\R^2)$. Note that the representation $\mathcal{U}$ is reducible as it leaves the following closed subspaces invariant $\{f \in \mathbb{L}_{2}(\R^2)\; |\; \textrm{supp}\mathcal{F}[f] \subset B_{\ul{0},\varrho}\}$, $\varrho>0$, where $B_{\ul{0},\varrho}$ denotes the ball with center $\ul{0}\in \R^2$ and radius $\varrho>0$ and where $\mathcal{F}:\mathbb{L}_{2}(\R^2) \to \mathbb{L}_{2}(\R^2)$ denotes the Fourier transform
given by
\[
\mathcal{F}f(\www)= \frac{1}{2\pi} \int \limits_{\R^2} f(\ul{x}) e^{-i(\www,\ul{x})}\, {\rm d}\ul{x},
\]
for almost every $\www \in \R^2$ and all $f\in \mathbb{L}_{2}(\R^2)$.

This differs from standard continuous wavelet theory, see for example \cite{LeeTS} and \cite{Antoine99}, where the wavelet transform is constructed by means of a quasi-regular representation of the similitude group \mbox{$\R^{d} \rtimes \mathbb{T} \times \R^{+}$}, which is unitary, \emph{irreducible} and square integrable (admitting the application of the more general results in \cite{Grossmann85}). For the image analysis this means that we do allow a stable reconstruction already at a \emph {single scale} orientation score for a proper choice of $\psi$. In standard wavelet reconstruction schemes, however, it is \emph{not} possible to obtain an image~$f$ in a well-posed manner from a ``fixed scale layer'', that is from $\mathcal{W}_{\psi}f(\cdot,\cdot,\sigma)\in
\mathbb{L}_{2}(\R^2\rtimes \mathbb{T})$, for fixed scale $\sigma>0$.

Moreover, the general wavelet reconstruction results \cite{Grossmann85} do not apply to the transform $f \mapsto U_f$, since our representation $\mathcal{U}$ is reducible. In earlier work we therefore provided a general theory \cite{DuitsRThesis}, \cite{DuitsMMaster}, \cite{DuitsM2004}, to construct wavelet transforms associated with admissible vectors/ distributions.\footnote{Depending whether images are assumed to be band-limited or not, for full details see \cite{DuitsMandRart2004}.} With these wavelet transforms we construct orientation scores $U_f :
\mathbb{R}^2 \rtimes \mathbb{T} \to \mathbb{C}$ by means of admissible line detecting vectors\footnote{Or rather admissible distributions~$\psi \in \mathbb{H}^{-(1+\epsilon),2}(\R^2)$, $\epsilon>0$ if one does not want a restriction to bandlimited images.} $\psi \in \mathbb{L}_{2} (\R^2)$ such that the transform $\mathcal{W}_{\psi}$ is \emph{unitary} onto the unique reproducing kernel Hilbert space $\mathbb{C}^{SE(2)}_{K}$ of functions on $SE(2)$ with reproducing kernel $K(g,h)=(\mathcal{U}_{g}\psi,\mathcal{U}_{h}\psi)$, which is a closed vector subspace of $\mathbb{L}_{2}(SE(2))$. For the abstract construction of the unique reproducing kernel space $\mathbb{C}^{\mathbb{I}}_{K}$ on a set $\mathbb{I}$ (not necessarily a group) from a function of positive type $K:\mathbb{I}\times \mathbb{I} \to \mathbb{C}$, we refer to the early work of Aronszajn \cite{Aronszajn}. Here we only provide the essential Plancherel formula, which can also be found in a slightly different way in the work of F\"{u}hr
\cite{Fuehrbook}, for the wavelet transform $\mathcal{W}_{\psi}$ and which provides a more tangible description of the norm on $\mathbb{C}^{SE(2)}_{K}$ rather than the abstract one in \cite{Aronszajn}.
To this end we note that we can write
\[
\begin{array}{l}
(\mathcal{W}_{\psi}f)(\ul{x},e^{i\theta}) =(\mathcal{U}_{(\ul{x},e^{i\theta})}\psi,f)_{\mathbb{L}_{2}(\R^2)} =(\mathcal{F}\mathcal{T}_{\ul{x}}\mathcal{R}_{\theta}\psi,\mathcal{F}f)_{\mathbb{L}_{2}(\R^2)} 
 = \mathcal{F}^{-1}(\overline{\mathcal{R}_{\theta}\mathcal{F}\psi} \cdot \mathcal{F}f)(\ul{x})
\end{array}
\]
where the rotation and translation operators on $\mathbb{L}_{2}(\R^2)$ are defined by $\mathcal{R}_{\theta}f(\ul{y})=f(R_{\theta}^{-1}\ul{y})$ and $\mathcal{T}_{\ul{x}}f(\ul{y})=f(\ul{y}-\ul{x})$. Consequently, we find that
{\small
\begin{equation} \label{Parceval}
\begin{array}{ll}
\|\mathcal{W}_{\psi}f\|^{2}_{\mathbb{C}_{K}^{SE(2)}} &=  \int \limits_{\R^2} \int \limits_{\mathbb{T}} |(\mathcal{F} \mathcal{W}_{\psi}f)(\www, e^{i\theta})|^2 {\rm d}\theta \, \frac{1}{M_{\psi}(\www)}\, {\rm d}\www \\ &= \int \limits_{\R^2} \int \limits_{\mathbb{T}} |(\mathcal{F}f)(\www)|^2 |\mathcal{F} \psi(R_{\theta}^{T}\www)|^2 {\rm d}\theta \,  \frac{1}{M_{\psi}(\www)}\, {\rm d}\www \\
& = \int_{\R^2} |(\mathcal{F}f)(\www)|^2 {\rm d}\www= \|f\|_{\mathbb{L}_{2}(\R^2)}^{2},
\end{array}
\end{equation}
}
where $M_{\psi} \in C(\R^2,\R)$ is given by $M_{\psi}(\www)= \int_{0}^{2\pi} |\mathcal{F}\psi(R_{\theta}^{T}\www)|^2{\rm d}\theta$. If $\psi$ is chosen such that $M_{\psi}=1$ then we gain $\mathbb{L}_{2}$-norm preservation. However, this is not possible as $\psi \in \mathbb{L}_{2}(\R^2)\cap \mathbb{L}_{1}(\R^2)$ implies that $M_{\psi}$ is a continuous function vanishing at infinity. Now theoretically speaking one can use a Gelfand-triple structure generated by $\sqrt{1+|\Delta|}$ to allow distributional wavelets\footnote{Just like the Fourier transform on $\mathbb{L}_{2}(\R^2)$, where $\ul{x} \mapsto e^{i \www \cdot \ul{x}}$ is not within $\mathbb{L}_{2}(\R^2)$.}
$\psi \in \mathbb{H}^{-k}(\R^2)$, $k>1$ with the property $M_{\psi}=1$, so that $\psi$ has equal length in each irreducible subspace (which uniquely correspond to the dual orbits of $SO(2)$ on $\R^2$), for details and generalizations see \cite{DuitsMandRart2004}. In practice, however, because of finite grid sampling, we can as well restrict $\mathcal{U}$ (which is well-defined) to the space of \mbox{bandlimited} images.

\begin{figure}[t]
\centering
\begin{tabular}{cccc}
{\small (a)} & {\small (b)} & {\small(c)} & {\small (d)} \\
\includegraphics[width=.23\linewidth]{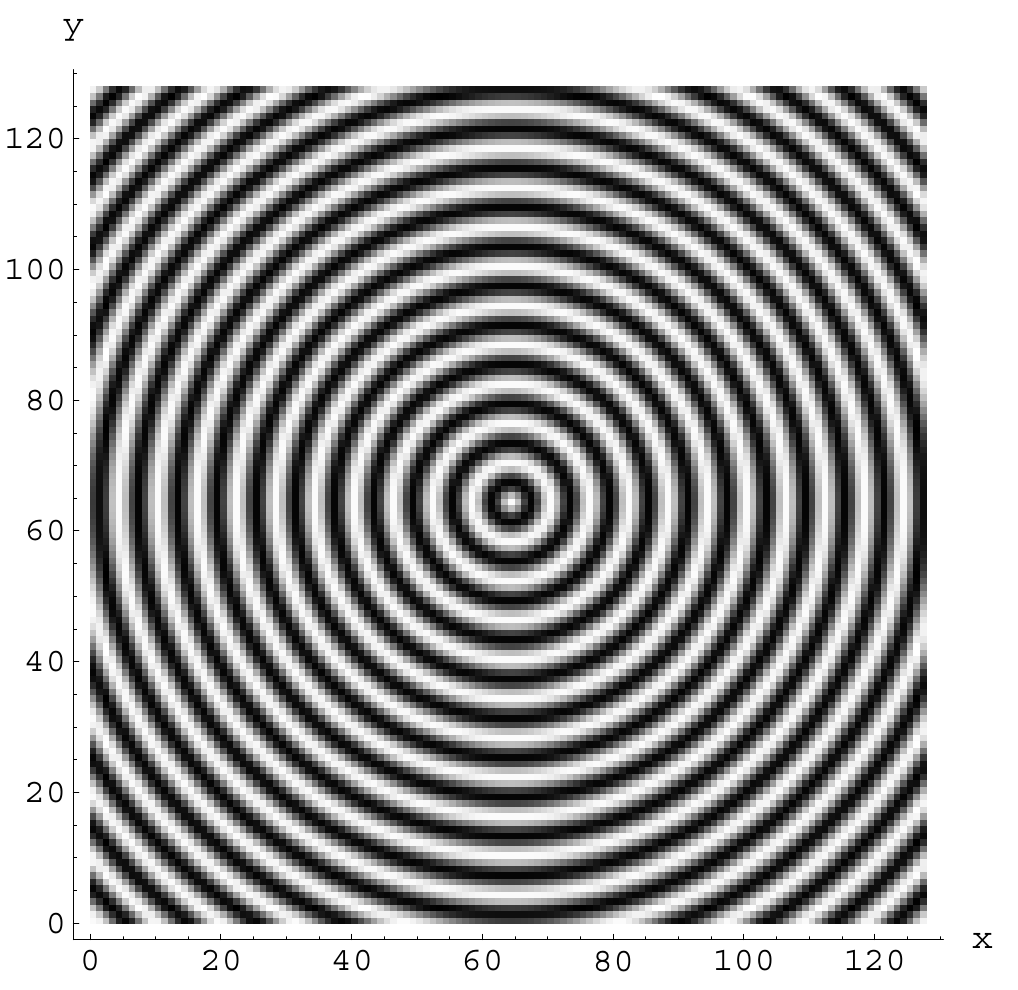} &
 \hspace{-4mm}\includegraphics[width=.19\linewidth]{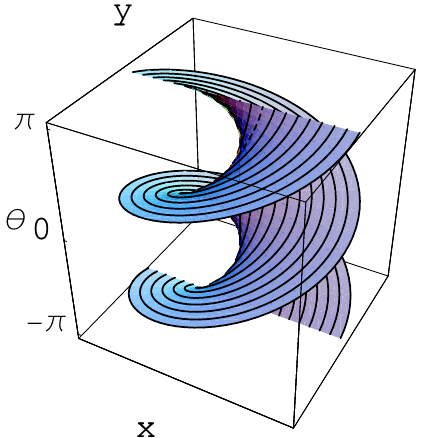}  \hspace{2mm}& \hspace{-2mm}\includegraphics[width=.245\linewidth]{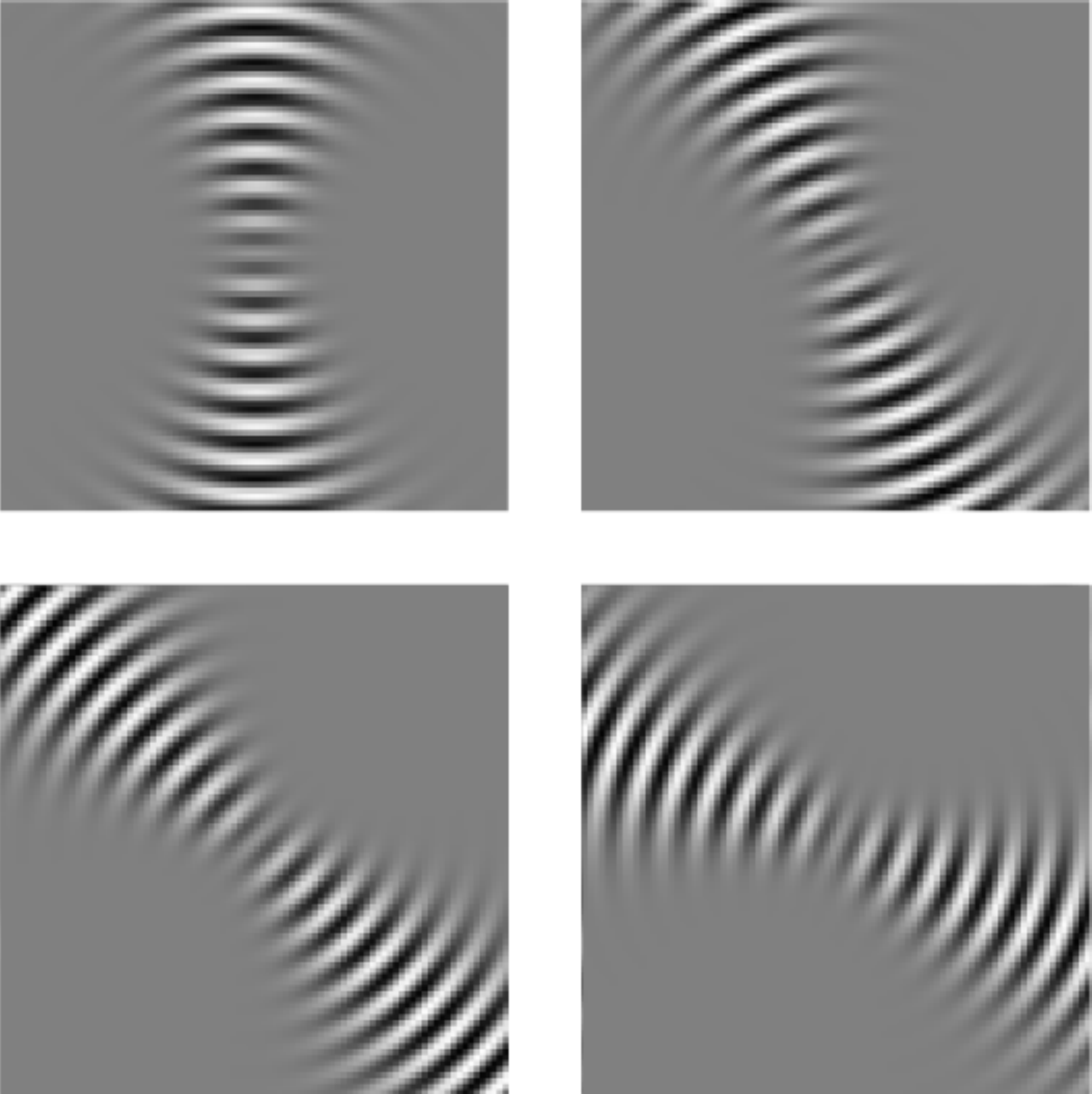} \hspace{-2mm} & \includegraphics[width=.245\linewidth]{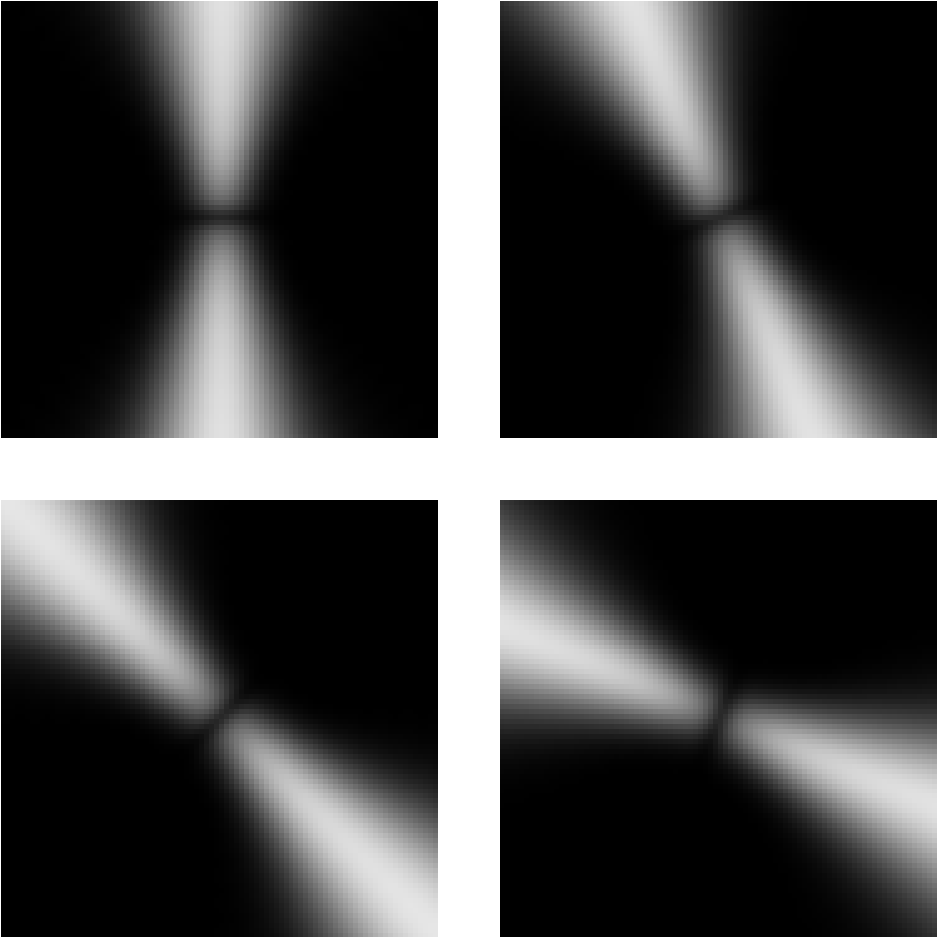} \\
\end{tabular}
\begin{tabular}{cccc}
(e) & (f) & (g) & (h) \\
\includegraphics[width=.2\linewidth]{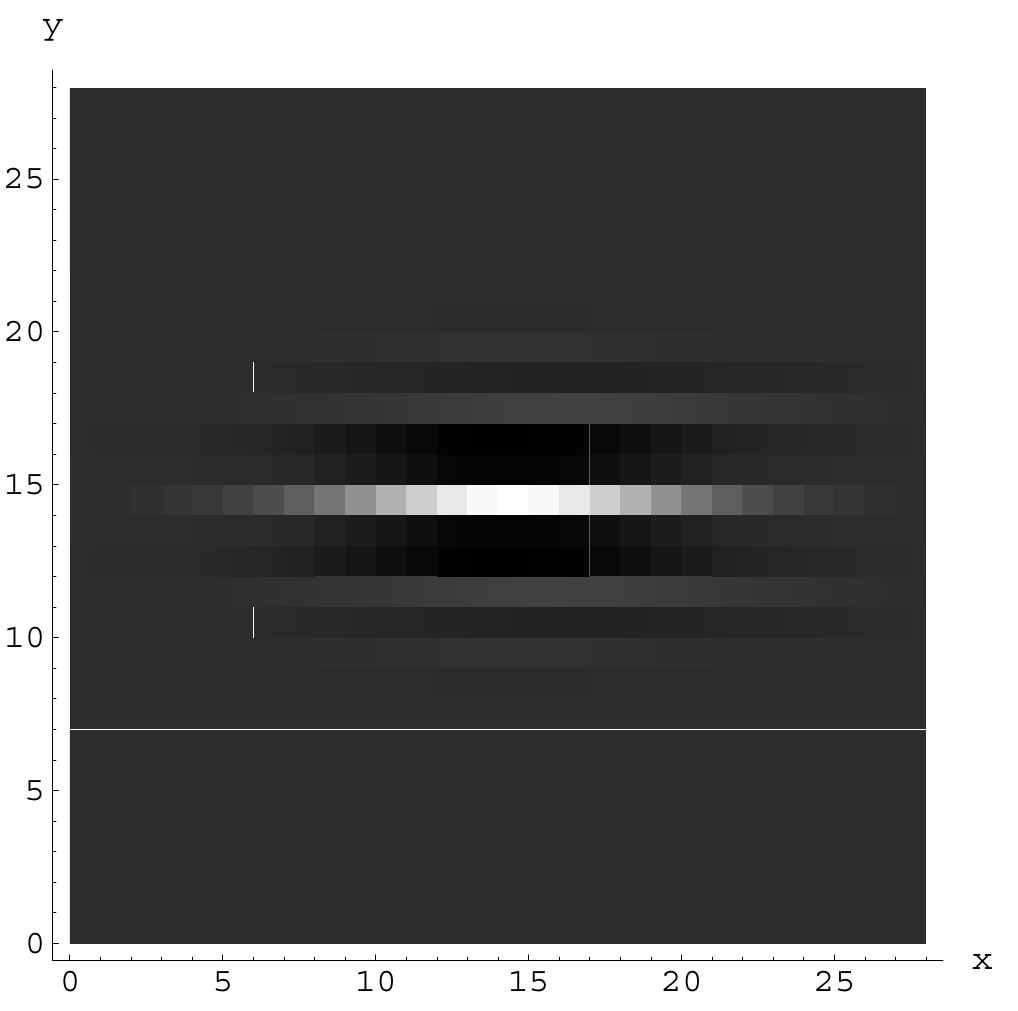} &\includegraphics[width=.2\linewidth]{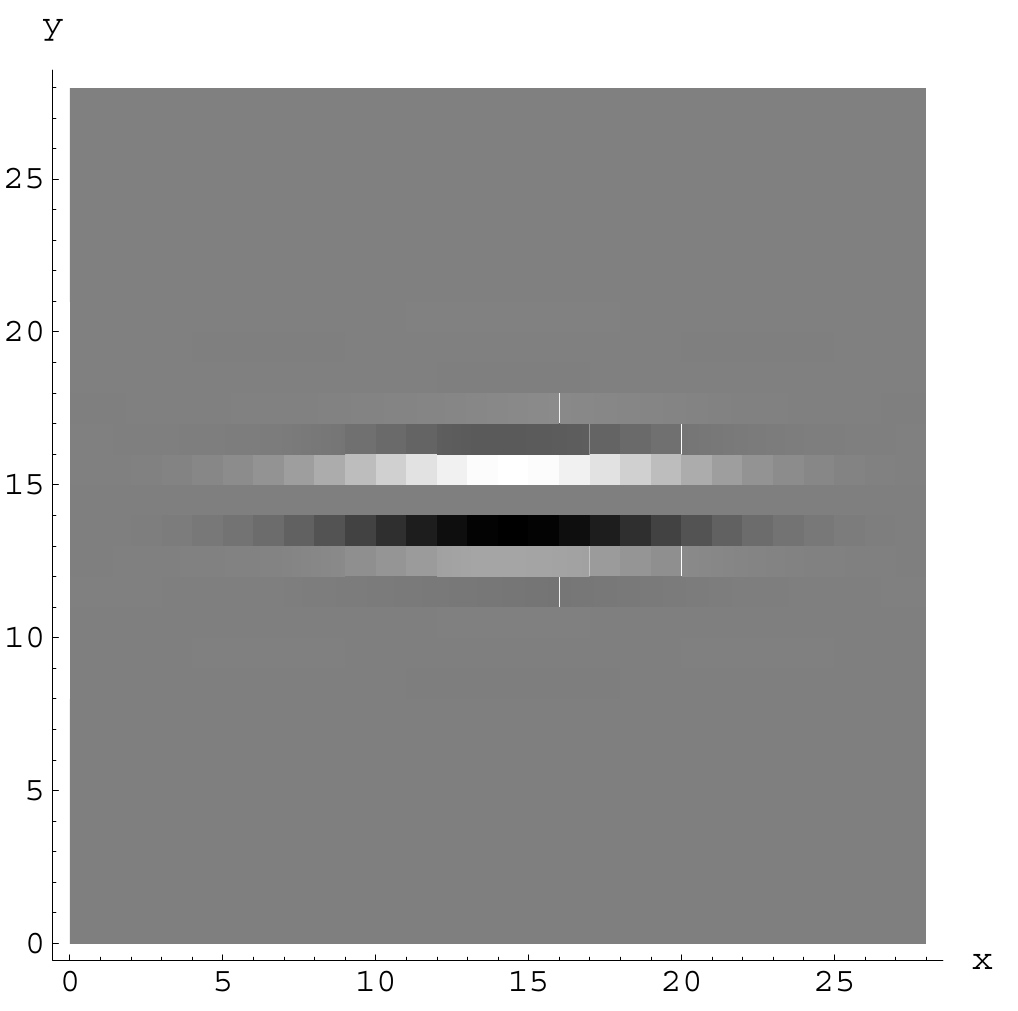} & \includegraphics[width=.2\linewidth]{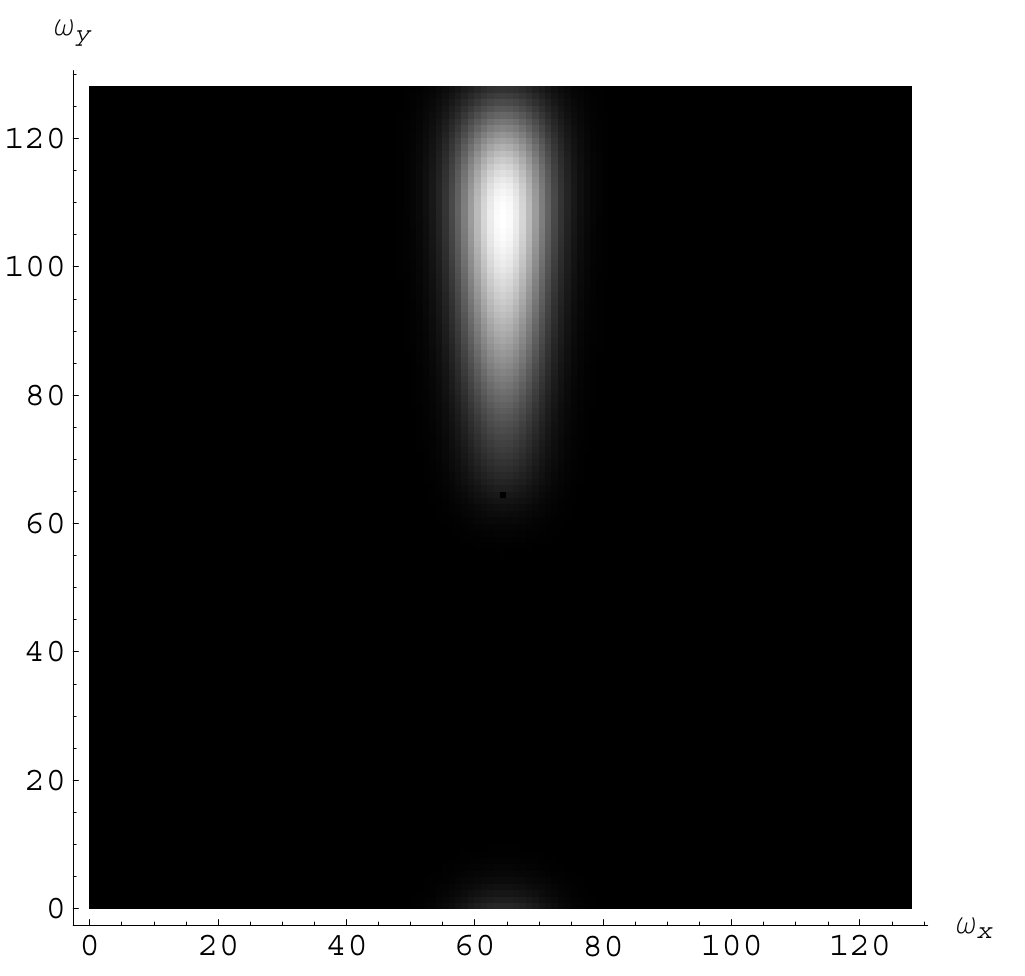} &\includegraphics[width=.2\linewidth]{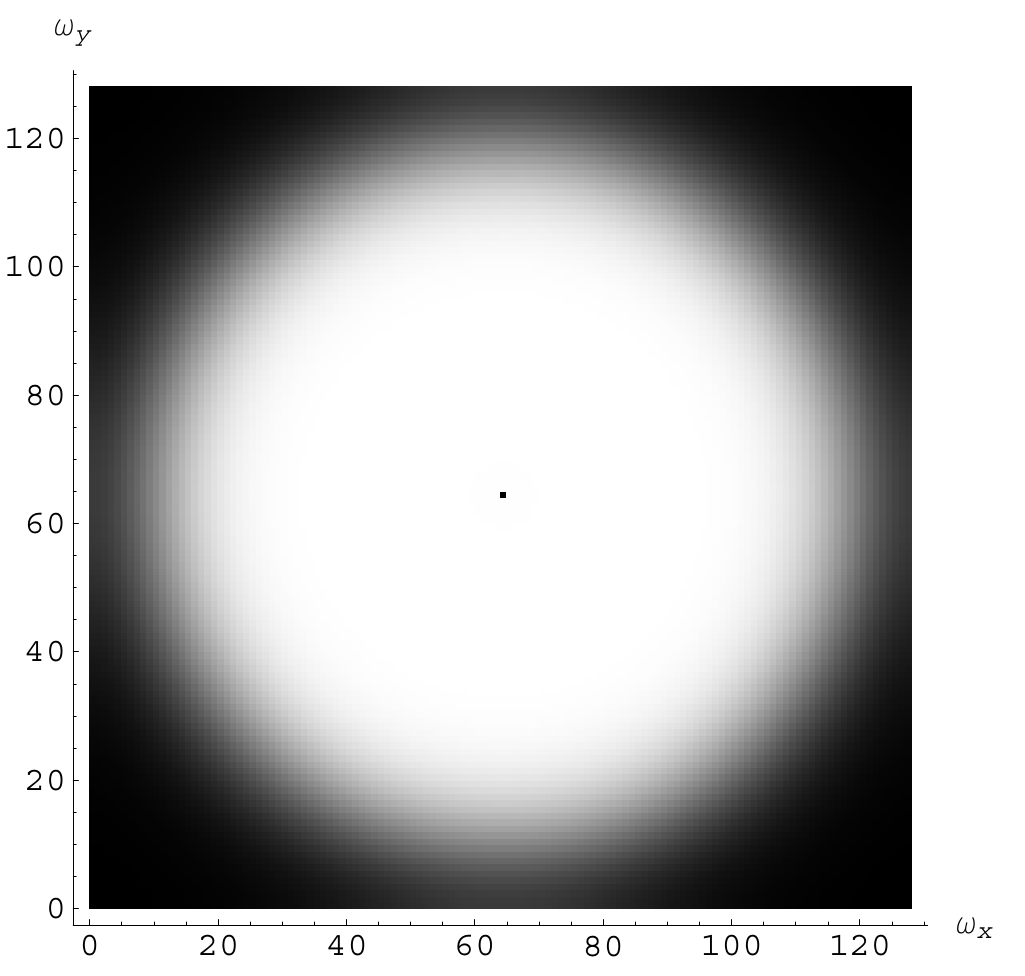} \\
\end{tabular}
\caption{(a) Example image $(x,y) \mapsto f(x,y)$. (b) The structure of the corresponding orientation score $U_f:=\mathcal{W}_{\psi}[f]$. The circles become spirals and all spirals are situated in the same helicoid-shaped plane.
(c) Real part of orientation score $(x,y)\mapsto U_f(x,y,e^{i\theta})$ displayed for 4 different fixed orientations. (d) The absolute vale $(x,y) \mapsto |U_f(x,y,e^{i\theta})|$ yields a \emph{phase-invariant and positive response} displayed for 4 fixed orientations.
(e) Real part of the wavelet {\small $\psi(\ul{x})=\frac{e^{-\frac{\|\ul{x}\|^{2}}{4s}}}{\sqrt{4\pi s}}\, \mathcal{F}^{-1}[\www \mapsto B^{k}\left( \frac{n_{\theta}((\phi\textrm{mod}\, 2\pi \!) - \! \frac{\pi}{2})}{2\pi}\right) \mathcal{M}(\rho)](\ul{x})$}, where {\small
$
\mathcal{M}(\rho)= \frac{e^{-\frac{\rho^2}{2 \sigma^2}}}{\sum_{k=0}^{q} (-1)^{k} \left(2^{-1} \sigma^{-2}\rho^{2} \right)^{k}},
$} with $\sigma=\frac{\varrho}{2}$ and Nyquist frequency $\varrho$ and $k$-th order $B$-spline $B_{k}=B_{0}*^{k}B_{0}$ and $B_{0}(x)=1_{[-\frac{1}{2},\frac{1}{2}]}(x)$ and parameter values $k=2$, $q=4$, $\frac{1}{2}\sigma^2=400$, $s=10$, $n_\theta=64$. (f) Imaginary part of $\psi$. (g) The function $|\mathcal{F}\psi|^{2}$ (h) The function $M_{\psi}$. }\label{fig:cakeexample}
\end{figure}

Finally, since the wavelet transform $\mathcal{W}_{\psi}$ maps the space of images $\mathbb{L}_{2}(\R^2)$ unitarily onto the space of orientation scores $\mathbb{C}^{SE(2)}_{K}$ (provided that $M_{\psi}>0$) we can reconstruct the original image $f:\R^2 \to \R$ from its orientation score $U_f: SE(2) \to \mathbb{C}$ by means of the adjoint
\begin{equation} \label{trueinvOBF}
\hspace*{-0.1cm}
\begin{array}{ll}
f &=\mathcal{W}_{\psi}^{*}\mathcal{W}_{\psi}[f]=\mathcal{F}^{-1}\left[
\www \mapsto \int_{0}^{2\pi} \mathcal{F}
[U_f(\cdot,e^{i\theta})](\www) \;  \mathcal{F}[
\mathcal{R}_{e^{i\theta}}\psi](\www) \; {\rm d \theta}
 \; M^{-1}_{\psi}(\www) \right]
 \end{array}
\end{equation}
For typical examples (and different classes) of wavelets $\psi$ such that $M_{\psi}=1$ and details on fast approximative reconstructions see \cite{Fran06b}, \cite{Duits2005IJCV},\cite{Duits2007PRIA}. For an illustration of a typical proper wavelet $\psi$ (i.e. $M_{\psi}\approx 1$) with corresponding transformation $\mathcal{W}_{\psi}f$ and corresponding $M_{\psi}:\R^2 \to \R^{+}$ usually looks like in our relatively fast algorithms, working with discrete subgroups of the torus, see Figure \ref{fig:cakeexample}.
With this well-posed, unitary transformation between the space of images and the space of orientation scores at hand, we can perform image processing via orientation scores, see \cite{Duits2007PRIA}, \cite{Duits2005IJCV}, \cite{Duits2007PRIA}, \cite{Duits2004PRIA}, \cite{Kali99a}.
However, for the remainder of the article we assume that the object~$U_f$ is some given function in $\mathbb{L}_{2}(SE(2))$ and we write $U \in \mathbb{L}_{2}(SE(2))$ rather than $U_f \in \mathbb{C}^{SE(2)}_{K}$.
For all image analysis applications where an object $U_f \in \mathbb
{L}_{2}(SE(2))$ is constructed from an image $f \in \mathbb{L}_{2}(\R^2)$, operators on the object~$U \in \mathbb{L}_{2}
(SE(2))$ must be left-invariant to ensure Euclidean invariant image processing \cite{DuitsRThesis}{p.153}. This applies also to the cases where the original image cannot be reconstructed in a stable manner as in channel representations \cite{Forssenchannel} and steerable tensor voting \cite{Fran04}.

\section{Left-invariant Diffusion on the Euclidean Motion Group \label{ch:2}}

\begin{figure}
\centerline{%
\includegraphics[width=0.8\hsize]{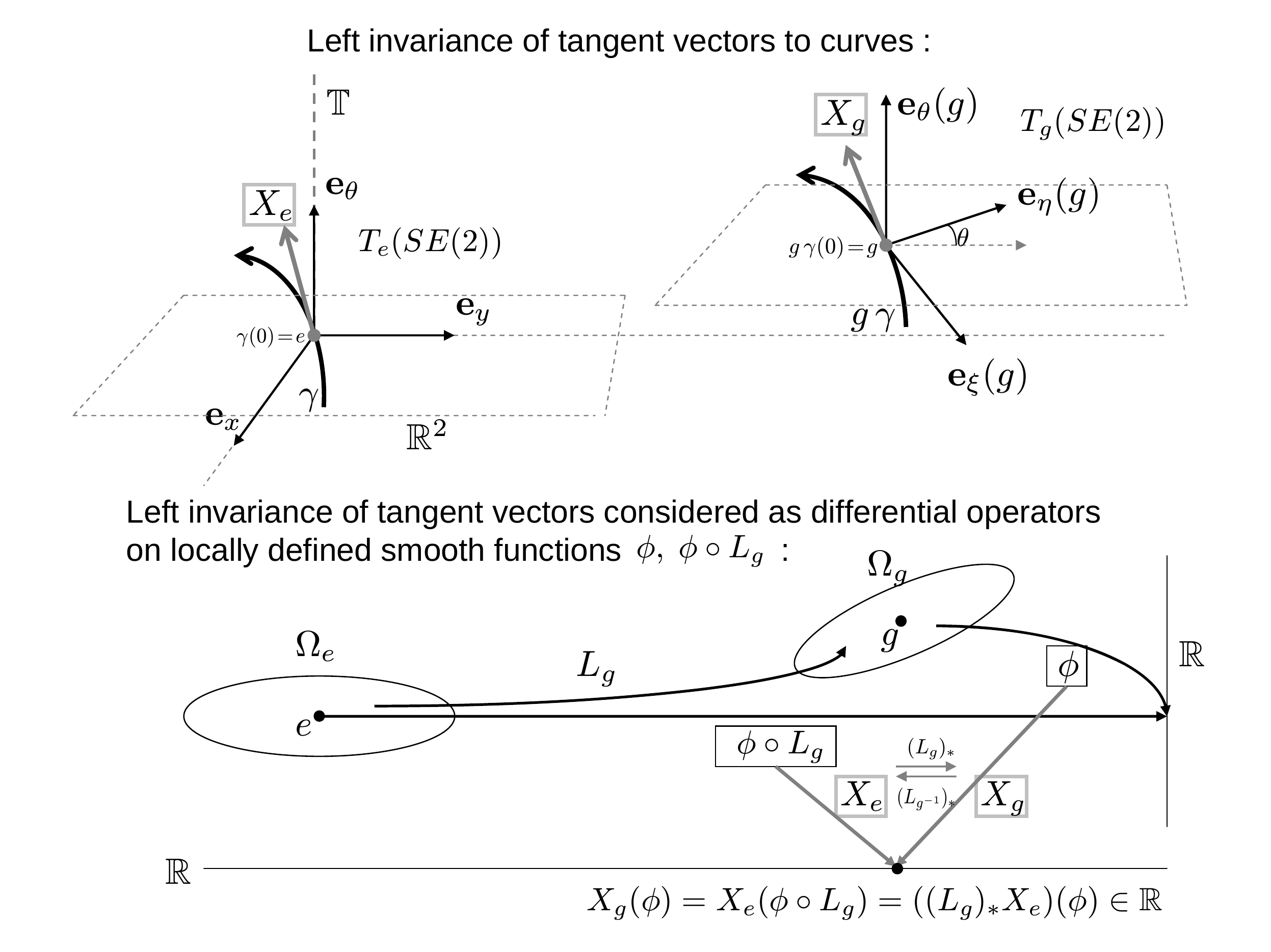}
}
\caption{Left invariant vector fields on $SE(2)$, where we both consider the tangent  vectors tangent to curves, that is $X_{g}=c^{1}\ul{e}_{\theta}(g) +c^{2} \ul{e}_{\xi}(g) +c^{3} \ul{e}_{\eta}(g)$ for all $g \in SE(2)$, and as differential operators on locally defined smooth functions, that is $X_{g}=c^{1}\left. \partial_{\theta} \right|_{g} +c^{2} \left. \partial_{\xi} \right|_{g} +c^{3} \left.\partial_{\eta}\right|_{g}$ for all $g \in SE(2)$. We see that the push forward of the left multiplication connects the tangent space $T_{e}(SE(2))$ to all tangent spaces $T_{g}(SE(2))$. Conversely,
the Cartan connection $D=d+\omega$ (\ref{Cartancomp}) on the vector bundle $(SE(2),T(SE(2)))$ connects all tangent spaces to $T_{e}(SE(2))$. Here we note that
$\omega_{g}(X_{g}):=(L_{g^{-1}})_{*} X_g= X_{e}$
for all left-invariant vector fields $X$. }\label{fig:leftinvariant}
\end{figure}

The group product within the group $SE(2)$ of planar translations and rotations is given by
\[
g g'=(\ul{x},e^{i \theta}) (\ul{x}',e^{i \theta'})= (\ul{x}+R_{\theta}\ul{x}',e^{i (\theta+\theta')}), \ \  g=(\ul{x},e^{i\theta}), g'=(\ul{x}',e^{i \theta'}) \in SE(2),
\]
with $R_{\theta}= \left(\begin{array}{cc}
\cos \theta & -\sin \theta \\
\sin \theta & \cos \theta
\end{array} \right) \in SO(2)$. The tangent space at the unity element $e=(0,0,e^{i0})$, $T_{e}(SE(2))$, is a 3D Lie algebra equipped with Lie product
{\small
$[A,B]= \lim_{t \downarrow 0}
t^{-2}\left(a(t)b(t)(a(t))^{-1}(b(t))^{-1}-e \right) ,$}
where $t \mapsto a(t)$ resp. $t\mapsto b(t)$ are \emph{any} smooth
curves in $G$ with $a(0)=b(0)=e$ and $a'(0)=A$ and $b'(0)=B$. Define $\{A_{1},A_{2},A_{3}\}:=\{\ul{e}_{\theta},\ul{e}_{x},\ul{e}_{y}\}$. Then $\{A_{1},A_{2},A_{3}\}$ form a basis of $T_{e}(SE(2))$ and their Lie-products are
\begin{equation} \label{Lieproduct}
\ [A_{1},A_{2}]=A_3, \
\ [A_{1},A_{3}]=-A_{2},\
\ [A_{2},A_{3}]=0\ .
\end{equation}
A vector field on $SE(2)$ 
is called left-invariant if for all $g \in G$ the push-forward of $(L_{g})_{*}X_{e}$ by left multiplication $L_g h=gh$ equals $X_{g}$, that is
\begin{equation} \label{leftinv}
(X_{g})=(L_{g})_{*}(X_{e}) \desda X_{g}f= X_{e}(f\circ L_{g}), \textrm{ for all }f \in C^{\infty}:\Omega_{g} \to \R,
\end{equation}
where $\Omega_{g}$ is some open set around $g \in SE(2)$. Recall that the tangent space at the unity element $e=(0,0,e^{i0})$, $T_{e}(G)$, is spanned by
{\small $T_{e}(G)=\textrm{span}\{\ul{e}_{\theta},\ul{e}_{x},\ul{e}_{y}\}= \textrm{span}\{(1,0,0),(0,1,0),(0,0,1)\}$}. By the general recipe of constructing left-invariant vector fields from elements in the Lie-algebra $T_{e}(G)$ (via the derivative of the right regular representation) we get the following basis for the space $\mathcal{L}(\textrm{SE}(2))$ of left-invariant vector fields :
\begin{equation}\label{leftinvSE2}
\{\mathcal{A}_{1},\mathcal{A}_{2},\mathcal{A}_{3}\}=\{\partial_{\theta},\partial_{\xi},\partial_{\eta}\}= \{\partial_{\theta},\cos \theta \, \partial_{x} +\sin \theta \, \partial_{y},-\sin \theta \, \partial_{x} +\cos \theta \, \partial_{y}\},
\end{equation}
with $\xi = x \, \cos \theta +y \, \sin \theta $, $\eta=-x\, \sin \theta  +y \, \cos \theta $. More precisely, the left-invariant vector-fields are
given by
{\small
\begin{equation} \label{baseonG}
\begin{array}{l}
\ul{e}_{\theta}(\ul{x},e^{i\theta})=\ul{e}_{\theta}, \ \;
\ul{e}_{\xi}(\ul{x},e^{i\theta}) = \cos \theta \, \ul{e}_{x} \!+\! \sin \theta \, \ul{e}_{y},  \ \;
\ul{e}_{\eta}(\ul{x},e^{i\theta}) =\!-\!\sin \theta \, \ul{e}_{x}\!+\! \cos \theta \, \ul{e}_{y},
\end{array}
\end{equation}
}
where we identified
$T_{g=(\ul{x},e^{i\theta})}(\R^2,e^{i\theta})$ with $T_{e}(\R^2,e^{i 0})$ and $T_{g=(\ul{x},e^{i\theta})}(\ul{x},\mathbb{T})$ with $T_{e}(\ul{0},\mathbb{T})$, by parallel transport (on $\R^2$ respectively $\mathbb{T}$). We can always consider these vector fields as differential operators (i.e. replace $\ul{e}_{i}$ by $\partial_{i}$, $i=\theta,\xi,\eta$), which yields (\ref{leftinvSE2}).
Summarizing, we see that for left-invariant vector fields the tangent vector at $g$ is related to
the tangent vector at $e$ by (\ref{leftinv}). In fact, the push forward $(L_{g})_{*}$ of the left multiplication puts a Cartan-connection\footnote{This Cartan connection can be related to a left-invariant metric induced by the Killing-form, which is degenerate on $SE(2)$. This can be resolved by pertubing the vectorfields into $\mathcal{L}(SO(3))\equiv so(3)$ by
{\small
$\{- \beta^2 y \cos \theta \, \partial_{\theta}\! +\! \cos \theta \, \sqrt{1\!+\!\beta^2 \, y^2}\, \partial_{x}\!+\! \sin \theta \,(1 \!+\!\beta^2 \, y^2) \, \partial_{y}, \beta^2 y \sin \theta \, \partial_{\theta} \!-\! \sin \theta \, \sqrt{1 \!+\!\beta^2 \,y^2}\, \partial_{x} \!+\! \cos \theta \, (1\!+\!\beta^2 \, y^2) \,\partial_{y}, \partial_{\theta}\}$, $0<\beta<<1$. See subsection \ref{ch:vectorbundles}. }.
} between tangent spaces, $T_{e}(SE(2))$ and $T_{g}(SE(2))$. 
Equality (\ref{leftinv}) sets the isomorphism between $T_{e}(SE(2))$ and $\mathcal{L}(SE(2))$, as $A_{i} \leftrightarrow \mathcal{A}_{i}$, $i=1,2,3$ implies
$[A_{i},A_{j}] \leftrightarrow [\mathcal{A}_{i},\mathcal{A}_{j}]$, $j=1,2,3$,  recall (\ref{Lieproduct}).
Moreover it is easily verified that
\[
\ [\mathcal{A}_{1},\mathcal{A}_{2}]=\mathcal{A}_1 \mathcal{A}_2-\mathcal{A}_2 \mathcal{A}_{1}=\mathcal{A}_3, \
\ [\mathcal{A}_{1},\mathcal{A}_{3}]=-\mathcal{A}_{2},\
\ [\mathcal{A}_{2},\mathcal{A}_{3}]=0\ .
\]
See Figure \ref{fig:leftinvariant} for a geometric explanation of left invariant vector fields, both considered as tangent vectors to curves in $SE(2)$ and as differential operators on locally defined smooth functions.
\textbf{Example: } \\
Consider $\ul{e}_{x} \in T_{e}(G)$, then the derivative of the right-regular representation
gives us
\begin{equation} \label{leftinv3}
\begin{array}{ll}
({\rm d}\mathcal{R}(A_{1})\Phi)(g)=({\rm d}\mathcal{R}(\ul{e}_{x})\Phi)(g)  &= \lim \limits_{h \downarrow 0}
\frac{\Phi(g \; e^{h \ul{e}_{x}})-\Phi(g)}{h} \\
 &= \lim \limits_{h \downarrow 0} \frac{\Phi(g \; (h,0,e^{i0}))-\Phi(g)}{h}\\ &=\lim \limits_{h \downarrow 0}
\frac{\Phi(\ul{x} + h (\cos \theta, \sin \theta),e^{i\theta})-\Phi(\ul{x},e^{i\theta})}{h}=\mathcal{A}_{1}\Phi(g) \\
 & =(\cos \theta \, \partial_{x} + \sin \theta \, \partial_{y}) \Phi(g)=\partial_{\xi} \Phi(g),
\end{array}
\end{equation}
for all $\Phi$ smooth and defined on some open environment around $g=(\ul{x},e^{i\theta}) \in G$.

%
Next we follow our general theory for left-invariant scale spaces on Lie-groups, see \cite{DuitsR2006SS2}, and set the following quadratic form on $\mathcal{L}(SE(2))$
{\small
\begin{equation} \label{QF}
Q^{\ul{D},\ul{a}}(\mathcal{A}_{1},\mathcal{A}_{2},\mathcal{A}_{3})= \sum_{i=1}^{3}\left( -a_{i} \mathcal{A}_{i} +\sum_{j=1}^{3}
D_{ij}\mathcal{A}_{i}\mathcal{A}_{j}\right),\  a_{i}, D_{ij} \in \R, D:=[D_{ij}]>0, D^{T}=D
\end{equation}
}
and consider the only 
linear left-invariant 2nd-order evolution equations
\begin{equation} \label{evolutionuptoorder2}
\boxed{
\left\{
\begin{array}{l}
\partial_{s}W=Q^{\ul{D},\ul{a}}(\mathcal{A}_{1},\mathcal{A}_{2},\mathcal{A}_{3}) \; W \ ,\\
\lim \limits_{s\downarrow 0}W(\cdot,s) =U_f(\cdot)\ .
\end{array}
\right.
}
\end{equation}
with corresponding resolvent equations (obtained by Laplace transform over $s$):
\begin{equation}\label{resolventequations}
\boxed{
P=\alpha(Q^{\ul{D},\ul{a}}(\mathcal{A}_{1},\mathcal{A}_{2},\mathcal{A}_{3})-\alpha I)^{-1}U_f.
}
\end{equation}
These resolvent equations are highly relevant as (for the cases $\ul{a}=\ul{0}$) they correspond to first order Tikhonov regularizations on $SE(2)$, \cite{DuitsR2006SS2}, \cite{Citti}.
They also have an important probabilistic interpretation, as we will explain next.

By the results in \cite{TerElst3}, \cite{Hebisch}, \cite{DuitsR2006AMS}, the solutions of these left-invariant evolution equations are given by $SE(2)$-convolution with the corresponding Green's function:
\begin{equation}\label{Gconv}
\begin{array}{l}
W(g,s)= (G_{s}^{\ul{D},\ul{a}}*_{SE(2)}U)(g) =\int \limits_{SE(2)} G_{s}^{\ul{D},\ul{a}}(h^{-1}g) U(h) \; {\rm d }\mu_{SE(2)}(h), \qquad g(\ul{x},e^{i\theta}),  \\
=\int \limits_{\R^2} \int \limits_{0}^{2\pi} G_{s}^{\ul{D},\ul{a}}(R_{\theta'}^{-1}(\ul{x}-\ul{x}'),e^{i(\theta-\theta')})\, U_{f}(\ul{x}',e^{i\theta'}) {\rm d}\theta' {\rm d}\ul{x}' \\[8pt]
P_{\alpha}(g)= (R_{\alpha}^{\ul{D},\ul{a}}*_{SE(2)}U)(g), \qquad R_{\alpha}^{\ul{D},\ul{a}}=\alpha \int \limits_{0}^{\infty}G_{s}^{\ul{D},\ul{a}} e^{-\alpha s} {\rm d}s.
\end{array}
\end{equation}
For Gaussian estimates of the Green's functions see the general results in \cite{TerElst3} and \cite{Hebisch}. See Appendix \ref{ch:terElst} for details on sharp Gaussian estimates for the Green's functions and formal proof of (\ref{Gconv}) in the particular case $D_{11}=D_{22}>0$, $D_{33}=0$ and $\ul{a}=\ul{0}$ (which is the Forward Kolmogorov equation of the contour enhancement process which we will explain next).

In the special case $D_{ij}=\delta_{i1}\delta_{j1}$, $\ul{a}=(\kappa_{0},1,0)$ our evolution equation (\ref{evolutionuptoorder2}) is the Kolmogorov equation
\begin{equation} \label{ForwardKolmogorovMumford}
\left\{
\begin{array}{l}
\partial_{s} W(g,s) = (\partial_{\xi} + D_{11}\partial_{\theta}^{2}) W(g,s), \qquad g \in SE(2),s>0 \\
W(g,0)=U(g)
\end{array}
\right.
\end{equation}
of Mumford's direction process, \cite{Mumford},
\begin{equation} \label{Mumford}
\left\{
\begin{array}{l}
\ul{X}(s)=X(s)\, \ul{e}_{x}+ Y(s)\, \ul{e}_{y} =\ul{X}(0)+ \int_{0}^{s}\cos \Theta(\tau)\, \ul{e}_{x}+\sin \Theta(\tau)\, \ul{e}_{y} \; {\rm d}\tau,  \\
\Theta(s)=\Theta(0)+ \sqrt{s} \, \epsilon_{\theta}+ s \, \kappa_{0} \, \qquad \epsilon_{\theta} \sim \mathcal{N}(0,2 D_{11}),
\end{array}
\right.
\end{equation}
for \emph{contour completion}. The explicit solutions of which we have derived in \cite{DuitsR2006AMS}.

However, within this article we will mainly focus on stochastic processes for \emph{contour enhancement}. For \emph{contour enhancement} we consider the particular case $\ul{a}=\ul{0}$. $D_{13}=D_{31}=D_{23}=D_{32}=0$. In particular we consider the case $D_{ij}=\delta_{ij}$, $D_{33}=0$, $\ul{a}=\ul{0}$ so that our evolution equation (\ref{evolutionuptoorder2}), becomes
\begin{equation} \label{ForwardKolmogorovCitti}
\left\{
\begin{array}{l}
\partial_{s} W(g,s) = (D_{11}(\partial_{\theta}^{2}) + D_{22}(\partial_{\xi})^{2}) W(g,s) \\
W(g,0)=U(g)
\end{array}
\right.
\end{equation}
which is the Kolmogorov equation of the following stochastic process for contour enhancement:
\begin{equation} \label{Cittiproc}
\left\{
\begin{array}{l}
\ul{X}(s)=\ul{X}(0)+\sqrt{s} \, \epsilon_{\xi} \; \int_{0}^{s}\cos \Theta(\tau) \; \ul{e}_{x}+\sin \Theta(\tau)\; \ul{e}_{y} \; {\rm d}\tau ,  \\
\Theta(s)=\Theta(0)+ \sqrt{s} \, \epsilon_{\theta}  ,
\end{array}
\right.
\end{equation}
with $\epsilon_{\xi} \sim \mathcal{N}(0,2 D_{22})$ and $\epsilon_{\theta} \sim \mathcal{N}(0,2 D_{11})$, $D_{11},D_{22}>0$.

In general the evolution equations (\ref{evolutionuptoorder2}) are the forward Kolmogorov equations of all linear left-invariant stochastic processes on $SE(2)$, as explained in \cite{DuitsR2006AMS}, \cite{vAlmsick2005}.
\begin{figure}
\centerline{%
\includegraphics[width=0.25\hsize]{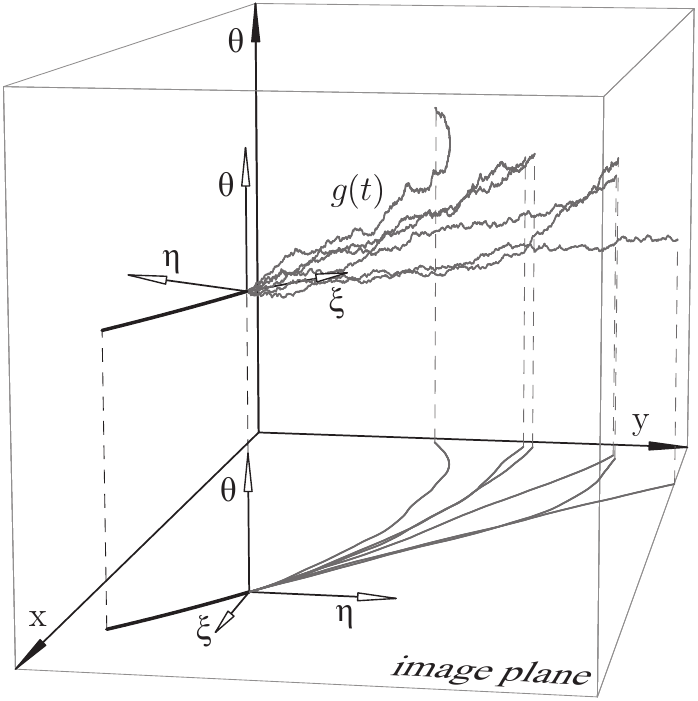}
\hfill
\includegraphics[width=0.33\hsize]{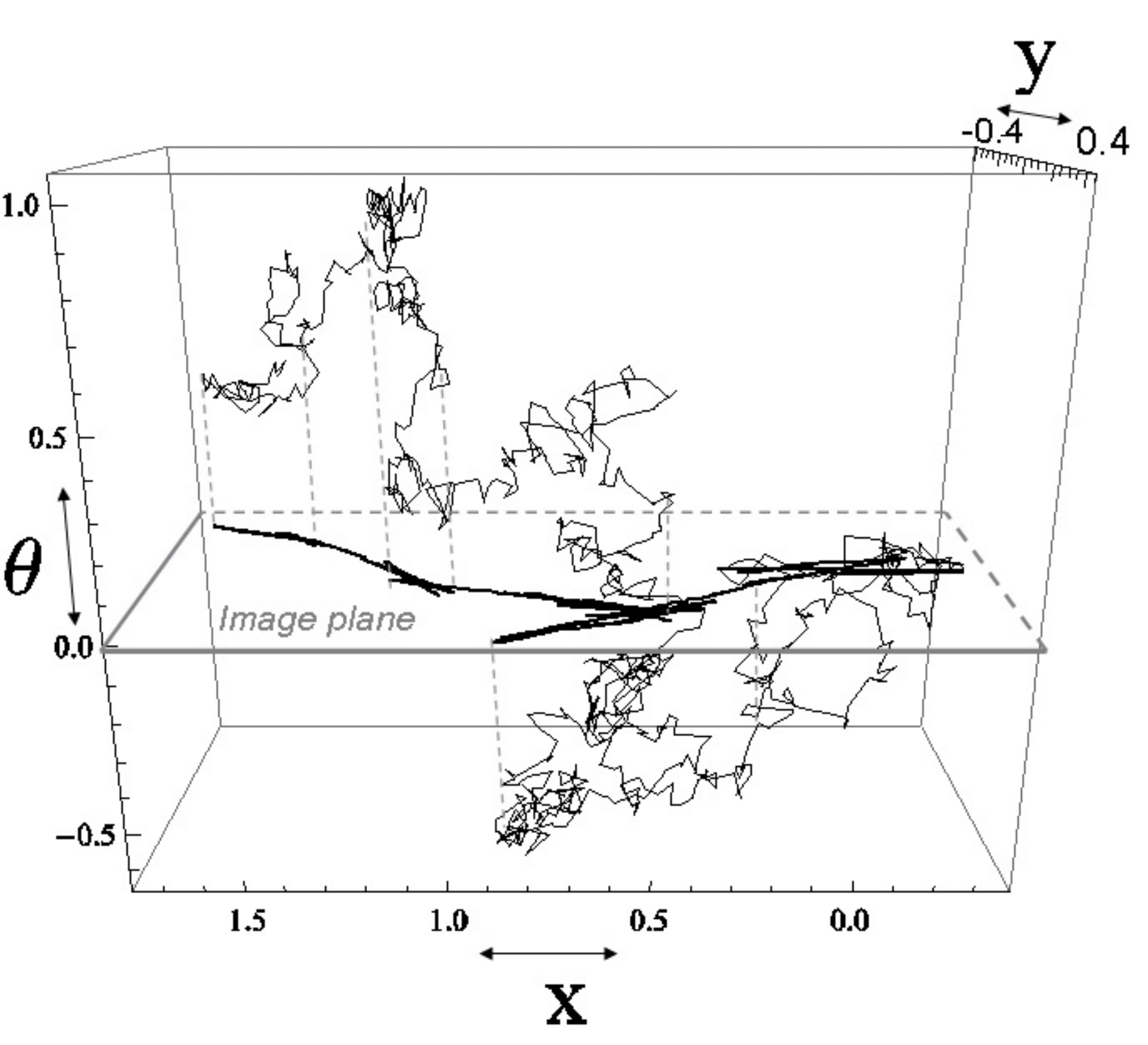}
\includegraphics[width=0.36\hsize]{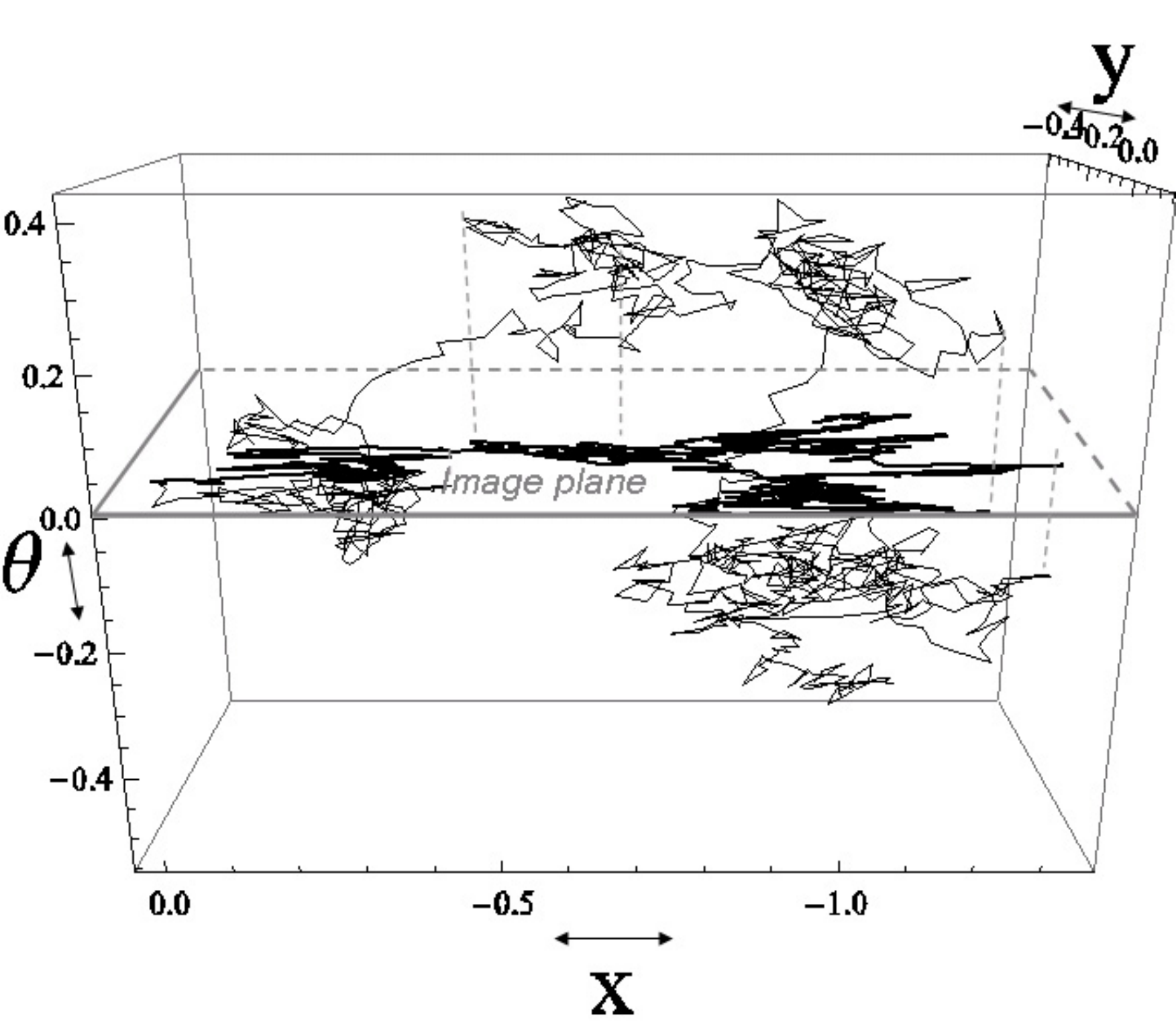}
}
\hbox{}\vspace{0.25cm}\hbox{}
\centerline{
\includegraphics[width=0.2\hsize]{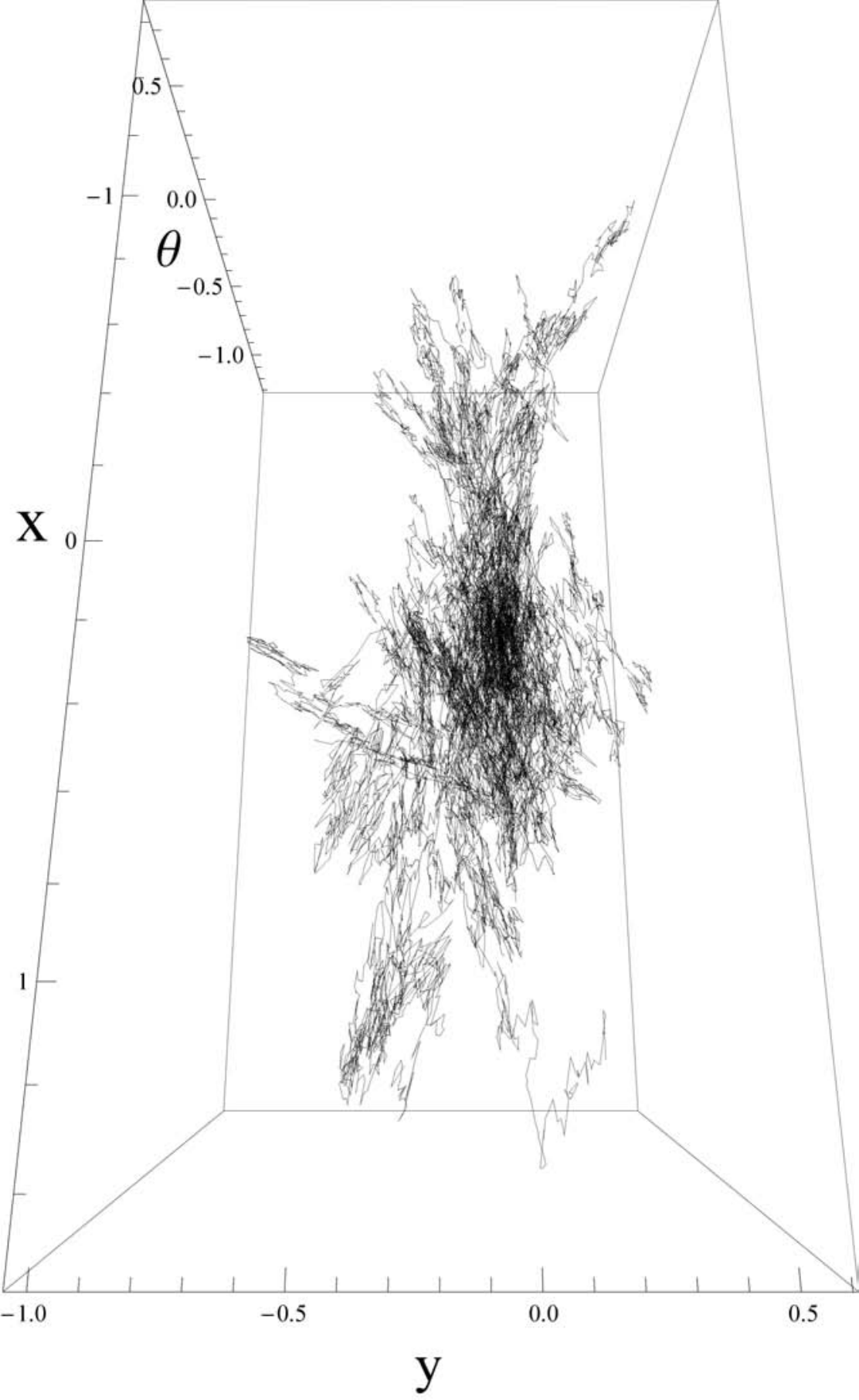}
\includegraphics[width=0.25\hsize]{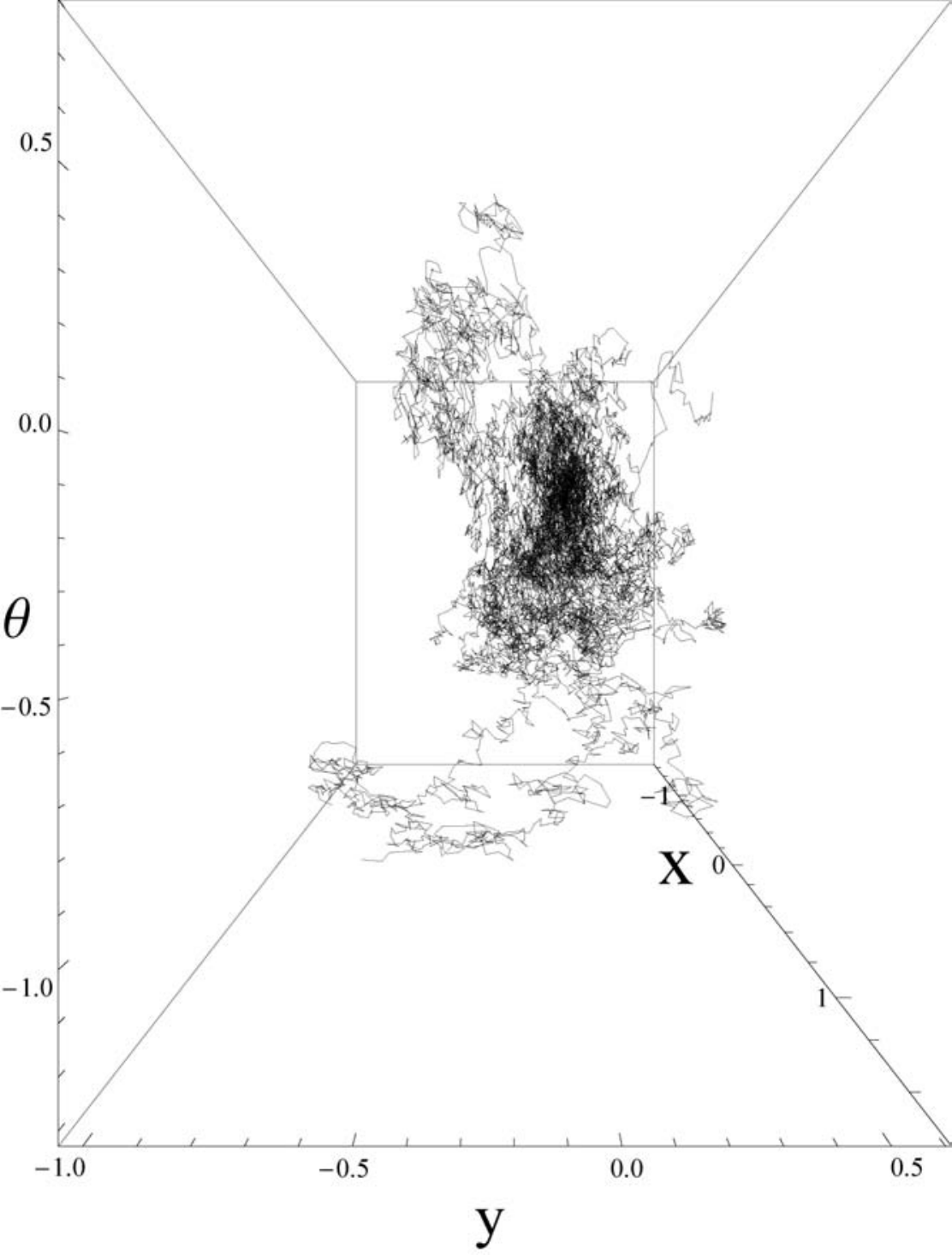}
\includegraphics[width=0.5\hsize]{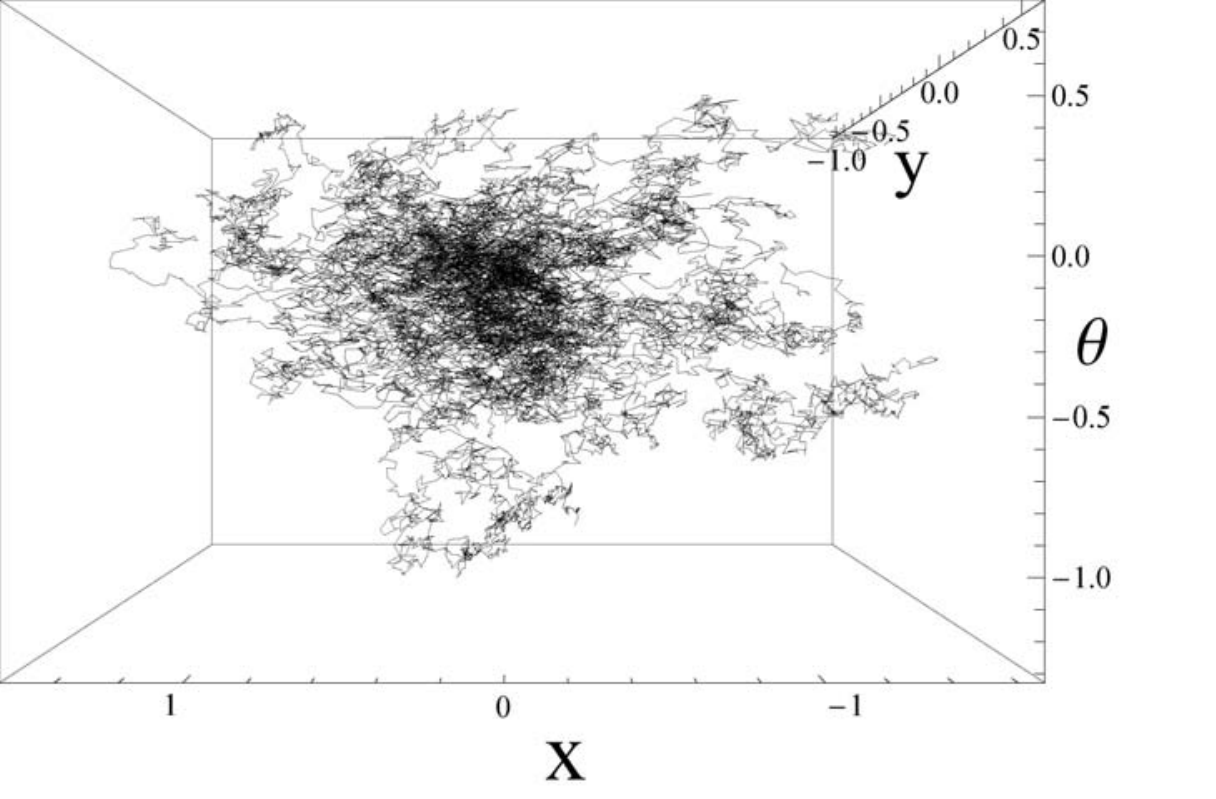}
}
\caption{Top left: six random walks in $SE(2)=\R^2 \rtimes \mathbb{T}$ (and their projection on $\R^2$) of direction processes for contour-completion by Mumford \cite{Mumford} with $\ul{a}=(\kappa_{0},1,0)$, $D=\textrm{diag}\{D_{11},D_{22},D_{33}\}$ for various parameter settings of $\kappa_{0}\geq 0$ and $D_{ii}>0$. Top middle one random walk ($500$ steps, with step-size $0.005$) and its projection to the image plane of the linear left-invariant stochastic process for contour enhancement within $SE(2)$ with parameter settings $D_{11}=D_{22}=\frac{1}{2}$ and $D_{33}=0$ (corresponding to Citti and Sarti's cortical model for contour enhancement, \cite{Citti}). Top right; one random walk ($800$ steps, with step-size $0.005$) of the stochastic process with parameter settings $D_{11}=\frac{1}{2}\sigma_{\theta}^{2}$, $D_{22}=\frac{1}{2}\sigma_{\xi}^{2}$, $D_{33}=\frac{1}{2}\sigma_{\eta}^{2}$, with $\sigma_{\theta}=0.75$, $\sigma_{\xi}=1$, $\sigma_{\eta}=0.5$ (other parameters have been set to zero).
Bottom: 30 random walks in $SE(2)=\R^2 \rtimes \mathbb{T}$ again with $D_{ij}=\delta_{ij}\sigma_{i}^2, \sigma_{\theta}=0.75$, $\sigma_{\xi}=1$, $\sigma_{\eta}=0.5$, viewed along $\theta$-axis (left) along $x$-axis (middle) along $y$-axis (right). Appropriate averaging of infinitely many of these sample paths yields the Green's functions, see Figure \ref{fig:comparisonspatialandFD}, of the forward Kolmogorov equations (\ref{evolutionuptoorder2}). Furthermore we note that Mumford's direction process is the only linear left-invariant stochastic process on SE(2) whose sample path projections on the image plane are differentiable. For contour completion this may be a reason to discard the other linear left-invariant stochastic processes, see \cite{MarkusThesis}. However, the Green's function of all linear left-invariant processes (so also the ones for contour-enhancement) are infinitely differentiable on $SE(2)\setminus \{e\}$ iff the H\"{o}rmander condition as we will discuss in section \ref{ch:hoermander}, see  (\ref{hoermandercond}), is satisfied.
}\label{fig:orientationbundle2}
\end{figure} 
With respect to this connection to probability theory we note that $W(g,s)$ represents the probability density of finding \emph{oriented} random walker\footnote{That is a random walker in the space $SE(2)$ where it is only allowed to move along horizontal curves which are curves whose tangent vectors always lie in $\textrm{span}\{\partial_{\theta}, \partial_{\xi}\}$ which is the horizontal subspace if we apply the Cartan connection on $P_{Y}=(SE(2),SE(2)/Y, \pi, R)$ see section \ref{ch:fiberbundles}. In previous work in the field of image analysis, \cite{Duits2005IJCV}, \cite{Duits2007PRIA}, we called these random walkers ``oriented gray value particles''. } (traveling with unit speed, which allows us to identify traveling time with arc-length $s$) at position $g$ given the initial distribution $W(\cdot,0)=U_f$ a traveling time $s>0$, whereas $P(g)$ represents the unconditional probability density of finding an \emph{oriented} random walker
at position $g$ given the initial distribution $W(\cdot,0)=U_f$ regardless its traveling time. To this end we note that traveling time $T$ in a Markov process is negatively exponentially distributed
\[
P(T=s)=\alpha e^{-\alpha s},
\]
since this is the only continuous memoryless distribution and indeed a simple calculation yields:
\begin{equation}\label{uncond}
\begin{array}{l}
P(x,y,\theta \; |\; U \textrm{ and } T=s)=(G_{s}^{D_{11}}*_{SE(2)}U)(x,y,\theta) 
\\[7pt]
P(x,y,\theta \; |\; U)= \int_{0}^{\infty} P(x,y,\theta \; |\; U \textrm{ and } T=s) P(T=s) {\rm d}s  
 =(R_{\alpha}^{D_{11}}*_{SE(2)}U)(x,y,\theta) \\[7pt]
 \textrm{ with } R_{s}^{D_{11}}=\alpha \int_{\R^{+}} G_{s}^{D_{11}}e^{-\alpha s}{\rm d}s,
\end{array}
\end{equation}
For exact solutions for the resolvent equations (\ref{resolventequations})(in the special case of Mumford's direction process), approximations and their relation to fast numerical algorithms, see \cite{DuitsR2006AMS}.


\section{Image Enhancement via left-invariant Evolution Equations on Invertible Orientation Scores \label{ch:leftinvimageproc}}

\begin{figure}
\centerline{\includegraphics[width=0.9\hsize]{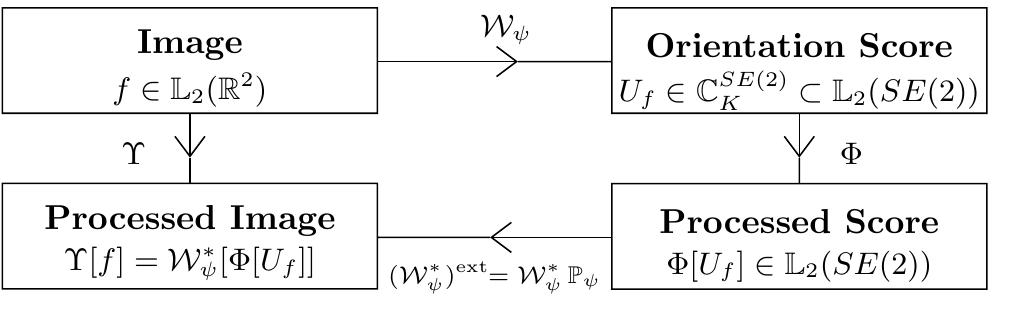}}
\centerline{%
\includegraphics[width=1.02\hsize]{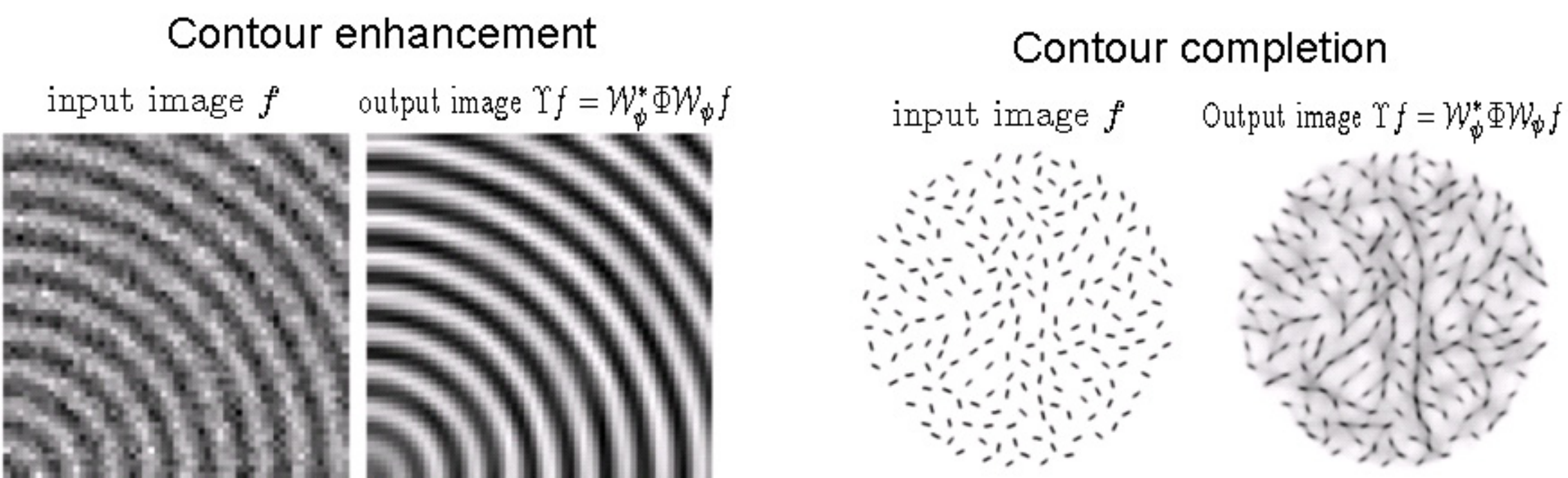}
}
\caption{
{\small Top Row: The complete scheme; for admissible vectors $\psi$ the linear map $\mathcal{W}_{\psi}$ is unitary from $\mathbb{L}_{2}(\R^2)$ onto a closed subspace {\small $\mathbb{C}_{K}^{SE(2)}$} of $\mathbb{L}_{2}(SE(2))$. So we can uniquely relate a transformation $\Phi:\mathbb{C}_{K}^{SE(2)} \to \mathbb{C}_{K}^{SE(2)}$ on an orientation score to a transformation on an image $\Upsilon=(\mathcal{W}_{\psi}^{*})^{ext }\circ \Phi \circ \mathcal{W}_{\psi}$, where $(\mathcal{W}_{\psi}^{*})^{ext}$ is given by (\ref{Wext}). 
Here we take $\Phi$ as a concatenation of \emph{non-linear} invertible greyvalue transforms and linear left-invariant evolutions (\ref{completionfield}), although we stress that for most practical applications it is better to replace the operator $\Phi$ by the adaptive evolution operator $\mathcal{W}_{\psi}f \mapsto u(x,y,e^{i\theta},t)$ defined by the non-linear adaptive left-invariant evolution equation (\ref{coherencesimple}) with certain stopping time $t>0$. Bottom row: automated contour enhancement (left) and completion (right). Parameter settings, left;
$D_{ij}=D_{ii}\delta_{ij}$, {\small $\alpha^{2}\frac{D_{11}}{D_{22}}=\frac{1}{4}$}, 
 $D_{33}=0$, $\ul{a}=\ul{0}$, right $D_{ij}=\delta_{i1}\delta_{j1}$, $\ul{a}=(0,1,0)$, $\frac{D_{11}}{\alpha}=0.1$.}
}\label{fig:leftinv}
\end{figure}

Now that we have constructed a stable transformation between images $f$ and corresponding orientation scores $U_f$, in Section \ref{ch:OS} we can relate operators $\Upsilon$ on images to operators $\Phi$ on
orientation scores in a robust manner, see Figure \ref{fig:leftinv}.
It is easily verified that $\mathcal{W}_{\psi} \circ \mathcal{U}_{g}=
\mathcal{L}_{g} \circ \mathcal{W}_{\psi}$ for all $g \in SE(2)$, where the left-representation $\mathcal{L}:G \to \mathcal{B}(\mathbb{L}_{2}(SE(2)))$ is given by $\mathcal{L}_{g}\Phi(h)=\Phi(g^{-1}h)$. Consequently, the net operator on the image $\Upsilon$ is Euclidean invariant if and only if the operator on the orientation score is left-invariant, i.e.
\begin{equation} \label{leftinvimportant}
\boxed{
\Upsilon \circ \mathcal{U}_{g}=\mathcal{U}_{g} \circ \Upsilon \textrm{ for all }g \in SE(2) \; \desda \; \Phi \circ \mathcal{L}_{g}=\mathcal{L}_{g} \circ \Phi \textrm{ for all }g \in SE(2),}
\end{equation}
see \cite{DuitsRThesis}{Thm. 21 p.153}.

Here the diffusions discussed in the previous section, section \ref{ch:2}, can be used to construct suitable operator $\Phi$ on the orientation scores.
At first glance the diffusions themselves (with certain stopping time $t>0$) or their resolvents (with parameter $\alpha>0$) seem suitable candidates for operators on orientation scores, as they follow from stochastic processes for contour enhancement and contour completion and they even map the space of orientation scores $\mathbb{C}^{SE(2)}_{K}$ into the space of orientation scores $\mathbb{C}^{SE(2)}_K$ again. But appearances are deceptive since if the operator $\Phi$ is left-invariant (which must be required, see Figure \ref{fig:leftinv}) and linear then the netto operator $\Upsilon$ is translation and rotation invariant boiling down to an isotropic convolution on the original image, which is of course not desirable.

So our operator $\Phi$ must be left-invariant and non-linear and still we would like to directly relate such operator to stochastic processes on $SE(2)$ discussed in the previous section.  Therefor we consider the operators
\begin{equation} \label{completionfield}
\begin{array}{l}
\Phi(U, V)=\alpha^{\frac{1}{p}}((Q^{D,\ul{a}}(\underline{A})-\alpha I)^{-1}(U)^p \; ((Q^{D,\ul{a}}(\underline{A}))^{*}-\alpha I)^{-1}(V)^{p})^{\frac{1}{2p}}, \\
         = \alpha^{\frac{1}{p}}\left((R_{\alpha}^{D,-\ul{a}}*_{SE(2)}(U)^p) \; \cdot \; (R_{\alpha}^{D,-\ul{a}}*_{SE(2)}(V)^{p})\right)^{\frac{1}{2p}}\, \qquad  p>0.
\end{array}
\end{equation}
where $U$ (the source distribution) and $V$ (the sink distribution) denote two initial distributions on $SE(2)$ and where we take the $p$-th power of both real part $\Re(U)$ and imaginary
part $\Im(U)$ separately in a sign-preserving manner, i.e. $(U)^{p}$ means $\textrm{sign}\{ \Re(U)\} \; |\Re(U)|^{p} + i \, \textrm{sign} \{\Im(U)\} \; |\Im(U)|^{p}$.
Here the function $\Phi(U, V) \in \mathbb{L}_{2}(SE(2))$ can be considered as the completion distribution\footnote{In image analysis these distributions are called ''completion fields``, where the word field is inappropriate. } obtained from collision of the forwardly evolving source distribution $U$ and backwardly evolving sink distribution $V$, similar to \cite{August2003}.

Within this manuscript we shall restrict ourselves to the case where both source and sink equal the orientation score of original image $f$, i.e. $U=V=U_f:=\mathcal{W}_{\psi}f$ and only occasionally (for example section \ref{ch:elastica}) we shall study the case where $U=\delta_{g_0}$ and $V=\delta_{g_1}$, where $g_0$ and $g_{1}$ are some given elements in $SE(2)$.

In section \ref{ch:CED} we shall consider more sophisticated and more practical alternatives to the operator given by (\ref{completionfield}). But for the moment we restrict ourselves to the case (\ref{completionfield}) as this is much easier to analyse and also much easier to implement as it requires two group convolutions (recall (\ref{Gconv})) with the corresponding Green's functions which we shall explicitly derive in the next section.

The relation between image and orientation score remains 1-to 1 if we ensure  that the operator on the orientation score again provides an orientation score of an image: Let {\small $\mathbb{C}_{K}^{SE(2)}$} denote\footnote{We use this notation since the space of orientation scores generated by proper wavelet $\psi$ is the unique reproducing kernel space on $SE(2)$ with reproducing kernel $K(g,h)=(\mathcal{U}_{g}\psi,\mathcal{U}_{h}\psi)$, \cite{DuitsRThesis}{p.221-222, p.120-122}} the space of orientation scores within $\mathbb{L}_{2}(SE(2))$, then the relation is 1-to 1 iff $\Phi$ maps {\small $\mathbb{C}_{K}^{SE(2)}$} into {\small $\mathbb{C}_{K}^{SE(2)}$}. However, we naturally extend the reconstruction to $\mathbb{L}_{2}(SE(2))$:
{\small
\begin{equation}\label{Wext}
\begin{array}{l}
(\mathcal{W}_{\psi}^{*})^{ext}U(g)=\mathcal{F}^{-1}\left[
\www \mapsto \int_{0}^{2\pi} \mathcal{F}
[U(\cdot,e^{i\theta})](\www)\; \mathcal{F}[
\mathcal{R}_{e^{i\theta}}\psi](\www) \; {\rm d \theta}
 \, M^{-1}_{\psi}(\www) \right], \
 \end{array}
 \end{equation}
 }
for all $U \in \mathbb{L}_{2}(SE(2))$.
So the effective part of a operator $\Phi$ on an orientation score is in fact $\mathbb{P}_{\psi}\Phi$ where $\mathbb{P}_{\psi}=\mathcal{W}_{\psi}(\mathcal{W}_{\psi}^{*})^{ext}$ is the orthogonal projection of $\mathbb{L}_{2}(SE(2))$ onto {\small $\mathbb{C}^{SE(2)}_{K}$}. Recall that $\Phi$ must be left-invariant because of (\ref{leftinvimportant}).

It is not difficult to show that the only linear left-invariant kernel operators on $\mathbb{L}_{2}(SE(2))$ are given by $SE(2)$-convolutions. Recall that these kernel operators are given by (\ref{Gconv}). Even these $SE(2)$-convolutions do not leave the space of orientation scores $\mathbb{C}^{SE(2)}_{K}$ invariant. Although,
\[
\begin{array}{ll}
(K*_{SE(2)} \mathcal{W}_{\psi}f)(g) &= \int \limits_{SE(2)} (\mathcal{U}_{h}\psi,f)_{\mathbb{L}_{2}(\R^2)} K(h^{-1}g) {\rm d}\mu_{SE(2)}(h)  \\
                                    &=  (\int \limits_{SE(2)} \mathcal{U}_{h}\psi \, K(h^{-1}g) {\rm d}\mu_{SE(2)}(h),f)_{\mathbb{L}_{2}(\R^2)} \\
                                    &=  (\int \limits_{SE(2)} \mathcal{U}_{g\tilde{h}^{-1}}\psi \, K(\tilde{h}) {\rm d}\mu_{SE(2)}(\tilde{h})\, ,f)_{\mathbb{L}_{2}(\R^2)}= (\mathcal{U}_{g}\tilde{\psi},f)_{\mathbb{L}_{2}(\R^2)}=\mathcal{W}_{\tilde{\psi}}f(g).
\end{array}
\]
for all $f \in \mathbb{L}_{2}(\R^2)$, $g \in SE(2)$, where $\tilde{\psi}=\int_{SE(2)} \mathcal{U}_{\tilde{h}^{-1}}\psi \, K(\tilde{h}) \; {\rm d}\mu_{SE(2)}(\tilde{h})$,
the reproducing kernel space associated to $\tilde{\psi}$ will in general not coincide with the reproducing kernel space associated to $\psi$. Here we recall from \cite{DuitsM2004},\cite{DuitsRThesis}, that $\psi$ determines the reproducing kernel $K(g,h)=(\mathcal{U}_{g}\psi,\mathcal{U}_{h}\psi)_{\mathbb{L}_{2}(\R^2)}$.

\section{The Heat-Kernels on $SE(2)$.\label{ch:diffSE2}}

In section \ref{ch:exactheatkernels} we present the exact formulae, which do not seem to appear in literature, of the Green's functions and their resolvents for linear anisotropic diffusion on the group $SE(2)$. Although the exact resolvent diffusion kernels (which take care of Tikhonov regularization on SE(2), \cite{DuitsR2006SS2}) are expressed in only 4 Mathieu functions, we also derive, in section \ref{ch:heisapprox}, the corresponding Heisenberg approximation resolvent diffusion kernels (which are rather Green's functions on the space of positions and velocities rather than Green's functions on the space of positions and orientations) which arise by replacing $\cos \theta$ by $1$ and $\sin \theta$ by $\theta$. Although these approximation Green's functions are not as simple as in the contour-completion case, \cite{DuitsR2006AMS}{ch:4.3}, they are more suitable if it comes to fast implementations, in particular for the Green's functions of the time processes. For comparison between the exact resolvent heat kernels and their approximations, see figure \ref{fig:comparisonspatialandFD}.

\subsection{The Exact Heat-Kernels on $SE(2)=\R^2 \rtimes SO(2)$. \label{ch:exactheatkernels}}

In this section we will derive the heat-kernels $K_{s}^{D}:SE(2) \to \R^{+}$ and the corresponding resolvent kernels $R_{\alpha,D}:SE(2) \to \R^{+}$ on $SE(2)$. Recall that $SE(2)$-convolution with these kernels, see (\ref{Gconv}), provide the solutions of the Forward Kolmogorov equations (\ref{ForwardKolmogorovCitti}) and recall that $R_{\alpha,D}=\alpha \, \int \limits_{0}^{\infty} K_{s}^{D}e^{-\alpha s}\, {\rm d}s$. During this chapter we set $D$ as a constant diagonal matrix. Although $D_{33}=0$ (as in (\ref{ForwardKolmogorovCitti})) has our main interest we also consider the more general case where $D_{33}\geq 0$.

The kernels $K_{s}^{D}$ and $R_{\alpha,D}$ are the unique solutions of the respectively the following  problems
{\small
\[
\left\{
\begin{array}{l}
\left(- D_{11} (\partial_{\theta})^2-D_{22}(\partial_{\xi})^2- D_{33}(\partial_{\eta})^2  + \alpha \right) R_{\alpha,D}^{\infty}= \alpha \delta_{e},  \\
R_{\alpha,D}^{\infty}(\cdot,\cdot,0)=R_{\alpha,D}^{\infty}(\cdot,\cdot,2\pi) \\
R_{\alpha,D}^{\infty} \in \mathbb{L}_{1}(SE(2)).
\end{array}
\right. \!
\textrm{ ,}
\left\{
\begin{array}{l}
\partial_{s} K_{s}^{D} = \left(D_{11} (\partial_{\theta})^2+D_{22}(\partial_{\xi})^2+ D_{33}(\partial_{\eta})^2\right)  K_{s}^{D} \\
\lim \limits_{s\downarrow 0} K_{s}^{D}=\delta_{e} \\
K_{s}^{D} \in \mathbb{L}_{1}(SE(2))
\end{array}
\right. \ .
\]
}
The first step here is to perform a Fourier transform with respect to the spatial part $\equiv \R^2$ of $SE(2)=\R^2 \rtimes \mathbb{T}$ , so that we obtain $\hat{R}_{\alpha,D}, \hat{K}_{s}^{D} \in \mathbb{L}_{2}(SE(2)) \cap C(SE(2))$ given by
\[
\begin{array}{l}
\hat{K}_{s}^{D}(\omega_{1},\omega_{2},\theta)=\mathcal{F}[K_{s}^{D}(\cdot,\cdot,\theta)](\omega_{1},\omega_{2}). \\
\hat{R}_{\alpha,D}(\omega_{1},\omega_{2},\theta)=\mathcal{F}[R_{\alpha,D}(\cdot,\cdot,\theta)](\omega_{1},\omega_{2}).
\end{array}
\]
Then $\hat{R}_{\alpha,D}$ and $\hat{K}_{s}^{D}$ satisfy
\begin{equation} \label{FDeq}
\begin{array}{ll}
(\alpha I-\mathcal{B}_{\www}) \hat{R}_{\alpha,D}= \frac{\alpha}{2\pi}\delta_{0} \  \  \ &
\textrm{ and }
\partial_{s} K_{s}^{D}=\mathcal{B}_{\www} K_{s}^{D}, \ \ \lim \limits_{s \downarrow 0} K_{s}^{D}(\www,\theta)=\delta_{e} \\
\end{array}
\end{equation}
where we define the operator
\[
\mathcal{B}_{\www}= -D_{22} \rho^2 \cos^{2}(\varphi-\theta) -D_{33} \rho^2 \sin^{2}(\varphi-\theta) + D_{11} (\partial_{\theta})^{2}
\]
where we expressed $\www \in \R^2$ in polar coordinates
\[
\www=(\rho \cos \varphi,\rho \sin \varphi) \in \R^2
\]
and where we note that $\mathcal{F}(\delta_e)=\frac{1}{2\pi} 1_{\R^2 } \otimes \delta_0^{\theta}$.
By means of the basic identities $\cos^{2}(\varphi-\theta)+\sin^{2}(\varphi-\theta)=1$ and $\cos(2(\varphi-\theta))=2 \cos^{2}(\varphi-\theta)-1$ we can rewrite operator $\mathcal{B}_{\www}$ in a (second order) Mathieu operator (corresponding to the well-known Mathieu equation $y''(z)+[(a-2q) \cos (2z)]y(z)=0$, \cite{Schaefke},\cite{Abra65})
\[
\mathcal{B}_{\www}= D_{11}\left((\partial_{\theta})^2 + aI - 2q \cos (2 (\varphi-\theta))\right),
\]
where $a=-\frac{\alpha + (\rho^2/2) (D_{22}+D_{33})}{D_{11}}$ and $q=\rho^2 \left( \frac{D_{22}-D_{33}}{4D_{11}} \right) \in \R$. Clearly, this unbounded operator
(with domain $\mathcal{D}(\mathcal{B}_{\www})=\mathbb{H}^{2}(\mathbb{T})$) is for each fixed $\www \in \R^2$ a symmetric operator of Sturm-Liouville type on $\mathbb{L}_{2}(\mathbb{T})$:
\[
\mathcal{B}_{\www}^{*}=
\mathcal{B}_{\www}.
\]
Its right inverse extends to a compact self-adjoint operator on $\mathbb{L}_{2}(\mathbb{T})$ and thereby $B_{\www}$ has the following complete orthogonal basis
of eigen functions
\[
\begin{array}{l}
\Theta_{n}^{\www}(\theta)=\textrm{me}_{n}(\varphi-\theta,q), \qquad n \in \mathbb{Z}, q=\rho^2 \left( \frac{D_{22}-D_{33}}{4D_{11}} \right) \in \R , \\
\mathcal{B}_{\www} \Theta_{n}^{\www}= \lambda_{n}^{\varrho} \, \Theta_{n}^{\www},
\end{array}
\]
whose eigen-values equal $\lambda_{n}^{\rho}= - a_{n}(q) D_{11}-\frac{\rho^{2}}{2}(D_{22}+D_{33})\leq -n^2 D_{11} \leq 0$, where $\textrm{me}_{n}(z,q)=\textrm{ce}_{n}(z,q)+i \, \textrm{se}_{n}(z,q) $ denotes the well-known Mathieu function (with discrete Floquet exponent $\nu=n$), \cite{Schaefke},\cite{Abra65}, and characteristic values $a_n(q)$ which are countable solutions of the corresponding characteristic equations \cite{Schaefke},\cite{Abra65}{p.723}, containing continued fractions. Note that at $\www=\ul{0}$, i.e. $\rho=0$, we have $\textrm{me}_{n}(z,0)=e^{inz}$, $\lambda_{n}^{0}=n^2$.

The functions $q \mapsto a_{n}(q)$ are analytic on the real line. Here we note that in contrast with the eigen function decomposition of the generator of the Forward Kolmogorov equation (\ref{ForwardKolmogorovMumford}) of Mumford's direction process \cite{DuitsR2006AMS}
Green's functions of the contour completion case \cite{DuitsR2006AMS}, we have $q \in \R$ rather than $q \in i \mathbb{R}$ and therefor we will not meet any nasty branching points of $a_n$.
The Taylor expansion of $a_{n}(q)$ for $n \neq 1,2,3$ (for the cases $n=1,2,3$ see \cite{Abra65}{p.730}) is given by
\[
a_n(q)= n^2 + \frac{1}{2(n^2-1)} q^2 + \frac{5 n^2 +7}{32(n^2-1)^3(n^2-4)}q^4 +
\frac{9 n^2 + 58 n^2 +29}{64(n^2-1)^5(n^2-4)(n^2-9)q^6} +O(q^8).
\]
For each fixed $\www \in \R^2$ the set $\{\Theta_{n}^{\www}\}_{n \in \mathbb{Z}}$ is a complete orthogonal basis for $\mathbb{L}_{2}(\mathbb{T})$ and moreover we have
\[
\langle \delta_{0},\phi \rangle= \phi(0)= \sum \limits_{n=-\infty}^{\infty} (\Theta_{n}^{\www},\phi)\Theta_{n}^{\www}(0)
\]
for all test functions $\phi \in \mathcal{D}(\mathbb{T})$. Consequently, the unique solutions of (\ref{FDeq}) are given by
\begin{equation} \label{soldiff}
\begin{array}{l}
\hat{K}_{s}(\www,\theta)= \sum \limits_{n=\infty}^{\infty} \Theta_{n}^{\www}(\theta) \Theta_{n}^{\www}(0) e^{\lambda_{n}^{\varrho}s}, \\
\hat{R}_{\alpha,D}(\www,\theta)= \alpha \sum \limits_{n=-\infty}^{\infty} \frac{\Theta_{n}^{\www}(\theta) \Theta_{n}^{\www}(0)}{\alpha-\lambda_{n}^{\rho}}  \\
\end{array}
\end{equation}
Or more explicitly formulated:
\begin{theorem} \label{th:exaxtdiff}
Let $D_{11}, D_{22}, D_{33}>0$, then
the heat kernels $\mathcal{K}_{t}^{D_{11}, D_{22}, D_{33}}$ on the Euclidean motion group which satisfy
\begin{equation} \label{timeprob2}
\left\{
\begin{array}{l}
\partial_{t }\mathcal{K} =\left(D_{11} (\partial_{\theta})^2 +D_{22} (\partial_{\xi})^2+D_{33}(\partial_{\eta})^2\right) \mathcal{K}  \\
\mathcal{K}(\cdot,\cdot,0,t)=\mathcal{K}(\cdot,\cdot,2\pi,t) \textrm{ for all }t>0. \\
\mathcal{K} (\cdot,\cdot,\cdot,0)=\delta_{e} \\
\mathcal{K}(\cdot,t) \in \mathbb{L}_{1}(G), \textrm{ for all }t>0.
\end{array}
\right.
\end{equation}
are given by
\[
\mathcal{K}_{t}^{D_{11},D_{22},D_{33}}(x,y,e^{i\theta}):=\mathcal{K}(x,y,e^{i\theta},t)=
\mathcal{F}^{-1}[\www \mapsto \hat{\mathcal{K}}_{t}^{D_{11},D_{22},D_{33}}(\www,e^{i\theta})](b_{1},b_{2})
\]
where
\[
\hat{\mathcal{K}}_{t}^{D_{11},D_{22},D_{33}}(\www,e^{i\theta})= e^{-t(1/2)(D_{22}+D_{33})\rho^2}
\left(
\sum \limits_{n=-\infty}^{\infty} \frac{\textrm{\textrm{me}}_{n}(\varphi,q)\textrm{\textrm{me}}_{n}(\varphi-\theta,q)}{2\pi}
e^{-t a_{n}(q)D_{11}}
\right)
\]
with $q=\frac{\rho^2(D_{22}-D_{33})}{4\, D_{11}}$ and $a_{n}(q)$ the Mathieu Characteristic (with Floquet exponent $n$) and with the property that $\mathcal{K}^{D_{11},D_{22},D_{33}}_{t}>0$ and
\[
\|\mathcal{K}^{D_{11},D_{22},D_{33}}_{t}\|_{\mathbb{L}_{1}(SE(2))}= \int \limits_{0}^{2\pi} \hat{\mathcal{K}}_{t}^{ D_{11},D_{22},D_{33}}(\ul{0},e^{i\theta}) \, {\rm d}\theta=
\sum \limits_{n=-\infty}^{\infty} (2\pi)^{-1}\int \limits_{0}^{2\pi} e^{in\theta} {\rm d}\theta e^{-t \, n^2 D_{11}}=1.
\]
\end{theorem}
Consider the case where $D_{11} \downarrow 0$, then
$a_{n}(q) \sim -2q $ as $q\to \infty$ and we have
\[
\begin{array}{ll}
\lim \limits_{D_{11}\downarrow 0} \hat{\mathcal{K}}_{t}^{D_{11},D_{22},D_{33}}(\www,e^{i\theta}) &=
e^{-\frac{t}{2}(D_{22}+D_{33})(\omega_x^2+\omega_{y}^2)} e^{-\frac{t}{2}(D_{22}-D_{33})(\omega_x^2-\omega_{y}^2)} \delta_{0}^{\theta} \\
 &=
e^{-t(D_{22}\omega_{x}^2 +D_{33}\omega_{y}^2)} \delta^{\theta}_{0}= \hat{\mathcal{K}}_{t}^{0,D_{22},D_{33}}(\www,e^{i\theta})  \delta^{\theta}_{0} \ .
\end{array}
\]
Finally we notice that the case $D_{11}=0$ yields the following operation on $\mathbb{L}_{2}(G)$:
\[
\! \!(\hat{\mathcal{K}}_{t}^{0,D_{22},D_{33}}*_{SE(2)}U)(g)=  \int \limits_{\R^{2}} G_{t}^{D_{22},D_{33}}(R_{\theta}^{-1}(\ul{x}-\ul{x}'))U(\ul{x}',e^{i\theta}){\rm d\ul{x}'}= \!
(\mathcal{R}_{e^{i\theta}}G_{t}^{D_{22},D_{33}}*_{\R^2}f)(\ul{x})
\]
$g=(\ul{x},e^{i\theta}) \in SE(2)$, where $G_{t}^{D_{22},D_{33}}(x,y)=G_{t\, D_{22}}^{d=1}(x)\; G_{t\, D_{33}}^{d=1}(y)$ equals the well-known anisotropic Gaussian kernel or heat-kernel on $\R^n$, and where $\mathcal{R}_{e^{i\theta}}\phi(\ul{x})=\phi(R_{\theta}^{-1}\ul{x})$ is the left regular action of $SO(2)$ in $\mathbb{L}_{2}(\R^2)$, which corresponds to anisotropic diffusion in each fixed orientation layer $U(\cdot,\cdot,\theta)$ where the axes of anisotropy coincide with the $\xi$ and $\eta$-axis. This operation is for example used in image analysis in the framework of channel smoothing \cite{FelsbergBspline}, \cite{Duits2005IJCV}. We stress that also the diffusion kernels with $D_{11}>0$ are interesting for computer vision applications such as the frameworks of tensor voting, channel representations and invertible orientation scores as they allow different orientation layers $\{U(\cdot,\cdot,\theta)\}_{\theta \in [0,2\pi)}$ to interfere. See Figure \ref{fig:GFs} (with inclusion of curvature as we will explain in subsection \ref{ch:curvest}) dependent heat-kernel $\mathcal{K}_{t}^{D_{11},D_{22},D_{33}}$, $D_{22}\gg D_{33} >0$ on $SE(2)$. For illustration of the
corresponding resolvent kernel $R_{\alpha,D}$ (with comparison to approximations we shall derive in section \ref{ch:heisapprox}) see Figure \ref{fig:comparisonspatialandFD}.

Next we shall derive a more suitable expression than (\ref{soldiff}) for the resolvent kernel $R_{\alpha,D}$. To this end we will unwrap the torus to $\R$ and replace the periodic boundary condition
in $\theta$ by an absorbing boundary condition at infinity. Afterwards we shall construct the true periodic solution by explicitly computing (using the Floquet theorem) the series consisting of (rapidly decreasing) $2\pi$-shifts of the solution with absorbing condition at infinity.

In our explicit formulae for the resolvent kernel $R_{\alpha,D}$ we shall make use of the non-periodic complex-valued Mathieu function which is a solution of the Mathieu equation
\begin{equation} \label{mathieu}
y''(z) +[(a-2q)\cos (2z)]y(z)=0, \qquad a,q \in \R
\end{equation}
and which is by definition\footnote{There exist several definitions of Mathieu solutions, for an overview see \cite{Abra65}{p.744, Table 20.10} each with different normalizations. In this article we always follow the consistent conventions by Meixner and Schaefke \cite{Schaefke}. However, for example \emph{Mathematica 5.2} chooses an unspecified convention. This requires slight modification of (\ref{menu}), see \cite{MarkusThesis}}, \cite{Schaefke}{p.115}, \cite{Abra65}{p.732}, given by
\begin{equation}\label{menu}
\begin{array}{l}
\textrm{me}_{\pm \nu}(z, q)=\textrm{ce}_{\nu}(z,q) \pm i \textrm{se}_{\nu}(z,q). \\ 
\end{array}
\end{equation}
Here $\nu=\nu(a,q)$ equals the Floquet exponent (due to the Floquet Theorem \cite{Schaefke}{ p.101}) of the solution, which means that
\begin{equation} \label{Floquet}
\textrm{me}_{\pm \nu}(z+\pi, q)=e^{i \nu z}\textrm{me}_{\pm \nu}(z, q),
\end{equation}
for all $z,q \in \mathbb{R}$.
\begin{theorem}\label{heatkernelsonSE2}
Let $\alpha>0$, $D_{22}\geq D_{33}>0$, $D_{11}>0$.
The solution $R_{\alpha,D}^{\infty}: \R^3 \setminus \{0,0,0\} \to \R $ of the problem
\[ 
\left\{
\begin{array}{l}
\left(- D_{11} (\partial_{\theta})^2-D_{22}(\partial_{\xi})^2- D_{33}(\partial_{\eta})^2  + \alpha \right) R_{\alpha,D}^{\infty}= \alpha \delta_{e},  \\
R_{\alpha,D}^{\infty}(\cdot,\cdot,\theta) \to 0 \textrm{ uniformly on compacta as }|\theta|\to \infty \\
R_{\alpha,D}^{\infty} \in \mathbb{L}_{1}(\R^3),
\end{array}
\right.
\] 
is given by
\[
R_{\alpha,D}^{\infty}(x,y,\theta)=
\mathcal{F}^{-1}[(\omega_x,\omega_y) \mapsto \hat{R}_{\alpha,D}^{\infty}(\omega_x,\omega_y,\theta)](x,y).
\]
In case $D_{33}< D_{22}$ we have
\begin{equation}\label{heatKinf}
\begin{array}{l}
\hat{R}_{\alpha,D}^{\infty}(\omega_x,\omega_y,\theta)
=
 \frac{-\alpha}{4\pi D_{11} W_{a,q}} \\
\left[
\textrm{me}_{\nu}\left(\varphi, \frac{(D_{22}-D_{33})\rho^2}{4\, D_{11}}\right) \textrm{me}_{-\nu}\left(\varphi-\theta,  \frac{(D_{22}-D_{33})\rho^2}{4\, D_{11}}\right)\,
\rm{u}(\theta) \right. \\
 \qquad \qquad + \left.
\textrm{me}_{-\nu}\left(\varphi, \frac{(D_{22}-D_{33})\rho^2}{4\, D_{11}} \right) \textrm{me}_{\nu}\left(\varphi-\theta, \frac{(D_{22}-D_{33})\rho^2}{4\, D_{11}}\right)\,
\rm{u}(-\theta)
\right].
\end{array}
\end{equation}
with $\www=(\rho \cos \phi, \rho \sin \phi)$, where $\theta \mapsto {\rm u}(\theta)$ denotes the unit step function, which is given by ${\rm u}(\theta)=1$  if $\theta>0$, ${\rm u}(\theta)=0$ if $\theta<0$ and where the Floquet exponent equals $\nu\left(
\frac{-(\alpha+(1/2)(D_{22}+D_{33})\rho^2)}{D_{11}}, \frac{(D_{22}-D_{33})\rho^2}{4\, D_{11}}
\right)$ and where $W_{a,q}= \textrm{ce}_{\nu}(0,q)\textrm{se}_{\nu}'(0,q)$ equals the Wronskian of $\textrm{ce}(\cdot,q)$ and $\textrm{se}(\cdot,q)$ with $a=\frac{-(\alpha+(1/2)(D_{22}+D_{33})\rho^2)}{D_{11}}$ and $q=\frac{(D_{22}-D_{33})\rho^2}{4\, D_{11}}$.

In case $D_{22}=D_{33}$ (which follows by taking the limit $D_{22}\to D_{33}$ in (\ref{heatKinf})) we have
\[
\hat{R}_{\alpha, D}^{\infty}(\www,\theta)= \frac{\alpha \, e^{- \sqrt{\frac{\alpha+D_{22}\rho^2}{D_{11}}}|\theta|}}{4\pi \sqrt{D_{11}}\sqrt{D_{22}\rho^2+\alpha}}, \qquad \rho=\|\www\|, D_{22}=D_{33},
\]
which yields for $D_{33}=D_{22}$:
\begin{equation} \label{D22isD33}
\begin{array}{l}
K_{s}^{D; \infty}(\ul{x},\theta)=\frac{1}{\sqrt{D_{11}}D_{22}}\frac{1}{(4\pi s)^{\frac{3}{2}}} e^{-\frac{\frac{\theta^{2}}{D_{11}}+\frac{r^{2}}{D_{22}}}{4s}}, \qquad r=\|\ul{x}\|, D_{22}=D_{33}, \\
R_{\alpha, D}^{\infty}(\ul{x},\theta)=\frac{\alpha}{4\pi} \frac{1}{\sqrt{D_{11}}D_{22}} \frac{e^{-\sqrt{\alpha}\sqrt{\frac{\theta^{2}}{D_{11}}+\frac{r^{2}}{D_{22}}} }}{\sqrt{\frac{\theta^{2}}{D_{11}}+\frac{r^{2}}{D_{22}}}}.
\end{array}
\end{equation}
\end{theorem}
\textbf{Proof }
Again we apply Fourier transform with respect to $\R^2$ only, this yields
\begin{equation} \label{eqFDkern}
\begin{array}{ll}
(D_{22}(\partial_{\xi})^2+D_{33}(\partial_{\eta})^2+D_{11}(\partial_{\theta})^2-\alpha I) R_{\alpha,D}^{\infty} = -\alpha \delta_{e} & \desda \\
(-D_{22}\rho^2 \cos^{2}(\varphi-\theta)-D_{33}\rho^2 \sin^{2}(\varphi-\theta)+D_{11}(\partial_{\theta})^2-\alpha I) \hat{R}_{\alpha,D}^{\infty} = -\frac{\alpha}{2\pi} \delta^{\theta}_{0} &\desda \\
(-D_{33}\rho^2 + (D_{33}-D_{22})\rho^2\cos^{2}(\varphi-\theta)+D_{11}(\partial_{\theta})^2-\alpha I) \hat{R}_{\alpha,D}^{\infty} = -\frac{\alpha}{2\pi} \delta^{\theta}_{0} & \desda \\
((\partial_{\theta})^2+ a I -2q \cos(2(\phi-\theta))) \hat{R}_{\alpha,D}^{\infty} = -\frac{\alpha}{2\pi D_{11}} \delta^{\theta}_{0} & \
\end{array}
\end{equation}
where $a=- \left( \frac{\alpha +(\rho^2/2)(D_{22}+D_{33})}{D_{11}}\right) $ and $q=\rho^2\left(\frac{D_{22}-D_{33}}{4\, D_{11}}\right)$.

We shall first deal with the cases $D_{33}<D_{22}$ and return to the case $D_{22}=D_{33}$ later.
In order to solve the last equation of (\ref{eqFDkern}), we first find the solutions $F,G$ of the equations
\[
\begin{array}{ll}
\left\{
\begin{array}{l}
((\partial_{\theta})^2+ a I -2q \cos(2(\phi-\theta))) F(\theta) = 0  \\
F(\theta) \to 0 \textrm{ as }\theta \to +\infty
\end{array}
\right.
&
\textrm{ and }
\left\{
\begin{array}{l}
((\partial_{\theta})^2+ a I -2q \cos(2(\phi-\theta))) G(\theta) = 0  \\
G(\theta) \to 0 \textrm{ as }\theta \to -\infty
\end{array}
\right.
\end{array}
\]
and then we make a continuous (but not differentiable) fit of these solutions.
Now for $a<0<q$ and $a<-2q$ we have\footnote{Floquet exponents always exponents come in conjugate pairs, therefor throughout this paper we set the imaginary part of the Floquet-exponent to a positive value.} $\textrm{Im}(\nu(a,q))>0$. We indeed have $a<0$ and $q>0$ since $\alpha>0, D_{22}-D_{33} >0$ and moreover we have
{\small
\[
a=- \left( \frac{\alpha +(\rho^2/2)(D_{22}+D_{33})}{D_{11}}\right)< -\left( \frac{(\rho^2/2)(D_{22}-D_{33})}{D_{11}}\right)=-2q.
\]
}
So consequently (recall (\ref{Floquet})) we find
\[
F(\theta)= C_{1} \textrm{me}_{-\nu}(\varphi-\theta,  q) \textrm{ for }\theta>0 \textrm{ and }G(\theta)= C_{2}\textrm{me}_{\nu}(\varphi-\theta,  q) \textrm{ for }\theta<0,
\]
now in order to make a continuous fit at $\theta=0$ we set $C_{1}=\lambda \, \textrm{me}_{\nu}(\varphi,  q)$ and $C_{2}=\lambda \, \textrm{me}_{-\nu}(\varphi,q)$ for a constant $\lambda \in \R$ yet to be determined.
\begin{equation} \label{lambdcalc}
\begin{array}{ll}
\frac{-\alpha}{2\pi D_{11}}\delta_{0} & =(\partial_{\theta}^{2} -2q \cos(2(\varphi-\theta))+a\, I )\hat{R}_{\alpha,D_{11}}^{\infty}(\www,\cdot) \\
           & = \lambda(\partial_{\theta}^{2} -2q \cos(2(\varphi-\theta))+a\, I )(\textrm{me}_{\nu}(\varphi,q)\textrm{me}_{-\nu}(\varphi-\theta,q)\rm{u}(\theta)+ \textrm{me}_{-\nu}(\varphi,q)\textrm{me}_{\nu}(\varphi-\theta,q)\rm{u}(-\theta) ) \\[8pt]
&= -\lambda( \textrm{me}_{-\nu}(\varphi,q)\textrm{me}_{\nu}'(\varphi,q)-\textrm{me}_{\nu}(\varphi,q)\textrm{me}_{-\nu}'(\varphi,q)) \delta_{0} + 0+0+ \lambda \hat{R}^{\infty}_{\alpha,D_{11}}(\delta_0'-\delta_{0}') \desda \\[8pt]
 & \lambda = -\frac{\alpha}{2\pi D_{11}} (W[\textrm{me}_{\nu}(\cdot,q),\textrm{me}_{-\nu}(\cdot,q)])^{-1} \desda \\
 & \lambda=  -\frac{\alpha}{4\pi D_{11}} (W[\textrm{ce}_{\nu}(\cdot,q),\textrm{se}_{\nu}(\cdot,q)])^{-1} \desda\\
 & \lambda= -\frac{\alpha}{4\pi D_{11}} (\textrm{ce}_{\nu}(0,q)\textrm{se}_{\nu}'(0,q)-0)= -\frac{\alpha}{4 \pi D_{11} W_{a,q}}
\end{array}
\end{equation}
where the Wronskian is given by $W[f,g](z)=f(z)g'(z)-g(z)f'(z)$, which is for solutions of the Mathieu-equation independent of $z$ so substitute $z=0$.

Now that we have explicitly derived the solution for the case $D_{22}>D_{33}$. We can take the limit $D_{33} \uparrow D_{22}$ and consequently $q \downarrow 0$. It directly follows from the Mathieu equation (\ref{mathieu}) that $\lim \limits_{q \downarrow 0} \textrm{me}_{\nu}(a,q)(\theta)=\textrm{me}_{\nu}(a,0)=e^{i \sqrt{a} \theta}$ for all $\theta \in [0,2\pi)$. Thereby we have
{\small
\[
\lim \limits_{D_{33} \uparrow D_{22}} \hat{R}^{\infty}_{\alpha,D_{11}}(\www,\theta)= \frac{-\alpha}{ 4\pi D_{11} i \sqrt{a}} e^{i \sqrt{a} |\theta|}=\frac{\alpha \, e^{- \sqrt{\frac{\alpha+D_{22}\rho^2}{D_{11}}}|\theta|}}{4\pi \sqrt{D_{11}}\sqrt{D_{22}\rho^2+\alpha}}
.
\]
}
Now the results (\ref{D22isD33}) follow by direct computation. $\hfill \Box$ \\
\\
We note that if $D_{22}=D_{33}$ the diffusion in the spatial part is isotropic and $\Delta=\partial_{\xi}^{2}+\partial_{\eta}^{2}=\partial_{x}^{2}+\partial_{y}^{2}$ commutes $\partial_{\theta}^{2}$ with so \emph{in case $D_{22}=D_{33}$
left-invariant diffusion on $\R^2 \times \mathbb{T}$ (with direct product) left-invariant diffusion on $\R^{2} \rtimes \mathbb{T}$ (with semi-direct product) and the kernels (\ref{D22isD33}) indeed coincide with the Green's-functions for anisotropic diffusion on $\R^3$}. We have employed this fact in \cite{Fran06b} in order to generalize fast Gaussian derivatives on images
\begin{equation} \label{GD}
f \in \mathbb{L}_{2}(\R^2) \mapsto \frac{d^{n+m}}{d x^{m} dy^{n}}(G_{\ul{s}}*_{\R^2}f)(\ul{x})=(G_{\ul{s}}^{(m,n)}*_{\R^2}f)(\ul{x}) \in \R, \qquad \ul{x} \in \R^2, \ul{s}=(s_1,s_{2}) \in \R^{+} \times \R^{+},
\end{equation}
with separable Gaussian kernels $G_{\ul{s}}^{d=2}(x,y)=G_{s_{1}}^{d=1}(x)G_{s_{2}}^{d=1}(y)$ (a property which is very useful to reduce the computation time) 
to fast Gaussian derivatives on orientation scores. To this end we note that (\ref{GD}) can at least formally be written as
\[
\frac{d^{m+n}}{dx^m dy^n} e^{s\Delta} f = \frac{d}{dx^m} e^{s \partial_{x}^2} \left( \frac{d}{dy^n}  e^{s \partial_{y}^2} f \right).
\]
Now since $[\Delta,\partial_{\theta}]=0$ we can perform a similar trick for left-invariant Gaussian derivatives on orientation scores:
\begin{equation} \label{GDOS1}
\frac{d^{m+n+l}}{{\rm d}\xi^{m}{\rm d}\eta^{n}{\rm d}\theta^l} e^{s(D_{11} \partial_{\theta}^{2}+D_{22}\Delta)}U=
\frac{d^{m+n}}{{\rm d}\xi^{m}{\rm d}\eta^{n}}e^{s D_{22} \Delta} \left( \frac{d^{l}}{{\rm d}\theta^l} e^{s D_{11} \partial_{\theta}^{2}} U\right),
\end{equation}
which can again be used to reduce computation time:
\begin{equation}\label{GDOS2}
\begin{array}{l}
\frac{d^{m+n}}{d\xi^m d\eta^n} (K_{s}^{D_{33}=D_{22}}*_{SE(2)} U)(\ul{x},e^{i\theta})= \int \limits_{\R^2} (\cos \theta \partial_{x}+ \sin \theta \partial_{y})^{m}
(-\sin \theta \partial_{x}+\cos \partial_{y})^{n} G_{s D_{22}}^{d=2}(\ul{x}-\ul{x}') \times \\
\left(\int_{-\pi}^{\pi} \frac{d^{l}}{{\rm d}\theta^l} G_{s,D_{11}}^{d=1}(\theta-\theta') U(\ul{x}',\theta') {\rm d}\theta'\right) {\rm d}\ul{x}',
\end{array}
\end{equation}
where we stress that the order of the derivatives matters.

Finally, we stress that we can expand the exact Green's function $R_{\alpha, D_{11}}$ as an infinite sum over $2\pi$-shifts of the solution $R^{\infty}_{\alpha, D_{11}}$ for the unbounded case:
\begin{equation} \label{homotopy}
R_{\alpha, D}(x,y,e^{i\theta})= \lim \limits_{N \to \infty} \sum \limits_{k=-N}^{N}
R^{\infty}_{\alpha, D}(x,y,\theta-2k \pi).
\end{equation}
Note that this splits the probability-density of finding a random walker (whose traveling time is negatively exponentially distributed $s \sim NE(\alpha)$) in $SE(2)$ at position (regardless its traveling time) $(x,y,e^{i\theta})$ given its starting position and orientation $e=(0,0,e^{0i})$ into the probability density of finding a random walker in $SE(2)$ at position $(x,y,e^{i\theta})$ given it started at $e=(0,0,e^{0i})$ \emph{ and }given the fact that the homotopy number of its path equals $k$, for $k \in \mathbb{Z}$.

The nice thing is that the sum in (\ref{homotopy}) (which decays rather rapidly) can be computed explicitly by means of the Floquet theorem, i.e. (\ref{Floquet}), and the geometrical series $\sum \limits_{n=0}^{\infty} r^{k}= \frac{1}{1-r}$ for $r=e^{i\nu}$ with $r=|e^{i\nu}|<1$ since the imaginary part of $\nu=\nu(a,q)$ is positive. By straightforward computations this yields the following result.
\begin{theorem}\label{th:thesolutionresolvent}
Let $\alpha, D_{11},D_{22}>0$ and $D_{33}\geq 0$. Then
the solution $R_{\alpha,D}: SE(2) \to \R$ of the problem
\[ 
\left\{
\begin{array}{l}
\left(- D_{11} (\partial_{\theta})^2-D_{22}(\partial_{\xi})^2- D_{33}(\partial_{\eta})^2  + \alpha \right) R_{\alpha,D}= \alpha \delta_{e},  \\
R_{\alpha,D}(\cdot,\cdot,\theta+ 2k \pi)=R_{\alpha,D}(\cdot,\cdot,\theta)  \textrm{ for all }k \in \mathbb{Z}, \\
R_{\alpha,D} \in \mathbb{L}_{1}(SE(2)),
\end{array}
\right.
\] 
is given by
\[
R_{\alpha,D}(\ul{x},\theta)= \sum \limits_{k \in \mathbb{Z}} R_{\alpha,D}^{\infty}(\ul{x},\theta +2 k \pi)
\]
the righthand side of which can be calculated using Floquet's theorem and (\ref{heatKinf}) yielding for $D_{33}<D_{22}$:
 {\small
\begin{equation} \label{final}
\begin{array}{l}
[\mathcal{F}R_{\alpha,D}(\cdot,\theta)](\www)= \frac{\alpha}{4 \pi D_{11}\textrm{ce}_{\nu}(0,q)\, \textrm{se}_{\nu}'(0,q) } \left\{ \right. \\
\left.\left(-\cot(\nu \pi)\left(\textrm{ce}_{\nu}(\varphi,q)\, \textrm{se}_{\nu}(\varphi-\theta,q)+\textrm{se}_{\nu}(\varphi,q)\, \textrm{se}_{\nu}(\varphi-\theta,q)\right)+ \right. \right.\\
\left. \left. \qquad
\textrm{ce}_{\nu}(\varphi,q)\, \textrm{se}_{\nu}(\varphi-\theta,q)-
\textrm{se}_{\nu}(\varphi,q)\, \textrm{ce}_{\nu}(\varphi-\theta,q)\right){\rm u}(\theta)\qquad+ \qquad \right. \\
\left.\left(-\cot(\nu \pi)\left(\textrm{ce}_{\nu}(\varphi,q)\, \textrm{ce}_{\nu}(\varphi-\theta,q)-\textrm{se}_{\nu}(\varphi,q)\,
\textrm{se}_{\nu}(\varphi-\theta,q)\right) + \right. \right.\\
\left. \qquad
\textrm{ce}_{\nu}(\varphi,q)\, \textrm{se}_{\nu}(\varphi-\theta,q)+
\textrm{se}_{\nu}(\varphi,q)\, \textrm{ce}_{\nu}(\varphi-\theta,q)\right){\rm u}(-\theta)
\ \ \}
\end{array}
\end{equation}
}
with $q=\frac{(D_{22}-D_{33})\rho^2}{4 D_{11}}$, $\www=(\rho \cos \varphi, \rho \sin \varphi)$ and Floquet exponent $\nu=\nu(a,q)$, $a=-\frac{\alpha+(1/2)(D_{22}-D_{33})\rho^2}{D_{11}}$ and where $\theta \mapsto {\rm u}(\theta)$ denotes the unit step function, which is given by ${\rm u}(\theta)=1$  if $\theta>0$, ${\rm u}(\theta)=0$ if $\theta<0$.
\end{theorem}
The results in the preceding theory on the resolvent Green's function of the contour enhancement process can be set in a variational formulation, like the variational formulation in \cite{Citti} (where $D_{33}=0$).
\begin{corollary}
Let $U \in \mathbb{L}_{2}(SE(2))$ and $\alpha, D_{11}, D_{22}>0$, $D_{33}\geq 0$. Then the unique solution of the variational problem
{\small
\begin{equation}\label{varprobSE2}
\arg \mbox{}\hspace{-0.6cm}\mbox{} \min \limits_{W \in \mathbb{H}^{1}(SE(2))} \int \limits_{SE(2)}
\frac{\alpha}{2}(W(g)-U(g))^{2} + D_{11}(\partial_{\theta} W(g))^{2} + D_{22}(\partial_{\xi} W(g))^{2} +D_{33}(\partial_{\xi} W(g))^{2} {\rm d}\mu_{SE(2)}(g)
\end{equation}
}
is given by
\[
W(g)=(R_{\alpha,D}*_{SE(2)}U)(g)= \int \limits_{SE(2)} R_{\alpha,D}(h^{-1}g) U(h) \; {\rm d}\mu_{SE(2)}(h)
\]
where the Green's function $R_{\alpha,D}:SE(2) \to \R^{+}$ is explicitly given in Theorem \ref{th:thesolutionresolvent}.
\end{corollary}
\textbf{Proof }By convexity of the energy
\[
\mathcal{E}(W):=\int \limits_{SE(2)}
\frac{\alpha}{2}(W(g)-U(g))^{2} + D_{11}(\partial_{\theta} W(g))^{2} + D_{22}(\partial_{\xi} W(g))^{2} +D_{33}(\partial_{\eta} W(g))^{2} {\rm d}\mu_{SE(2)}(g)
\]
the solution of the variational problem (\ref{varprobSE2}) is unique. Along the minimizer we have
\[
\lim \limits_{h \downarrow 0} \frac{\mathcal{E}(W+h \delta)-\mathcal{E}(W)}{h}=0
\]
for all pertubations $\delta \in \mathbb{H}^{1}(SE(2))$. So by integration by parts we find
\[
\left(\alpha (W-U) - D_{11}\partial_{\theta}^{2}W-D_{22}\partial_{\xi}^{2}W-D_{33}\partial_{\eta}^{2}W, \delta\right)_{\mathbb{L}_{2}(SE(2))}=0
\]
for all $\delta \in \mathbb{H}^{1}(SE(2))$. Now $\mathbb{H}^{1}(SE(2))$ is dense in
$\mathbb{L}_{2}(SE(2))$ and therefore
\[
\alpha \, U=\left(\alpha I - ( D_{11}\partial_{\theta}^{2}+D_{22}\partial_{\xi}^{2}+D_{33}\partial_{\eta}^{2})\right)W
\]
so $W=\alpha \left(\alpha I - ( D_{11}\partial_{\theta}^{2}+D_{22}\partial_{\xi}^{2}+D_{33}\partial_{\eta}^{2})\right)^{-1}U$
and by left-invariance and linearity this resolvent equation is solved by a $SE(2)$-convolution with the smooth Green's function $R_{\alpha,D}:SE(2) \backslash \{e\} \to \R^{+}$ from Theorem \ref{th:thesolutionresolvent}.$\hfill \Box$ \\
 \\
\textbf{Remark:} We looked for a variational formulation of the contour \emph{completion} process as well, but in vain.

\subsection{The Heisenberg Approximations of the heat-kernels on $SE(2)$ \label{ch:heisapprox}}

If we approximate $\cos \theta \approx 1$ and $\sin \theta \approx \theta$ the left-invariant vector fields are approximated by
\begin{equation} \label{approxVF}
\hat{A}_{1}=\partial_{\theta}, \hat{A}_{2}=\partial x + \theta \partial_{y}, \hat{A}_{3}= -\theta \partial_{x} + \partial_{y}
\end{equation}
which are left-invariant vector fields in a 5 dimensional Nilpotent Lie-algebra of Heisenberg type. In our previous related work, we used this
replacement to explicitly derive more tangible Green's functions which are (surprisingly) good approximations\footnote{In fact in the field of image analysis the approximative Green's functions is often mistaken for the exact Green's functions.} of the exact Green's functions of the direction process with the goal of \emph{contour completion} (i.e. $a_{2}, D_{11} \neq 0$, other parameters are set to zero) for reasonable parameter settings, see \cite{DuitsR2006AMS}. In fact this replacement will provide Green's functions on the group of positions and velocities rather than Green's functions on the group of positions and orientations, see \cite{Thornber2}{ App. C}.

\begin{figure}
\centerline{
\includegraphics[width=0.55\hsize]{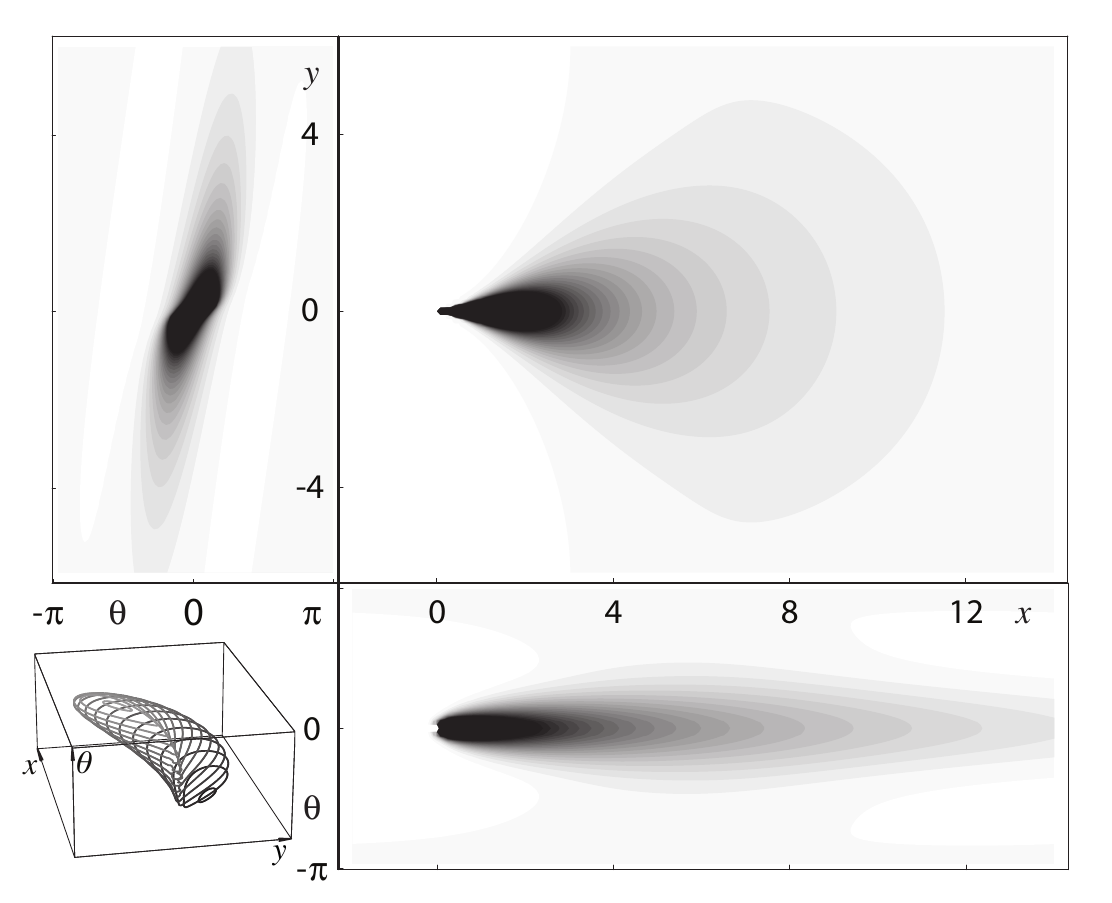}
\mbox{}\hspace{-0.6cm} \mbox{}
\includegraphics[width=0.55\hsize]{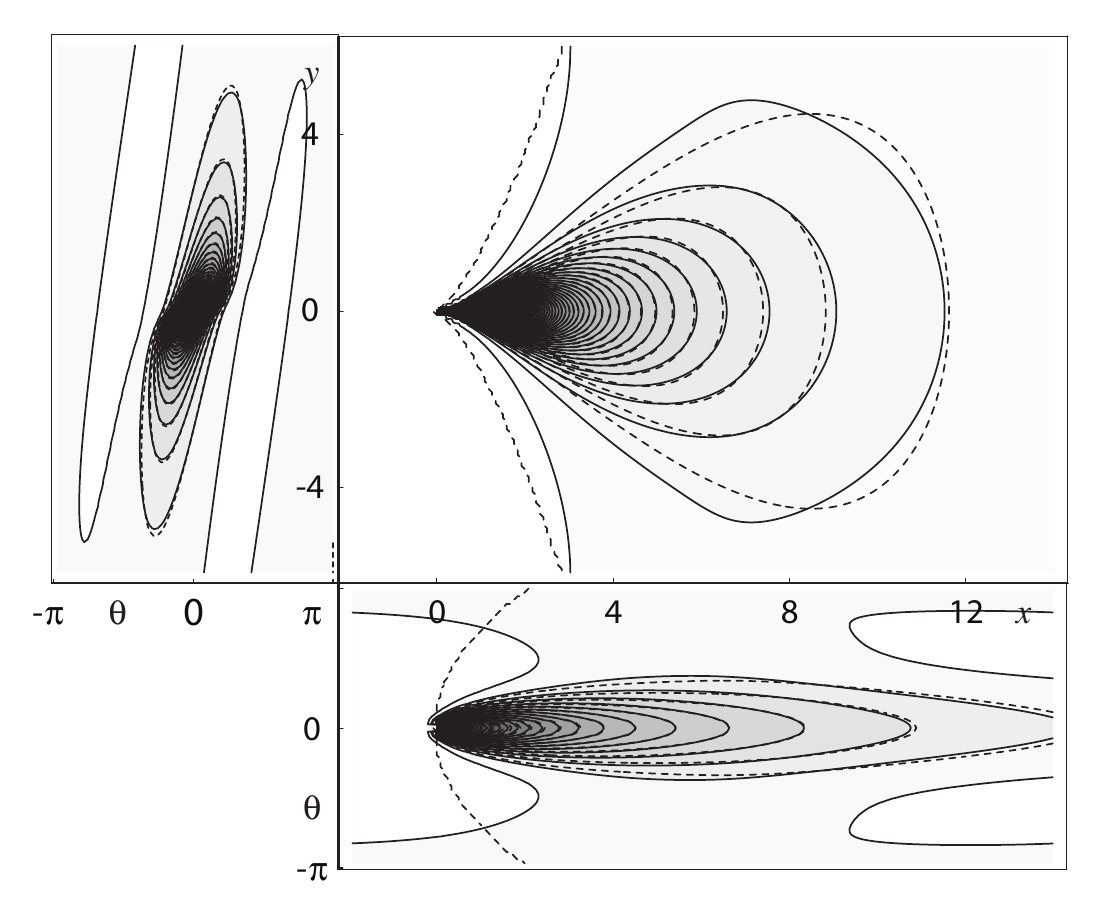}
\vspace{-0.1cm}\mbox{}
}
\centerline{
\includegraphics[width=0.3\hsize]{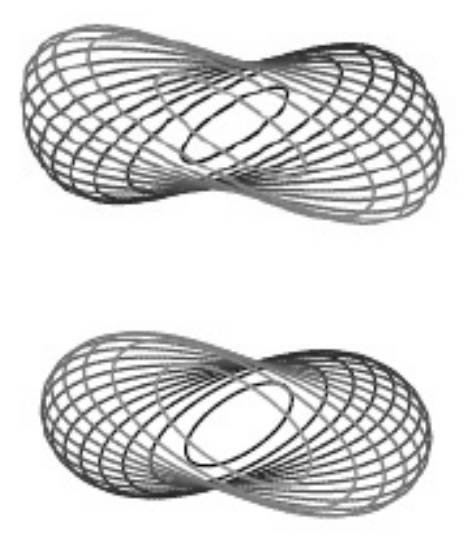}
\includegraphics[width=0.5\hsize]{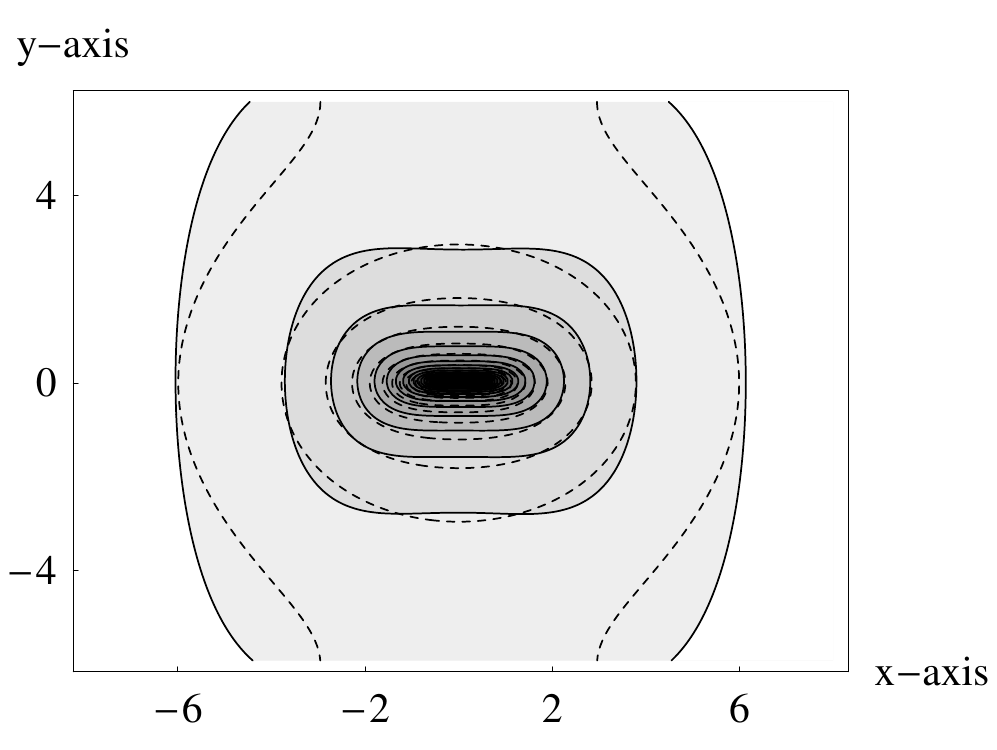}
}
\caption{
{\small Top Row:
Left: isocontour-plots of the marginals of the exact resolvent Green's function of the direction process $\overline{R}_{\alpha}^{D_{11}, a_{2}=1}$ given in \cite{DuitsR2006AMS}, $\kappa_0=0$, $\alpha=\frac{1}{10}$, $D_{11}=\frac{1}{32}$. 
Right: a comparison of the level curves of the marginals of $\overline{R}_{\alpha}^{D_{11}}$ (given by (\ref{approx})) and $R_{\alpha}^{ D_{11}, a_{2}=1}$. Dashed lines denote the level sets of the approximation $\overline{R}_{\alpha}^{D_{11}}$. The small difference is best seen in the iso-contours close to zero. In the exact case oriented photon may loop, whereas in the approximate case oriented particles must move forward.
Bottom Row, right:
A comparison between the exact Green's function of the resolvent diffusion process {\small $\alpha=\frac{1}{30}$, $D_{11}=0.1$, $D_{22}=0.5$} in Theorem \ref{th:thesolutionresolvent} and the approximate Green's function (in dashed lines) of the resolvent process with infinite lifetime $\lim_{\alpha \to 0}\alpha^{-1}(D_{11}\hat {A}_{1}^{2}+D_{22}\hat{A}_{2}^{2}-\alpha I)^{-1}\delta_{e}$, {\tiny $D_{11}=0.1$, $D_{22}=0.5$} given by (\ref{Happrox}). Bottom left: 3D-view on a stack of iso-contours (top: approximation, bottom exact) viewed along $\theta$-direction.
}}\label{fig:comparisonspatialandFD}
\end{figure}

Here we will derive the Green's functions for \emph{contour enhancement}, which are the heat-kernels on $SE(2)$.
In the case of \emph{contour-completion}, however, one has the interesting situation that the approximative left-invariant vector field $\hat{A}_{2}=\partial_{x}+\theta \partial_{y}$ together with the diffusion generator $(\partial_{\theta})^{2}$ and the identity operator $I$ and all commutators form an $8$-dimensional nil-potent Lie-algebra spanned by
\mbox{$
\{I, \partial_x,\partial_{\theta},\partial_{y},\theta \partial_{y}, \partial_{\theta}^{2}, \partial_{\theta}\partial_{y},\partial_{y}^2\}
$}.
From this observation and \cite{Varadarajan}{Theorem 3.18.11 p.243} it follows that the approximations of the Green's functions (which are again Green's functions but of a different Heisenberg type of group of dimension 5)
\begin{equation} \label{approx}
\begin{array}{ll}
\overline{K}_{s}^{D_{11}, a_{2}=1}(x,y,\theta) &=\delta(x-s) \frac{\sqrt{3}}{2\, D_{11}\pi x^2} e^{-\frac{3(x \theta -2y)^2+x^2(\theta-\kappa_{0} x )^{2}}{4 x^3 D_{11}}}\\
\overline{R}_{\alpha}^{D_{11}, a_{2}=1}(x,y,\theta) &= \alpha \,
\frac{\sqrt{3}}{2\, D_{11}\pi x^2} e^{-\alpha x} e^{-\frac{3(x \theta -2y)^2+x^2(\theta-\kappa_{0} x )^{2}}{4 x^3 D_{11}}} \rm{u}(x),
\end{array}
\end{equation}
where $\rm{u}$ denotes the 1D-Heavy-side/unit step function.
This technique can not be applied to the diffusion case, as the commutators of the separate diffusion generators provide infinitely many directions.
Here we follow \cite{Citti} and apply a coordinates transformation 
\begin{equation} \label{coordinatetrafo}
\overline{K}^{D_{11},D_{22}}_{s}(x,y,\theta)=\tilde{K}_{s}(x',\omega',t')=\tilde{K}_{s}\left(\frac{x}{\sqrt{2D_{22}}}, \frac{\theta}{\sqrt{2D_{11}}}, \frac{2(y-\frac{x\theta}{2})}{\sqrt{D_{11}D_{22}}}\right)
\end{equation}
where we note
{\small
\[
\begin{array}{l}
\partial_{s} \overline{K}^{D_{11},D_{22}}_{s}= \left(D_{11} \partial_{\theta}^{2} +D_{22} (\partial_{x}+\theta \partial_{y})^2 \right) \hat{K}^{D_{11},D_{22}}_{s} \desda \\[8pt]
\partial_{s} \tilde{K}^{D_{11},D_{22}}_{s}= \frac{1}{2}\left((\partial_{\omega'}-2 x' \partial_{t'})^{2} + (\partial_{x'}+2 \omega' \partial_{t'})^{2} \right) \tilde{K}^{D_{11},D_{22}}_{s}= \frac{1}{2}\Delta_{K} \tilde{K}^{D_{11},D_{22}}_{s},
\end{array}
\]
}
which provides us the left-invariant evolution equation on the usual Heisenberg group $H_{3}$ generated by Kohn's Laplacian. The Heat-kernel on $H(3)$ is well-known, for explicit derivations see \cite{Gaveau},\cite{Faraut84},\cite{DuitsR2006SS2}, and is given by
\begin{equation} \label{heatkernelonH3}
K_{s}^{\ul{D}}(x,\omega,t)=\frac{1}{(2\pi s)^2} \int_{\R} \frac{2\tau}{\sinh (2\tau)} \cos \left(\! \frac{2\, \tau \, t}{s \sqrt{D_{11} D_{22}}} \!\right)\, e^{-\frac{\left(\frac{x^2}{D_{11}s} +\frac{\omega^2}{D_{22}s}\right) \tau}{ 2\tanh(2\tau)}} \; {\rm d}\tau,
\end{equation}
and as a result by (\ref{coordinatetrafo}) we obtain\footnote{
Note that our approximation of the Green's function on the Euclidean motion group does not coincide with the formula by Citti 
in \cite{Citti}.}
the following Heisenberg-type approximation of the Green's function
\begin{equation} \label{Happrox}
\begin{array}{l}
\overline{K}_{s}^{D_{11},D_{22}}(x,y,\theta)=\frac{1}{2 D_{11}D_{22}}\tilde{K}_{s}^{H_3}\left(\frac{x}{\sqrt{2D_{11}}}, \frac{\theta}{\sqrt{2D_{11}}}, \frac{2(y-\frac{x\theta}{2})}{\sqrt{D_{11}D_{22}}}\right) \\
 \ \ \ =  \frac{1}{8 D_{11}D_{22}\pi^2 s^2} \int \limits_{\R} \frac{2\tau}{\sinh (2\tau)} \cos \left(\frac{2 \tau(y-\frac{x\theta}{2})}{s\sqrt{D_{11}D_{22}}}\right)e^{-\frac{\left( \frac{x^2}{s\, D_{22}}+ \frac{\theta^2}{s\, D_{11}}\right) \tau}{\tanh(2\tau)}} \; {\rm d}\tau \\
\lim \limits_{\alpha \to 0} \alpha^{-1}\hat{R}_{\alpha}(x,y,\theta)= \frac{1}{4 \pi D_{11} D_{22}}
\frac{1}{\sqrt{\frac{1}{16}\left( \frac{x^2}{D_{22}} +\frac{\theta^2}{D_{11}}\right)^2 + \frac{(y-\frac{1}{2}x\theta)^2}{D_{11}D_{22}}}}.
\end{array}
\end{equation}
See Figure \ref{fig:comparisonspatialandFD} for illustrations of both the exact resolvent Green's function $R^{D}_{\alpha}$ and its approximation

\subsubsection{The H\"{o}rmander condition and the Underlying Stochastics of the Heisenberg approximation of the Diffusion process on $SE(2)$ \label{ch:hoermander}}

A differential operator $L$ defined on a manifold $M$ of dimension $n \in \mathbb{N}, n<\infty$ is called hypo-elliptic if for all distributions $f$ defined on an open subset of $M$ such that $Lf$ is $C^{\infty}$ (smooth), $f$ must also be $C^{\infty}$. In his paper \cite{Hoermander}, H\"{o}rmander presented a sufficient and essentially necessary condition for an operator of the type
\[
L= c+ X_0 + \sum \limits_{i=1}^{r} (X_i)^2, \qquad r \leq n
\]
where $\{X_{i}\}$ are vector fields on $M$, to be hypo-elliptic. This condition, which we shall refer to as the H\"{o}rmander condition is that among the set
\begin{equation} \label{hoermandercond}
\{X_{j_1}, [X_{j_1},X_{j_2}], [X_{j_1},[X_{j_2},X_{j_3}]], \ldots , [X_{j_1},[X_{j_2},[X_{j_3}, \ldots, X_{j_k}]]]\ldots \; |\; j_i \in \{0,1,\ldots,r\}\}
\end{equation}
there exist $n$ which are linearly independent at any given point in $M$. Note that if $M$ is a Lie-group and we restrict ourselves to left-invariant vector fields than it is sufficient to check whether the vector fields span the tangent space at the unity element.

If we apply this theorem to the Forward Kolmogorov equation of the direction process than we see that the H\"{o}rmander condition is satisfied since we have $M=SE(2)\times \R^{+}$, $X_{0}=-\partial_{s}-\partial_{\xi}$, $X_{1}=\partial_{\theta}$ and we have
\[
\textrm{dim}\; \textrm{span} \{-\partial_{s}-\partial_{\xi},\partial_{\theta}, [\partial_{\theta}, -\partial_{s}-\partial_{\xi}], [\partial_{\theta},[\partial_{\theta}, -\partial_{s}-\partial_{\xi}]]\}= \textrm{dim} \textrm{span}\{\partial_{s},\partial_{\theta},\partial_{\xi},\partial_{\eta}\}= 4
\]
and indeed the Green's function of Mumford's direction process is infinitely differentiable on $SE(2)$, see \cite{DuitsR2006AMS}. Similarly the Green's function of the resolvent direction process determined by $L R=\delta_{e}$, with $L=-\partial_{\xi}+ D_{11}(\partial_{\theta})^{2}-\gamma I$ is infinitely differentiable on $SE(2)\backslash\{e\}$, for explicit formulae see \cite{DuitsR2006AMS}. To this end we set $M=SE(2)$ and we note that
\[
\textrm{span}\{\partial_{\theta}, [\partial_{\theta},\partial_{\xi}], [\partial_{\theta}, [\partial_{\theta},\partial_{\xi}]]\}= \textrm{span}\{\partial_{
\theta},\partial_{\xi}, \partial_{\eta}\}=\mathcal{L}(SE(2)).
\]
However, in the case of the direction process the Heisenberg approximation of the time dependent Green's function (\ref{approx}) is singular. This is in contrast to its  resolvent kernel where the Laplace transform takes care of the missing direction direction in the tangent space:
\[
\begin{array}{l}
\textrm{dim} \; \textrm{span} \{-\partial_{s}-(\partial_{x}+\theta \partial_y),\partial_{\theta}, [\partial_{\theta}, -\partial_{s}-(\partial_{x}+\theta \partial_y)], [\partial_{\theta},[\partial_{\theta}, -\partial_{s}-(\partial_{x}+\theta \partial_y)]]\}\\
= \textrm{dim} \; \textrm{span}\{\partial_{s}+\partial_{x},\partial_{y},\partial_{s}+\partial_{x}+\theta \partial_{y}\}= 3
\end{array}
\]
By the preceding it follows that this deficit does not occur in the Heisenberg approximation (\ref{Happrox}) of the pure diffusion (contour completion) case. This can be understood by the H\"{o}rmander condition, since set $S=SE(2)\times \R^{+}$ then we have
\[
\textrm{dim}\; \textrm{span} \{\partial_{s},\partial_{x}+\theta \partial_{y},\partial_{\theta},[\partial_{\theta},\partial_{x}+\theta \partial_{y}]=\partial_{y} \}=4.
\]
This puts us to the following question: \\
``Can we get physical insight in the induced smoothing in the remaining directions in the diffusion processes on $SE(2)$ generated by hypo-elliptic operators which are not elliptic ?" \\
Before we provide an affirmative answer to this question we get inspiration from the heat kernel on the 3D-Heisenberg group $H_{3}$, recall (\ref{heatkernelonH3}), which is smooth in all directions, despite the fact that diffusion is only done in $\partial_{x}+\omega \partial_{t}$ and $\partial_{\omega}-x\partial_{t}$-direction. Here,
the induced smoothness in $t$ direction, has an elegant stochastic interpretation. As shown in \cite{Gaveau}, the underlying stochastic process (with the diffusion equation on $H_{3}$ as the forward Kolmogorov equation) is given by
\begin{equation} \label{Stochprocess}
\left\{
\begin{array}{l}
Z(s)=X(s) + i \, W(s)=Z_{0} +\varepsilon \sqrt{s},\ \varepsilon \sim \mathcal{N}(0,1) \\
T(s)= 2 \int \limits_{0}^{s} W{\rm d}X-X{\rm d}W, \ s>0
\end{array}
\right.
\end{equation}
so the random variable $Z$ is a Brownian motion in the complex plane 
and the random variable $T(s)$ measures the deviation from a sample path with respect to a straight path $Z(s)=Z_{0}+ s(Z(s)-Z_{0})$ by means of the \emph{stochastic} integral $T(s)=2  \int \limits_{0}^{s} W{\rm d}X-X{\rm d}W$. To this end we note that for\footnote{A Brownian motion is a.e. not differentiable in the classical sense, nor does the integral in (\ref{Stochprocess}) make sense in classical integration theory.} $s \mapsto (x(s),\omega(s)) \in C^{\infty}(\R^{+},\R^2)$ such that the straight-line from $X_{0}$ to $X(s)$ followed by the inverse path encloses an oriented surface $\Omega \in \R^2$, we have by Stokes' theorem
that
$
2 \mu(\Omega)=-\int_{0}^{s} (-X'(t)W(t)+X(t)W'(t))\, {\rm d}t+0= \int_{0}^{s} W{\rm d}X-X{\rm d}W.
$

Now by the coordinate transformation in (\ref{coordinatetrafo}) we directly deduce that the underlying
stochastic process of the Heisenberg approximation of the diffusion process on $SE(2)$ is given by
\[
\left\{
\begin{array}{l}
X(s)+i\,\Theta(s)= X(0)+i\, \Theta(0) + \sqrt{s}(\epsilon_x+i \, \epsilon_{\theta}), \textrm{ where }\epsilon_{x} \sim \mathcal{N}(0,2D_{11}), \epsilon_{\theta} \sim \mathcal{N}(0,2D_{22}) \\
Y(s)=\frac{X(s)\Theta(s)}{2} +\frac{1}{2}\int_{0}^{s} \Theta {\rm d}X-X{\rm d}\Theta=\int_{0}^{s}\Theta(t)-\Theta(0){\rm d}t,
\end{array}
\right.
\]
which provides a better understanding of the ``implicit smoothing'' (by means of the commutators) within the H\"{o}rmander condition of the Heisenberg approximation of the diffusion process on $SE(2)$.

\section{Modes\label{ch:elastica}}

The concept of a completion distribution is well-known in image analysis, see for example \cite{Thornber}, \cite{Zweck}, \cite{August2003}, \cite{Duits2005IJCV}.
The idea is simple: Consider two left-invariant stochastic processes on the Euclidean motion group, one with forward convection say its forward Kolmogorov equation is generated by $A$ and one with the same stochastic behavior but with backward convection, i.e. its forward Kolmogorov equation is generated by the adjoint of $A^{*}$ of $A$. Then we want to compute the probability that random walker from both stochastic processes collide. This collision probability density is given by
\[
C^{(U,V)}= (A-\alpha I)^{-1}U (A^{*}-\alpha I)^{-1}W \qquad U,W \in \mathbb{L}_{2}(SE(2)) \cap \mathbb{L}_{1}(SE(2))
\]
where $U, W$ are initial distributions. This collision probability is called a completion field as it serves as a model for perceptual organization in the sense that elongated local image fragments are completed in a more global coherent structure.
These initial distributions can for example be obtained from an image by means of a well-posed invertible wavelet transform constructed by a reducible representation of the Euclidean motion group as explained in \cite{Duits2005IJCV}. Alternatives are lifting using the interesting framework of curve indicator random fields \cite{August} or (more ad-hoc) by putting a limited set of delta distributions after tresholding some end-point detector or putting them simply by hand \cite{Zweck}. Here we do not go into detail on how these initial distributions can be obtained, but only consider the case $U= \delta_{(0,0,\theta_{0})}$ and $W=\delta_{(x_1,y_{1},\theta_{1})}$, $x_{1},y_{1} \in \R$. In this case we obtain by means of (\ref{approx}) the following approximations of the completion fields:
{\small
\begin{equation} \label{complfield}
\begin{array}{ll}
\hat{C}_{g_0,g_1}^{\alpha,D_{11},\kappa_{0},\kappa_{1}}(\ul{x},\theta)&=\!\alpha^2\left( (\hat{A}\!-\!\alpha I)^{-1} \delta_{\ul{x}_{0},\theta_{0}}\right)(\ul{x},\theta)  \left((\hat{A}^{*}\!-\!\alpha I)^{-1} \delta_{\ul{x}_1,-\theta_{1}}\right)(\ul{x},\theta) \\
 &=
T_{\alpha, \kappa_0, D_{11} \; ; \ul{x}_0,\theta_0}(x,y,\theta)\,
T_{\alpha, \kappa_{1}, D_{11} \; ; -x_1,y_1,\theta_1}(-x,y,-\theta),
\end{array}
\end{equation}
}
with corresponding modes, obtained by solving for
\[
\left\{
\begin{array}{l}
\partial_{y}\{T_{\alpha, \kappa_0, D_{11} \; ; \ul{0},\theta_0}(x,y,\theta)
\, T_{\alpha, \kappa_{1}, D_{11} \; ; -x_1,y_1,\theta_1}(-x,y,-\theta)\}=0 \\
\partial_{\theta}\{T_{\alpha, \kappa_0, D_{11} \; ; \ul{0},\theta_0}(x,y,\theta)
\, T_{\alpha, \kappa_{1}, D_{11} \; ; -x_1,y_1,\theta_1}(-x,y,-\theta)\}=0
\end{array}
\right.
\]
These modes depend on only on the difference $\kappa_0-\kappa_{1}$ but not on $D_{11}$ nor on $\alpha$:
\begin{equation} \label{modes}
\begin{array}{ll}
y(x) &=  x \theta_{0} + \frac{x^3}{x_1^3}(-2 y_1+x_{1}(\theta_{0}-\theta_{1})) + \frac{x^2}{x_1^2}(3 y_1 +x_1(\theta_1-2 \theta_0)) \\
 &+\frac{(\kappa_{0}-\kappa_1)}{ x_1^{3}} (x-x_1)^2(\frac{x_1}{2}-x)x^2 \\
\theta(x) &= \theta_{0} + 2 \frac{x}{x_1^2}(3 y_1 + x_1(\theta_1-2\theta_0))
-3 \frac{x^2}{x_{1}^{3}}(2y_1+x_1(\theta_{1}-\theta_{0})) \\
&+\frac{(\kappa_{0}-\kappa_1)}{x_1^{3}}
x(x-x_1)(-3x^2+3x_1 x -x_1^{2})\ ,
\end{array}
\end{equation}
where $x \in [0,x_{1}]$ and $y(0)=0$, $\theta(0)=\theta_{0}$ and $y(x_1)=y_1$, $\theta(x_1)=-\theta_{1}$ and
$\frac{dy}{dx}(0)=\theta_{0} \textrm{ and } \frac{dy}{dx}(x_1)=-\theta_{1}$, see Figure \ref{fig:CompletionFields}. These modes are the unique minimizers of the following variational problem
{\small
\begin{equation} \label{varprob}
\textrm{argmin} \{\mathcal{E}(y)=\int \limits_{0}^{x_1} \left(y''(x)-\frac{(\kappa_{1}-\kappa_{0})\, c(x)}{x_1^{3}}\right)^2 {\rm d}x\; |\;  y(0)\!=\!0,y(x_1)\!=\!y_1, y'(0)\!=\!\theta_{0}, y'(x_1)\!=\!-\theta_{1}\}
\end{equation}
}
where \mbox{$c(x)=20(x-x_1)(x-\frac{x_1}{2})x$} and
$y'(x)=\theta(x)+\frac{(\kappa_{1}-\kappa_{0})\, d(x)}{x_1^{3}}$ with \\ \mbox{$d(x)=-2 x^2(x-x_1)^2$.}
The variational problem (\ref{varprob}), for the case $\kappa_{0}=\kappa_{1}$ is indeed the corresponding (with arclength replaced by $x$) approximation of the elastica functional in \cite{Mumford} and indeed $\frac{\partial E}{\partial v}(y)=0$ for all $v \in \mathcal{D}((0,x_1))$ if and only if $y^{(4)}(x)=(\kappa_{1}-\kappa_{0})c^{(2)}(x)$ under the conditions $y(0)=0,y(x_1)=y_1, y'(0)=\theta_{0}, y'(x_1)=-\theta_{1}$.

We note that because of left-invariance with respect to the 5-dimensional Heisenberg type of group we have
\[
\hat{S}_{\alpha,\kappa_0,D_{11}; \ul{x}',\theta'}(x,y,\theta)=
\hat{S}_{\alpha, \kappa_0,D_{11}; e}(x-x',y-y'-\theta'(x-x'),\theta-\theta').
\]
As a result the approximate completion field (and thereby its mode) is not left-invariant on $\R^2 \rtimes \mathbb{T}$ and thereby its marginal is not Euclidean invariant. As a result the formulas do depend\footnote{If \{x,y\} is aligned with $g_0$ the result is different then if it is aligned with $g_{1}$.} on the choice of coordinate system $\{x,y\}$.

However, this problem does not arise for the exact completion field
\[
C^{g_0,g_1,\alpha,D_{11},\kappa_0,\kappa_1}=\alpha^2\left( (A\!-\!\alpha I)^{-1} \delta_{\ul{x}_{0},\theta_{0}}\right)  \left((A^{*}\!-\!\alpha I)^{-1} \delta_{\ul{x}_1,\theta_{1}}\right),
\]
since by left-invariance of the generator $A$ we have
\[
\begin{array}{l}
(A-\alpha I)^{-1} \delta_{g_0}=(A-\alpha I)^{-1} \mathcal{L}_{g_0}\delta_{e}=\mathcal{L}_{g_0}(A-\alpha I)^{-1} \delta_{e} \textrm{ and therefore } \\
C^{h\, g_0,h\, g_1,\alpha,D_{11},\kappa_0,\kappa_{1}}=\mathcal{L}_{h} C^{g_0,g_1,\alpha, D_{11},\kappa_0,\kappa_{1}} \textrm{ for all }h \in \R^2\rtimes \mathbb{T}.
\end{array}
\]
Throughout this paper we shall often use the following convention
\begin{definition} \label{def:horizontal}
A curve $s \mapsto \gamma(s)=(x(s),y(s),e^{i\theta(s)})$ in $SE(2)$ is called horizontal iff $\theta(s)= \arg (x'(s)+i \, y'(s))$. Then $\gamma$ is called the lifted curve in $SE(2)$ of the curve $s \mapsto \ul{x}(s)=(x(s),y(s))$ in $\R^2$.
\end{definition}
Note that this sets a bijection between horizontal curves in $SE(2)$ and curves in $\R^2$.

\subsection{Elastica curves, geodesics, modes and zero-crossings of completion fields}

In his paper Mumford \cite{Mumford}{p.496} showed that the modes of the direction process are given by elastica curves which are by definition curves
$t \mapsto \ul{x}(t)$ in $\R^{2}$, with length $L$ and prescribed boundary conditions
\begin{equation}\label{BC}
\begin{array}{l}
\ul{x}(0)=\ul{x}_{0}, \ul{x}(L)=\ul{x}_{1} \textrm{ with prescribed directions }\\
\arg (x'(0)+i\, y'(0))=\theta_{0} \textrm{ and }
\arg (x'(L)+i\, y'(L))=\theta_{1}
\end{array}
\end{equation}
which minimize the functional
\[
\mathcal{E}_{\epsilon}(\ul{x})=\int_{0}^{L} \kappa^{2}(s) +\epsilon \; {\rm d}s, \qquad \textrm{ with }\epsilon=4 \alpha D_{11}, \textrm{ with } \kappa(s)=\|\ddot{\ul{x}}(s)\|, \textrm{ and }s>0 \textrm{ arclength.}
\]
Here he uses the following discrete version (with $N$ steps) of the stochastic process
\begin{equation} \label{discretedirectionP}
\left\{
\begin{array}{l}
e^{i\theta(s_k+\Delta s)}=e^{i\left(\theta(s_k) +\sqrt{\Delta s} \, \ve_k\right)}, \qquad \ve_k \sim \mathcal{N}(0,\sigma^2), \qquad D_{11}=\frac{1}{2}\sigma^2 \\
\ul{x}(s_k+\Delta s)=\ul{x}(s_k) + \Delta s \, \left(
\begin{array}{l}
\cos \theta(s_k) \\
\sin \theta(s_k)
\end{array}
\right) \\
\Delta s=\frac{L}{N} \textrm{, with }L \sim NE(\alpha), k=0,\ldots,N-1, \\
\end{array}
\right.
\end{equation}
where the physical dimension of $\ve_k$ and $\sigma=\sqrt{2\, D_{11}}$ is {\small $[LENGTH]^{-1/2}$}, for the definition of the mode:
\begin{definition}
The mode of the direction process with parameters $\alpha,D_{11}$ is a curve in $\R^2$ which is the point-wise limit of the maximum likelihood  curve of the discrete direction process with given boundary conditions (\ref{BC}).
\end{definition}
In his paper, \cite{Mumford}, Mumford states that the probability density of a discrete realization, a polygon of length $L$, $\Gamma= \bigcup \limits_{i=1}^{n} \overline{\ul{x}(s_i), \ul{x}(s_{i+1})}$ whose sides have length $\frac{L}{N}$ and the $\theta_i$ are discrete Brownian motion scaled down by $\Delta s = \frac{L}{N}$:
\[
\theta_{i+1}=\theta_{i}+ \sqrt{\frac{L}{N}} \ve_{i+1}, \qquad i \in \{0,\ldots,N-1\},
\]
$\ve_i$ independent normal random variables with mean $0$, standard deviation $\sigma=\sqrt{2\, D_{11}}$, equals
\[
\begin{array}{ll}
P(\Gamma) &= \frac{\alpha}{(\sqrt{2\pi} \sigma)^{N-1}}e^{-\sum \limits_{i=0}^{N-1} \ve_{i}^{2}/(2\sigma^2) -\alpha L}  \\
 &\equiv e^{-\sum \limits_{i=0}^{N-1} \frac{L}{N}\left(\frac{\theta_{i+1}-\theta_{i}}{L/N}\right)/(2\sigma^2) -\alpha L}
\end{array}
\]
which converges as $N\to \infty$ to
\begin{equation} \label{epsilonMumford}
e^{-\frac{1}{D_{11}}\int_{0}^{L} \kappa^{2}(s) +\epsilon \; {\rm d}s}, \textrm{ with }\epsilon=4\alpha D_{11}, \qquad D_{11}=\frac{1}{2}\sigma^2,
\end{equation}
from which he deduced that the modes (or maximum likelihood curves) are elastica curves.

The drawback of this construction is that the definition of the mode is obtained by means of a discrete approximation. However, Olaf Wittich brought to our attention that the above construction is quite similar to the minimization of the Onsager-Machlup functional $\int \left[ \frac{1}{2}\|\dot{\gamma}\|_{M}^{2}-\frac{1}{12}\textrm{Scal}(\gamma)\right]$ which under sensible conditions yields the asymptotically most probable path in a Brownian motion on a manifold $M$, \cite{Takahashi}, which does not require the jump from the continuous to the discrete setting and back. This is a point for future investigation.

\begin{figure}
\centerline{
\includegraphics[width=0.33\hsize]{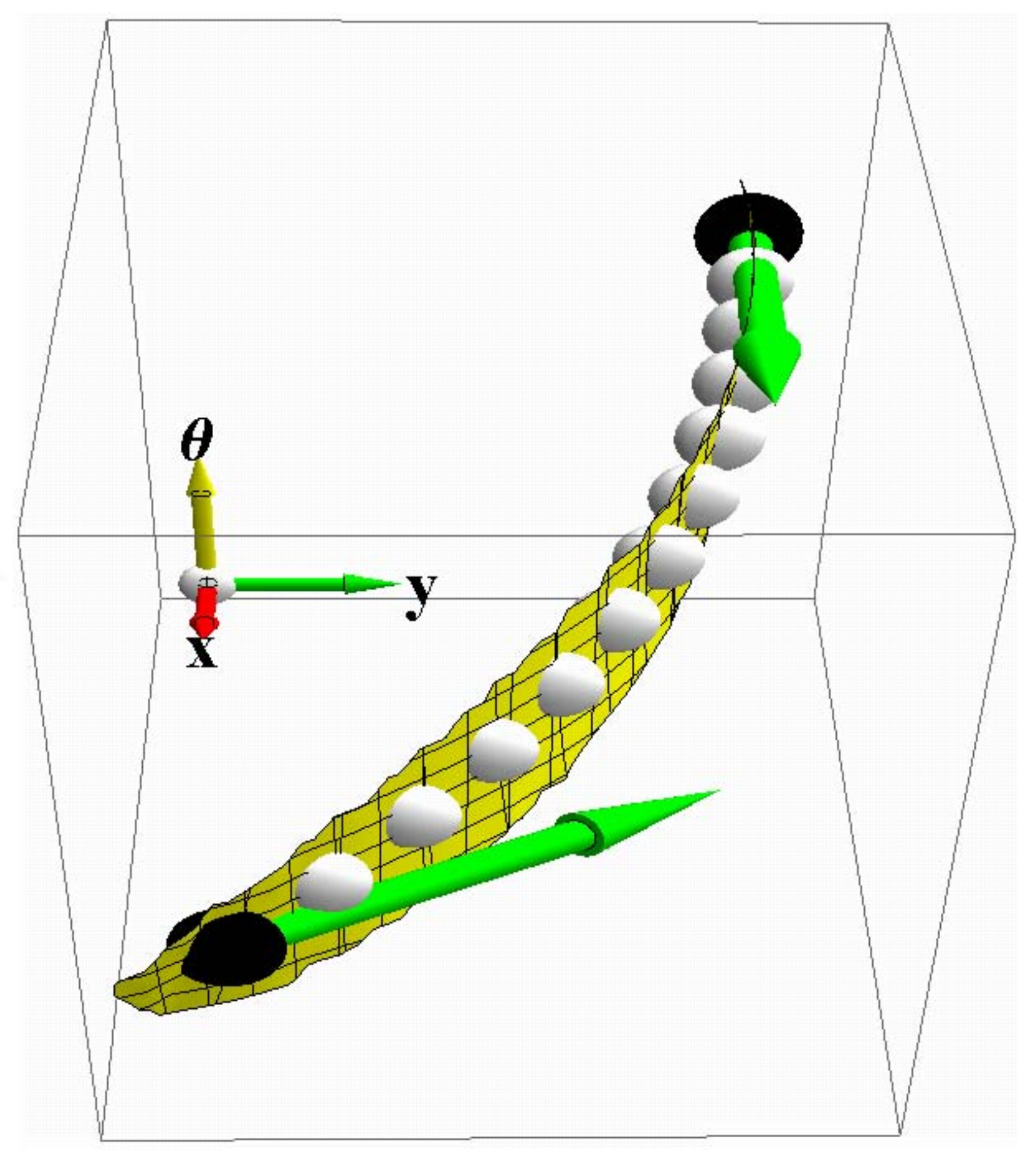}
\includegraphics[width=0.33\hsize]{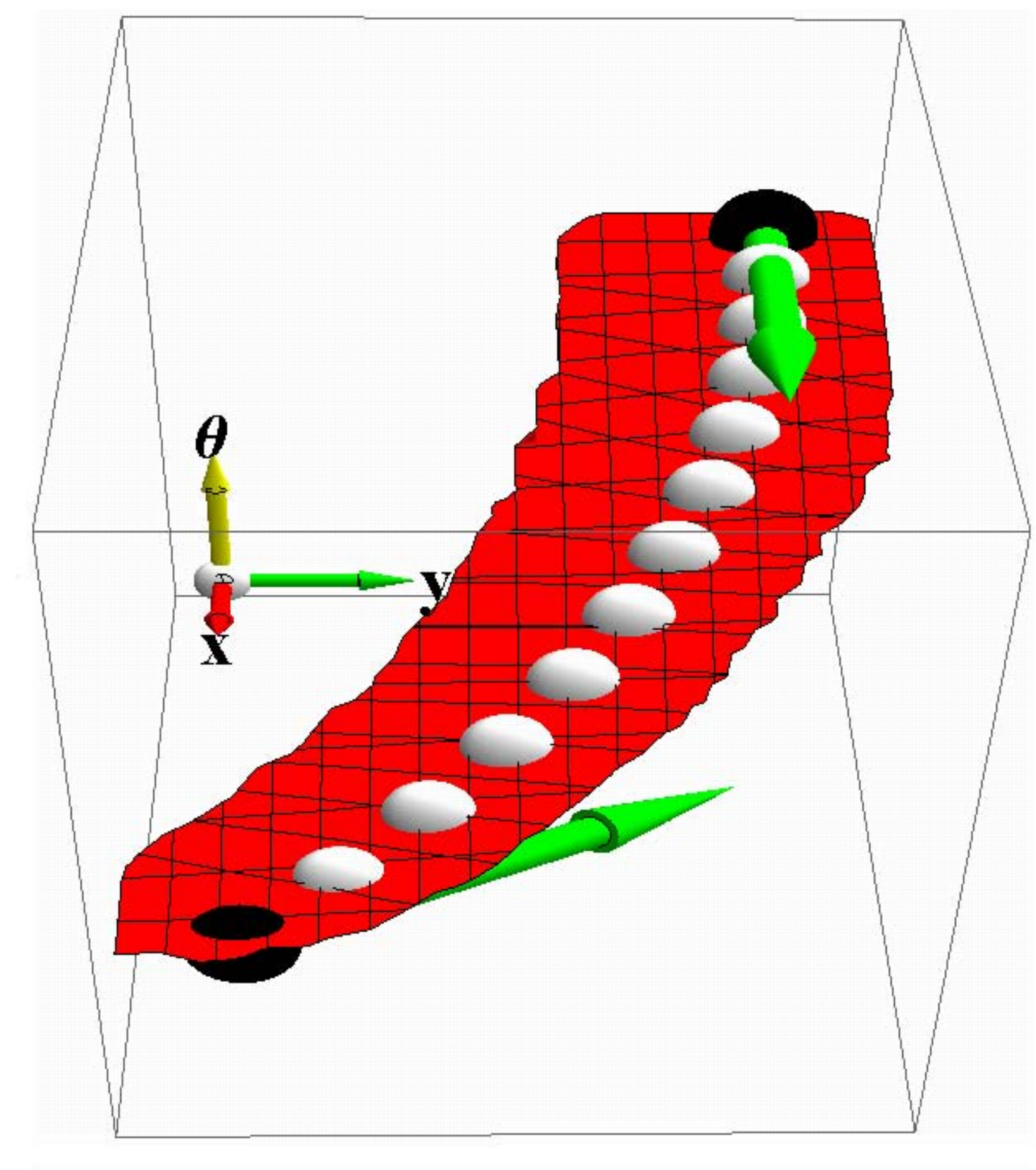}
\includegraphics[width=0.33\hsize]{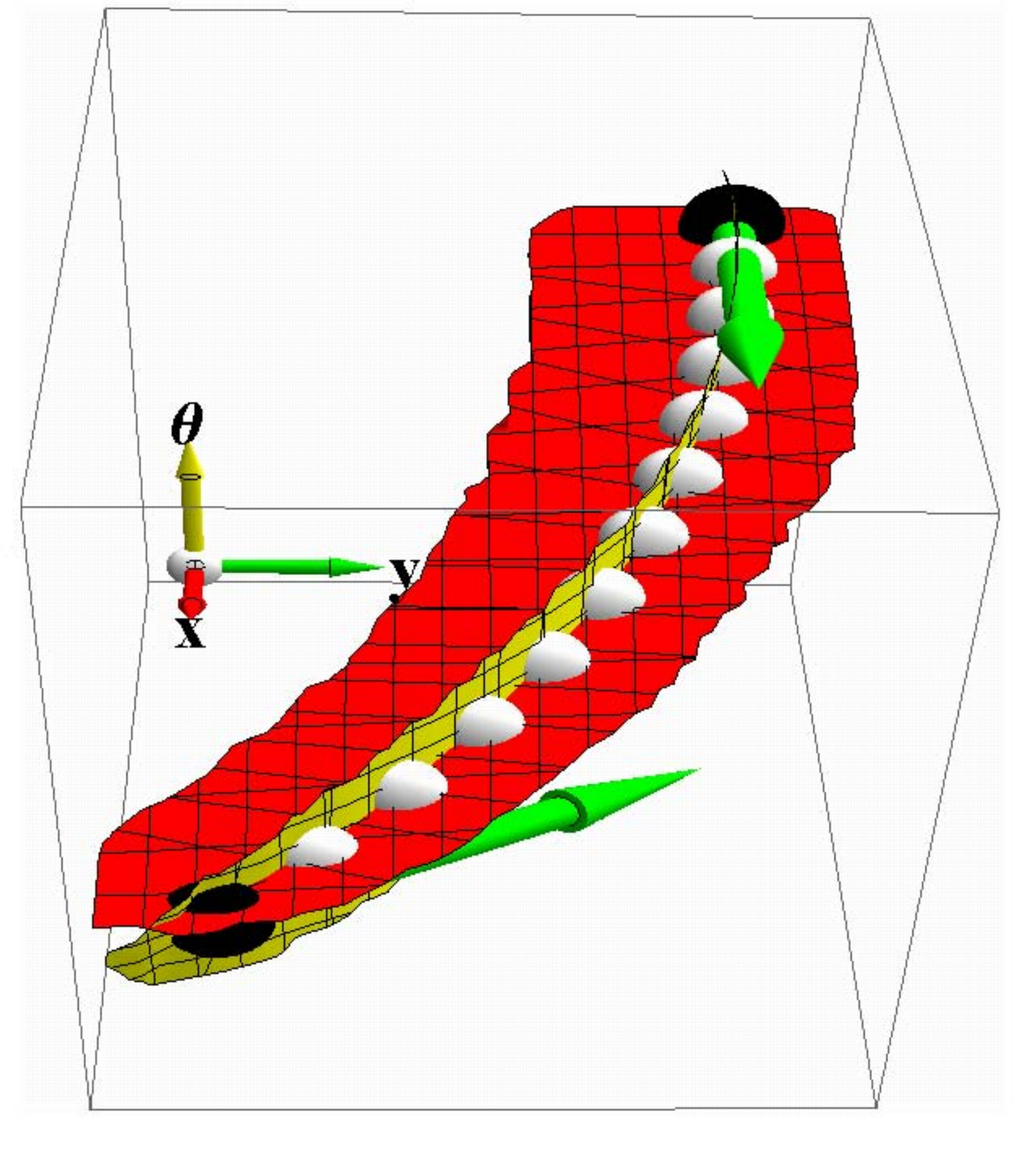}
}
\centerline{
\includegraphics[width=0.33\hsize]{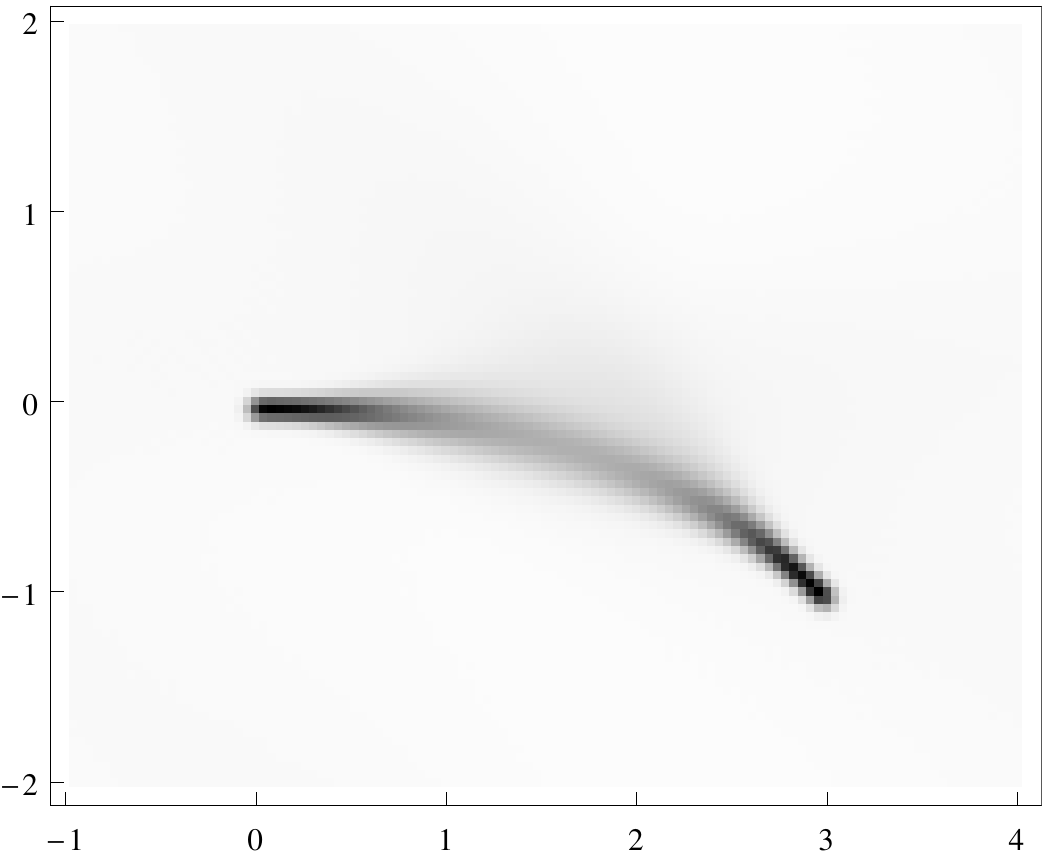}
\includegraphics[width=0.33\hsize]{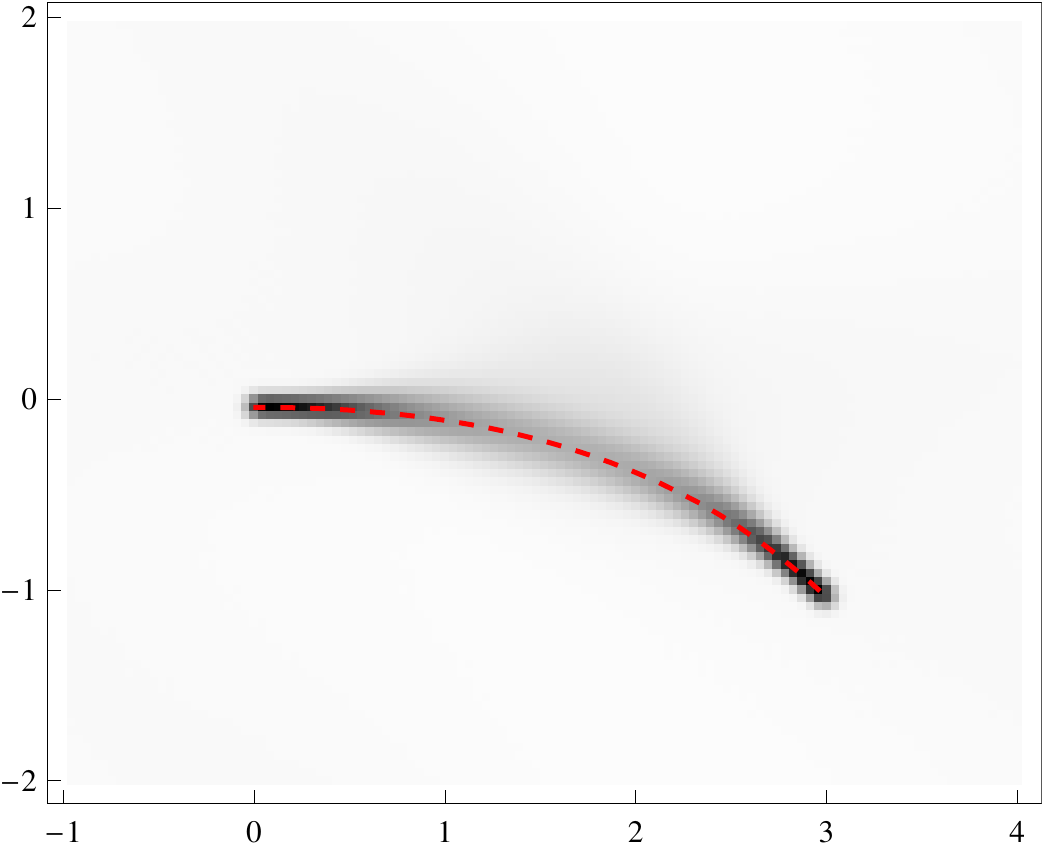}}
\caption{We computed an elastica curve that passes two given points $g_0=(\ul{x}_{0},e^{i\theta_{0}})=(0,0,1)$ and $g_{1}=(\ul{x}_{1},e^{i\theta_{1}})=(3,-1,\frac{7}{4\pi})$ via the shooting algorithm (\ref{NR}), $D_{11}=\frac{1}{32}$, $\alpha=\frac{1}{10}$,  and we lifted this curve in the Euclidean motion group via $\theta(s)=\angle (\dot{\ul{x}}(s),\ul{e}_{x})$, with parameter $\epsilon=4 \alpha D_{11}$ illustrated by a line of white balls centered around equidistant points on the elastica. Furthermore we computed the zero-crossing of the planes where the exact completion field  $C_{\alpha, D_{11}, \kappa_{0}=\kappa_{1}=0}^{g_{0},g_{1}}=\alpha^{2}(\alpha I +\partial_{\xi}-\partial_{\theta}^{2} )^{-1} \delta_{g_{0}}(\alpha I -\partial_{\xi}-\partial_{\theta}^{2} )^{-1} \delta_{g_{1}}$ (whose $xy$-marginal we depicted at the bottom) has  zero respectively $\theta$ (red-plane) and $\eta$ (yellow-plane) derivative. Note that the zero-crossing of these two planes is extremely close to the elastica curve. The green arrows reflect $\left. \ul{e}_{\xi} \right|_{g_{0}}$ and $\left. \ul{e}_{\xi} \right|_{g_{1}}$.}\label{fig:zerocross}
\end{figure}

Mumford's observation that the modes of the direction process (with parameters $D_{11}>0$ and $\alpha>0$) coincide with elastica curves (with $\epsilon=4\alpha D_{11}$) raises the following two challenging questions:
\begin{enumerate}
\item 
Is there a connection between the unique curve determined by the zero crossings of
\mbox{$
\partial_{\theta}C^{g_0, g_1, \alpha, D_{11}}$} and \mbox{$\partial_{\eta}C^{g_0, g_1, \alpha, D_{11}}$}, i.e. the unique curve determined by the intersection of the planes
\[
\{g \in SE(2)\; |\; \partial_{\theta}C^{g_0, g_1, \alpha, D_{11}}(g)=0 \} \textrm{ and } \{g \in SE(2) \; |\; \partial_{\eta}C^{g_0,g_1, \alpha, D_{11}}(g)=0\}
\]
and the (lifted) modes/elastica curves ? Do they coincide, likewise their respective Heisenberg approximations (\ref{varprob}) and (\ref{modes}) ?
\item What is the connection between this result and the well-known Onsager-Machlup functional which describes the asymptotic probability of a diffusion particle on a complete Riemannian manifold staying in a small ball around a given trajectory, \cite{Takahashi}.
\end{enumerate}
Numerical computations seem to indicate that the unique curve induced by the zero crossings of $\partial_{\theta}C=0$ and $\partial_{\eta}C=0$ closely approximate the elastica curves ! In fact they even seem to coincide, see figure \ref{fig:zerocross}.  The main problem with mathematically underpinning this numerical observation here is that we only have elegant formulae for the exact Green's functions in the Fourier domain. This problem did not occur in the case of the Heisenberg approximation.

\begin{figure}
\centerline{%
\includegraphics{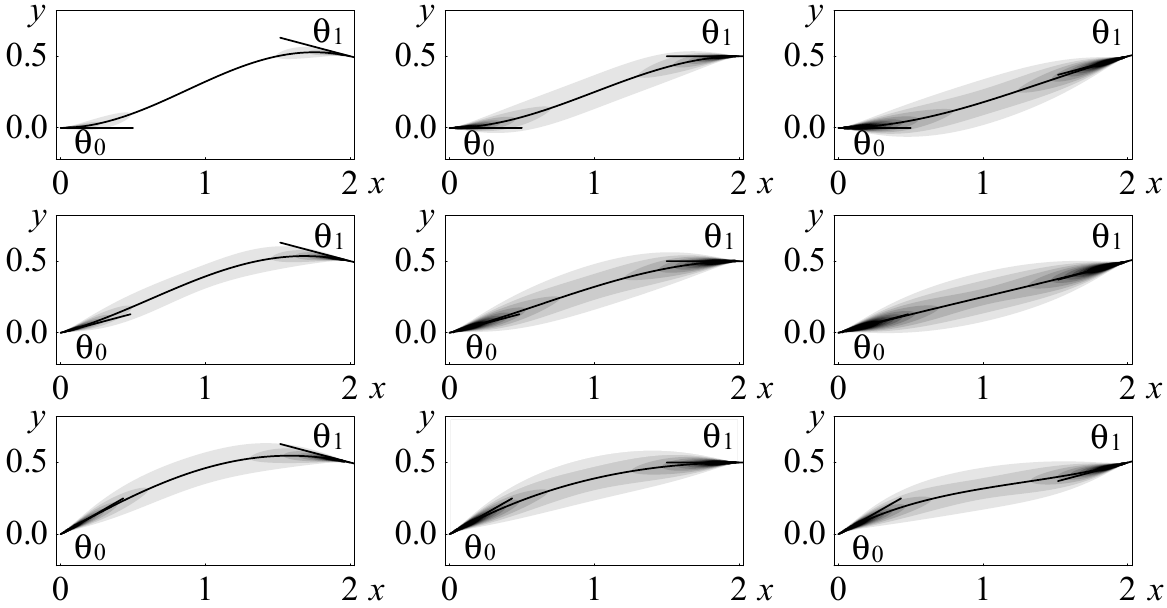}
}\caption{The shading in these plots denotes the marginal of the analytic completion field approximation (\ref{complfield}) obtained via integration over $\theta$ for $x \in (0,2)$, $y \in (-0.2,0.8)$  i.e. $\int_{\R} T_{D_{11}=0.5, \theta_0=2, \kappa_0=0}^{x_0=0, y_0=0, \theta_0}(x,y,\theta) T_{D_{11}=0.5, \theta_1=2, \kappa_1=0}^{x_1=2, y_1=0.5 ,\theta_1}(-x,y,-\theta) {\rm d\theta}$ for $\theta_{0}=0^\circ, 15^\circ, 30^\circ$ from top to bottom, and for $\theta_{1}=-15^\circ, 0^\circ, 15^\circ$ from left to right.
The lines drawn on top of these completion fields are the modes~(\ref{modes}), the optimal connecting lines.
}
\label{fig:CompletionFields}
\end{figure}

In the Heisenberg approximation case the intersection of the planes given by 
\[
\partial_{y}\hat{C}^{g_0,g_1, \alpha, D_{11}}(x,y,\theta)=0 \textrm{ and } \partial_{\theta}\hat{C}^{g_0,g_1, \alpha, D_{11}}(x,y,\theta)=0
\]
yield the $B$-spline solution (\ref{modes}) (minimizing the approximate elastica functional \ref{varprob} where the role of arc-length $s>0$ is replaced by $x$), which does not depend on $\alpha$, illustrated in Figure \ref{fig:CompletionFields}. This is due to the fact that the approximate resolvent Green's function
\[
\overline{R}_{\alpha}^{D_{11}, a_{2}=1}(x,y,\theta) = \alpha \,
\frac{\sqrt{3}}{2\, D_{11}\pi x^2} e^{-\alpha x} e^{-\frac{3(x \theta -2y)^2+x^2(\theta-\kappa_{0} x )^{2}}{4 x^3 D_{11}}} \rm{u}(x),
\]
 satisfies
\begin{equation} \label{alphaindep}
\overline{R}_{\alpha}^{D_{11}, a_{2}=1}(x,y,\theta)= \alpha e^{-\alpha x}\lim \limits_{\tilde{\alpha} \to 0} \frac{\overline{R}_{\tilde{\alpha}}^{D_{11}, a_{2}=1}(x,y,\theta)}{\tilde{\alpha}}
\end{equation}
which coincides with the fact that the random walker in the approximate case is not allowed to turn (it should always move forward in $x$-direction).

Regarding the exact case both the intersecting curve of the planes
\begin{equation}\label{planes}\{g \in SE(2) \; |\; \partial_{\theta}C^{g_0, g_1, \alpha, D_{11}}(g)=0\} \textrm{ and } \{g \in SE(2) \; |\; \partial_{\eta}C^{g_0, g_1, \alpha, D_{11}}(g)=0\}
\end{equation}
and the elastica curves will depend on $\alpha$. However, the Green's functions\footnote{for explicit formula for the exact resolvent Green functions of the direction proces similar to (\ref{final}) we refer to \cite{DuitsR2006AMS}.} $R_{\alpha}^{D_{11}, a_{2}=1}$ and thereby the intersection of the planes (\ref{planes}) , only depend on the \emph{quotient} $\frac{D_{11}}{\alpha}$ (after rescaling position variables by $\ul{x} \mapsto \frac{1}{\sqrt{D_{11}}} \ul{x}$) whereas the elastica curves only depend on the \emph{product} $\epsilon= 4\alpha*D_{11}$. So the curves will certainly not coincide for all parameter settings $(\alpha,D_{11})$. Furthermore the resolvent Green's function of the exact direction process $R_{\alpha}^{D_{11}, a_{2}=1}$ does not satisfy (\ref{alphaindep}). However for the special case $\alpha \to \infty$ and thereby $\frac{D_{11}}{\alpha} \to 0$ the approximate Green's function converges to the exact Green's function, see \cite{DuitsR2006AMS}, \cite{MarkusThesis}. Moreover as $\alpha \to \infty$ the arc-length $s>0$ of the projection of the path of an exact random walker of the direction process on the spatial plane (which is differentiable,  recall from Figure \ref{fig:orientationbundle2}) tends to the initial direction which is along the $x$-direction and as a result the energy (\ref{varprob}) (where we recall $\theta(x)=y'(x)$) tends to the elastica energy $\int (\kappa(s))^2\;{\rm d}s$ for fixed length curves, since for horizontal curves we have $\kappa(s)=\theta'(s)$, as $\alpha \to \infty$. So by the different dependence on the parameters $D_{11}$ and $\alpha$ we can only expect the intersection curve of the planes $\partial_{\theta}C=0$ and $\partial_{\eta}C=0$ to coincide with the elastica curve in the limiting case $\alpha \to \infty$.

In image analysis applications we typically have that $0<\frac{D_{11}}{\alpha}<< 1$ in which case the Heisenberg approximation is a good approximation, as a result for these parameter settings the intersections of the planes $\partial_{\theta}C^{g_0,g_1, \alpha, D_{11}}=0$ and $\partial_{\eta}C^{g_0,g_1, \alpha, D_{11}}=0$ depend very little on the parameter $\frac{D_{11}}{\alpha}$ and as a result they turn out to be close approximations of the elastica curves.

This serves as a theoretical motivation for our curve extraction from completion fields $C:SE(2)\to \R^{+}$ of orientation scores \ref{complfield} via the zero crossings of $\partial_{\theta}C$ and $\partial_{\eta}C$ which is useful for \emph{detecting} noisy elongated structures (such as catheters) in many medical image analysis applications.


Finally, we note that the above observations are quite similar to the result in large deviation theory, \cite{wittich}, where the time integrated unconditional Brownian bridge measure on a manifold $SE(2)$ uniformly tends to the geodesics which for $\alpha \to \infty$, see Appendix \ref{ch:app2}. Here, we stress that the Brownian bridge measure on the manifold $SE(2)$ is related to the completion fields of the contour \emph{enhancement} process on $SE(2)$ (with $D_{33}\geq 0$) and in limiting case $D_{33} \to 0$ the corresponding geodesics (where we restrict ourselves to horizontal curves) coincide with the horizontal minimizers of $\int \sqrt{\kappa^{2}+\epsilon} {\rm d}s$. These horizontal curves were also reported by Citti and Sarti \cite{Citti} as geodesics on $SE(2)$ and note that these curves are, in contrast to the elastica curves, coordinate independent on the space $SE(2)$ !

In respectively Section \ref{ch:elastica} and Section \ref{ch:geodesics} we will derive exact formulae for the curvature of geodesics and the curvature of elastica curves which are well-known and we provide a numeric comparison between the curves. For reasonable parameter settings these curves turn out to be quite close to eachother. The well-known problem with elastica curves is that their curvature involves Elliptic functions. Nevertheless, it is still possible to integrate the curvature twice and to provide analytic formulae for the curves themselves, \cite{MarkusThesis}. The geodesics however do not suffer from this problem, but they cause numerical problems in shooting algorithms, because of singularities outside the boundary conditions.

Therefore we also point to Appendix \ref{ch:app} where we explicitly compute the geodesics. Here we will order our results in much more abstract and structured way by means of Pfaffian systems. Moreover, by applying the Bryant and Griffiths approach \cite{Bryant} on the Marsden-Weinstein reduction for Hamiltonian systems admitting a Lie group of symmetries (developed for elastica curves) to the geodesics we are able to get nice analytic formulae, which do not seem to appear in literature, for the geodesics which have the advantage that they do not involve special functions.

Moreover, we refer to Appendix \ref{ch:app3} for the computation of snakes/actice shape models in $SE(2)$ based on completion fields of orientation scores and elastica curves.
See figure \ref{fig:zerocross}.

\section{The Underlying Differential Geometry: The Cartan Connection on the principal fiber bundle $P_H=(SE(2),SE(2)/H,\pi,R)$ \label{ch:Cartan}}

The goal is to obtain a connection on $SE(2)$ such that the exponential curves are geodesics and compute its curvature and torsion. Moreover we would like to relate the connection to a left-invariant Riemannian metric where the left-invariant vector fields $\{\mathcal{A}_{i}\}_{i=1}^{3}=\{\partial_{\theta},\partial_{\xi},\partial_{\eta}\}$ serve as a moving frame of reference and we want to compute the covariant derivatives of these left-invariant vector fields. To this end we recall the general Cartan-connection construction. In contrast to the more familiar Levy-Cevita connection this connection does not require a metric. Nevertheless, it is possible to relate this connection to a left invariant metric constructed from the Killing form on the Lie algebra $T_{e}(SE(2))$. Unfortunately this killing form is degenerate on $SE(2)$ therefore we will embed the Lie-algebra of $SE(2)$ in the Lie-algebra of $SO(3)$ where the Killing form is non-degenerate. Throughout this section we will use the Einstein summation convention.

Let $G$ be a Lie group of finite dimension $n$ with unit element $e$ and subgroup $H$. By setting the equivalence relation $a\sim b \desda a^{-1}b \in H$ on $G$ we get the left cosets as equivalence classes $[g]=g H$. Let $G/H$ denote the partition of left cosets on $G$. Let $\pi: G \to G/H$ be the projection of $G$ onto $G/H$ given by $\pi(g)=[g]$ and let $R$ be the right multiplication given by $R_h g=gh$. Note that $\pi(R_h g)=\pi(g)$. This yields a principal fiber bundle $P_H=(G,G/H, \pi, R)$ with structure group $H$. A Cartan/Ehresmann connection\footnote{In the common case of Riemannian geometry, with Riemannian connection $\nabla$, one can create a Lie-algebra valued one-form by means of $\omega(a^{l} \partial_{l})=\Gamma^{k}_{ij} {\rm d}x^{i} (a^{l} \partial_{l}) \partial_{k}=w^{k}_{j} (a^{i} \partial_{i})$, where the 1-forms $\omega^{k}_{j}$ are given by $\omega^{k}_{j}=\Gamma^{k}_{ij}{\rm d}x^{i}$, where the Christoffel symbols are given by $\Gamma^{k}_{ij}=({\rm d }x^{i},\nabla_{X_i} X_{j})$, with $({\rm d }x^{i},X_{j})=\delta^{i}_{j}$. Necessary and sufficient conditions for $\nabla$ the map $\nabla$ given by $\nabla_{X}X_{j}=\omega^{k}_{j}(X) X_{k}$ to be a Riemannian connection are $\omega^{j} \wedge \omega^{i}_{j}=0$ and ${\rm d}g_{ij}=g_{kj}\omega^{k}_{i}+g_{ik}\omega^{k}_{j}$. Note however that a Cartan connection in contrast to the Riemannian connection does not require a metric. } $\omega :TP \to T_{e}(G)$ is a Lie-algebra valued 1-form on $P_H$ such that
\begin{equation}\label{CartanEhresmann}
\begin{array}{l}
1) \   \omega((R_{h})_{*}Y)  = \textrm{Ad}(h^{-1})\omega(Y)  \textrm{ for all }h \in H \textrm{ and all vector fields }Y { on }P\\
2) \  \omega( {\rm d}\mathcal{R}(X))=X \textrm{ for all } X \in T_{e}(H),
\end{array}
\end{equation}
where $(R_{g})_{*}$ denotes the push-forward of the right-multiplication and
\begin{equation} \label{adjointrep}
\textrm{Ad}(g)=(R_{g^{-1}} L_{g})_{*}
\end{equation}
 which equals the derivative (at the unity element) of the conjugation automorphism on $h \mapsto hgh^{-1}$, $h,g \in H$. Note that requirement 1) means $\omega_{gh}((R_{h})_{*}Y_{g})=Ad(h^{-1})\omega_{g} (Y_{g})$, for all vector fields $Y$ and all $g, h \in H$.

In particular we take the Cartan-Maurer form $\omega_{h}=(L_h^{-1})^{*} \omega_{e} :T_{g}(P_H)\to T_{e}(H) \subset T_{e}(G)$, $g \in G$, with $\omega_{e} = I$.
By using the restrictions $\{\left. \mathcal{A}_{i} \right|_{g}\}_{i=1}^{n}$ of the left-invariant vector fields $\{\mathcal{A}_{i}\}_{i=1}^{n}$, with corresponding covectors $\{{\rm d}\mathcal{A}^{i}\}_{i=1}^{n}$ with $\langle {\rm d}\mathcal{A}^{i}, \mathcal{A}_{j}\rangle =\delta^{i}_{j}$ (also known as the Maurer Cartan co-frame),to $g$ as a local basis for $T_{g}(G)$ for all $g \in G$, where we assume that $\{\mathcal{A}_{i}\}_{i=1}^{n}$ are ordered such that the first $m \leq n$, $\{\mathcal{A}_{1},\ldots,\mathcal{A}_{m}\}$ elements generate the subgroup $H$ we can express the Cartan-Maurer form on $P_H$ as follows
\begin{equation} \label{MaurerCartan}
\omega_{g}(X_{g})= \sum \limits_{i=1}^{m} \langle \left. {\rm d}\mathcal{A}^{i}\right|_{g}, X_{g} \rangle_{T_{g}(G)} A_{i},
\end{equation}
for all vector fields on $P_H$ and where we recall that $\left. \mathcal{A}_{i}\right|_{g=e}=A_{i} \in T_{e}(G)$.

Next we give a brief derivation of (\ref{MaurerCartan}). First recall that the left-invariant vector fields
$\{\mathcal{A}_{i}\}_{i=1}^{n}$ satisfy $\left.\mathcal{A}_{i}\right|_{g}= (L_{g})_{*} A_{i}$, i.e. they are obtained from $T_{e}(G)$ by push forward of the left multiplication and therefor the dual elements (the corresponding co-vector fields) are obtained by the pull-back from $T_{e}(G)$
\[
\left. {\rm d} \mathcal{A}^{i} \right|_{g}= (L_{g})^{*} {\rm d}A^{i}
\]
since we have $\langle (L_{g})^{*} {\rm d }A^{i}, (L_{g})_{*} A_{i}\rangle= \langle {\rm d }A^{i}, A_{j}\rangle =\delta^{i}_{j}$. Now direct computation yields
\[
\begin{array}{ll}
(L_{h^{-1}})_{*} X_{g}  &=   \sum \limits_{i=1}^{m} \langle A^{i}, (L_{h^{-1}})_{*} X_{g}\rangle A_{i}          \\
                        &=    \sum \limits_{i=1}^{m} \langle (L_{h})^{*} {\rm d} A^{i},  X_{g}\rangle A_{i}        \\
                        &=   \sum \limits_{i=1}^{m}  \langle \left. {\rm d}\mathcal{A}^{i}\right|_{g},  X_{g}\rangle A_{i}.
\end{array}
\]
In case of the Maurer-Cartan form we see that 2) is satisfied since left-invariant vector fields are obtained by the derivative of the right representation and satisfy $X_{e}=(L_{h^{-1}})_{*}(L_{h})_{*}X_{e}=(L_{h^{-1}})_{*} X_{h}= \omega(X_{h})$.

Finally we note that in the case of the Maurer-Cartan connection 2) is also satisfied as
for all $h,g \in H$ and all vector fields $Y$ on $P_H$ we get
\[
\begin{array}{ll}
\omega_{gh}((R_{h})_{*} Y_{g})= (L_{(gh)^{-1}})_{*} ((R_{h})_{*} Y_{g})&=
(L_{h^{-1}} \circ L_{g^{-1}})_{*} ((R_{h})_{*} Y_{g})= (L_{h^{-1}})_{*}\circ  (L_{g^{-1}})_{*} \circ (R_{h})_{*} Y_{g} \\
 &= (L_{h^{-1}} \circ R_{h})_{*} (L_{g^{-1}})_{*} Y_{g} = \textrm{Ad}(h^{-1}) \omega_{g} Y_{g},
\end{array}
\]
where we note that $R_{h}$ and $L_{g}$ commute for any pair of elements $g,h \in G$.

The horizontal subspace of $T_{g}(G)$ is defined as $\mathcal{H}_{g}= \textrm{Ker}(\omega_{g})= \textrm{span}_{i=m+1}^{n} \{A_{i}\}$. A smooth curve within $\gamma: [0,1] \to P_H$ is horizontal if all tangent vectors $c'(t)$ are horizontal, that is within $\mathcal{H}_{c(t)}$. A horizontal lift $c^{*}:[0,1]\to P_H$ of a curve $c:[0,1]\to P_H$ is a horizontal curve with $\pi(c^{*})=c$. It can be shown that a horizontal lift $c^{*}$ of $c$ is uniquely determined by $c^{*}=g_{0}$ and $\pi(g_{0})=c_{0}$ for some given point $g_{0} \in G$. From the first property (\ref{CartanEhresmann}) of the Cartan form it follows that $\mathcal{H}_{hg}=(R_{h})_{*} H_{g}$ and horizontal lifts are uniquely determined by right action of $H$ in the principal fiber bundle $P_H$, where $T_{g}(G)=\mathcal{V}_{g} \oplus \mathcal{H}_{g}$, with $\mathcal{V}_{g}$ the space of vertical tangent vectors. Consequently, the dimension of  $\mathcal{H}_{g}$ equals the dimension of $(G/H)$. \\ \\
\textbf{Example: } $G=SE(2)=\R^2 \rtimes \mathbb{T}$, $H=\mathbb{T}$ and $\omega= (L_{(0,0,e^{i\theta})^{-1}})_{*}=(L_{(0,0,e^{-i \theta})})_{*}$ we have $\mathcal{V}_{g}=\textrm{span}\{\left. \partial_{\theta} \right|_{g}\}$ and $\mathcal{H}_{g}=\textrm{span}\{\left. \partial_{\xi} \right|_{g},\left.\partial_{\eta} \right|_{g}\}$ and horizontal lifts are obtained by multiplication with $(0,0,e^{-i \theta})$ from the right. \\
 \\
\textbf{Example: } $G=SE(2)=\R^2 \rtimes \mathbb{T}$, $H=Y:=\{(0,h,e^{i0})\; |\; h \in \R\}$ and $\omega= (L_{(0,h,0)^{-1}})_{*}$, so in components the Cartan-Maurer form reads
\begin{equation} \label{Maurer}
\omega_{g}(X_{g})= \langle \left. {\rm d}\eta \right|_{g} , X_{g}\rangle \partial_{y}
\end{equation}
we have $\mathcal{V}_{g}=\textrm{span}\{\left. \partial_{\eta} \right|_{g}\}$ and $\mathcal{H}_{g}=\textrm{span}\{\left. \partial_{\xi} \right|_{g},\left.\partial_{\theta} \right|_{g}\}$ and horizontal lifts are obtained by multiplication with $(0,y,0)$ from the right.
\ \\ \\
Now that horizontal lifts are determined by the right action of $H$ on $G$ we can introduce the concept of parallel transport. To this end we will use the left-invariant vector fields as a frame of reference in $T(G)$, i.e. we use their restrictions to $g \in G$, $\left. \{\mathcal{A}_{i}\}_{i=1}^{n}\right|_{g}$ as a basis for $T_{e}(G)$ for all $g \in G$. Now the Parallel transport of a tangent vector {\small $X_{[g]}= \xi^{i} \left. \mathcal{A}_{i} \right|_{g}$} on $G/H$ along a curve $c:[0,1]\to G/H$ is
\[
\tau_{t}(X_{g})= c_{i}^{*}(t) \xi^{i}, \textrm{ where } c^{*} \textrm{ is a horizontal lift of }c \textrm{ and }c_{i}^{*}(t)=\langle {\rm d}\mathcal{A}^{i}, \frac{d c^{*}(t)}{dt}\rangle \mathcal{A}_{i}.
\]
This definition is independent on the choice of horizontal lift and $\tau_{t}$ is an isomorphism between the tangent spaces $T_{c(0)}(G/H)$ and $T_{c(t)}(G/H)$.
Now the covariant derivative of the vector field $Y$ on $G/H$ along the curve $c:[0,1]\to G/H$ in the point $[g]=c(0)$ is defined as
\[
\nabla_{X_[g]} Y = \lim \limits_{h \downarrow 0} \frac{1}{h}((\tau^{h})^{-1}Y_{c(h)}-\tau^{0}Y_{c(0)})=\lim \limits_{h \downarrow 0} \frac{1}{h}((\tau^{h})^{-1}Y_{c(h)}-Y_{[g]}),
\]
with $X_{[g]}=\frac{dc}{dt}(0)$. The vector field $Y$ is called parallel along the curve $c$ if $\nabla_{X_{c(t)}}Y=0$, with $X_{c(t)}=\dot{c}(t)=\frac{dc}{dt}(t)$, $t\in [0,1)$.
A curve which is covariantly constant, i.e.
\begin{equation}
\nabla_{\dot{c}} \dot{c}=0
\end{equation}
is called an auto-parallel. In Riemannian geometry, such curves coincide with geodesics, i.e. the unique smooth curves $c$, with $c(a)=p, c(b)=q$ which minimize
\begin{equation}
\int_{a}^{b} \sqrt{g_{ij}(\ul{x}(c(t))) \dot{x}^{i}(t) \dot{x}^{j}(t)} {\rm dt}.
\end{equation}
However, due to the torsion in the Cartan connection auto-parallels and geodesics no longer coincide. Moreover, the meaning of geodesics as paths with minimal arc-length can not be used as we did not consider a connection induced by a Riemannian metric (yet). 
%

If we want to express the connection in a (possibly degenerate) metric we note that this metric must be left-invariant because of 2) in (\ref{CartanEhresmann}). Moreover, because of 1) in (\ref{CartanEhresmann}) this metric must be invariant under the adjoint representation $\textrm{Ad}: G \to T_{e}(G)$ of $G$ on its Lie-algebra. This brings us to the (possibly degenerate) left-invariant metric $m_{K}$ induced by the Killing form $K$:
\begin{equation} \label{leftinvmetric}
\begin{array}{l}
m_{K}(\mathcal{A}_g,\mathcal{B}_g)=K(\mathcal{A}_e, \mathcal{B}_e)\ , \textrm{ for all }\mathcal{A}, \mathcal{B} \in \mathcal{L}(G) \\
K(A,B)= \textrm{trace}\,(\textrm{ad}\, A \circ \textrm{ad}\, B),
\end{array}
\end{equation}
where the adjoint representation $\textrm{ad}:T_{e}(G) \to \mathcal{B}(T_{e}(G))$ of the Lie-algebra on itself is given by $(\textrm{ad}(A))(B)=[-A,B]$, which is the derivative of the $Ad$ representation mentioned before.
For the moment we will assume that the killing form is non-degenerate (which is not the case for $G=SE(2)$).
The matrix elements with respect to our moving frame of reference {\small $\left. \tilde{A}_{i}\right|_{g}=\left. {\rm d}\mathcal{R}(A_{i})\right|_{g}$} the components of $m_{K}$
are given by 
\begin{equation} \label{compmetric}
\begin{array}{l}
g_{ij}  \equiv \textrm{trace}((\textrm{ad}\,  \mathcal{A}_{i} \circ \textrm{ad}\, \mathcal{A}_j))=\langle {\rm d}\mathcal{A}^{k} ,  c_{ik}^{l} c^{q}_{il} \mathcal{A}_{q} \rangle=
\, c_{ik}^{l} c^{k}_{il},
\end{array}
\end{equation}
where the structure constants $c_{ik}^{l}$ are 
defined by $[A_{i} , A_j]=c^{k}_{ij} A_{k}$.

\subsection{Vector bundles \label{ch:vectorbundles}}

If we consider the trivial case $H=e$, in which case we have $G/H \equiv G$, $\omega=0$ and thereby every tangent vector is horizontal, so it does not make a lot of sense to consider a principal fiber bundle $P_{H}$ with structure group $H$. In such a situation one rather considers the action of the group $G$ onto itself. The Cartan form on for example $SE(2)$ would now be given by
\[
\omega_{g}(X_{g})=\langle \left. {\rm d}\theta \right|_{g}, X_{g} \rangle \partial_{\theta}+
\langle \left. {\rm d}\xi \right|_{g}, X_{g} \rangle \partial_{x} +
\langle \left. {\rm d}\eta \right|_{g}, X_{g} \rangle \partial_{y},
\]
which corresponds to (\ref{MaurerCartan}) the connection in case $H=G$, but now defined on $G$ rather than
$G/G \equiv \{e\}$. This means that we shall consider the vector bundle $E=(G,T(G))$, which we shall consider next.

Let $t \mapsto c(t)$ be a smooth curve in $G$ with $X(t)=c'^{i}(t)\mathcal{A}_{i}$. Let $\mu : G \to E$ be a section in $E$. Let $\{\mu_{k}\}$ be the sections in $E$ aligned with the left-invariant vector fields $\{\mathcal{A}_{k}\}$ on $G$. Then
the Cartan connection in components reads
\begin{equation}\label{Cartancomp}
(D \mu)X(t):= D_{X(t)} \mu(c(t))=  \dot{a}^{k}(c(t)) \mu_{k}(c(t)) + \dot{c}^{i} a^{k}(t) \Gamma^{j}_{ik}(c(t)) \mu_{j}(c(t)),
\end{equation}
for all sections $\mu(c(t))= a^{k}(c(t)) \mu_{k}(c(t))$ and
where $D_{\mathcal{A}_{i}} \mu_{j}=\Gamma^{k}_{ij} \mu_{k}$. The Cartan connection $D$ is given by $D= {\rm d}+\omega$ since
\[
D a^{k} \mu_{k} = {\rm d}(a^{k}) \mu_{k} +a^{k} \omega(\mu_k),
\]
where we note that by the chain rule we have
\begin{equation} \label{dpart}
({\rm d} (a^{k})\mu_{k})X(t)= \frac{\partial a^{k}}{\partial \xi_k} \mu_{k}(c(t)) \dot{c}^{l}(t) \delta^{k}_{l}=  \mu_{k}(c(t)) \frac{\partial a^{k}}{\partial \xi_k} \dot{c}^{k}(t)= \mu_{k}(c(t)) \dot{a}^{k}(t),
\end{equation}
where we used short notation $a(t)=a(c(t))$ and moreover
\begin{equation}\label{wpart}
(a^{k}\omega \mu_{j})(\dot{c}^{i}\mathcal{A}_{i})=a^{k}(t) \mu_{j}(c(t))\omega(\mathcal{A}_{i})= \dot{c}^{i}(t) a^{k}(t) \mu_{j}(c(t)) \Gamma^{j}_{ik}(c(t))
\end{equation}
and the terms in (\ref{dpart}) and (\ref{wpart}) indeed add up to the right hand side of (\ref{Cartancomp}). Now the Cartan connection can be split up in a symmetric and anti-symmetric part:
\[
\Gamma^{i}_{kl}=\overline{\Gamma}^{i}_{kl}+K^{i}_{kl}, \ \ \overline{\Gamma}^{i}_{kl}=\frac{1}{2}(\Gamma^{i}_{kl}+\Gamma^{i}_{lk}),
K^{i}_{kl}=\frac{1}{2}(\Gamma^{i}_{kl}-\Gamma^{i}_{lk})
\] 
where $\overline{\Gamma}^{i}_{kl}$ are the components of a Levy-Cevita connection induced by the metric $m_{K}$ given by (\ref{leftinvmetric}). Now by left invariance we have with respect to our frame of reference that \\
\mbox{$\overline{\Gamma}^{i}_{kl}=\frac{1}{2}g^{im}(g_{mk,l}+g_{ml,k}-g_{kl,m})=0$  }, where we use the convention: indices after the comma to index left-invariant differentiation (i.e. $g_{mk,l}=\mathcal{A}_{l}g_{mk}$).
So in our case the components of the Cartan tensor coincide with the components of the contorsion tensor $K$:
\begin{equation} \label{christcartan}
\Gamma^{i}_{kl}=K^{i}_{kl}=\frac{1}{2} g^{im}(c_{mkl}+c_{mlk}-c_{klm})
\end{equation}
where $c_{mkl}=g_{mp} c^{p}_{kl}$, where $c^{p}_{kl}$ are the structure constants of the Lie-algebra 
\[
[\mathcal{A}_{k},\mathcal{A}_{l}]=c^{p}_{k,l} \mathcal{A}_{p}, p,k,l=1,\ldots,n,
\]
where the components of the left-invariant metric tensor are given by
\begin{equation}\label{compleftinvmetrictensor}
g_{ij}(g)=m_{K}(\left.\mathcal{A}_{i}\right|_{g},\left.\mathcal{A}_{j}\right|_{g})=\frac{1}{I_{\textrm{ad}}} c_{jk}^{l} c^{k}_{il} =\frac{1}{2} c_{jk}^{l} c^{k}_{il}.
\end{equation}
where the Dynkin index $I_{\textrm{ad}}$ of the adjoint representation, coincides with the dual Coxeter number \cite{DiFrancesco}, which in case of
$SO(3)\equiv SU(2) \equiv A_{1}$ equals $2$, \cite{Hiller}.

Now the curvature of the Cartan connection $D$ is
\begin{equation} \label{curvatureofCartan}
D^2=D \circ D = ({\rm d}+\omega) \circ ({\rm d}+\omega) = {\rm d}\omega + \omega \wedge \omega,
\end{equation}
where we note that since $\omega$ is a one-form we have $D^2\mu ={\rm d}\omega \mu -\omega {\rm d}\mu +\omega {\rm d}\mu + \omega \wedge \omega \mu$ for all section $\mu$ in $E$. In components
(for details see \cite{Jost}{p.111-112}) (\ref{curvatureofCartan}) reads
\[
D^{2}=\frac{1}{2}(w_{j\, , i}- w_{i,\, , j}+[w_i,w_{j}]) {\rm d}\xi^{i} \wedge {\rm d}\xi^{j}
\]
This provides the following formula for the curvature tensor:
\[
\begin{array}{ll}
R^{j}_{i,kl} &= \Gamma^{j}_{li,\ k} -  \Gamma^{j}_{ki,\ l} + \sum_{\lambda=1}^{n} \Gamma^{j}_{k\lambda}\Gamma^{\lambda}_{li}-\Gamma^{j}_{l\lambda}\Gamma^{\lambda}_{ki}\\
 &= \sum_{\lambda=1}^{n} \Gamma^{j}_{k\lambda}\Gamma^{\lambda}_{li}-\Gamma^{j}_{l\lambda}\Gamma^{\lambda}_{ki},
\end{array}
\]
which after some computation using the symmetries of the curvature tensor gives
\begin{equation}\label{Cartan2}
R^{j}_{i,kl} = \frac{1}{4}  c_{\lambda i}^{j} c^{\lambda}_{kl},
\end{equation}
which is the formula by \cite{Cartan}{p.187}.

We would like to apply (\ref{Cartan2}) and (\ref{compmetric}) to the case of the Euclidean motion group. But here a problem arises as the metric $m_{K}$
induced by the killingform $K$ on $SE(2)$ is degenerate. Therefore we embed,
by means of a small parameter $\beta$, the Lie-algebra of $SE(2)$ spanned by $\{\partial_{\theta},\partial_{\xi},\partial_{\eta}\}$ into the Lie-algebra of $SO(3)$ which is $so(3)=\{X \in \R^{3\times 3}\;|\; X^{T}=-X\}$ whose Killing form is non-degenerate. We simply obtain the appropriate sectional curvatures and covariant derivatives by taking the limit $\beta \to 0$ afterwards.

\begin{figure}
\centerline{\includegraphics[width=0.75\hsize]{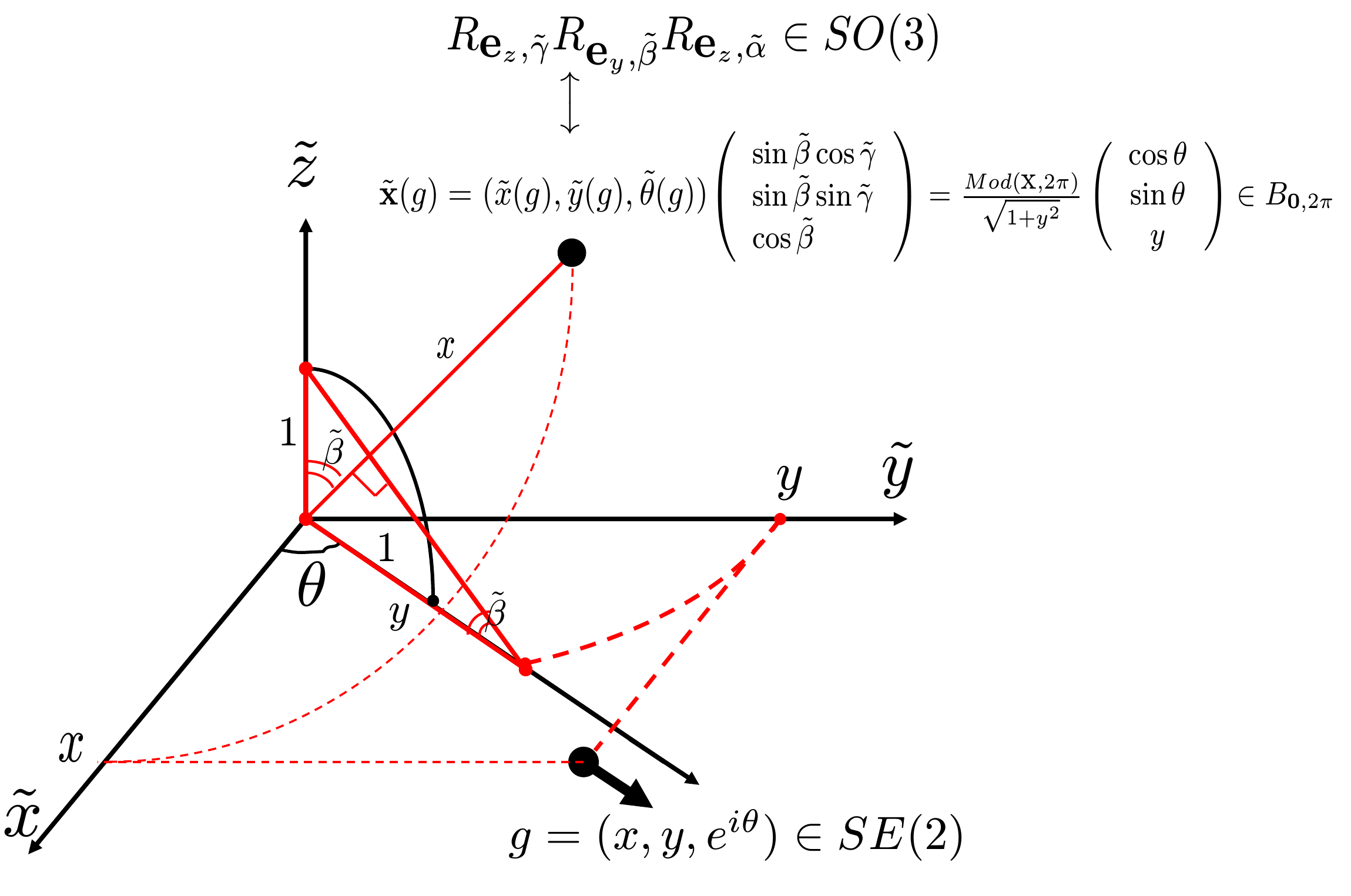}}
\caption{The embedding of $SE(2)$ in $SO(3)$ by (\ref{trafo}) for $\beta=1$. The group $SO(3)$ can be identified with a unit ball $B_{\ul{0},2\pi}$ with radius $2\pi$ by means of the Euler angle parametrization. 
Here all points on the sphere $\{\tilde{\ul{x}} \in \R^3\; | \; \|\tilde{\ul{x}}\|=\tilde{\alpha}=2\pi\}$ are identified with the origin. Now given a fixed member of $(x,y,e^{i\theta}) \in SE(2)$ we can obtain the corresponding element in $SO(3)$ as follows. First consider the point $(\tilde{x},\tilde{y},\tilde{z})=(x,y,0)$ with attached direction $\theta=\angle ((\tilde{x},\tilde{y}),(1,0))$ and construct the unique half-line $\ell$ trough the origin with direction $\theta$. Then project $(x,y,0)$ on the $\tilde{y}$-axis and rotate the point so that it ends up at point $P$ at the line $\ell$. Then find the unique point on the unit sphere $Q$ such that $\vec{OQ} \bot \vec{PN}$ where $N$ is the north-pole. Finally, scale $Q$ with $x$ modulo $2\pi$.
}\label{fig:trafo}
\end{figure}

The Euler angle parametrization of $SO(3)$ is given by $R_{\ul{e}_{z},\tilde{\gamma}} R_{\ul{e}_{y},\tilde{\beta}}
R_{\ul{e}_{z},\tilde{\alpha}}$. Then a basis of left-invariant vector fields on $SO(3)$ (in Euler-angles) is given by
\begin{equation} \label{VFSO3}
\begin{array}{l}
\mathcal{B}_{1}= \cot \tilde{\beta} \cos \tilde{\gamma} \, \partial_{\tilde{\gamma}} -
\frac{\cos \tilde{\gamma}}{\sin \tilde{\beta}} \, \partial_{\tilde{\alpha}} +\sin \tilde{\gamma} \, \partial_{\tilde{\beta}}  \\
\mathcal{B}_{2}= -\cot \tilde{\beta} \, \sin \tilde{\gamma} \, \partial_{\tilde{\gamma}} -
\frac{\cos \tilde{\gamma}}{\sin \tilde{\beta}}\, \partial_{\tilde{\alpha}} +\sin \tilde{\gamma} \, \partial_{\tilde{\beta}} \\
\mathcal{B}_{3}= \partial_{\tilde{\gamma}},
\end{array}
\end{equation}
where we note that $[\mathcal{B}_{1},\mathcal{B}_{2}]=\mathcal{B}_{3}, [\mathcal{B}_{2},\mathcal{B}_{3}]=\mathcal{B}_{1}$, $[\mathcal{B}_{3},\mathcal{B}_{1}]=\mathcal{B}_{2}$. Now apply the coordinate transformation
\begin{equation} \label{trafo}
\begin{array}{l}
\tilde{\alpha}=\beta x \\
\tilde{\beta}=\frac{\pi}{2}-\arctan(\beta y)\\
 \tilde{\gamma}=\theta
\end{array}
\end{equation}
and multiply $\mathcal{B}_{1}$ and $\mathcal{B}_{2}$ with $\beta$ then we obtain the following vector fields
\begin{equation} \label{VFmod}
\begin{array}{ll}
\mathcal{A}_{2}^{\beta} &:=\beta \mathcal{B}_{1} =- \beta^2 y \cos \theta \, \partial_{\theta} + \cos \theta \, \sqrt{1+\beta^2 \, y^2}\, \partial_{x}+ \sin \theta \,(1+\beta^2 \, y^2) \, \partial_{y}  \\
\mathcal{A}_{3}^{\beta} &:= \beta \mathcal{B}_{2}= \beta^2 y \sin \theta \, \partial_{\theta} - \sin \theta \, \sqrt{1+\beta^2 \,y^2}\, \partial_{x}+ \cos \theta \, (1+\beta^2 \, y^2) \,\partial_{y}  \\
\mathcal{A}_{1}^{\beta} &:=\mathcal{B}_{3}=\partial_{\theta}.
\end{array}
\end{equation}
These vector fields again form a three dimensional Lie algebra:
\[
[\mathcal{A}_{1}^{\beta},\mathcal{A}_{2}^{\beta}]=\mathcal{A}_{3}^{\beta}, \ \ \ [\mathcal{A}_{1}^{\beta},\mathcal{A}_{3}^{\beta}]=-\mathcal{A}_{2}^{\beta}, \ \ \
[\mathcal{A}_{2}^{\beta},\mathcal{A}_{3}^{\beta}]=\beta^{2} \; \mathcal{A}_{1}^{\beta}.
\]
and which converges to $\{\mathcal{A}_{1},\mathcal{A}_{2},\mathcal{A}_{3}\}=\{\partial_{\theta},\partial_{\xi}, \partial_{\eta}\}$ for $\beta
\to 0$.  To get a geometrical understanding of the embedding of $SE(2)$ into $SE(3)$ by means of the coordinate transformation (\ref{trafo}) see Figure \ref{fig:trafo}.
The components of the left-invariant metric tensor are given by
\begin{equation}
g_{ij}(g)=m_{K}(\left.\mathcal{A}_{i}\right|_{g},\left.\mathcal{A}_{j}\right|_{g})= \frac{1}{2}
\textrm{trace}({\rm ad}(A_{i})\circ {\rm ad}(A_{j}))= \frac{1}{2} c_{jk}^{l} c^{k}_{il},
\end{equation}
where $c_{jk}^{l}$ denote the structure constants of the perturbed Lie algebra, which are related to the structure constants of the Lie algebra $so(3)$ of $SO(3)$ properly scaled by $\beta$.
As a result the components of the left-invariant Killing-form-metric with respect to this perturbed Lie-algebra, recall (\ref{compleftinvmetrictensor}), are given by\footnote{Note that the physical dimension of $\epsilon $ is $[length]^{-1}$.}
\begin{equation} \label{metric}
G=[g_{ij}]= \left(\begin{array}{ccc}
1 & 0 & 0 \\
0 & \beta^{2} & 0 \\
0 & 0 &  \beta^{2}
\end{array}\right)
\end{equation}
and the non-zero components of the constant Riemann-curvature tensor, recall (\ref{Cartan2}), are now given by
{\small
\[
\begin{array}{l}
R^{1}_{212}=-R^{1}_{221}=-R^{1}_{331}=R^{1}_{313}=R^{2}_{332}=-R^{2}_{323}=-R^{3}_{232}=R^{3}_{223}=\frac{\beta^{2}}{4} \to 0 \textrm{ as }\beta \to 0 \\
R^{2}_{112}=-R^{2}_{121}=R^{3}_{113}=-R^{3}_{131}=\frac{1}{4}
\end{array}
\]
}
From which we deduce that in the limiting case $\beta \to 0$ the non-normalized sectional curvature $|\partial_{\theta} \wedge \partial_{\xi}|K(\partial_{\theta} \wedge \partial_{\xi})=(R(\partial_{\theta},\partial_{\xi})\partial_{\xi},\partial_{\theta})=R_{1212}=g_{1m}R^{m}_{212}$ of the plane spanned by $\{\partial_{\xi},\partial_{\theta}\}$  is constant $=\frac{1}{4}$. Similarly the non-normalized curvature of the plane spanned by $\{\partial_{\eta},\partial_{\theta}\}$ is constant $=\frac{1}{4}$ in contrast to the spatial plane spanned by $\{\partial_{\xi},\partial_{\eta}\}$, which is of course flat. This explains the presence of curvature in orientation scores. Moreover, the curvedness of the Cartan connection on the space $SE(2)$ is important for application of differential geometrical operators on orientation scores.
The covariant derivative of a vector field $\ul{v}$ on $SE(2)$ is a (1,1)-tensor field whose components are given by
\[
\nabla_{j'}v^{i}= \partial_{j'}v^{i} + \Gamma^{i}_{j'k'} v^{k'},
\]
The covariant derivative of a co-vector field $\ul{a}$ on $SE(2)$ is (0,2)-tensor field with components:
\[
\nabla_{j} a_i = \partial_{j} a_i - \Gamma^{k}_{ji} a_k.
\]
In particular for the gradient ${\rm d} U= \sum \limits_{i=1}^{3} U_{,i} {\rm d }\mathcal{A}^{i}$ of an orientation score $U:SE(2) \to \mathbb{C}$
and the corresponding vector field
\[
\mathcal{G}^{-1}{\rm d U}=\sum \limits_{i=1}^{3} U_{,i} \mathcal{A}_{i} =\frac{\partial U}{\partial \theta} \partial_{\theta}+\frac{\partial U}{\partial \xi} \partial_{\xi}+
\frac{\partial U}{\partial \eta} \partial_{\eta}
\]
This yields the following covariant second order derivations\footnote{Note that a rescaling of $\theta$, say $\tilde{\theta} \to \lambda \theta$ directly results in a rescaling of the covariant derivatives: $\nabla_{\theta}\nabla_{\xi} \to  \lambda \nabla_{\theta}\nabla_{\xi}$, $\nabla_{\theta}\nabla_{\eta} \to  \lambda \nabla_{\theta}\nabla_{\eta}$, $\nabla_{\theta}\nabla_{\theta} \to  \lambda^{2} \nabla_{\theta}\nabla_{\theta}$ . } on orientation scores $U$ (in the limiting case $\varepsilon \to 0$):
{\small
\begin{equation}\label{covH}
\begin{array}{ll}
[\nabla_{i} \nabla_{j} U] &=[\nabla_{i} U_{,j}]=[\partial_{i} U_{,j} -\Gamma^{\lambda}_{ij}\, U_{,\lambda}]=
\left(
\begin{array}{ccc}
\partial_{\theta} \partial_{\theta }U &  \partial_{\xi} \partial_{\theta }U + \frac{1}{2}\partial_{\eta}U & \partial_{\eta} \partial_{\theta} U - \frac{1}{2}\partial_{\xi}U  \\
\partial_{\theta} \partial_{\xi }U - \frac{1}{2}\partial_{\eta} & U_{\xi \xi} & U_{\xi \eta}  \\
 \partial_{\theta} \partial_{\eta} U + \frac{1}{2}\partial_{\xi}& U_{\eta \xi} & U_{\eta \eta}
\end{array}
\right) \\
 &=
\left(
\begin{array}{ccc}
U_{\theta\theta} & \frac{1}{2}(U_{\theta \xi} + U_{\xi \theta}) & \frac{1}{2}(U_{\theta \eta} + U_{\eta \theta}) \\
\frac{1}{2}(U_{\xi \theta} +U_{\theta \xi}) & U_{\xi \xi} & U_{\xi \eta}  \\
\frac{1}{2}(U_{\eta \theta} + U_{\theta \eta}) & U_{\eta \xi} & U_{\eta \eta}
\end{array}
\right).
\end{array}
\end{equation}
}
with $U_{ij}=\partial_{j} (\partial_{i} U)$, $i,j=1,\ldots,3$,
where we recall (\ref{christcartan}) and $\partial_{\eta}=[\partial_{\theta},\partial_{\xi}]$ and $\partial_{\xi}=-[\partial_{\theta},\partial_{\eta}]$. \\
 \\
So by anti-symmetry of the Christoffel symbols $\Gamma^{i}_{jk}=-\Gamma^{i}_{kj}$ we have \begin{equation}
\nabla_{ij}U+ \nabla_{ji} U= \mathcal{A}_{i}\mathcal{A}_{j}U+ \mathcal{A}_{j}\mathcal{A}_{i}U.
\end{equation}
and therefore
our linear \emph{and non-linear} (that is the conductivity $(D_{ij}(W))(g,s)$ explicitly depends on $W$) evolution equations on $SE(2)$ can be straightforwardly expressed in covariant derivatives:
\[
\left\{
\begin{array}{ll}
\partial_{s} W(g,s) &= \sum \limits_{i=1}^{3} \mathcal{A}_{i} (D_{ij}(W))(g,s) \mathcal{A}_{j} W= \sum \limits_{i=1}^{3} \nabla_{i} (D_{ij}(W))(g,s) \nabla_{j} W \\
W(g,0)&=|U_{f}(g,s)|
\end{array}
\right.
\]
where $|U_{f}(g,s)|$ denotes the absolute value of the orientation score $U_f \in \mathbb{L}_{2}(SE(2))$ of image $f \in \mathbb{L}_{2}(\R^2)$
where we note that $\sum_{ij} \nabla_{i}(D_{ij} \nabla_{j} U)= \sum_{i,j}D_{ij}\nabla_{i} \nabla_j U + \sum_{ij}(\nabla_{i}D_{ij}) \nabla_{j}U= \sum_{i,j}D_{ij}\mathcal{A}_{i} \mathcal{A}_j U + \sum_{ij}(\mathcal{A}_{i}D_{ij}) \mathcal{A}_{j}U$.

\subsection{Auto-parallels}

Recall from Riemannian differential geometry that curves which are covariantly constant (or auto parallel) that is
\[
\nabla_{\dot{\ul{x}}} \dot{\ul{x}}= \ul{0} \desda  \ddot{x}^{i} = -\Gamma_{kl}^{i} \dot{x}^{k} \dot{x}^{l},
\]
where $\ul{x}=x^{i}\partial_{i}$ and where the well-known Christoffel symbols read $\Gamma^{i}_{kl}=
\frac{1}{2}g^{im}(g_{mk,l}+g_{ml,k}-g_{kl,m})$ coincide with the path-length minimizers, i.e. geodesics, if and only if the connection is torsion free.

The Cartan connection, however, is not torsion free. Therefore the auto-parallels on $SE(2)$ do not coincide with the geodesics. In fact the auto-parallels coincide with the exponential curves. To this end we note that the Christoffel symbols (\ref{christcartan}) with respect to our basis of left-invariant vector fields are anti-symmetric and as a result we have for auto parallel curves
\[
\frac{d}{ds} \langle {\rm d}\mathcal{A}^{k}, c'(s) \rangle = \Gamma^{k}_{ij} \langle {\rm d}\mathcal{A}^{i}, c'(s) \rangle \langle {\rm d}\mathcal{A}^{j}, c'(s) \rangle=0
\]
and thereby we have that
\[
\langle {\rm d}\mathcal{A}^{k}, c'(s) \rangle = \langle {\rm d} \mathcal{A}^{k},c'(0)\rangle=\textrm{constant}= c^{k} \textrm{ for }k=1,\ldots, n,
\]
so $c'(s)=c^{k} \left. \mathcal{A}_{k} \right|_{c(s)}$, for some constants $c_{k} \in \R$. Now these curves
exactly coincide with the exponential curves $c(s)=\textrm{exp}(s c^{k} \left. \mathcal{A}_{k} \right|_{c(0)})c(0) $ within the Lie group.  \\
\\
\textbf{Example :} \\
Recall that the left-invariant vector fields in $SE(2)$ were given by (\ref{leftinvSE2}) as a result
auto-parallels $\gamma$ in the Euclidean motion group  are given by the following set of equations
\[
\frac{d}{ds} \gamma(s)= c^{k} \left. \mathcal{A}_{k} \right|_{\gamma(s)} \gamma(s), k=1,2,3,
\]
or explicitly in \emph{$(x,y,\theta)$-coordinates} (not to be mistaken with the $(\xi,\eta,\theta)$-coordinates)
\[
\left\{
\begin{array}{l}
\dot{\gamma_{3}}(s)=1 \\
\dot{\gamma_{1}}(s)= c^{2} \cos \gamma_{1}(s) - c^{3} \sin \gamma_{1}(s) \qquad \gamma(0)=g_{0}\\
\dot{\gamma_{2}}(s)= c^{2} \cos \gamma_{1}(s) + c^{3} \sin \gamma_{1}(s)
\end{array}
\right.
\]
the unique solution of which is given by
\begin{equation} \label{circspiral}
\begin{array}{l}
\gamma(t)=\exp(t(\sum \limits_{i=1}^{3}c^{i} \mathcal{A}_{i})) g_{0}=
(x_{0}+ \frac{c^3}{c^1}(\cos (c^1 t+\theta_{0})-\cos \theta_{0})
+ \frac{c^2}{c^{1}}(\sin(c_1 t\!+\!\theta_{0})-\sin \theta_{0}), \\ y_{0}\!+\! \frac{c^3}{c^1}(\sin (c^1 t\!+\!\theta_{0})-\sin \theta_{0})\!-\! \frac{c^2}{c^{1}}(\cos (c^1 t\!+\!\theta_{0})\!-\!\cos \theta_{0}),e^{i(c^1 t \!+\!\theta_{0})}),
\end{array}
\end{equation}
for $c^{1}\neq 0$, which is a circular spiral with radius $\frac{\sqrt{(c^2)^2+(c^3)^2}}{c_{1}}$ and central point
\[
(-\frac{c^3}{c^{1}}\cos \theta_{0}-\frac{c^2}{c^{1}}\sin \theta_{0}+x_{0}, \frac{c^2}{c^{1}}\cos \theta_{0}-\frac{c^3}{c^{1}}\sin \theta_{0}+y_{0}).
\]
This result is easily deduced by the method of characteristics for first order PDE's. For $c^{1}=0$ we get a straight line in the plane $\theta=\theta_{0}$:
\[
\gamma(t)= (x_{0}+t \, c^{2} \cos \theta_{0}- t \, c^{3}\sin \theta_{0}, y_{0}+t \, c^{2} \sin \theta_{0}+ t \, c^{3}\cos \theta_{0} , e^{i\theta_{0}}),
\]
which coincides with (\ref{circspiral}) by taking the limit $c^{1}\to 0$.

\subsection{Fiber bundles and the concept of horizontal curves \label{ch:fiberbundles}}

Recall Definition \ref{def:horizontal}, where we provide the definition of a horizontal curve in $SE(2)$. Next we will set this definition in a differential geometrical context (which justifies the word horizontal).

In case $G=SE(2)=\R^2 \rtimes \mathbb{T}$, $Y=\{(0,h,e^{i0})\; |\; h \in \R\}$ and $\omega= (L_{(0,y,0)^{-1}})_{*}$  we have $\mathcal{V}_{g}=\textrm{span}\{\left. \partial_{\eta} \right|_{g}\}$ and $\mathcal{H}_{g}=\textrm{span}\{\left. \partial_{\xi} \right|_{g},\left.\partial_{\theta} \right|_{g}\}$ and horizontal lifts are obtained by multiplication with $(0,h,0)$ from the right
, where we note that
\[
g(0,h,0)=(x,y,e^{i\theta})(0,h,0)=(x -h \, \sin \theta ,y+ h \, \cos \theta ,e^{i\theta}) =g+ h \ul{e}_{\eta}.
\]
This particular choice is important in image analysis as it is the differential geometrical description of ``lifting'' of curves. That is to each smooth planar curve in $C(\R^{+},\R^2)$ we can create a curve in $C(\R^{+},SE(2))$ by setting
\begin{equation} \label{lifting}
\begin{array}{l}
s \mapsto \ul{x}(s) \in C(\R^{+},\R^2) \qquad \leftrightarrow \\
 s \mapsto (\ul{x}(s),e^{i\theta(s)}) \in C(SE(2),\R) \textrm{ with }\theta(s)= \textrm{arg}(x'(s)+iy'(s))= \angle \ul{x}'(s).
\end{array}
\end{equation}
\emph{Convention:
By $\angle \ul{x}'(s)$ we mean the angle that $\ul{x}'(s)$ makes with the fixed $\ul{e}_{x}$-axis, i.e.
\[
\theta(s)= \angle \ul{x}'(s) \desda \cos \theta(s) = \frac{\ul{x}'(s) \cdot \ul{e}_{x}}{\|\ul{x}'(s)\|}.
\]
}
In this setting such a curve is horizontal, since its tangent field is a horizontal vector field. Notice that right multiplication with a fixed element $(0,h,0)$ provides a horizontal lift of the curve:
\[
\begin{array}{l}
g(s) \textrm{ is horizontal }\Rightarrow g (0,h,0)(s) =g(s)(0,h,0) \textrm{ is horizontal }. \\
\pi(g(s))=\pi(g(s)(0,h,0)), \textrm{ for all }s>0.
\end{array}
\]
where we note
$\angle \ul{b}'(t)=\theta(t) \Rightarrow
\angle (h \theta'(t)(
-\cos \theta(t),-\sin \theta(t))+\ul{b}'(t)
)=\theta(t)$.
We note that in general right multiplication of a horizontal curve with a constant element in $SE(2)$ does not preserve the horizontality. Actually this only works for the subgroup $Y$. This is in contrast with left multiplication with a fixed element:
\[
g(s) \textrm{ is horizontal }\Rightarrow (\tilde{g} \, g)(s)=\tilde{g} \, g(s) \textrm{ is horizontal } \textrm{ for all }\tilde{g} \in SE(2).
\]
The setting in this example is important to relate elastica curves to geodesics which minimize
\begin{equation}\label{distanceonfiber}
\begin{array}{ll}
d_{SE(2)}(g,g_0) &=
\inf \left\{ \int_{0}^{1} \sqrt{ (\theta'(t))^{2} + \epsilon \|\ul{x}'(t)\|^2 }\; {\rm d}t \; |\;  \right. \\
 &\left. \gamma \textrm{ is a smooth horizontal curve connecting }g \textrm{ and }g_{0}, \gamma(0)=g_0, \gamma(1)=g\right\},
\end{array}
\end{equation}
since only for horizontal curves we have $\kappa(s)=\theta'(s)$, so we may write
\[
d_{SE(2)}(g,g_0)= \inf \{ \int_{0}^{L} \sqrt{ (\kappa(s))^{2} + \epsilon }\; {\rm d}s \; |\; \gamma(0)=g_0, \gamma(L)=g, \gamma \textrm{ smooth and horizontal }\}.
\]
The auto-parallels in the fiber bundle are now given by
\begin{equation}\label{hor}
\gamma(t)=\exp(t(\sum \limits_{i=1}^{2}c^{i} \mathcal{A}_{i})) g_{0}=
(x_{0}+ \frac{c^2}{c^{1}}(\sin(c_1 t\!+\!\theta_{0})-\sin \theta_{0}), \\ y_{0}\- \frac{c^2}{c^{1}}(\cos (c^1 t\!+\!\theta_{0})\!-\!\cos \theta_{0}),e^{i(c^1 t \!+\!\theta_{0})}),
\end{equation}
where we note that they follow from the vector bundle case by omitting the vertical direction $\partial_{\eta}$, so we get them from (\ref{circspiral}) by setting $c_{3}=0$. Note that these auto-parallels are indeed horizontal as we have
\[
\textrm{arg}\,\{ \gamma_{1}'(t) + i\gamma_{2}'(t)\} = \textrm{arg}\, \{c_{2}\cos(c_{1}t + \theta_{0}) + i \, \sin(c_{1}t + \theta_{0})\}= c_{1}t+\theta_{0}.
\]
Also the covariant derivatives are again blind for the vertical direction $\partial_{\eta}$ and they are given by
\[
[\nabla_{j} \nabla_{i}U]=
\left(
\begin{array}{cc}
\partial_{\theta} \partial_{\theta} U & \partial_{\xi} \partial_{\theta} U  \\
\partial_{\theta} \partial_{\xi} U & \partial_{\xi} \partial_{\xi} U
\end{array}
\right)= \left(
\begin{array}{cc}
\partial_{\theta} \partial_{\theta} U & \partial_{\theta}  \partial_{\xi} U  \\
 \partial_{\xi} \partial_{\theta} U & \partial_{\xi} \partial_{\xi} U
\end{array}
\right)
\]
for horizontal gradients ${\rm d} U= \partial_{\theta} U {\rm d \theta} + \partial_{\xi} U {\rm d}\xi$ , i.e. $\partial_{\eta}U=0$.

\subsubsection{Horizontality and the extraction of spatial curvature from orientation scores \label{ch:curvest}}

Orientation scores $U$ and their absolute value $|U|=|U_f|=\sqrt{(\Re(U_f))^{2}+(\Im(U_f))^{2}}$ in general do not satisfy $\partial_{\eta}U=0$, $\partial_{\eta}|U|$. Nevertheless, in our linear and non-linear diffusion schemes (in section \ref{ch:leftinvimageproc} and in section \ref{ch:CED}) on orientation scores, we include the direction $\partial_{\xi} + \kappa \partial_{\theta}$, where $\kappa$ equals the horizontal curvature (i.e. the spatial curvature of projected curves on the spatial plane). 

\begin{figure}
\centerline{\includegraphics[width=0.7\hsize]{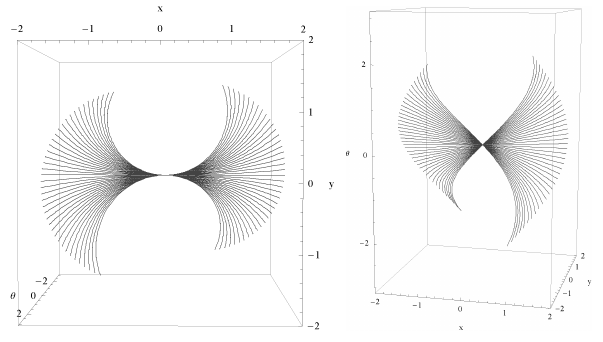}}
\caption{ All horizontal exponential curves through a fixed point $g \in SE(2)$ for different curvature values, shown from 2 different perspectives. The left-sided image shows that these curves correspond to circular arcs if projected onto te spatial plane.}
\label{fig:horizontal}
\end{figure}

Let $U:SE(2) \to \R^{+}$ be some positive smooth function on $SE(2)$. This could for example be the absolute value of a (processed) orientation score of an image, which is positive and phase invariant see Figure \ref{fig:cakeexample} (d). \\
Then by embedding $SE(2)$ into $\R^3$, the exponential curves through $g_{0}$ with direction $c^{i}\left. \mathcal{A}_{i} U \right|_{g_{0}}$ form ``tangent spirals'' to the orientation score $U:SE(2) \to \R^{+}$. In particular, horizontal exponential curves are $\textrm{exp}(s( \kappa \partial_{\theta}+\partial_{\xi})) g_{0}$ and given by (\ref{hor}). See Figure \ref{fig:horizontal}.
In this section we will obtain fast algorithms for curvature estimation at position, say $g_{0}\in SE(2)$, in the domain of $U$, by finding the tangent spiral (exponential curve) through $g_0$ that fits $U$ in an optimal way. For the exact definition of such an optimally fitting (horizontal) tangent spiral we first need a few preliminaries.

We introduce the following (left-invariant) first fundamental form on $T(SE(2))\times T(SE(2)$
\begin{equation} \label{fff}
g_{ij} {\rm d}\mathcal{A}^{i} \otimes {\rm d}\mathcal{A}^{j}= {\rm d}\theta \otimes {\rm d}\theta + \beta^{2}{\rm d}\xi \otimes {\rm d}\xi + \beta^2{\rm d}\eta \otimes {\rm d}\eta
\end{equation}
where $\beta>0$, $g_{ij}$ is the diagonal matrix with $\{1,\beta^2,\beta^2\}$ as respective diagonal elements. To this end we recall recall (\ref{metric}), where we note that (\ref{fff}) and (\ref{metric}) coincide. Recall that the metric (\ref{fff}) does not coincide with the Cartan connection on $SE(2)$ since it is not right-invariant. However, the corresponding metric connection does correspond to the Cartan connection on $SO(3)$, recall (\ref{VFSO3}) and (\ref{VFmod}).

The physical dimension of $\beta$ equals $1/\textrm{length}$ and $\beta$ is the fundamental parameter that relates distance on the torus to the distance in the spatial plane. The inner-product between two left-invariant vector fields is now given by
\[
(c^{i}_{1}\mathcal{A}_{i},c^{j}_{2}\mathcal{A}_{j})_{\beta}= g_{ij} c^{i}_{1} c^{j}_{2} = c^{\theta}_{1}c^{\theta}_{2}+ \beta^2\, c^{\xi}_{1}c^{\xi}_{2}+ \beta^2\, c^{\eta}_{1}c^{\eta}_{2}, 
\]
where we use the convention $c^{1}_k=c^{\theta}_{k}, c^{2}_{k}=c^{\xi}_{k}, c^{3}_{k}=c^{\eta}_{k}, k=1,2$. The norm of a left-invariant vector field $c^{i}\mathcal{A}_{i}$ is now given by
\begin{equation} \label{betanorm}
|c^{i}\mathcal{A}_{i}|_{\beta}=\sqrt{(c^{i}\mathcal{A}_{i},c^{i}\mathcal{A}_{i})_{\beta}}=\sqrt{(c^{\theta})^2+ (\beta \, c^{\xi})^2+(\beta \, c^{\eta})^2}=:\|\ul{c}\|_{\beta},
\end{equation}
with $\ul{c}=(c^{1},c^{2},c^{3}) \in \R^3$. Here we stress that the norm $|\cdot|_{\beta}: \mathcal{L}(SE(2)) \to \R^+$ is defined on the space $\mathcal{L}(SE(2))$ of left-invariant vector fields on $SE(2)$, whereas the norm $\|\cdot \|_{\beta}: \R^3 \to \R^{+}$ is defined on $\R^3$.

The gradient ${\rm d} U$ of $U:SE(2)\to \mathbb{R}^{+}$ is given by
\[
{\rm d}U=\frac{\partial U}{\partial \theta} {\rm d}\theta + \frac{\partial U}{\partial \xi}{\rm d}\xi + \frac{\partial U}{\partial \eta} {\rm d}\eta.
\]
It is a co-vector field. The corresponding vector field equals
\begin{equation} \label{broertjegradient}
\mathcal{G}^{-1} {\rm dU}= \frac{\partial U}{\partial \theta}  \partial_{\theta} + \beta^{-2}\frac{\partial U}{\partial \xi} \partial_{\xi} + \beta^{-2}\frac{\partial U}{\partial \eta} \partial_{\eta},
\end{equation}
where $\mathcal{G}^{-1}:T(SE(2))^{*} \to T(SE(2))$ the inverse of the fundamental bijection between the tangent space and its dual. Note that $\mathcal{G}^{-1} {\rm d}\mathcal{A}_{k}=g^{ki} \mathcal{A}_{i}$, with $g^{ij}g_{kl}=\delta^{i}_{k}\delta^{j}_{l}$. The norm of a co-vector field is given by
\[
|a_{i}{\rm d}A^{i}|^2_{\beta}= g^{ij} a_i a_j = (a_{\theta})^{2}  + \beta^{-2} (a_{\xi})^2 +  \beta^{-2}(a_{\eta})^2= \|\ul{a}\|_{\beta^{-1}} \textrm{ with }\ul{a}=(a^{1},a^{2},a^{3}).
\]
Finally we stress that if we differentiate a smooth function $U: SE(2) \to \R^{+}$ along an exponential curve $\gamma(t)=g_{0}\textrm{exp}(t(\sum c^{i} A_{i})) $ passing $g_{0}$ we get (by application of the chain rule)
\begin{equation} \label{dalongexp}
\boxed{
\begin{array}{ll}
\frac{d}{dt} U (\gamma(t)) &= \langle {\rm d} U , \gamma'(t)\rangle =\sum_{i=1}^{3} c^{i} \left. \mathcal{A}_{i} U \right|_{\gamma(t)} \\
 &=
 c^{1} \, U_{\theta}(\gamma(t)) + c^{2}\, U_{\xi}(\gamma(t)) + c^{3} \, U_{\eta}(\gamma(t)). 
\end{array}
}
\end{equation}
Or in words: \emph{The exponential curves $\{g_0 e^{t c^{i} A_{i}}\}_{g_0 \in SE(2)}$ are the characteristics of the left-invariant vector field $c^{i}\mathcal{A}_{i}$}.

After these two preliminaries we return to our goal of finding the optimal tangent spiral at position $g_0 \in SE(2)$ given $U:SE(2) \to \R^{+}$.
\begin{definition}
The solution of the following minimization problem
{\small
\begin{equation} \label{minprob}
\ul{c}_{*}=\arg \min \limits_{\{c^{i}\}_{i=1}^3} \left\{ \left|\left. \frac{d}{dt}\, {\rm d} U(\gamma(t)) \right|_{t=0}\right|^{2}_{\beta} \; |\; \gamma(t)=g_{0}\textrm{exp}(t(\sum \limits_{i=1}^{3} c^{i} A_{i})) \; ;\;
(c^{\theta})^2+ \beta^2(c^{\xi})^2+\beta^2(c^{\eta})^2=1
 \right\},
\end{equation}
}
yields the optimal tangent spiral $\{g_0 e^{t c^{i}_{*} A_{i}}\}_{g_0 \in SE(2)}$ at position $g_0 \in SE(2)$ given $U:SE(2) \to \R^{+}$.
\end{definition}
By means of (\ref{dalongexp}) and the chain rule the energy in (\ref{minprob}) can be rewritten as
\begin{equation} \label{Hess}
\begin{array}{ll}
\left|\left.\frac{d}{dt}({\rm d} U)(\gamma(t))\right|_{t=0}\right|^{2}_{\beta} &=\left\|\nabla (\nabla U)^{T}(\gamma(0)) \cdot \gamma'(0)\right\|^{2}_{\beta^{-1}}  \\
 &= \left\|\left.
 \left(
 \begin{array}{ccc}
\partial_{\theta}(\partial_{\theta} U) &   \partial_{\xi}(\partial_{\theta} U)  &  \partial_{\eta}(\partial_{\theta} U) \\
 \partial_{\theta}(\partial_{\xi} U) &  \partial_{\xi}(\partial_{\xi} U)  &  \partial_{\eta}(\partial_{\xi} U) \\
 \partial_{\theta}(\partial_{\eta} U) &  \partial_{\xi}(\partial_{\eta} U)  & \partial_{\eta}(\partial_{\eta} U) \\
 \end{array}
 \right)\right|_{g_0}
 \left(
 \begin{array}{l}
 c^{1} \\
 c^{2} \\
 c^{3}
 \end{array}
 \right)
 \right\|^{2}_{\beta^{-1}}=: \|\left. H U \right|_{g_0} \, \ul{c}\|^{2}_{\beta^{-1}},
\end{array}
\end{equation}
where $\nabla U:=(\partial_{\theta}U,\partial_{\xi}U, \partial_{\eta}U)$ and
where the non-covariant Hessian $H U$  is not to be mistaken with the covariant Hessian form consisting of covariant derivatives of the Cartan connection (\ref{covH}). We return to this later.

Note that the minimization problem (\ref{minprob}) can now be rewritten as
\[
\arg \min \limits_{\ul{c}} \left\{ \|(H U)(g_0) \;  \ul{c}\|^{2}_{\beta^{-1}} \; |\; \|\ul{c}\|_{\beta}=1\right\}.
\]
Set $M_{\beta}:=\textrm{diag}\{1,\beta^{-1},\beta^{-1}\} \in GL(3,\R)$ and $H_{\beta}U= M_{\beta} \, HU \, M_{\beta}$, then by the Euler-Lagrange theory the gradient of $\nabla_{\ul{c}}\|(H U) \ul{c}\|^{2}_{\beta^{-1}}=\nabla_{\ul{c}}(\ul{c},(H U)^{T} M_{\beta}^{2}(H U)\ul{c})_{1}$ at the optimum $\ul{c}_{*}$ is linearly dependent on the gradient of the side condition $\nabla_{\ul{c}} (1-\|\ul{c}\|^{2}_{\beta})=\nabla_{\ul{c}} (\ul{c},M_{\beta}^{-2}\ul{c})_{1}$:
\[
\begin{array}{l}
(H U(g_0))^T M_{\beta}^{2} (H U(g_0)) \ul{c}_{*}=\lambda \; M_{\beta}^{-2} \ul{c}_{*} \; \desda \;
(H_{\beta}U)^{T}(H_{\beta}U)\tilde{\ul{c}}= \lambda \; \tilde{\ul{c}},
\end{array}
\]
for some Lagrange multiplier $\lambda \in \R$, where $\tilde{\ul{c}}= M_{\beta}^{-1} \ul{c}_{*}$.

So we have shown that the minimization problem (\ref{minprob}) requires eigensystem analysis of $(H_{\beta}U)^{T}H_{\beta}U$ rather than eigensystem analysis of the covariant Hessian given by (\ref{covH}). The eigensystem of the covariant Hessian, however, correspond to the Euler-Lagrange equation for the following minimization problem
(for simplicity we set $\beta=1$)
\begin{equation} \label{minprob2}
\arg \min \limits_{c^{i}} \left\{ \left|\frac{d^{2}}{dt^{2}}U(\gamma(t))\right|^{2} \; |\; \gamma(t)=g_{0}\textrm{exp}(t(\sum \limits_{i=1}^{3} c^{i} A_{i})) \; ;\; \sum \limits_{i=1}^{3} (c^{i})^2=
(c^{\theta})^{2} + (c^{\xi})^{2} +(c^{\eta})^{2}=1 \right\},
\end{equation}
which by means of (\ref{dalongexp}) and again the chain rule can be rewritten as
\[
\begin{array}{l}
\left|\frac{d^{2}}{dt^{2}}U(\gamma(t))\right|^{2} = \left| \frac{d}{dt} \left(\nabla U \cdot \gamma'(t)\right) \right|^2 = \ul{c}^{T} (H U) \ul{c}
\end{array}
\]
and as a result the Euler-Lagrange equations for the minimization problem (\ref{minprob2}) correspond to the eigensystem of $\frac{1}{2}(HU+ (HU)^{T})$, which coincides with covariant Hessian $\nabla \nabla^{T} U$ given by (\ref{covH}):
\[
\nabla \nabla^{T} U \ul{c} =\frac{1}{2}(HU+ (HU)^{T}) \ul{c}= \lambda \ul{c}.
\]
Experiments on images consisting of lines with ground truth curvatures show that minimization problem (\ref{minprob2}) is certainly not preferable over (\ref{minprob}) for spatial curvature estimation.

\textbf{Remarks :}
\begin{itemize}
\item On the commutative group $\R^2$ (i.e. the domain of images $f$ rather than the domain of the orientation scores $U_f$) we do not have this difference, since here the Hessian
{\small $
Hf=
\left(
\begin{array}{cc}
f_{xx}& f_{xy} \\
   f_{yx} & f_{yy}
\end{array}
\right)
$
}
is square symmetric and thereby $H f=\frac{1}{2}(Hf + (Hf)^{T})$ and $(Hf)^{T}(Hf)$ have the same eigenvectors with respective eigenvalues $\{\lambda_{n}\}$ and $\{(\lambda_{n})^{2}\}$. \\
\item
If the spatial gradient vanishes at $g_{0}$ then $\left.\partial_{\xi}U \right|_{g_{0}}=\left.\partial_{\eta}U \right|_{g_{0}}=0$, the problems (\ref{minprob}) and (\ref{minprob2}) have the same minimizer and in this case the covariant Hessian (\ref{covH}) and the non-covariant Hessian (\ref{Hess}) coincide. \\
\end{itemize}
Sofar we did not include the concept of horizontality. Formally, because of the shape of our admissible vectors/distributions $\psi$ in the wavelet transforms, the orientation scores $U_f=\mathcal{W}_{\psi}f$ and their absolute value $|U_f|$ usually do not have a horizontal gradient at locations $g_{0}$ of elongated structures, i.e. in general the gradient does not satisfy $
\partial_{\eta}|U_f| \bigr|_{g_{0}}=0$.
Nevertheless, our algorithm in section \ref{ch:CED} requires horizontal curvature estimates from the absolute value of a (processed) orientation score $|U|$.

Therefor we suggest the following 2 methods for curvature estimation: \\
\\
1. Compute the eigen vectors of $(\tilde{H}_{\beta}^{hor}|U|)^{T} (\tilde{H}_{\beta}^{hor}|U|)$ with horizontal Hessian
\begin{equation} \label{hbu}
\tilde{H}_{\beta}^{hor} |U|=
\left(
\begin{array}{cc}
\beta^{2}\partial_{\theta} \partial_{\theta} |U| & \beta \partial_{\xi} \partial_{\theta} |U|\\
\beta \partial_{\theta} \partial_{\xi} |U| & \partial_{\xi} \partial_{\xi} |U| \\
\beta \partial_{\theta} \partial_{\eta} |U| & \partial_{\xi} \partial_{\eta} |U|
\end{array}
\right)
\end{equation}
to this end we note/recall that the optimum $\ul{c}_{*}= \arg \min\{\|\tilde{H}_{\beta}^{hor} |U|(g_0) \, \ul{c}\|^{2}_{\beta^{-1}} \; |\; \|\ul{c}\|_{\beta}=1\}$ with $\ul{c}=(c^\theta,c^\xi)=c^{\theta}\ul{e}_{\theta}+ c^{\xi} \ul{e}_{\xi}$ satisfies $2 (\tilde{H}_{\beta}^{hor} |U|)^{T} \tilde{H}_{\beta}^{hor} |U|\, \tilde{\ul{c}}= 2 \lambda \tilde{\ul{c}}$, $\ul{c}_{*}=M_{\beta}\tilde{\ul{c}}$ for some Lagrange multiplier $\lambda$. Then we compute the curvature of the projection $\ul{x}(s(t))=\mathbb{P}_{\R^{2}}( g_{0}\textrm{exp}(t(\sum c^{i}_{*} A_{i})))$ of the exponential curve in $SE(2)$ on the ground plane from the eigenvector $\ul{c}_{*}=(c^{\theta}_{*},c^{\xi}_{*})$ with smallest eigen value:
\begin{equation} \label{approach1}
\kappa_{est}=\|\ddot{\ul{x}}(s)\|\textrm{sign}(\ddot{\ul{x}}(s)\cdot \ul{e}_{\eta}) =\frac{c^{\theta}_{*}}{c^{\xi}_{*}}
\end{equation}
2. An alternative approach, however, would be to compute the best exponential curve where we do not restrict ourselves to horizontal curves. In this case we compute the curvature of the projection $\ul{x}(s(t))=\mathbb{P}_{\R^{2}}\left( g_{0}\textrm{exp}(t(\sum_{i=1}^{3} c^{i}_{*} A_{i})) \right)$ of the optimal exponential curve in $SE(2)$ on the ground plane from an eigenvector $\ul{c}_{*}=(c^{\theta}_{*},c^{\xi}_{*},c^{\eta}_{*})$. This eigen vector of $(\tilde{H}_{\beta}|U|)^{T} (\tilde{H}_{\beta}|U|)$, where the $3\times 3$-Hessian is given by
\begin{equation} \label{hbu2}
\tilde{H}_{\beta} |U|=
\left(
\begin{array}{ccc}
\beta^{2}\partial_{\theta} \partial_{\theta} |U| & \beta \partial_{\xi} \partial_{\theta} |U| & \beta \partial_{\eta} \partial_{\theta} |U|\\
\beta \partial_{\theta} \partial_{\xi} |U| & \partial_{\xi} \partial_{\xi} |U| & \partial_{\eta} \partial_{\xi} |U|\\
\beta \partial_{\theta} \partial_{\eta} |U| & \partial_{\xi} \partial_{\eta} |U| & \partial_{\eta} \partial_{\eta} |U|
\end{array}
\right),
\end{equation}
belongs to the pair of eigen vectors closest to the plane $\{\left. \ul{e}_{\xi} \right|_{g_{0}},\left. \ul{e}_{\theta} \right|_{g_{0}}\}$ and has the smallest eigen value. The curvature estimation is now given by
\begin{equation} \label{approach2}
\kappa_{est}=\|\ddot{\ul{x}}(s)\|\textrm{sign}(\ddot{\ul{x}}(s)\cdot \ul{e}_{\eta}) =\frac{c^{\theta}_{*} \textrm{sign}(c_{*}^{\xi})}{\sqrt{(c^{\xi}_{*})^2+(c^{\eta}_{*})^2}}.
\end{equation}
Note that in this alternative approach, in contrast to the other, we do not include the concept of horizontality by restricting ourselves to fitting only horizontal exponential curves, but we simply discard the eigen value (which may be small) corresponding to the eigen vector which is most pointing out the ``correct'' horizontal plane in our selection of eigen vector with smallest eigen value. \\
3. In stead of the Hessian in approach 2. one can also use the eigen vectors of the so-called ``structure tensor'', given by $G_{t}*(\partial_{\theta}|U|, \beta \partial_{\xi}|U|, \beta\partial_{\eta}|U|)^{T} (\partial_{\theta}|U|, \beta \partial_{\xi}|U|, \beta \partial_{\eta}|U|)$, this corresponds to the method proposed by van Ginkel \cite{vanGinkel} who considered curvature estimation from  non-invertible orientation scores.

\begin{figure}
\centerline{
\includegraphics[width=0.7 \hsize]{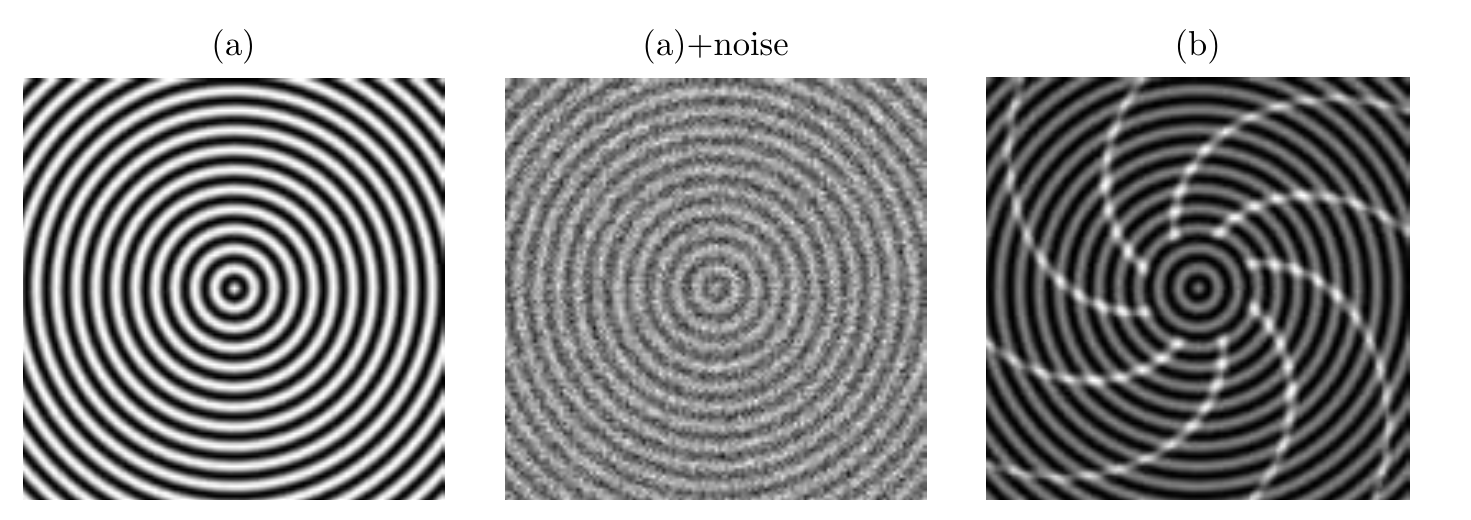}
}
\caption{Curvature estimation test images.}\label{fig:curvexptestimg}
\end{figure}

\begin{figure}
\centerline{\includegraphics[width=0.7\linewidth]{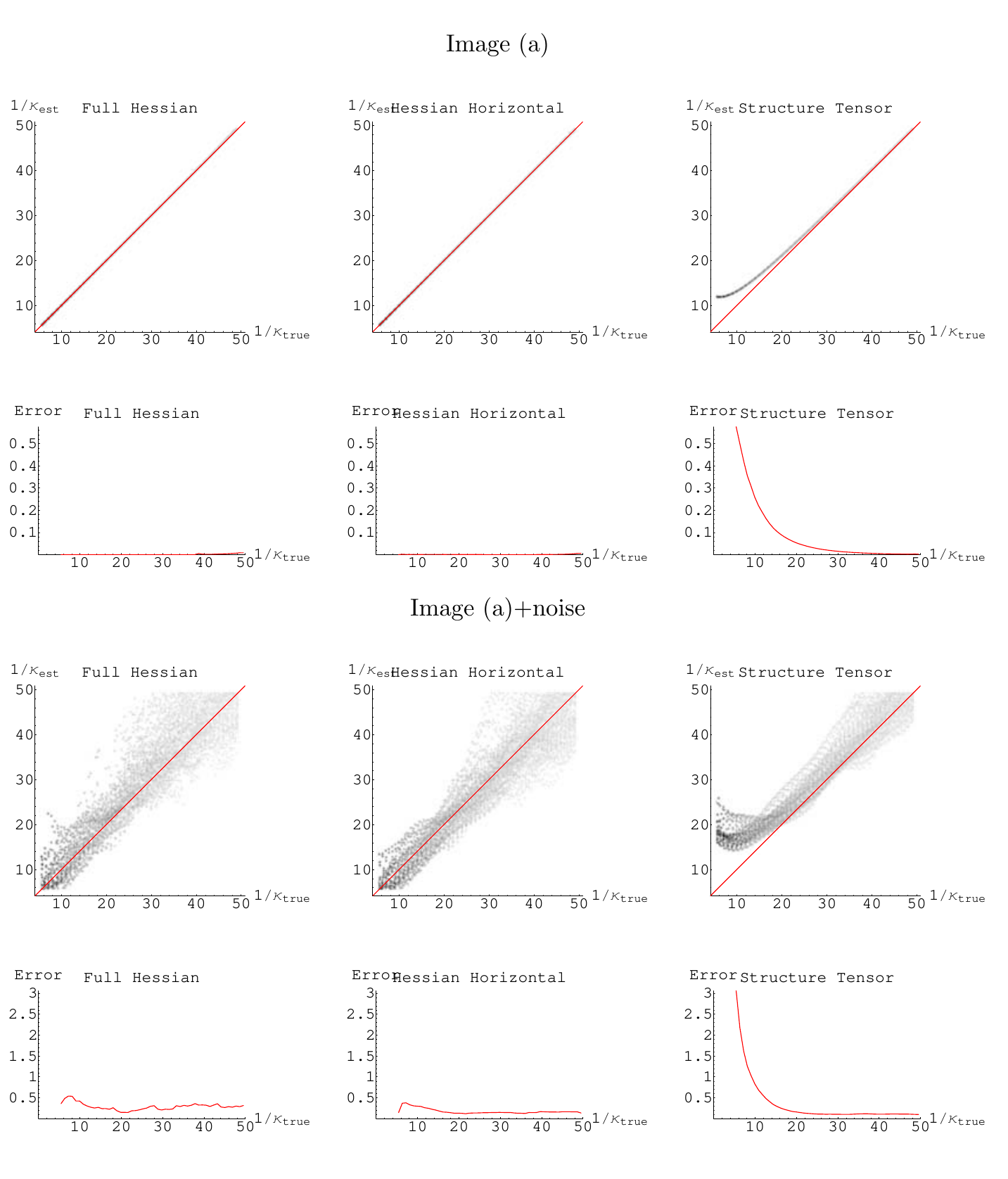}}
\centerline{{\tiny Image (b)}}
\centerline{\includegraphics[width=0.7\linewidth]{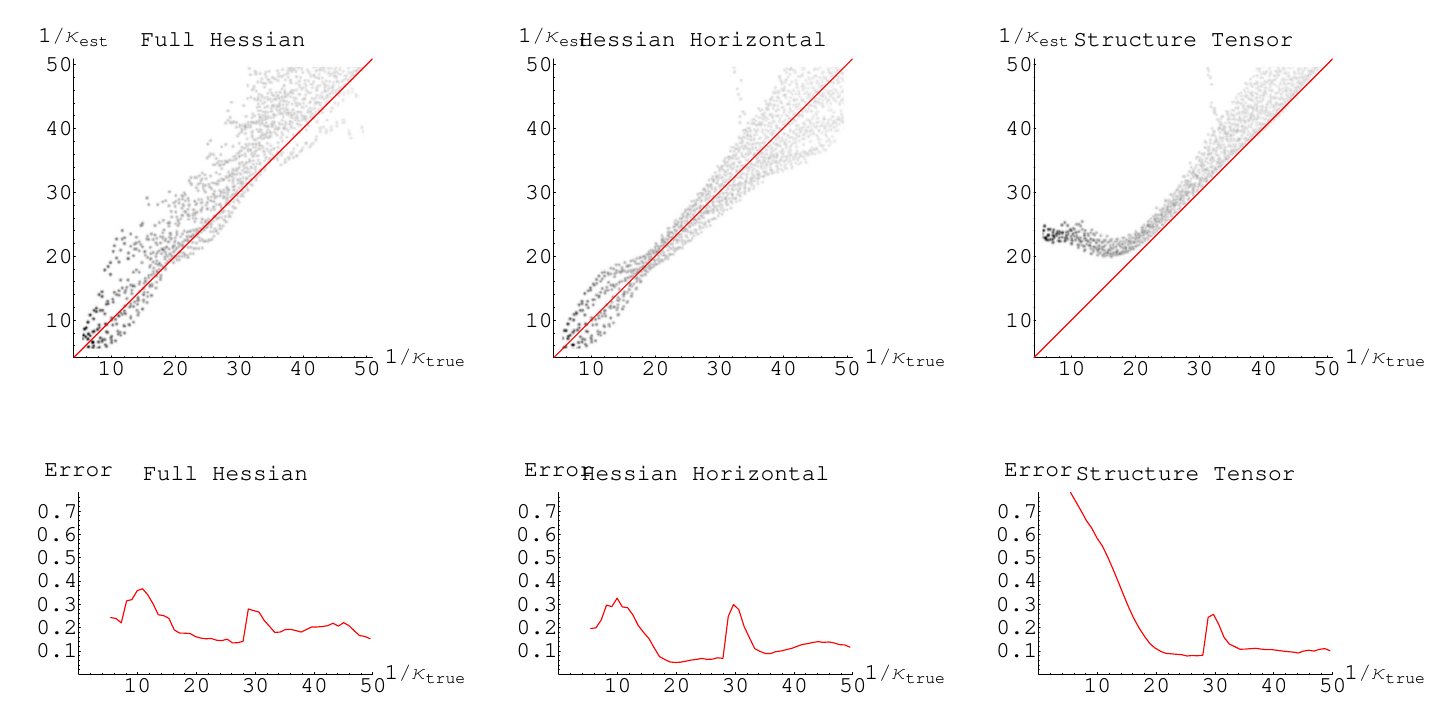}}
\caption{Top: Curvature estimation results (left column approach 1. in section \ref{ch:curvest}, middle column approach 2. in section \ref{ch:curvest}, right column approach 3. (by van Ginkel \cite{vanGinkel}) in section \ref{ch:curvest}) on image (a) and image (a)+noise for the three methods. For both images, the first row shows the density plot of true curvature against estimated curvature, and the second row shows the $\ell_2$-error as function of the different curvatures.
Bottom row: Shows curvature estimation results on image (b) (see Figure \ref{fig:curvexptestimg}) }\label{fig:concentricresults}
\end{figure}

For curvature estimation (comparing all three above methods) on orientation scores of noisy example images, see Figure \ref{fig:curvexptestimg}, and Figure \ref{fig:concentricresults}.

\section{Elastica \label{ch:elasticatrue}}

Let $\epsilon \geq 0$. Then for a smooth curve $t \mapsto \ul{x}(t)$ in $\R^{2}$, with length $L$  we define
\[
\mathcal{E}_{\epsilon}(\ul{x})=\int_{0}^{L} \kappa^{2}(s) +\epsilon \; {\rm d}s.
\]
where $s$ denotes the arc-length parameter\footnote{arc-length in $\R^2$, not arc-length in $SE(2)$.} defined by
\[
s(t)=\int_{0}^{t} \|\frac{\partial \ul{x}(\tau)}{\partial \tau}\| {\rm d \tau},
\]
and where curvature  of the planar curve is given by $\kappa(s)=\|\ddot{\ul{x}}(s)\|$.

Let $\mathcal{C} = \{s \mapsto \ul{x}(s) \in C^{\infty}(\R^{+},\R^{2}) \;|\; \ul{x}(0)=\ul{x}_{0}, \ul{x}(1)=\ul{x}_{1}, \langle \dot{\ul{x}}=\theta_{0}, \langle \dot{\ul{x}}=\theta_{1} \}$ be the space of smooth planar curves
which connect $\ul{x}_{0}$ and $\ul{x}_{1}$ such that the starting and ending direction are prescribed.

We sometimes also consider $\mathcal{C}_{L}=\{s \mapsto \ul{x}(s) \in \mathcal{C}\; |\;
 s \in [0,L]\}$ the space of smooth planar curves with \emph{fixed} total length $L$ which connect $\ul{x}_{0}$ and $\ul{x}_{1}$ such that the starting and ending direction are prescribed.

On $\mathcal{C}$ we can consider the following optimization problem, for $\epsilon>0$,
\begin{equation} \label{nolim}
\textrm{Find }s \mapsto \ul{x}(s) \in \mathcal{C}\textrm{ such that } \mathcal{E}_{\epsilon}(\ul{x}) \textrm{ is minimal.}
\end{equation}
Similarly we can consider the optimization problem on $\mathcal{C}_{L}$:
\begin{equation} \label{lim}
\textrm{Find }s \mapsto \ul{x}(s) \in \mathcal{C}_{L}\textrm{ such that } \mathcal{E}_{0}(\ul{x}) \textrm{ is minimal.}
\end{equation}

We first consider (\ref{nolim}), so let $\epsilon>0$. Let $\ul{x}$ be the optimal curve with tangent $\ul{t}=\dot{\ul{x}}$, normal $\ul{n}=\ddot{\ul{x}}$ and curvature  $\kappa= \ddot{\ul{x}} \cdot \ul{n}$. Then any infinitesimal deformation of this curve should yield lower energy. Since we can always re-parameterize our curves we only need to consider deformation of the curve in normal direction
\begin{equation}\label{deformation}
\ul{x}_{NEW}(s)=\ul{x}(s)+ h \delta(s)\, \ul{n}(s),          h>0
 \end{equation}
with $\delta$ twice differentiable and compactly supported within the open interval $(0,L)$.
Then we stress that the arc-length parameter $s_{NEW}$ of the pertubed curve $\ul{x}_{NEW}$ does not coincide with the arc-length parameter $s$ of the original curve
$\ul{x}$. In fact we have
\[
\dot{\ul{x}}_{New}(s)=(1-\delta h \kappa(s)) (\ul{t}(s) +h \delta'(s) \, \ul{n}(s)) + O(h^{2})\equiv (1-\delta h \kappa(s)) (\ul{t}(s) +h \delta'(s) \, \ul{n}(s))
\]
Then
\begin{equation} \label{kappanew}
\left\{
\begin{array}{l}
{\rm d}s_{NEW}\equiv(1-h \delta \,\kappa){\rm d}s  \\
\ul{t}_{NEW}=\frac{d\ul{x}_{NEW}}{ds_{NEW}}=\frac{ds}{ds_{NEW}}\frac{d \ul{x}_{New}}{ds}\equiv\ul{t} +h\delta' \, \ul{n}\\
\ul{n}_{NEW}= \frac{d^{2} \ul{x}_{New}}{ds_{NEW}^{2}} \equiv\ul{n} -h\delta' \ul{t}
\end{array}
\right.
\Rightarrow
\kappa_{NEW}=\frac{d \ul{t}_{NEW}}{ds_{NEW}} \cdot \ul{n}_{NEW}= \kappa +h\delta'' +h\delta \kappa^{2} +O(h^{2})
\end{equation}
From which it follows that
\[
\begin{array}{ll}
\lim \limits_{h\downarrow 0} \frac{\mathcal{E}(\ul{x}+h\delta \,\ul{n})- \mathcal{E}(\ul{x},\ul{n})}{h} &=
\lim \limits_{h\downarrow 0}\frac{1}{h}\left(\int_{0}^{L_{NEW}} \kappa_{NEW}^{2}(s) +\epsilon  \, {\rm d} s_{NEW}-\int_{0}^{L}\kappa^{2}(s) +\epsilon  \, {\rm d} s \right) \\
 &= \lim \limits_{h\downarrow 0}\frac{1}{h} \int_{0}^{L} ((\kappa +h\delta'' +h\delta \kappa^{2})^{2} +\epsilon)  \, (1-h\delta \kappa) -(\kappa^{2}+\epsilon){\rm d} s
\end{array}
\]
so by partial integration we find
$
\lim \limits_{h \to 0}
\frac{\mathcal{E}(\ul{x}+ h\delta \,\ul{n})- \mathcal{E}(\ul{x}+\delta \,\ul{n})}{h}=0 \textrm{ for all }\delta
$ if and only if
\begin{equation} \label{elastics}
2 \ddot{\kappa}+\kappa^{3}= \epsilon \kappa,
\end{equation}
however we stress that not all solutions of (\ref{elastics}) lead to global minimization of (\ref{nolim}). They can be local minima or even saddle points.

For problem (\ref{lim}) we note that the deformations must be length preserving, in this case we have
\[
\int_{0}^{L_{New}} {\rm d}s_{NEW} =\int_{0}^{L} {\rm d}s \desda \int_{0}^{L} \kappa(s) \delta(s) {\rm d}s=0  \]
Consequently the optimation is similar as above with the only difference that $\epsilon$ has to be replaced by an Euler lagrange multiplier
\begin{equation} \label{elastics2}
2 \ddot{\kappa}+\kappa^{3}= \lambda \kappa, \qquad \lambda>0.
\end{equation}

Now (\ref{elastics}) and (\ref{elastics2}) provide the curvature of elastica curves, which is unique if we set
\[
\kappa(0)=\kappa_0 \textrm{ and }\dot{\kappa}(0)=\kappa_{0}',
\]
for some positive constants $\kappa_{0}$ and $\kappa_{0}'$ which we will determine later.
To get the elastica curves themselves we have to integrate the Frenet formulas:
\[
\frac{d}{ds}\left(
\begin{array}{l}
\dot{\ul{x}}(s) \\
\ddot{\ul{x}}(s)
\end{array}
\right)=
\left(
\begin{array}{cc}
0 & \kappa(s)I_{2} \\
-\kappa(s) I_{2} & 0
\end{array}
\right)
\left(
\begin{array}{l}
\dot{\ul{x}}(s) \\
\ddot{\ul{x}}(s)
\end{array}
\right),
\]
with $I_2$ the identity matrix and where
the solution is uniquely determined by $s \mapsto \kappa(s)$ and $\angle \dot{\ul{x}}(0)=\theta_{0}$, to this end we note that the exponential of a skew symmetric matrix is orthogonal and therefor $\|\dot{\ul{x}}(s)\|=\|\ddot{\ul{x}}(s)\|=1$ and $\dot{\ul{x}} \cdot \ddot{\ul{x}}=0$ for all $s \geq 0$, so that $\angle \dot{\ul{x}}(0)=\theta_{0}$ sets the initial condition $\dot{\ul{x}}(0)= \cos \theta_{0} \ul{e}_x + \sin \theta_{0} \ul{e}_y, \ddot{x}(0)=-\sin \theta_0 \ul{e}_x +\cos \theta_{0} \ul{e}_{y}$ and thereby the full solution $s \mapsto (\dot{\ul{x}}(s),\ddot{\ul{x}}(s))$.

Now we get the elastica curve by integration over $s >0$:
\[
\ul{x}(s)=\ul{x}(0)+\int \limits_{0}^{s} \dot{\ul{x}}(t) \; {\rm d}t.
\]
The three free parameters $\kappa_0$, $\kappa_{0}'$ and the length of the elastic $L>0$ have to be set such that
\[
\ul{x}(L)=\ul{x}_{1} \textrm{ and }\angle \dot{\ul{x}}(L)= \theta_{1}.
\]
This can be done by means of a shooting algorithm where we use the $B$-spline solutions (\ref{modes}) (which correspond to the coordinate dependent mode-lines of a product of two Heisenberg approximations of the Green's functions) as an initial condition.

\begin{figure}
\centerline{\includegraphics[width=0.4\hsize]{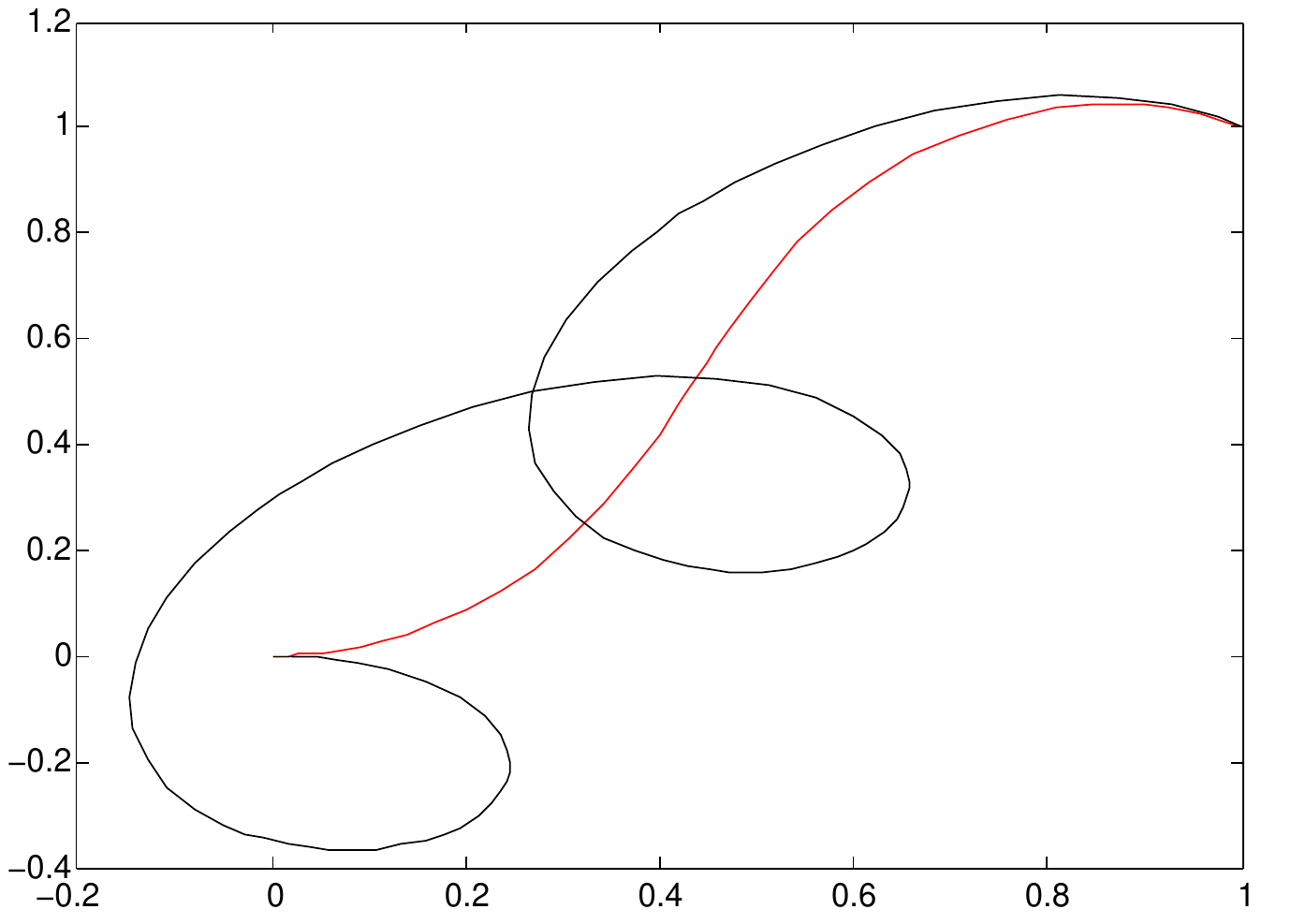}}
\caption{A good initial guess, such as (\ref{initpar}), for the shooting algorithm is crucial as we do not want to get stuck in local minima. Typically a dampened Newton-Raphson algorithm
$\ul{k}_{n+1}=\ul{k}_{n} -q\, D\Phi(\ul{x}_{n}) \Phi(\ul{x}_n)$, $0<q \leq 1$
in the shooting algorithm,
i.e. to find the zeros of (\ref{NR}) avoids wide jumps, so that the chance of getting stuck in a local minima is reduced.
\label{fig:localminima}}
\end{figure}

The shooting algorithm works as follows: First we write everything in one ODE-system
{\scriptsize
\[
\frac{d}{ds}
\left(
\begin{array}{c}
x(s) \\
y(s) \\
\theta(s) \\
\dot{x}(s) \\
\dot{y}(s) \\
\ddot{x}(s) \\
\ddot{y}(s) \\
\kappa(s) \\
\dot{\kappa}(s)
\end{array}
\right)
=\left(
\begin{array}{ccccccccc}
0 & 0 & 0 & 1 & 0 & 0 &0 &0 & 0 \\
0 & 0 & 0 & 0 & 1 & 0 &0 &0 & 0 \\
0 & 0 & 0 & 0 & 0 & 0 &0 &1 & 0 \\
0 & 0 & 0 & 0 & 0 & \kappa(s) &0 &0 & 0 \\
0 & 0 & 0 & 0 & 0 & 0 &\kappa(s) &0 & 0 \\
0 & 0 & 0 &  -\kappa(s) & 0 &0 &0 & 0 &0 \\
0 & 0 & 0 & 0 & -\kappa(s) & 0 &0 &0 & 0   \\
0 & 0 & 0 & 0 & 0 & 0 &0 &1  & 0 \\
0 & 0 & 0 & 0 & 0 & 0 &-\frac{\kappa^{2}(s)+\epsilon}{2} & 0 & 0
\end{array}
\right)
\left(
\begin{array}{c}
x(s) \\
y(s) \\
\theta(s) \\
\dot{x}(s) \\
\dot{y}(s) \\
\ddot{x}(s) \\
\ddot{y}(s) \\
\kappa(s) \\
\dot{\kappa}(s)
\end{array}
\right)
\]
} with initial condition
\begin{equation} \label{ic}
(x(0),y(0),\theta(0),\dot{x}(0),\dot{y}(0),\ddot{x}(0),\ddot{y}(0),\kappa(0),\dot{\kappa}(0))=
(0,0,0,1,0,0,1,\kappa_0,\kappa_0'),
\end{equation}
where we note by means of left-invariance we can assume that $g(0)=e$. Now this system of equations has a unique solution and can be numerically solved by a standard Runge-Kutta method yielding the numeric solution $(\overline{x}(s),\overline{y}(s),\theta(s)=\arg(\overline{x}'(s)+ i\overline{y}'(s)))$. This defines a function $\PPPP: \R \times \R \to C(\R^{+},\R^2)$ which maps $(\kappa_0,\kappa_0')$ (which determines the initial condition (\ref{ic})) to the spatial curve $s \mapsto (\overline{x}(s),\overline{y}(s))$ solution.
So finally, we apply a (dampened) Newton-Raphson scheme on the function $\Phi: \R \times \R \times \R^{+} \to \R^{+}$ given by
\begin{equation} \label{NR}
(\kappa_{0},\kappa_{0}',L) \mapsto \|(\PPPP(\kappa_0,\kappa_0'))(L)-\ul{x}_{1}\|^{2} + \gamma^{2}(\angle (\dot{\PPPP}(\kappa_0,\kappa_0'))(L) \; -\theta_{1})^{2},
\end{equation}
for suitable choice of $\gamma>0$, where we use finite difference approximations for the derivatives of $\Phi$. At this point we note that the initial guess for $\kappa_{0}$, $\kappa_{0}'$ and $L$, which can be derived from (\ref{modes}) is given by
\begin{equation} \label{initpar}
\begin{array}{lll}
\kappa(0)= \frac{y''(0)}{(1+(y'(0))^2)^{3/2}}, &
\kappa'(0)= \frac{y'''(0)}{(1+(y'(0))^2)^{3/2}}- 3 \frac{(y''(0))^{2}y'(0)}{(1+(y'(0))^2)^{5/2}},  &
L=\int \limits_{0}^{x_1} \sqrt{(y'(x))^{2}+1} {\rm d}x,
\end{array}
\end{equation}
with $y'(0)=\theta_{0}$, $y''(0)=\frac{2}{x_{1}^{2}}(3 y_1 +  x_1 (\theta_1 -2\theta_{0}))$ and $y'''(0)=\frac{6}{x_{1}^{3}}(-2 y_1 +  x_{1}(\theta_{0}-\theta_{1}))$, $y(x)=x\theta_{0}+\frac{x^3}{x_1^{3}} (- 2 y_1 + x_1 (\theta_0 - \theta_{1}))+\frac{x^2}{x_1^{2}}(3 y_1 +x_{1}(\theta_{1}-2 \theta_{0}))
$.
To this end we note that this good initial guess is highly relevant to avoid the shooting algorithm to get stuck at local minima, as illustrated in Figure \ref{fig:localminima}.

The horizontal curve in $SE(2)$ corresponding to the elastic is given by $s \mapsto g(s)=(\ul{x}(s), \angle \dot{\ul{x}}(s))$ and indeed satisfies
\[
g(0)=(\ul{x}_0,e^{i\theta_{0}}) \in SE(2) \textrm{ and }g(L)=(\ul{x}_{1},e^{i\theta_{1}}) \in SE(2).
\]
\subsubsection*{Exact derivation of the elastica}

In this subsection we will investigate well-known analytic formulae for (the curvature $\kappa(s)$ of) elastica curves, to get some analytic grip on the behavior of these curves. It turns out that exact formula for elastica curves involve special functions (Jakobi-elliptic or theta functions) with practical disadvantages due to the twice integration of their curvature.

Consider the ordinary differential system
\begin{equation} \label{objectiveeq}
\left\{
\begin{array}{l}
2 \kappa''+\kappa^{3}= \epsilon \kappa, \\
\kappa(0)=\kappa_{0}, \ \
\kappa'(0)=\kappa_{0}'
\end{array}
\right.
\end{equation}
A multiplication of the ODE by $\kappa'$ and integration over arc-length yields
\begin{equation} \label{equafterint}
\left(2 \, \frac{ d \kappa}{ds}\right)^{2}= -\kappa^{4} + 2 \epsilon \, \kappa^{2} + 4 C_{1},
\end{equation}
where $C_{1}$ is an integration constant, related to the initial conditions by means of
\[
C_{1}= (\kappa_{0}')^{2}+\frac{1}{4}\kappa_{0}^{4} - \frac{\epsilon}{2} \kappa^{2}_{0} \geq (\kappa_{0}')^{2}-\frac{\epsilon^{2}}{4} \geq -\frac{\epsilon^{2}}{4}.
\]
From which it directly follows that $s$ is an elliptic integral in $\kappa$, \cite{Mumford},
\begin{equation}\label{sell}
s= \int_{\kappa_{0}}^{\kappa(s)} \frac{2 {\rm d}u }{{\sqrt{-u^{4}+ 2 \epsilon u^2 + 4 C_{1}}}},
\end{equation}
which only holds for $C_{1}\geq (\kappa_{0}')^{2}-\frac{\epsilon^{2}}{4}$.
Now for $C_{1}\neq 0$ this can be rewritten as follows
\[
\begin{array}{ll}
\frac{i s}{2}  &= \int_{\kappa_{0}}^{\kappa(s)} \frac{ {\rm d}u }{{\sqrt{u^{4}- 2 \epsilon u^2 - 4 C_{1}}}} = \int_{\kappa_{0}}^{\kappa(s)} \frac{ {\rm d}u }{{\sqrt{(u- \gamma_{+}^{\epsilon})(u+\gamma_{-}^{\epsilon})}}} 
     = \int_{\frac{\kappa_{0}}{\sqrt{\gamma_{-}}}}^{\frac{\kappa(s)}{\sqrt{\gamma_{-}}}}
     \frac{ {\rm d}v }{\sqrt{1-v^{2}}\sqrt{1- \frac{\gamma_{-}^{\epsilon}}{\gamma_{+}^{\epsilon}}\, v^{2}}}
\end{array}
\]
where $v=u/\sqrt{\gamma_{-}^{\epsilon}}$ and where $\gamma_{\pm}^{\epsilon}= \epsilon \pm \sqrt{\epsilon^{2}+4C_{1}}$ are the real zero's of
$u \mapsto u^{2} - 2 \epsilon u - 4 C_{1}$.
We have
\begin{equation} \label{eqn1}
\sqrt{\gamma_{+}} \frac{i\, s}{2}= \int_{0}^{\frac{\kappa(s)}{\sqrt{\gamma_{-}}}} \frac{ {\rm d}v }{\sqrt{1-v^{2}}\sqrt{1- \frac{\gamma_{-}^{\epsilon}}{\gamma_{+}^{\epsilon}}\, v^{2}}} -
i C_{2} \sqrt{\gamma_{+}},
\end{equation}
where the constant $C_{2}= \textrm{sign}\{ C_{1}^{\epsilon}\} \, \frac{1}{i \sqrt{\gamma_{+}}} \int_{0}^{\frac{\kappa_{0}}{\sqrt{\gamma_{-}}}} \frac{ {\rm d}v }{\sqrt{1-v^{2}}\sqrt{1- \frac{\gamma_{-}^{\epsilon}}{\gamma_{+}^{\epsilon}}\, v^{2}}}$.
As a result we can rewrite (\ref{eqn1}) \begin{equation} \label{kappaana}
\kappa(s)= \sqrt{\gamma_{-}^{\epsilon}} \textrm{sn}\left(\frac{\sqrt{\gamma_{+}^{\epsilon}}i(s+C_{2})}{2}, \frac{\gamma_{-}^{\epsilon}}{\gamma_{+}^{\epsilon}}\right),
\end{equation}
where $\textrm{sn}(\cdot,k)$ denotes the Jacobi elliptic function of the first kind,
which is the solution of (\ref{objectiveeq}) for $C_{1}\neq 0$ if and only if
\[
\begin{array}{l}
C_{1}=C_{1}^{\epsilon}:=(\kappa_{0}')^{2}+\frac{1}{4}\kappa_{0}^{4} - \frac{\epsilon}{2} \kappa^{2}_{0} \textrm{ and }
C_{2}= \textrm{sign}\{ C_{1}^{\epsilon}\} \left. \int \limits_{0}^{\frac{\kappa(s)}{\sqrt{\gamma_{-}}}} \frac{ {\rm d}v }{\sqrt{1-v^{2}}\sqrt{1- \frac{\gamma_{-}^{\epsilon}}{\gamma_{+}^{\epsilon}}\, v^{2}}} \right|_{\gamma_{\pm}=\epsilon \pm \sqrt{\epsilon^{2}+ 4\,C_{1}^{\epsilon}}}.
\end{array}
\]
So we see that the curvature is a periodic function with period
\begin{equation}
T_{C_1,\epsilon}=
-\frac{2}{\sqrt{\gamma^{\epsilon}_{+}}} \int_{0}^{\pi/2} \frac{{\rm d}\theta}{\sqrt{1-(1-\frac{\gamma_{-}^{\epsilon}}{\gamma_{+}^{\epsilon}}) \sin^{2} \theta}} ,
\end{equation}
note that $C_{1} \mapsto T_{C_{1},\epsilon}$ is a monotonically decreasing differentiable function with
$\lim \limits_{C_{1} \downarrow 0} \frac{1}{T_{c_{1},\epsilon}}=0$,
so for applications $C_{1}$ is typically small since the number of periods over a fixed interval of interest corresponds to the number of turns the elastica makes during this interval. For $C_{1}\to 0$
solutions $s \mapsto \ul{x}(s)$ tend to straight lines (they are straight lines if $\kappa_{0}=\kappa'_{0}=0$).
For $\kappa_{0}=\sqrt{\epsilon}$ and $\kappa'(0)=0$ the solutions are circles.

Finally we note that the elastica $s \mapsto (x(s),y(s))\equiv z(s):=x(s)+i y(s)$ follow by their curvature $\kappa(s)$ by means of
\begin{equation} \label{int}
\frac{dz}{ds}=e^{i
\theta(s)} \textrm{ and }
\frac{d\theta}{ds}=\kappa(s),
\end{equation}
now the primitive $\theta(s)$ can easily be derived analytically from (\ref{kappaana}),
but the second integration step which provides the actual curve $z(s)$ is a non trivial expansion in elliptic functions, for details and derivations see \cite{MarkusThesis} .
This problem is partially resolved in Mumford's approach \cite{Mumford}{p.502-505}, where the Jakobi-elliptic functions are replaced by theta functions and where the arc-length
parametrization $z(s)=x(s)+i y(s)$ does not involve an integration. But even in this approach the standard solutions
\[
s \mapsto c \, \frac{d}{ds} \log \theta_{\Lambda}(s-\eta)- a s, \textrm{ with }a,c \in \mathbb{C} \textrm{ and } \Lambda= \mathbb{Z}+ it \mathbb{Z}, \eta=-it/4 \textrm{ or }\Lambda= \mathbb{Z}+ (1/2)(it+1) \mathbb{Z}, \eta=0,
\]
involve several parameters which are not straightforwardly related (to $\epsilon>0$ and) the boundary conditions $\ul{x}(0)=\ul{x}_{0}$, $\dot{\ul{x}}(0)=(\cos \theta_{0},\sin \theta_{0})$, $\ul{x}(1)=\ul{x}_{1}$, $\dot{\ul{x}}(1)=(\cos \theta_{1},\sin \theta_{1})$. So for computation purposes a shooting algorithm of the type (\ref{NR}) is preferable over an entirely exact approach.

\subsection{The corresponding geodesics \label{ch:geodesics}}

Next we are going to repeat the proceeding with the ``only'' difference 
that we take a square root of the integrand so that we have a homogenous energy
\[
\mathcal{E}(\ul{x})=\int \limits_{0}^{L} \sqrt{\kappa^{2}(s)+\epsilon }\; {\rm d} s,
\]
$\ul{x} \in \mathcal{C}$, which is related to the Cartan connection, recall the second Example in section \ref{ch:Cartan}.

We stress that the corresponding lifted curve $s \mapsto (\ul{x}(s), \angle \dot{\ul{x}}(s))$, recall definition \ref{def:horizontal}, is a geodesic on $SE(2)$ in the classical sense, since by our restriction to horizontal curves we have $\kappa(s)=\dot{\theta}(s)$ and we can rewrite the energy as
\[
\boxed{
\int \limits^{L}_{0} \sqrt{\kappa^{2}(s) +\epsilon} \; {\rm d}s= \int \limits_{0}^{L} \sqrt{|\dot{\theta}(s)|^2 + |\langle {\rm d}\xi, \dot{\ul{x}}(s) \rangle|^{2}\beta^2}\; {\rm d}s=\int \limits_{0}^{L} |\dot{\gamma}(s)|_{\beta} \;{\rm d}s= \int \limits_{0}^{L} \sqrt{\sum_{i,j} \dot{\gamma}^{i}(s) \dot{\gamma}^{j}(s) g_{ij}} \;{\rm d}s, }
\]
with $\beta^{2}=\epsilon$, where we recall (\ref{metric}) and (\ref{betanorm}).

\begin{figure}
\centerline{
\includegraphics[width=0.25\hsize]{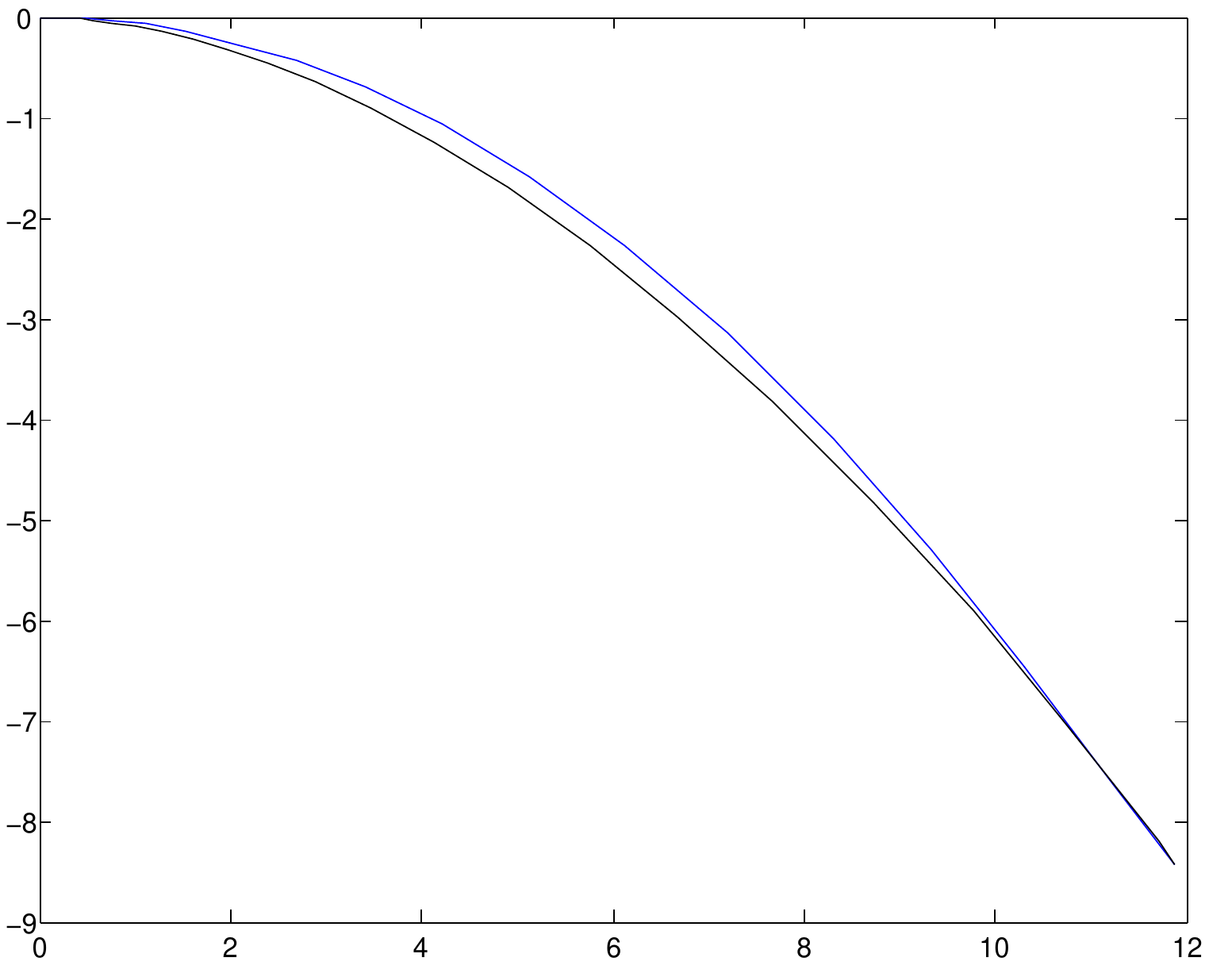}
\includegraphics[width=0.25\hsize]{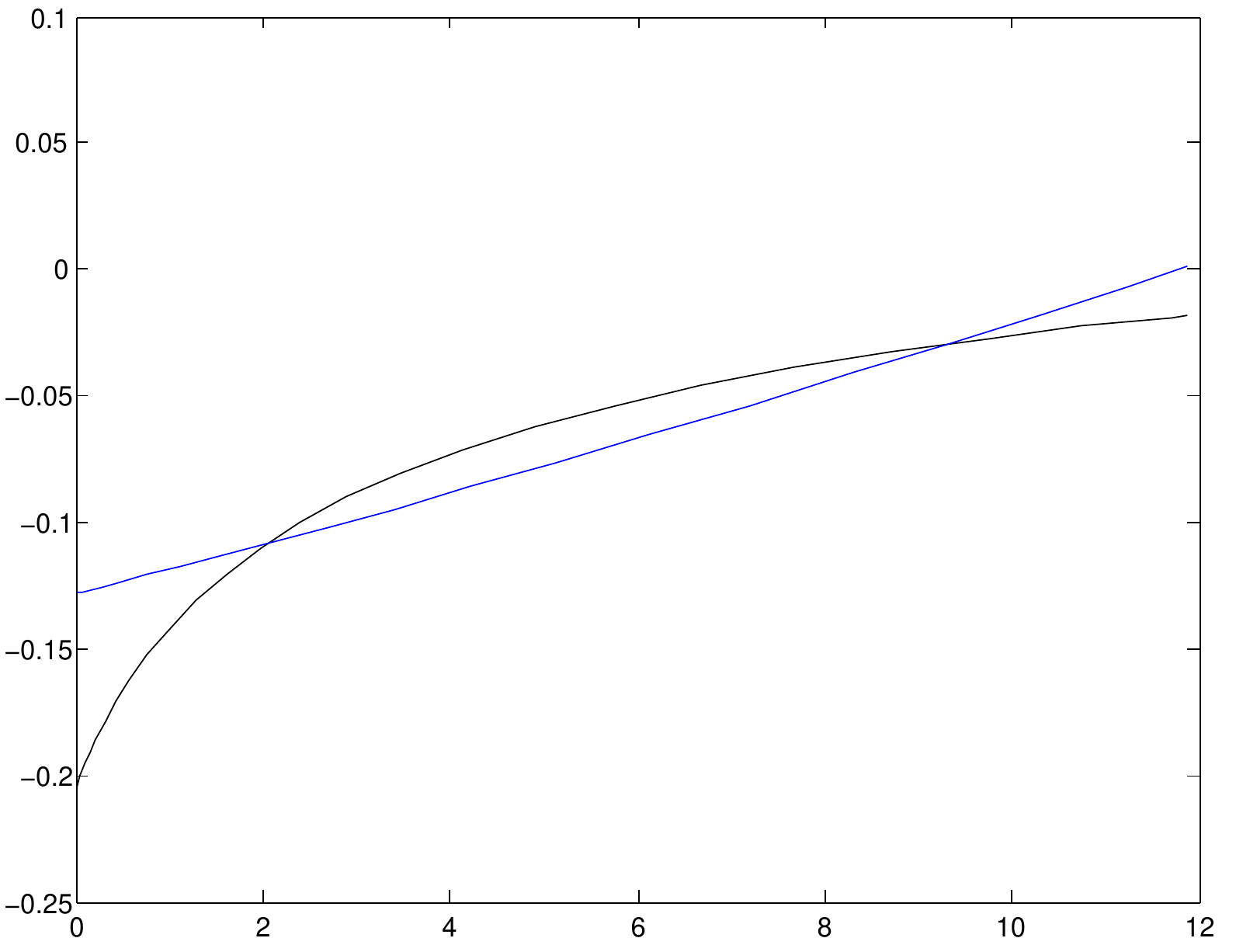}
}
\centerline{
\includegraphics[width=0.25\hsize]{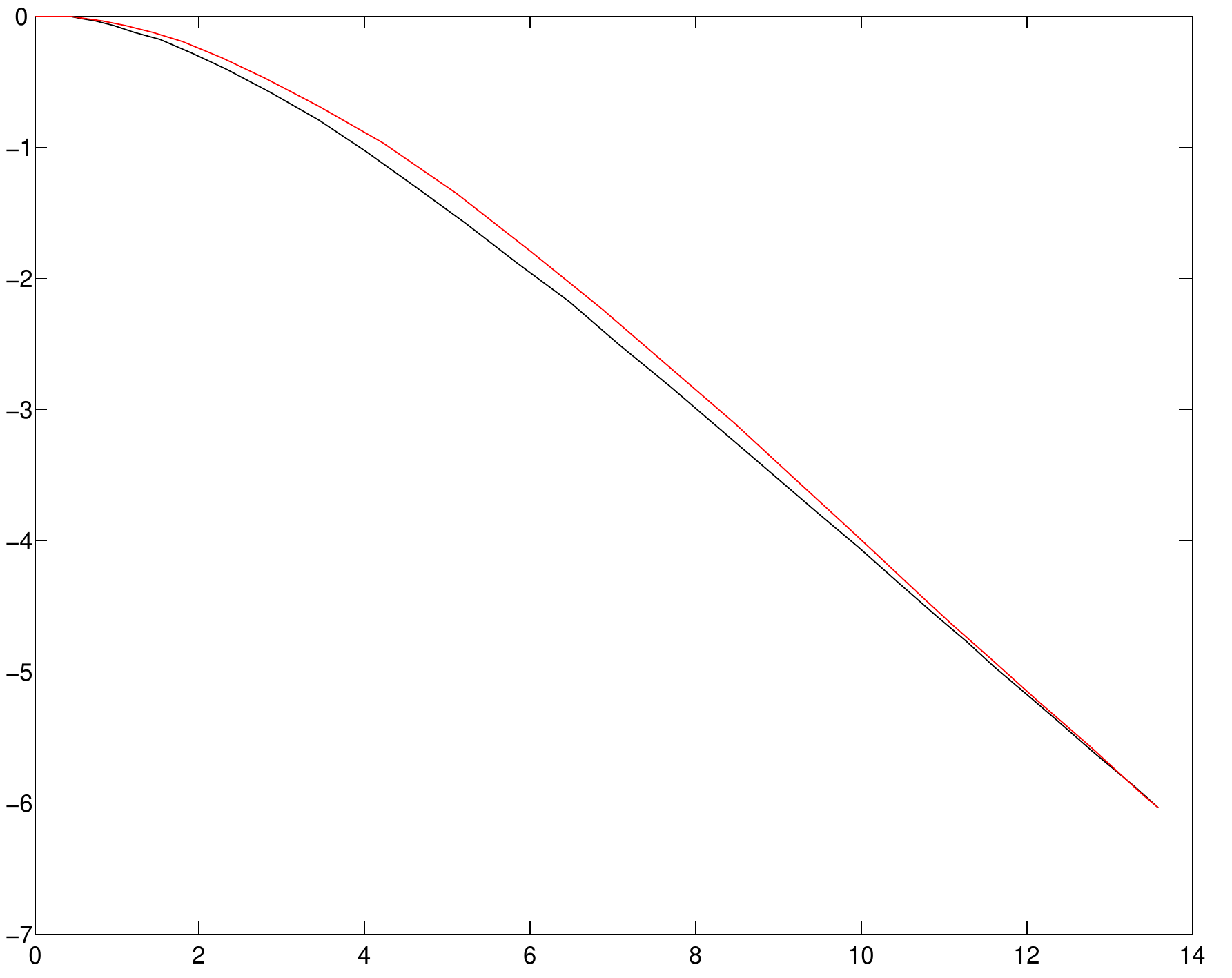}
\includegraphics[width=0.25\hsize]{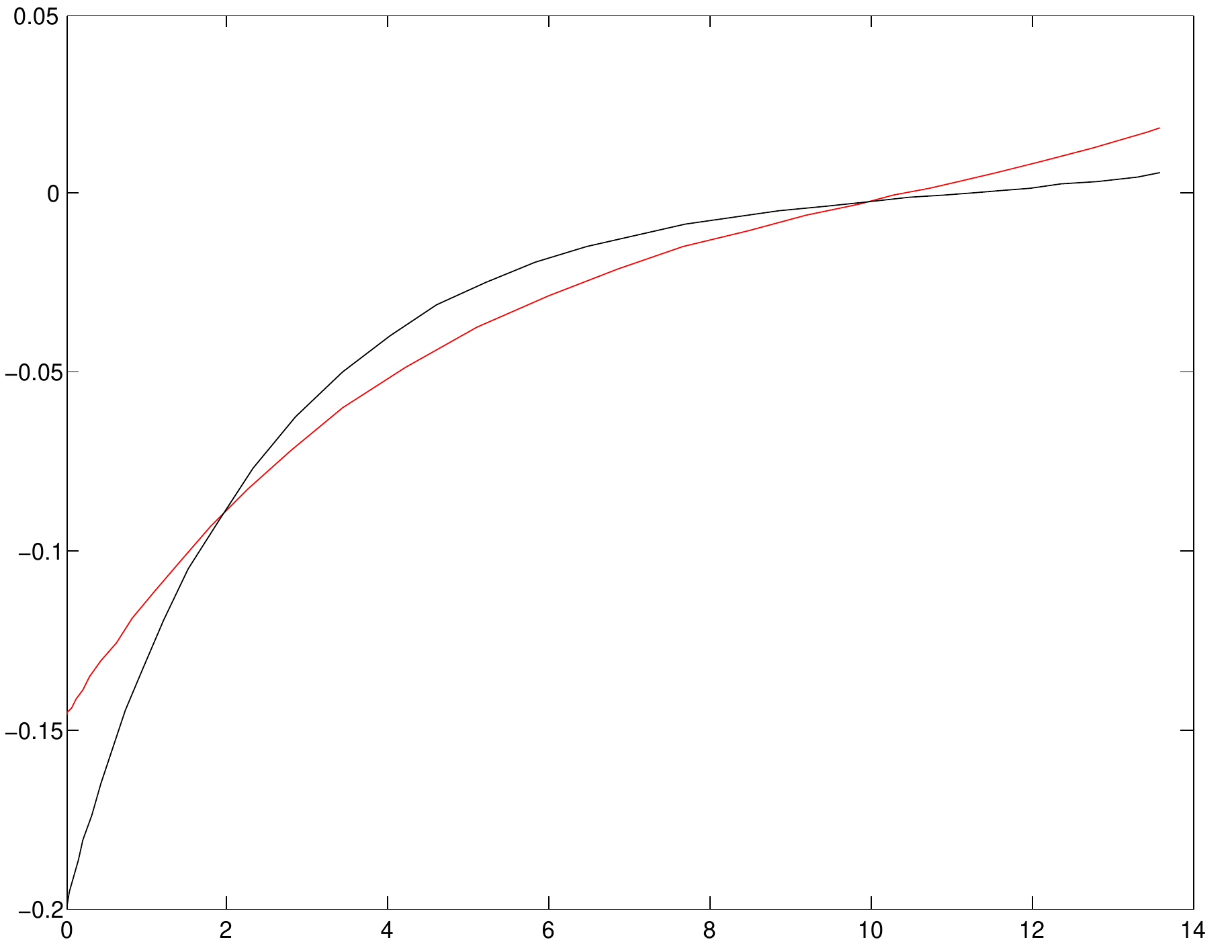}
}
\centerline{\includegraphics[width=0.25\hsize]{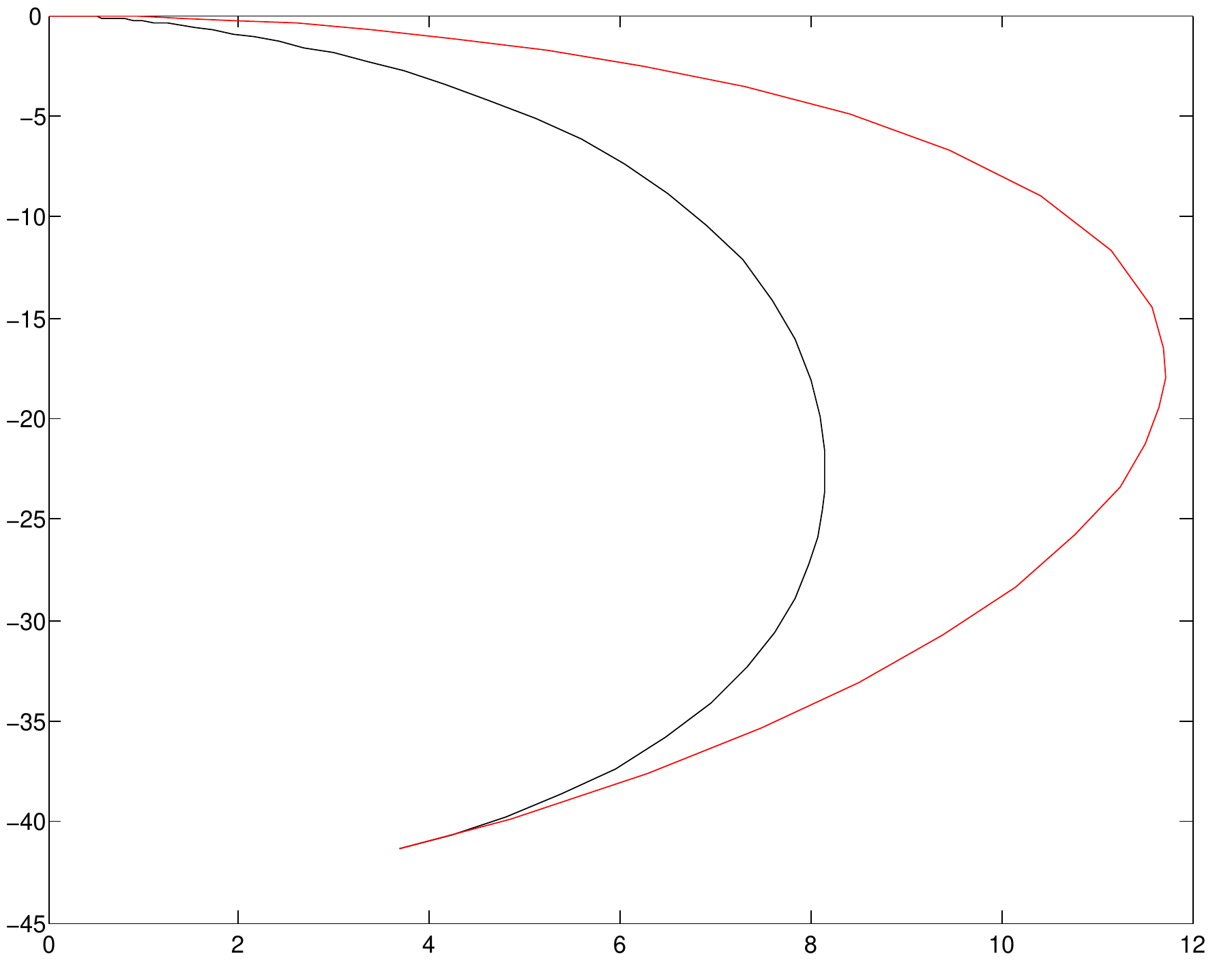}
\includegraphics[width=0.25\hsize]{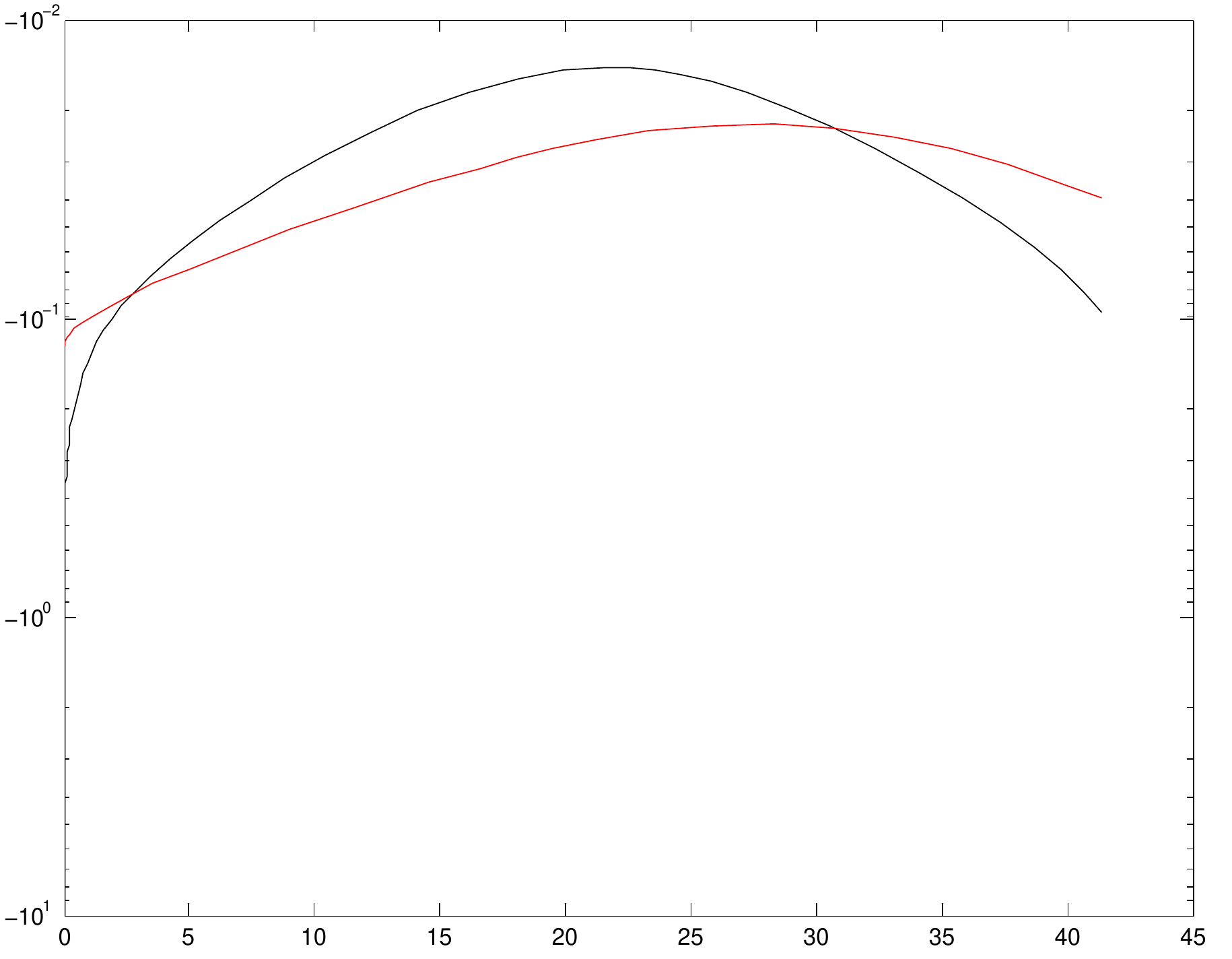}
}
\caption{ Left column; joint plots of the elastica and geodesics in $\R^2$. Right column;
the corresponding graphs of the curvature $\kappa(y)$ as a function of the $y$-coordinate of points along the curve. We use the same $\epsilon$ for both elastics and geodesics. This means we should set $\epsilon=4 \alpha D_{11}=\frac{D_{11}}{D_{22}}$, where $(\alpha,D_{11})$ are the parameters of Mumford's direction process and where $(D_{11},D_{22})$ are the parameters of the contour enhancement process, also proposed by Citti and Sarti\cite{Citti} as a cortical model for contour enhancement.
Top row parameter settings: $\epsilon=0.0125$, Length geodesic $s=15$, $\kappa(0)=-0.20502$, $\kappa'(0)= -\epsilon^{-\frac{1}{2}} (\kappa_0^{2}+\epsilon) \kappa_{0}$, so that $s_{max}=\infty, \ul{x}_{0}=\ul{0}, \theta_{0}=0, \ul{x}_{1} \approx (11.88, -8.43535), \theta_{1}=-51.8947^{\circ}$.
Middle row: Length geodesic $s=15$, $\epsilon=0.125$, $\kappa_{0}=-0.2$, $\kappa'(0)= -\epsilon^{-\frac{1}{2}} (\kappa_0^{2}+\epsilon) \kappa_{0} \Rightarrow s_{max}=\infty$,   $\ul{x}_{0}=\ul{0}, \theta_{0}=0, \ul{x}_{1} \approx (13.58, -6.09), \theta_{1}=-29.36^{\circ}$.
Bottom row: Illustration of an extreme case, Length geodesic $s=48.997$,
$\kappa'(0)= -0.0058, \kappa_{0}=-0.125$
 $\kappa(0)$ and $\kappa'(0)$ such that $\ul{x}^{\textrm{geodesic}}(s=48.997)=\ul{x}_{1}$ and $\angle (\dot{\ul{x}}^{\textrm{geodesic}}(s=45),\ul{e}_{x})=\theta_{1}$,
$\ul{x}_{0}= \ul{0}, \theta_{0}=0, \ul{x}_{1}\approx (-41.383) , \theta_{1}=50.63^{\circ}$.
Typically, the geodesics have more curvature at the boundaries whereas the elastica curves have more curvature in the middle of the curve. For reasonable parameter settings (such as top row where we set $\epsilon= 4 \alpha D_{11}=\frac{D_{11}}{D_{22}}=0.0125$) the geodesics are close to the elastica curves. This makes sense since the random walkers in the contour enhancement process are allowed to turn in negative $\xi$-direction in contrast to random walkers in the contour completion process
\label{fig:comparison}, recall Figure \ref{fig:orientationbundle2}.}
\end{figure}

The energy after deformation (\ref{deformation}) becomes
\[
\begin{array}{ll}
\mathcal{E}(\ul{x} +h\delta \ul{n})=\int \limits_{0}^{L_{NEW}} \sqrt{\kappa_{NEW}^{2}(s) +\epsilon } \;{\rm d}s_{NEW} &=  \int \limits_{0}^{L} \sqrt{\kappa^{2} +2 h\delta'' \kappa + 2 \delta h \kappa^{3} + \epsilon +O(h^{2})} (1-\delta h \kappa) {\rm d}s \\
  &= \int \limits_{0}^{L}\sqrt{\kappa^{2}+\epsilon}\sqrt{1+\frac{2h\delta'' \kappa +2 h\delta \kappa^3}{\kappa^{2} +\epsilon}}(1-h\delta \kappa){\rm d}s \\
  &= \int \limits_{0}^{L}\sqrt{\kappa^{2}+\epsilon}\left(1+ \frac{h\delta'' \kappa + h\delta \kappa^3}{\kappa^{2} +\epsilon} + O(h^{2}) \right)(1-h\delta \kappa){\rm d}s \\
  &= \int \limits_{0}^{L}\sqrt{\kappa^{2}+\epsilon}\left(1+ \frac{h\delta'' \kappa + h\delta \kappa^3}{\kappa^{2} +\epsilon} -\delta \kappa + O(h^{2}) \right){\rm d}s \\
  &= \mathcal{E}(\ul{x}) + h \int \limits_{0}^{L} \sqrt{\kappa^{2}+\epsilon }\left(
  \frac{h\delta'' \kappa +\delta \kappa^{3}}{\kappa^2 +\epsilon}-\delta \kappa\right) {\rm d}s +O(h^2)
 \end{array}
\]
where we used $\sqrt{1+x}=1 +\frac{1}{2} x + O(x^2)$ and (\ref{kappanew}). So in order to get a local minima the energy should increase under all possible deformations parameterized by $\delta$ and therefore we have
\begin{equation} \label{zdiv}
\begin{array}{l}
\left(\frac{\kappa}{\sqrt{\kappa^2+\epsilon}} \right)'' +\frac{\kappa^{3}}{\sqrt{\kappa^{2}+\epsilon}} - \kappa \sqrt{\kappa^{2}+\epsilon}=0 \desda
\left(\frac{\kappa}{\sqrt{\kappa^2+\epsilon}} \right)'' =\epsilon \, \frac{\kappa}{\sqrt{\kappa^2+\epsilon}} \desda \\[8pt]
\kappa''(s)-\kappa^{3} -\frac{3 \kappa (\kappa')^{2}}{\kappa^{2}+\epsilon}=\kappa \epsilon \ ,
\end{array}
\end{equation}
the solution of which is straightforwardly derived from (\ref{zdiv}) by substitution
\begin{equation} \label{zform}
z=\frac{\kappa}{\sqrt{\kappa^2+\epsilon}}
\end{equation}
which gives us
\[
\kappa^{2}(s)= \frac{\epsilon (z(s))^{2}}{1-(z(s))^{2}} \textrm{ for }(z(s))^{2} \leq 1
\]
with $z(s)=\frac{1}{2}(z_0 -\frac{1}{\sqrt{\epsilon}}z_{0}') e^{-\sqrt{\epsilon}s}
+\frac{1}{2}(z_0 +\frac{1}{\sqrt{\epsilon}}z_{0}') e^{\sqrt{\epsilon}s}$, i.e.
\begin{equation}\label{z}
z(s)= z_{0} \cosh(\sqrt{\epsilon}\, s) + \frac{z_0'}{\sqrt{\epsilon}} \sinh (\sqrt{\epsilon} \, s)
\end{equation}
with $z_0=\frac{\kappa_{0}}{\sqrt{\epsilon +\kappa_{0}^{2}}}$, $z_{0}'=\frac{\epsilon \kappa_{0}'}{(\epsilon + \kappa_{0}^{2})^{\frac{3}{2}}}$ which is only valid for
\begin{equation}\label{smax}
s \in [0,s_{\textrm{max}}):=[0, \frac{1}{\sqrt{\epsilon}} \log \left( \frac{1+\sqrt{1- (z_0^{2}-(\epsilon^{-\frac{1}{2}}z_{0}')^{2})}}{z_0 + \epsilon^{-\frac{1}{2}} z_{0}'}\right))
\end{equation}
For a comparison between the elastica and the geodesics derived in this section see Figure \ref{fig:comparison}. Here we recall that the $\epsilon$ of the elastica curves has to be set to
$\epsilon=\beta^2=4\alpha D_{11}$
(recall (\ref{epsilonMumford})), whereas the $\epsilon$ of the geodesics has to be set to
\[
\epsilon=\frac{D_{11}}{D_{22}}
\]
(See (\ref{epsilonCitti}) and see also (\ref{Lagrangianxiands})). Both parameters $\epsilon$ have the physical dimension $[LENGTH]^{-2}$. In our comparison in Figure \ref{fig:comparison} we have set $D_{22}=\frac{1}{4\alpha}$, so that $\epsilon=4 \alpha D_{11}=\frac{D_{11}}{D_{22}}$.

In appendix \ref{ch:app} we derive an exact tangible formula (\ref{solutiongeodesic}) (where the parameters are given by (\ref{pars}))
for the geodesics \emph{in the general case} by means of symplectic differential geometry and Noether's theorem. Moreover in appendix \ref{ch:app}, we will derive an important conservation law (the so-called co-adjoint orbit condition) along the geodesics and we re-derive (\ref{zdiv}) in a shorter and much more structured (but also more abstract) way.

\section{Non-linear adaptive diffusion on orientation scores for qualitative improvements of coherence enhancing diffusion schemes in image processing. \label{ch:CED}}

A scale space representation $u_{f}:\R^d \times \R^{+} \to \R$ of an image $f:\R^{d} \to \R$ is usually obtained by solving an evolution equation on the additive group $(\R^d,+)$. The most common evolution equation, in image analysis, is the diffusion equation,
\[
\left\{
\begin{array}{l}
\partial_{s} u_{f}(\ul{x},s) = \nabla_{\ul{x}} \cdot ( C(u_{f})(\ul{x},s) \, \nabla_{\ul{x}} u_{f})(\ul{x},s) \\
u_{f}(\ul{x},0)=f(\ul{x}),
\end{array}
\right.
\]
where $C: \mathbb{L}_{2}(\R^2 \times \R^{+}) \cap C^{2}(\R^2 \times \R^{+}) \to C^{1}(\R^2 \times \R^{+}))$ is a function which takes care of adaptive conductivity, that is conductivity depending on the local differential structure at $(\ul{x},s , u_f(\ul{x},s))$.
In case $C=1$ the solution is given by convolution $u_{f}(\ul{x},s)=(G_{s}*f)(\ul{x})$ with a Gaussian kernel $G_{s}(\ul{x})= \frac{1}{(4\pi s)^{\frac{d}{2}}} e^{-\frac{\|\ul{x}\|}{4s}}$ with scale, $s=\frac{1}{2} \sigma^{2}>0$.

As pointed out by Perona and Malik \cite{Perona}, non linear image adaptive anisotropic diffusion (diffuse less at locations with strong gradients in the image) are straightforwardly taken into account for by replacing the isotropic generator $\Delta=\nabla \cdot \nabla$ by $\nabla \cdot (c(\|\nabla u_f(\cdot,s)\|) \nabla)$, i.e. $C(u_{f})(\ul{x},s)=c(\|\nabla_{\ul{x}} u_{f}(\ul{x},s)\|)$, where $c:\R^{+} \to \R^{+}$ is some smooth strictly decaying positive function vanishing at infinity. This is based on the intuitive idea that if (locally) the gradient is large you do not want to diffuse too much.
By restricting ourselves to positively valued $c>0$ one ensures that the diffusion is always forward, and thereby ill-posed backward diffusion is avoided.
The most common choices are
\begin{equation} \label{choicec}
c(t)= e^{-\frac{c}{\left( \frac{\lambda}{t} \right)^{2p}}}, \ \  \ c(t)= \frac{1}{\left( \frac{t}{\lambda}\right)^{2p}+1} \textrm{ and } c(t)= \frac{1}{\sqrt{\left(\frac{t}{\lambda}\right)^{2}+1}},
\end{equation}
involving parameters $p>\frac{1}{2}, c, \lambda>0$.
The corresponding flux magnitude functions are given by
\[
\phi(t)=t \,c(t), \textrm{ with }t=\| \nabla u_{f}\|.
\]
The sign of
\begin{equation} \label{aap}
\phi'(t)=c(t)+t c'(t)
\end{equation}
is important, since if $\phi'(t)>0$ then the magnitude $\phi(t)$, $t=\| \nabla u_{f}\|$, of the flux
\begin{equation} \label{flux}
c(\|\nabla u_f\|) \nabla u_f
\end{equation}
(by Gauss Theorem) increases as $\|\nabla U_f\|$ increases, whereas if $\phi'(t)<0$ the magnitude $\phi(\| \nabla u_{f}\|)$ of the flux (\ref{flux}) decreases as $\|\nabla U_f\|$ increases. Typically, this introduces an extra ``sharpening effect'' of lines and edges. However, this sharpening effect (besides the decay of the conductivity function $c:\R^{+}\to \R^{+}$) should not be mistaken for ill-posed backward diffusion because in all cases $c(t)\geq 0$ for all $t>0$. To this end we note that the Perona and Malik equation can be rewritten in Gauge-coordinates $w$ along the normalized gradient $\ul{e}_{w}=\frac{1}{\|\nabla u_f\|}\nabla u_f$ and $v$ along the normalized vector $\ul{e}_{v}=\frac{1}{\|\nabla u_f\|}(-\partial_{y}u_f,\partial_{x}u_f)$ orthogonal to the gradient, using (\ref{aap}): 
\begin{equation} \label{PeronainGauge}
\begin{array}{ll}
\frac{\partial u_f}{\partial s}  =\diev (c(\|\nabla u_f\|)\nabla u_f)
 &= \frac{\partial }{\partial w}\left(c(\frac{\partial u_f}{\partial w}) \frac{\partial u_f}{\partial w}\right) +c(\frac{\partial u_f}{\partial w}) \frac{ \partial^{2} u_f}{\partial v^2} \ \  \desda \\
\frac{\partial u_f}{\partial s}  &= \phi'(\frac{\partial u_f}{\partial w}) \;  \frac{ \partial^{2} u_f}{\partial w^2} + c(\frac{\partial u_f}{\partial w}) \frac{ \partial^{2} u_f}{\partial v^2} \ ,
\end{array}
\end{equation}
with {\small $\frac{ \partial^{2} u_f}{\partial w^2} = \frac{1}{\|\nabla_{\ul{x}} u_f\|^{2} }(\nabla_{\ul{x}} u_f) H_{\ul{x}} [u_f] (\nabla_{\ul{x}} u_f)^{T}$} and {\small $\frac{\partial u_f}{\partial w}=\|\nabla_{\ul{x}} u_f\|$}. Now the coefficient in front of $\frac{ \partial^{2} u_f}{\partial w^2}$ in the righthand-side of (\ref{PeronainGauge}) is negative iff $\phi'(\|\nabla u_f \|)<0$, but this does \emph{not} correspond to inverse diffusion.

The intervals, respective with the choices of $c(t)$ in (\ref{choicec}), where $\phi'$ is positive are
\[
[0,\lambda (2 cp)^{-\frac{1}{2p}}), \  \  [0,\lambda(2p-1)^{-\frac{1}{2p}}) \textrm{ and }\R^{+}.
\]  So the ``sharpening effect'' does not occur in the case $c(t)= \frac{1}{\sqrt{\left(\frac{t}{\lambda}\right)^{2}+1}}$. For the other two choices there is always a danger that the ``sharpening effect'' due to switching sign of $\phi'(t)$ can cause ``staircasing effects'', \cite{Weic98}{ p.52}: That is step-edges will evolve as a staircase over time due to the fact that strong gradients will lead to an effective sharpening of the data whereas weak gradients will lead to relatively smoothing of the data.

A further improvement of the Perona and Malik scheme is introduced by Joachim Weickert, \cite{Weic99b}, who also uses the \emph{direction} of the gradient $\nabla_{\ul{x}}u_{f}$ of $u_{f}$, which is not used in the algorithms of Perona and Malik type. Therefor he proposed the so-called coherence enhancing diffusion schemes (CED-schemes) where the diffusion constant $c$ is replaced by a diffusion matrix:
\begin{equation}\label{CED}
\hspace{-0.5cm}\mbox{}
\begin{array}{l}
S(u_f)(\ul{x},s)= (G_{\sigma} * \nabla u_{f}(\cdot,s) (\nabla u_{f}(\cdot,s))^{T})(\ul{x})    \\
C(u_f)(\ul{x},s)= \alpha I + (1-\alpha) \, e^{-\frac{c}{(\lambda_{1}(S(u_f)(\ul{x},s))-\lambda_{2}(S(u_f)(\ul{x},s)))^{2}}} \,  \ul{e}_{2}(S(u_f)(\ul{x},s))\, \ul{e}_{2}^{T}(S(u_f)(\ul{x},s))
\end{array}
\end{equation}
where $\alpha \in (0,1)$, $c>0, \sigma>0$ are parameters and where the help-matrix $S$, with eigen values $\{\lambda_{i}(S(u_f)(\ul{x},s))\}_{i=1,2}$ is used to get a measure for local anisotropy $e^{-\frac{c}{(\lambda_{1}(S(u_f)(\ul{x},s))-\lambda_{2}(S(u_f)(\ul{x},s)))^{2}}}$ together with an orientation estimate $\ul{e}_{2}(S(u_f)(\ul{x},s))$ which is the eigen vector with smallest eigen value (orthogonal to the average gradient). In order to get robust/reliable orientation estimates it is essential to apply a componentwise smoothing on the so-called ``structure-tensor field'' $\nabla u_{f} \otimes \nabla u_{f}$. The amount of smoothing/averaging of the structure tensor field is determined by $\sigma>0$.

\begin{figure}
\centerline{
\includegraphics[width=0.9\linewidth]{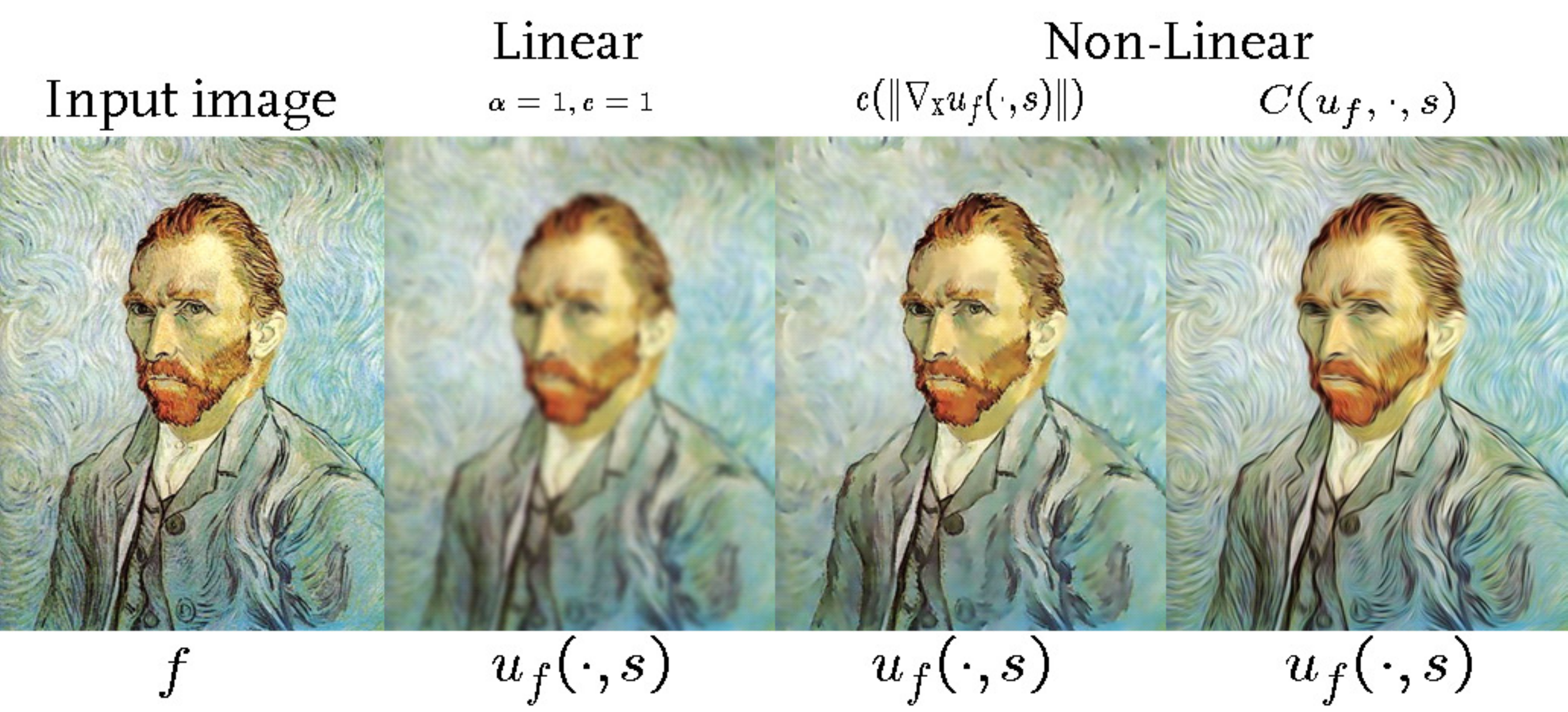}
}
\caption{From left to right: input image $f$ of the well-known portret of Van Gogh, computed on comparable slices $u_{f}(\cdot,s)$ in a linear scale space representation $C=1$, Perona en Malik non-linear scale space representation (left case in (\ref{choicec})) and coherence enhancing diffusion (CED) given by (\ref{CED}) by Weickert, \cite{Weic99b}. }\label{fig:vangogh}
\end{figure}

This lead to useful and visually appealing diffusions of the famous Van Gogh paintings and fingerprint images, see Figure \ref{fig:vangogh}.
Nevertheless, this elegant method fails in image analysis applications with (almost) crossing lines and contours as it starts to create strong artificial curvatures at crossing locations where the gradient is ill-defined.

As this is a major drawback in many (medical) imaging applications we are going to solve this problem by considering
similar non-linear adaptive evolution equations on invertible orientation scores. To this end we note that in invertible orientation scores crossing lines are nicely torn apart in the Euclidean motion group. Moreover, in our orientation scores we have full information on both local direction and local curvature (!) at hand which enables us to steer the diffusions in a left-invariant manner on the orientation scores (and thereby Euclidean invariant manner on images via the unitary wavelet transforms). See Figure \ref{fig:coherence}.

\subsection{Coherence Enhancing Diffusion on Orientation Scores}

\begin{figure}[t]
\centering
\begin{minipage}{0.28\linewidth}
\includegraphics[width=\linewidth]{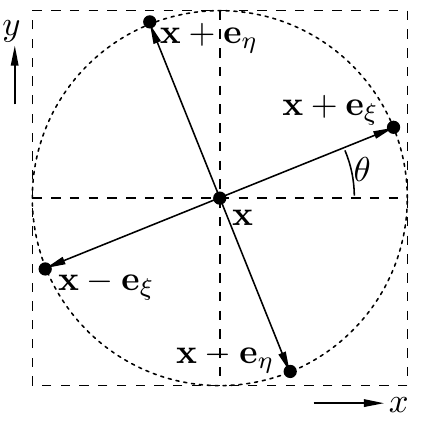}
\end{minipage}
\begin{minipage}{0.5\linewidth}
\begin{equation*}
\tiny
\begin{split}
\partial_\theta u &\approx \frac{1}{2 s_\theta} \left( u(\ul{x},l+1) - u(\ul{x},l-1) \right) \\
\partial_\theta^2 u &\approx \frac{1}{s_\theta^2} \left( u(\ul{x},l+1) - 2 u(\ul{x},l) + u(\ul{x},l-1) \right) \\
\partial_\xi u &\approx \frac{1}{2} \left( u(\ul{x}+\ul{e}_\xi^l, l) - u(\ul{x}-\ul{e}_\xi^l, l)\right) \\
\partial_\xi^2 u &\approx u(\ul{x}+\ul{e}_\xi^l, l) - 2 u(\ul{x}, l) + u(\ul{x}-\ul{e}_\xi^l, l) \\
\partial_\eta u &\approx \frac{1}{2} \left( u(\ul{x}+\ul{e}_\eta^l, l) - u(\ul{x}-\ul{e}_\eta^l, l)\right) \\
\partial_\eta^2 u &\approx u(\ul{x}+\ul{e}_\eta^l, l) - 2 u(\ul{x}, l) + u(\ul{x}-\ul{e}_\eta^l, l) \\
\end{split}
\end{equation*}
\end{minipage}
\begin{equation*}
\tiny
\begin{split}
\partial_{\xi}\partial_{\theta} u &\approx \frac{1}{4 s_\theta} \bigr (u(\ul{x}+\ul{e}_\xi^l, l+1) - u(\ul{x}+\ul{e}_\xi^l, l-1) - u(\ul{x}-\ul{e}_\xi^l ,l+1) + u(\ul{x}-\ul{e}_\xi^l, l-1) \bigr)\\
\partial_{\theta}\partial_{\xi} u &\approx \frac{1}{4 s_\theta} \bigr( u(\ul{x}+\ul{e}_\xi^{l+1}, l+1) - u(\ul{x}+\ul{e}_\xi^{l+1}, l-1) - u(\ul{x}-\ul{e}_\xi^{l-1} ,l+1) + u(\ul{x}-\ul{e}_\xi^{l-1}, l-1) \bigr)\\
\end{split}
\end{equation*}
\caption{Finite difference scheme of (\ref{coherencesimple}) where we use second order B-spline interpolation, \cite{Unser}, for sampling on the grid of our moving frame $\{\ul{e}_{\theta},\ul{e}_{\xi}=\cos \theta \, \ul{e}_{x} +\sin \theta \, \ul{e}_{y},\ul{e}_{\eta}=-\sin \theta \, \ul{e}_{x} +\cos \theta \, \ul{e}_{y}\}$.}\label{fig:Bsplineinterpolation}
\end{figure}

\begin{figure*}[t]
\centering
\begin{tabular}{cccc}
\scriptsize Original & \scriptsize CED-OS $t=30$ & \scriptsize CED $t=30$ \\
\includegraphics[width=0.2\linewidth]{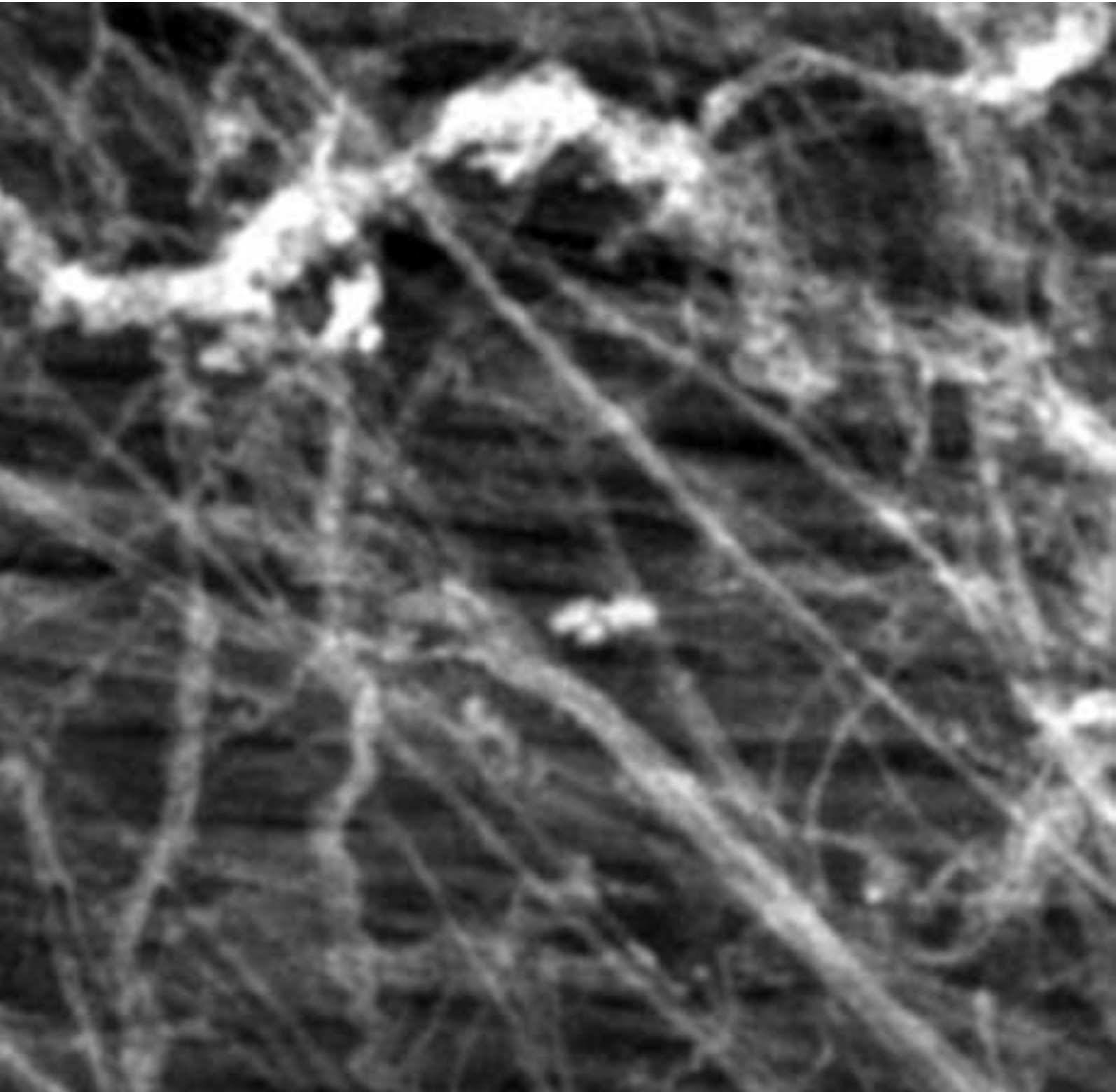} &
\includegraphics[width=0.2\linewidth]{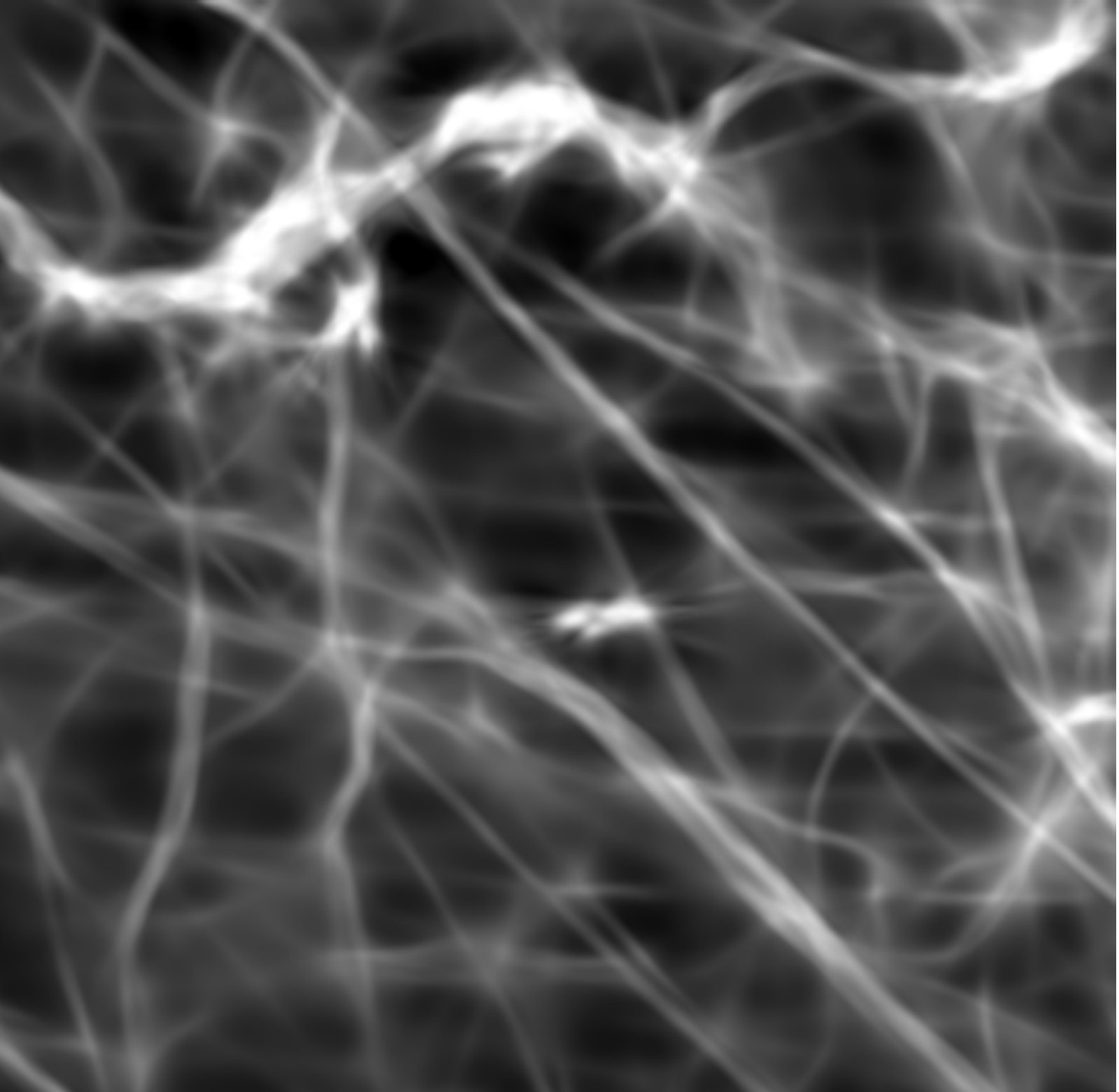} &
\includegraphics[width=0.2\linewidth]{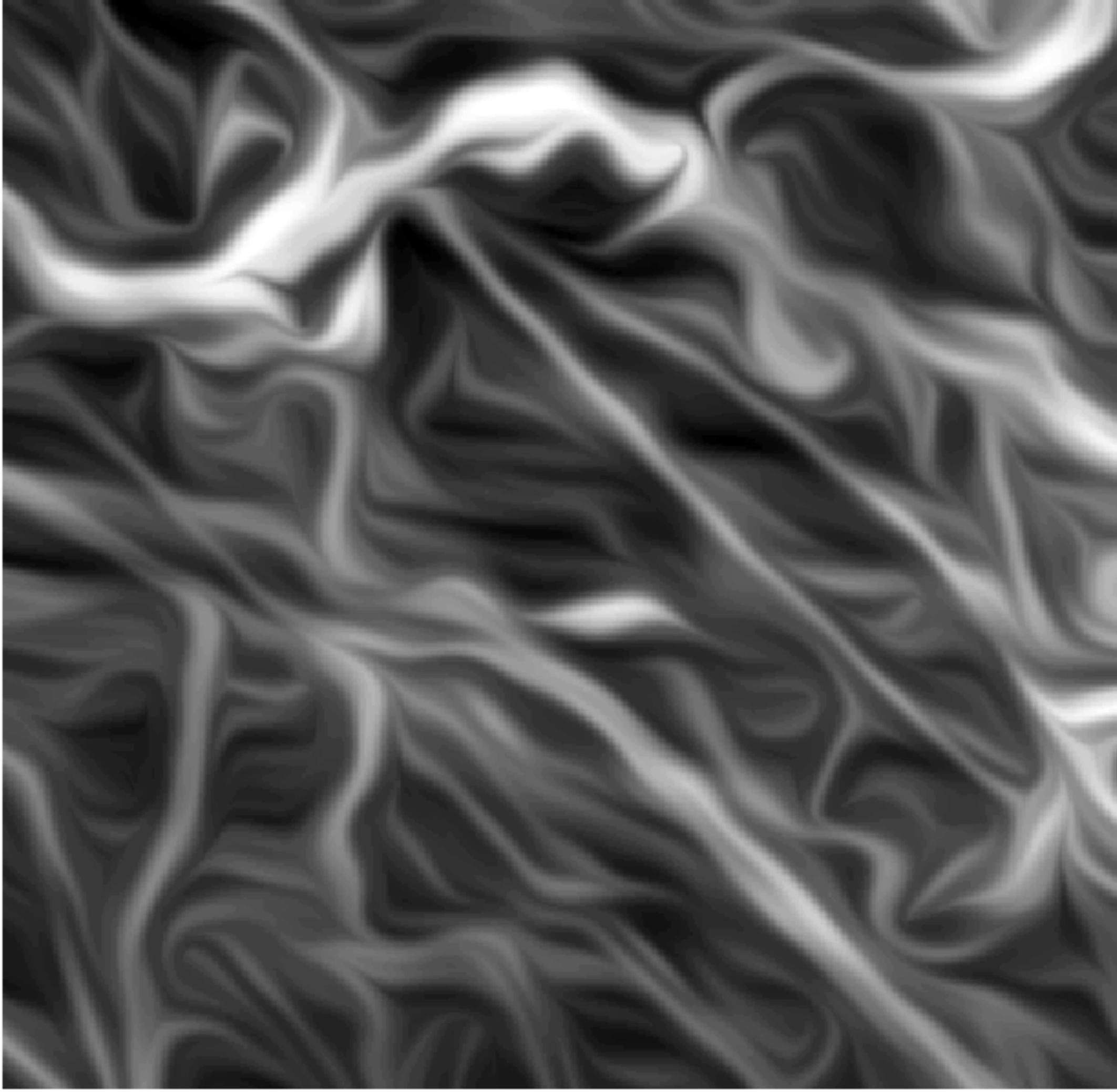}
\end{tabular}
\begin{tabular}{cccc}
\scriptsize Original & \scriptsize CED-OS $t=30$ & \scriptsize CED $t=30$ \\
\includegraphics[width=0.2\linewidth]{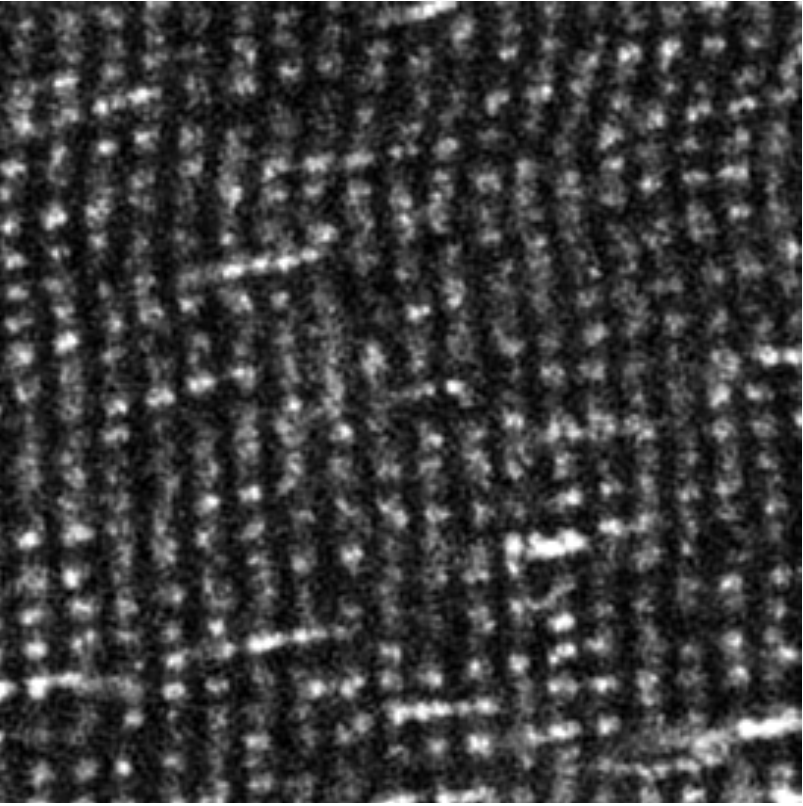} &
\includegraphics[width=0.2\linewidth]{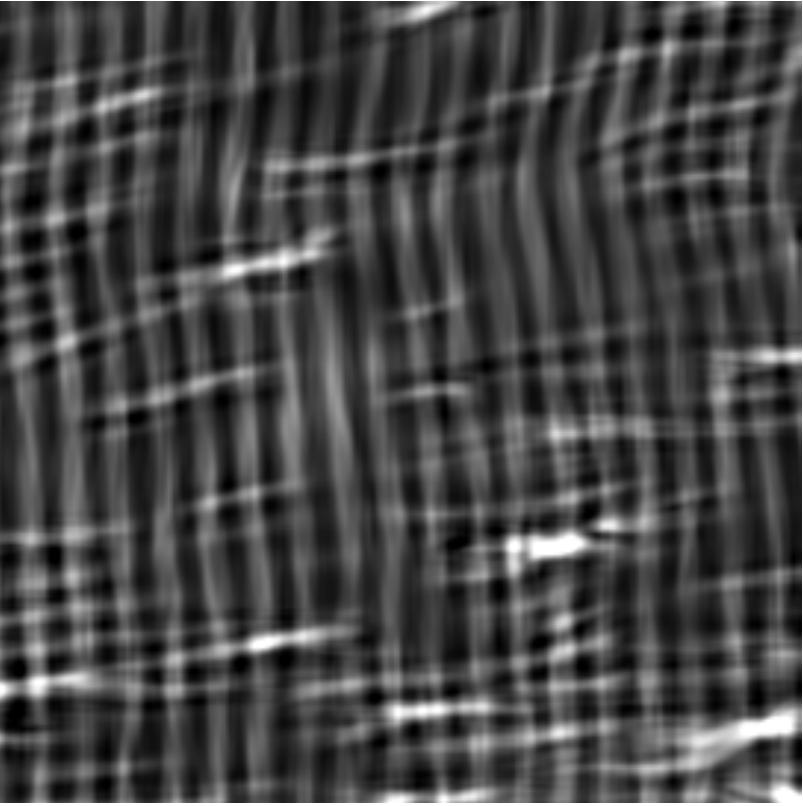} &
\includegraphics[width=0.2\linewidth]{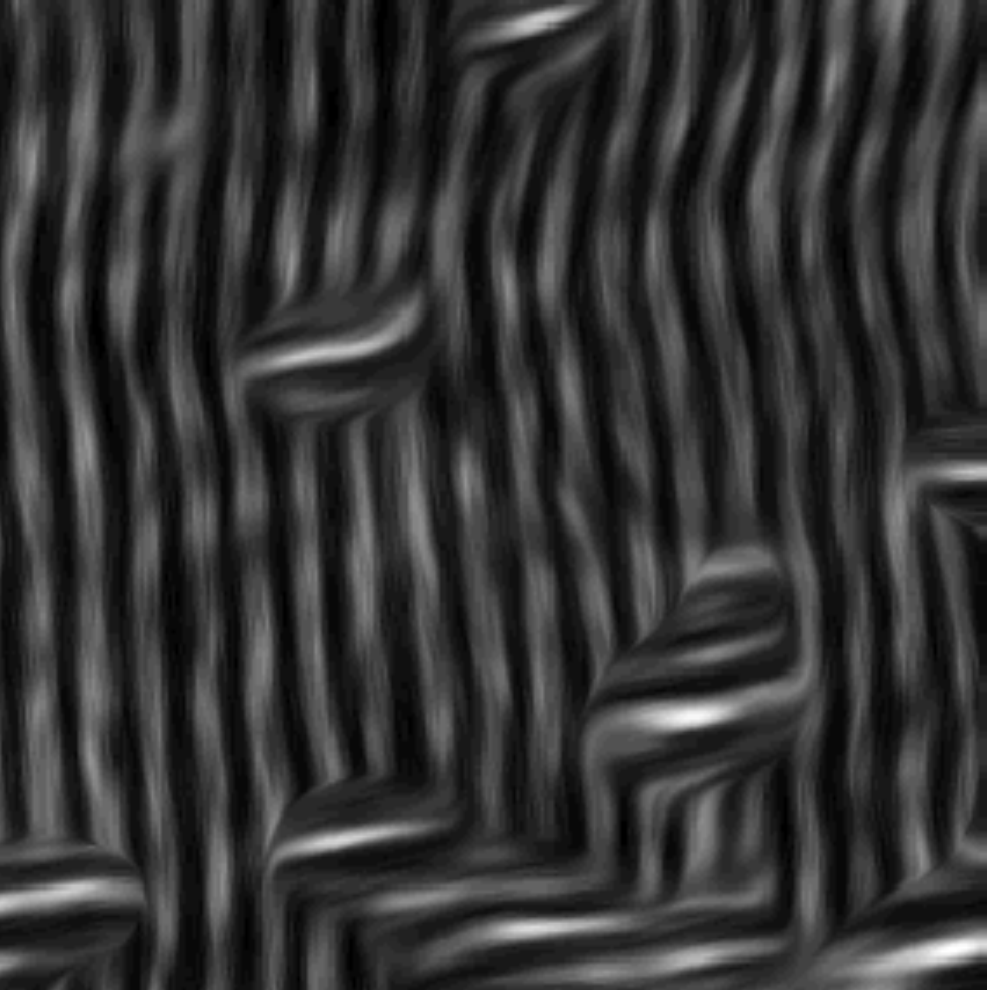}
\end{tabular}
\centerline{\includegraphics[width=0.8\linewidth]{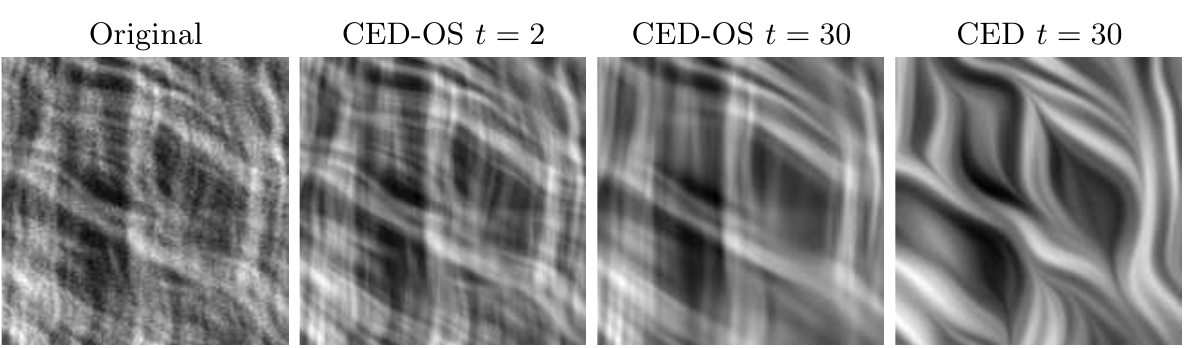}}
\caption{Medical image applications. Top row: Result of coherence enhancing diffusion on orientations scores (CED-OS), see (\ref{coherencesimple}) and (\ref{coherence2}), and standard coherence enhancing diffusion directly on the image, see (\ref{CED}), (CED) of bone-tissue. Middle row: Result of coherence enhancing diffusion on orientations scores (CED-OS) and standard coherence enhancing diffusion directly on the image (CED) of 2-photon microscopy images of a muscle cell. Bottom row coherence enhancing diffusion on orientation scores(CED-OS) and standard coherence diffusion (CED) on medical images of collageen fibers of the heart. All these applications clearly show that coherence enhancing diffusion on orientation scores (CEDOS) properly enhances crossing fibers whereas (CED) fails at crossings: CED creates a ``van Gogh'' type of painting out of the original image $f$. }\label{fig:coherence}
\end{figure*}

In order to obtain adaptive diffusion on orientation scores we will use the following basic non-linear left-invariant evolution equations on $SE(2)$ as a starting point
{\small
\begin{equation}\label{coherencesimple}
\hspace{-0.3cm}
\left\{ \!
\begin{array}{l}
\partial_t U(g,t) = (\!\begin{array}{lll} \beta \partial_\theta  \! \!&\! \! \partial_\xi \! \!&\! \! \partial_\eta \end{array} \!)
\begin{pmatrix}
(D_{11}(U))(g,t)  & 0  & 0 \\
0 & (D_{22}(U))(g,t) & 0 \\
0 & 0 & (D_{33}(U))(g,t)
\end{pmatrix}
\begin{pmatrix} \beta \partial_\theta \\ \partial_\xi \\ \partial_\eta \end{pmatrix} U(g,t), \\
\textrm{ for all } g \in SE(2),t>0,  \\[8pt]
U(g,t=0)=\mathcal{W}_{\psi}[f](g) \textrm{ for all }g \in SE(2),
\end{array}
\right.
\end{equation}
}
with $\beta>0$ (recall \ref{metric}) and
where the functions $D_{kk}:\mathbb{L}_{2}(SE(2) \times \R^{+}) \cap C^{2}(SE(2) \times \R^{+}) \to C^{1}(SE(2) \times \R^{+})$, $k=1,2,3$ given by
\[
(g,t) \mapsto (D_{kk}(U))(g,t) \geq 0, \qquad U \in \mathbb{L}_{2}(SE(2) \times \R^{+}),
\]
should be chosen dependent on the local Hessian $HU(\cdot,t)$ of $U(\cdot,t)$ (similarly as was done in the CED-scheme (\ref{CED})) such that at strong orientations $D_{33}$ should be small so that we have anisotropic diffusion in the spatial plane along the preferred direction $\partial_{\xi}$, while at weak directions $D_{33}$ and $D_{22}$ should be relatively large and isotropic $D_{22}\approx D_{33}$. Usually we set $D_{22}(U)(g,t)=1$, since in general there is no reason to make $D_{22}(U)(g,t)$ dependent on $g$ and in such cases a simple re-parametrization of time yields $D_{22}=1$. \\
\textbf{Example 1:} \\
For example one can take $D_{22}(U)(g,t)=1$, $D_{11}(U)(g,t)=D_{33}(U)(g,t)=e^{-\frac{(s(|U|)(g,t))^2}{c}}$, where $c>0$ is a standard (Perona \and Malik-)parameter where $s(U)(g,t)$ is a measure for orientation
strength like
\begin{equation}\label{orstrength}
s(U)(g,t)=\max (-\textrm{Re}(\lambda_{1}(H |U|(g,t))),0),
\end{equation}
where $\lambda_{1}(g,t)$ is the largest eigenvalue
of the Hessian $H |U(\cdot,t)|(g)= [\mathcal{A}_{j}\mathcal{A}_{i}|U(\cdot,t)|](g)$, $i=1,\ldots 3, j=1,\ldots,3$ where $i$ is the row index. Here we stress that we take the Hessian of the absolute value $|U(\cdot,t)|$, since the absolute value of an orientation score is phase invariant, i.e. it does not matter if you are on top of a line (large real part) or on the edge of a line (large imaginary part), recall Figure \ref{fig:cakeexample} (d) and recall subsection \ref{ch:curvest}. In practice we use Gaussian
derivatives (\ref{GD}), rather than usual derivatives, of $|U(\cdot,t)|$ (so isotropic with scale $s_{1}$ in the spatial part $\equiv\R^2$ and scale $s_{2}$ in the angular part $\equiv\mathbb{T}$ with periodic boundary conditions) at small scales $s_{1}=\beta^2 s_{1}, s_{2}>0$ (typically $\sqrt{2\, s_{1}}$ is in the order of say 2 pixels and $\sqrt{2 \, s_{2}}$ is in the order of say $\frac{2\pi}{16}$). See Figure \ref{fig:coherence}. \\
\textbf{Example 2:}\\
Another modification (or rather slight improvement) in the non-linear diffusion system (\ref{coherencesimple}) is obtained by replacing the orientation strength $s(U)(g,t)$ (\ref{orstrength}) in the first example by
\begin{equation} \label{orstrength2}
s(g,t)=\max \{- \sum \limits_{i=1}^{2} (\tilde{\ul{e}}_{i}^{O}(g,t))^{T}\; (H_{\beta}|U(\cdot,t)|(g))^{T}H_{\beta}|U(\cdot,t)|(g)\; \tilde{\ul{e}}_{i}^{O}(g,t)  \; , \;0\}
\end{equation}
where we recall that the symmetric Hessian $H_{\beta}|U(\cdot,t)|$ was given by (\ref{hbu}) and where $\tilde{e}_{1}^{O}(g,t)$ and $\tilde{e}_{2}^{O}(g,t)$ denote the remaining eigenvectors (with largest two eigen values) of the symmetric matrix $H_{\beta}|U(\cdot,t)|(g))^{T}H_{\beta}|U(\cdot,t)|(g)$:
\[
\begin{array}{l}
(H_{\beta}|U(\cdot,t)|(g))^{T}H_{\beta}|U(\cdot,t)|(g)\; \tilde{\ul{c}}(g,t)=\lambda_{0} \tilde{\ul{c}}(g,t) \\
(H_{\beta}|U(\cdot,t)|(g))^{T}H_{\beta}|U(\cdot,t)|(g)\; \tilde{\ul{e}}_{1}^{O}(g,t)=\lambda_{1} \tilde{\ul{e}}_{1}^{O}(g,t) \\
(H_{\beta}|U(\cdot,t)|(g))^{T}H_{\beta}|U(\cdot,t)|(g)\; \tilde{\ul{e}}_{2}^{O}(g,t)=\lambda_{2} \tilde{\ul{e}}_{1}^{O}(g,t), \qquad |\lambda_{0}|\leq |\lambda_{1}|\leq |\lambda_{2}|.
\end{array}
\]
So the righthand side of (\ref{orstrength2}) is to be considered as the Laplacian in the tangent plane orthogonal to the vector $\tilde{c}(g,t)$. Here we recall from section \ref{ch:curvest} that $\tilde{c}(g,t)$ corresponds to the best exponential curve fit $\gamma(s)= g \, e^{s \sum \limits_{i=1}^{3} \tilde{c}_{i}(g) A_i}$ to the data $|U(\cdot,t)|$. \\
\\
With respect to the numerics of (\ref{coherencesimple}) and (\ref{coherence2}), we implemented a forward finite difference scheme using central differences along the moving frame $\{\theta,\xi,\eta\}$ where we used 2nd order $B$-spline interpolation, \cite{Unser}, to get the equidistant samples on the $\{\xi,\eta,\theta\}$-grid from the given samples on the $\{x,y,\theta\}$-grid, see figure \ref{fig:Bsplineinterpolation}, thereby our method is second order accurate on $SE(2)$. As our algorithm (for details see \cite{Fran06b}) is only first order accurate in time, we took small time steps in our contour-enhancement experiments. With this respect we note that a Crank-Nickolson scheme for time integration is second order in time and can improve computation time.

Another issue for reduction of computation time is the time dependent conductivity matrix, which from a strict point of view needs to be updated at each time step of the evolution. In practice, usually the updating of the conductivity matrix in our finite difference scheme does not have to be done at every time step and even the linear case where
the conductivity matrix is not updated at all (so the matrix is determined only by the absolute value of the initial condition $|U(\cdot,0)|=|\mathcal{W}_{\psi}f|$) yields good results.

\subsubsection{Including adaptive curvatures in the diffusion scheme using Gauge-coordinates \label{ch:gauge}}

In subsection \ref{ch:curvest} we discussed two methods of how to obtain curvature estimates in orientation scores. This was done by finding the best exponential curve fit to the absolute value of the orientation score (which is phase invariant, recall Figure \ref{fig:cakeexample} (d)). We distinguished between two approaches. In the first approach we considered the best \emph{horizontal} exponential curve fit to the data (\ref{approach1}), whereas in the second approach (\ref{approach2}) we considered the best exponential curve fit to the absolute value of the orientation score. Both approaches yield a curvature estimate which in this paragraph we assume to be given. We shall write $(\kappa_{est}(|U|))(g,t)$ for the curvature estimate of the score $U$ via its absolute value $|U|$ at location $g\in SE(2)$ at time $t>0$. Since we only want to include curvature at strongly oriented structures we shall multiply it with a front factor:
\[
\kappa(|U|)(g,t)=\left(1-e^{-\left(\frac{d_{\kappa}}{D_{33}(g,t)}\right)^{4}}\right) \kappa_{est}(|U|)(g,t),
\]
where $d_{\kappa}$ controls the soft threshold on including the curvature estimate
$\kappa_{est}(|U|)(g,t)$.

\begin{figure}
\centerline{
\includegraphics[width=1\hsize]{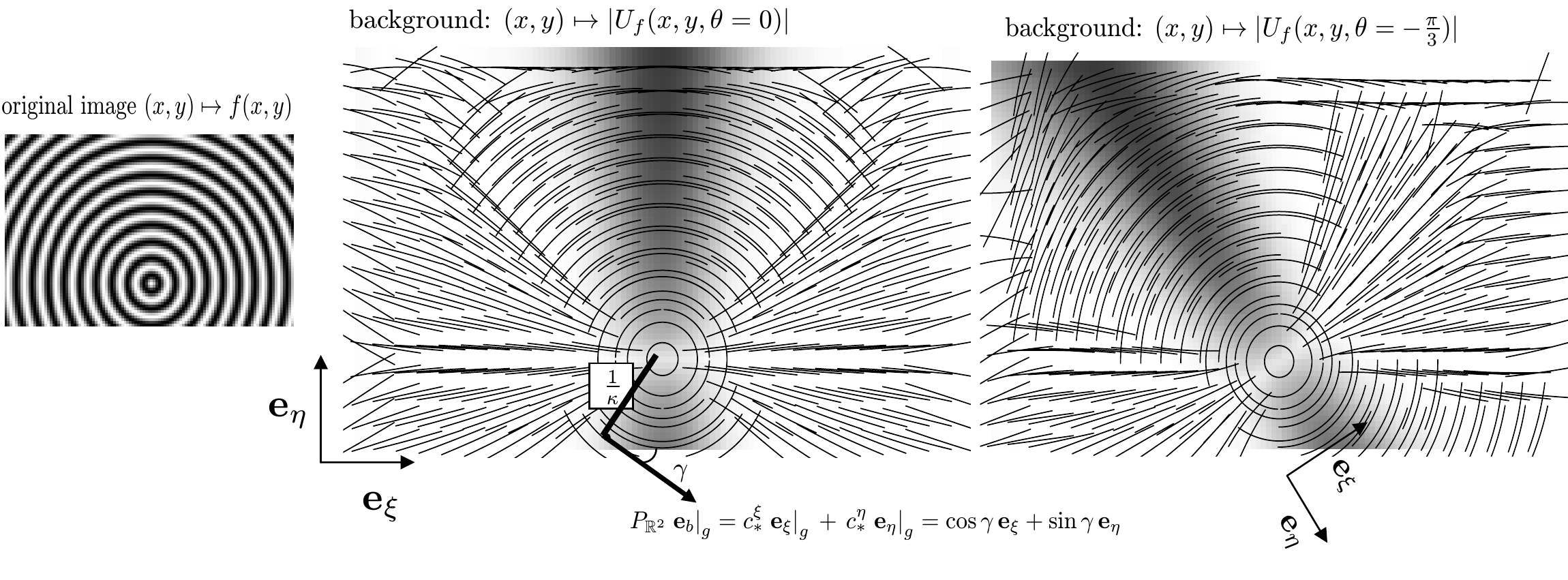}
}
\caption{Illustration of the projection $P_{\R^2}\ul{e}_{b}$ of the vector field $\ul{e}_{b}=\partial_{b}$ on the spatial plane (image-plane) plotted on fixed
orientation layers $|U|(\cdot,\theta)$ of the absolute value $|U|$ of the orientation score $U$. Along this vector field we plotted circular arcs to also include the curvature $\kappa$. These circular arcs correspond to the projections of best exponential curve fits to the data $|U|$ at each location $g$ on the image plane.
From left to right, the original image, and plots for $\theta=\frac{-\pi}{2}$ and $\theta=-\frac{\pi}{3}$. Note that at positions $g \in SE(2)$ near strongly oriented structures, the $\ul{e}_{b}$ vector is better aligned with the local image structure than the $\ul{e}_{\xi}$ vector.} \label{fig:dh}
\end{figure}

Now we can include curvature in our scheme (\ref{coherencesimple}) by replacing $\partial_{\xi} \mapsto \partial_{\xi} +\kappa \partial_{\theta}$. To this end we recall that the exponential curve $s \mapsto e^{{s \left. (\partial_{\xi}+\kappa \partial_{\theta})\right|_{e}}}=e^{s \partial_{x}+\kappa \partial_{\theta}}$ yields a circular spiral (\ref{circspiral}) whose projection on $\R^2$ is a circle with radius $|\kappa|^{-1}$ if $\kappa$ is constant. Moreover, along horizontal curves we have
\[
\frac{d}{ds}U(\gamma(s))=\kappa(s) \frac{\partial U}{\partial \theta}(\gamma(s))+ \frac{\partial U}{\partial \xi}(\gamma(s))
\]
where $\kappa(s)=\frac{d\theta}{ds}$ and $\langle {\rm d}\xi, \dot{\gamma}(s)\rangle=1$, see Appendix \ref{ch:app3}.
See Figure \ref{fig:GFs}.
\begin{figure}
\centerline{
\includegraphics[width=0.2\hsize]{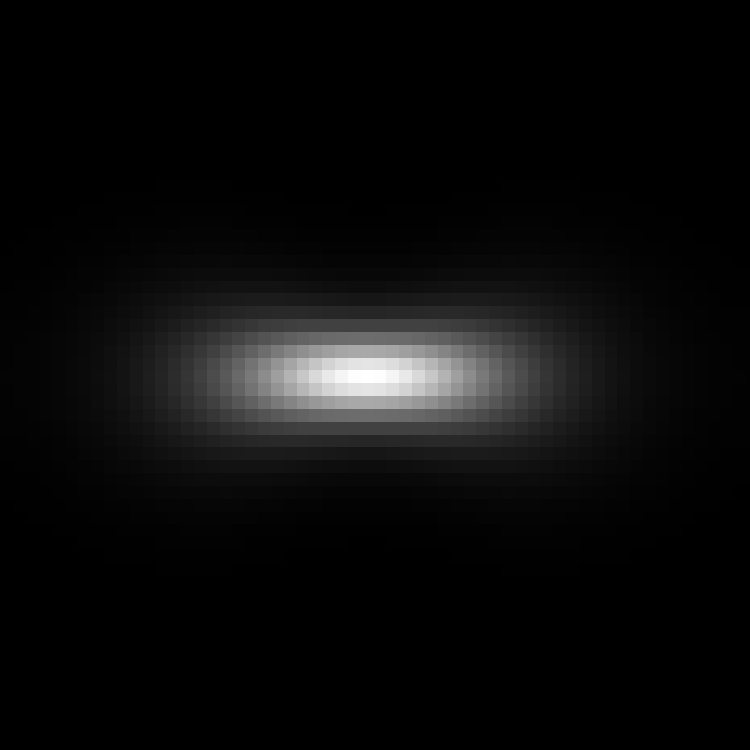}
\includegraphics[width=0.25\hsize]{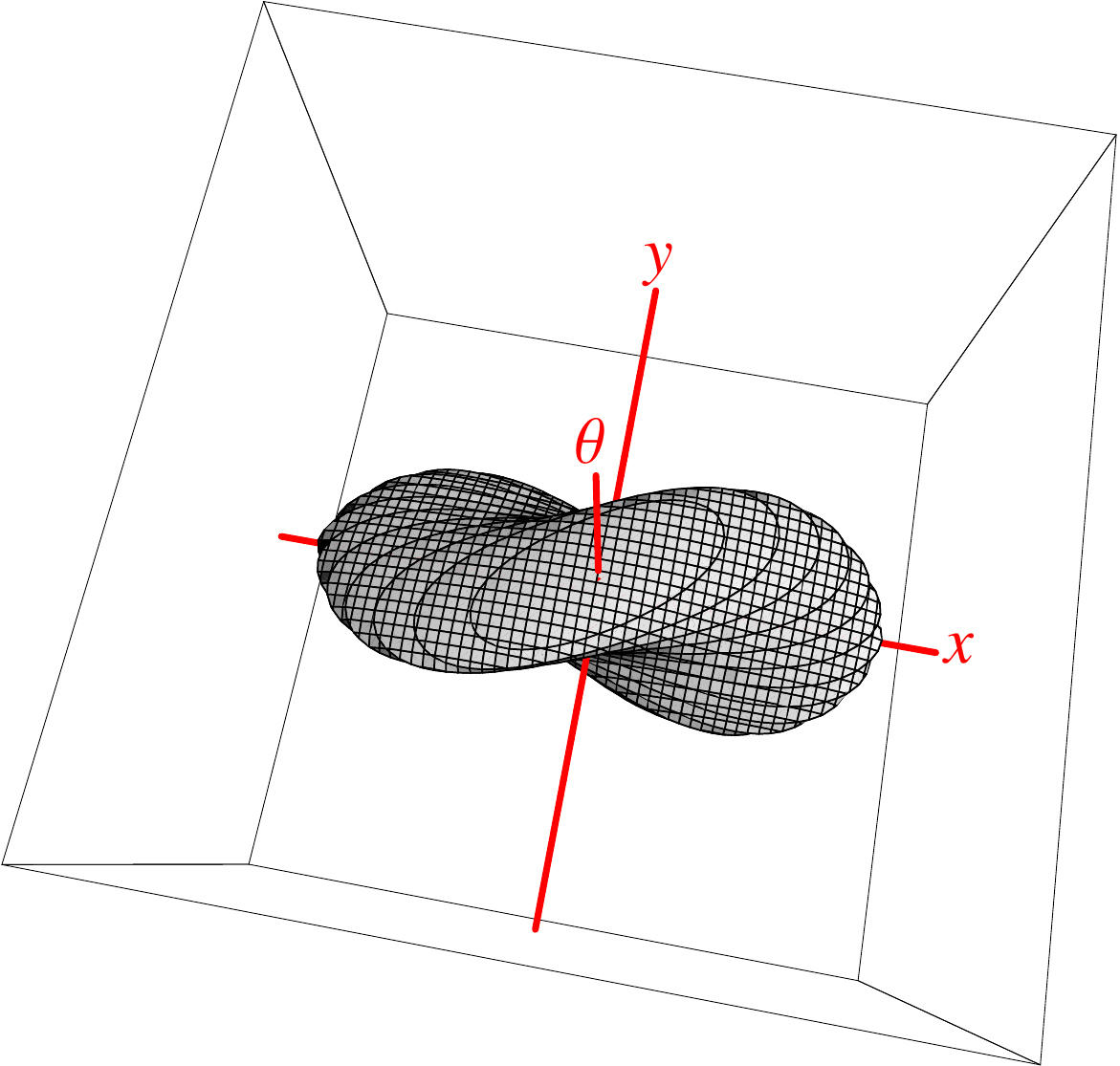}
\hfill
\includegraphics[width=0.2\hsize]{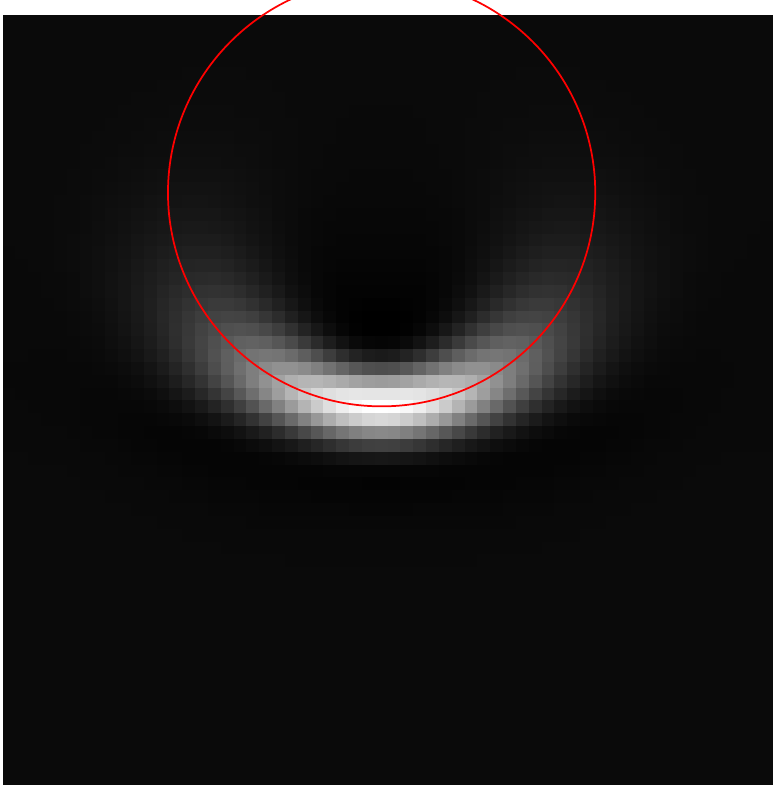}
\includegraphics[width=0.25\hsize]{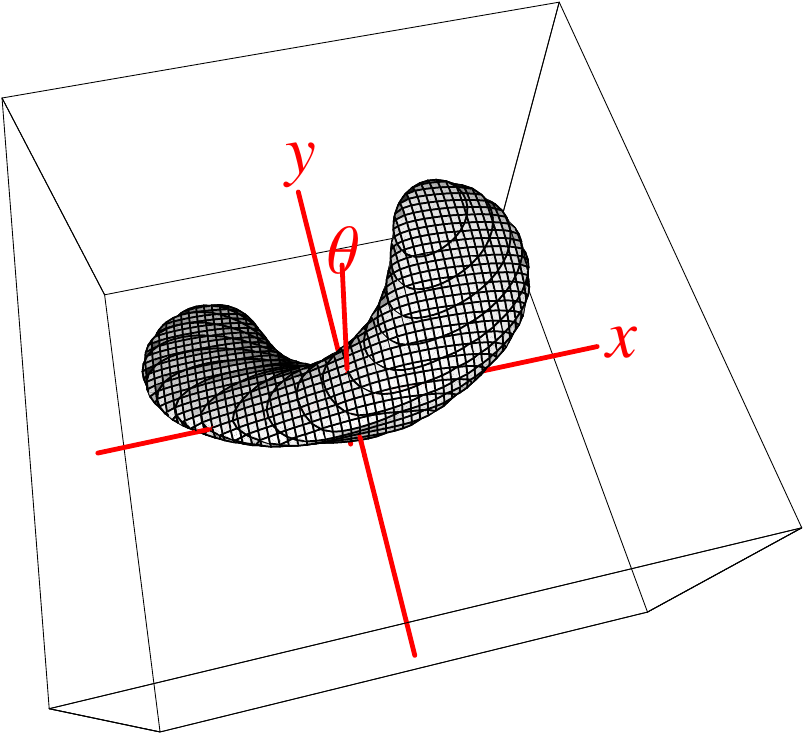}
}
\caption{Illustrations of the heat-kernels $K_{t}^{D}:SE(2) \to \R^{+}$ on $SE(2)$. Left: $D=\textrm{diag}\{D_{11},D_{22},0\}$ in left-invariant coordinate frame $\{\partial_{\theta},\partial_{\xi},\partial_{\eta}\}$. Right: $D=\textrm{diag}\{0,D_{bb},0\}$ in gauge-coordinate frame $\{\partial_{a},\partial_{b},\partial_{c}\}$ with $\gamma=0$, $\kappa=0.06$ , $\beta=1$, $D_{bb}=1.0036$ and $t=70$. }\label{fig:GFs}
\end{figure}
Here we should be careful since $\{\partial_{\theta}, \partial_{\xi} +\kappa \partial_{\theta}, \partial_{\eta}\}$ are (in contrast to $\frac{1}{\beta}\{\beta \partial_{\theta},\partial_{\xi},\partial_{\eta}\}$) no longer orthonormal with respect to the $(\cdot,\cdot)_{\beta}$ inner product. Therefore we are going to introduce the gauge coordinates, aligned with the optimally fitting exponential curve
\[
s \mapsto g \, \exp(s \sum \limits_{i=1}^{3} c_{*}^i(g,t) A_i), \qquad \ul{c}_{*}(g,t)=(c^{\theta}_{*}(g,t),c^{\xi}_{*}(g,t), c^{\eta}_{*}(g,t)) \in \R^3,
\]
with $\|\ul{c}_{*}\|_{\beta}=(c^{\theta}_{*})^2+\beta^{2}(c^{\xi}_{*})^2+\beta^{2}(c^{\eta}_{*})^2=1$,
to the orientation score data $|U(\cdot,t)|$ at position $g \in SE(2)$ at time $t>0$.
These Gauge coordinates are given by
\[
\left\{
\begin{array}{l}
\partial_{a}= \beta^{2}\sqrt{(c_{*}^{\xi})^2+(c_{*}^{\eta})^2} \partial_{\theta} -\frac{c^{\theta}_{*} c^{\xi}_{*}}{\sqrt{(c_{*}^{\xi})^2+(c_{*}^{\eta})^2}} \partial_{\xi}-\frac{c^{\theta}_{*} c^{\eta}_{*}}{\sqrt{(c_{*}^{\xi})^2+(c_{*}^{\eta})^2}} \partial_{\eta} \\
\partial_{b}=\beta (c_{*}^{\xi}\partial_{\xi}+c_{*}^{\eta}\partial_{\eta} + c_{*}^{\theta}\partial_{\theta}) \\
\partial_{c}= \frac{-c^{\eta}_{*}}{\sqrt{(c^{\xi}_{*})^2+(c^{\eta}_{*})^2}} \partial_{\xi} + \frac{c^{\xi}_{*}}{\sqrt{(c^{\xi}_{*})^2+(c^{\eta}_{*})^2}} \partial_{\eta}\ .
\end{array}
\right.
\]
Note that the gauge-vector is along the best exponential curve-fit direction, i.e. $\partial_{b} = \ul{c}_{*}$
and note that the span of the tangent vectors $\{\partial_{a},\partial_{b}\}$ corresponds with
$\textrm{span}\{\partial_{a},\partial_{c}\} \equiv (\ul{c}_{*})^{\bot}$.
\begin{figure}
\centerline{
\includegraphics[width=0.45\hsize]{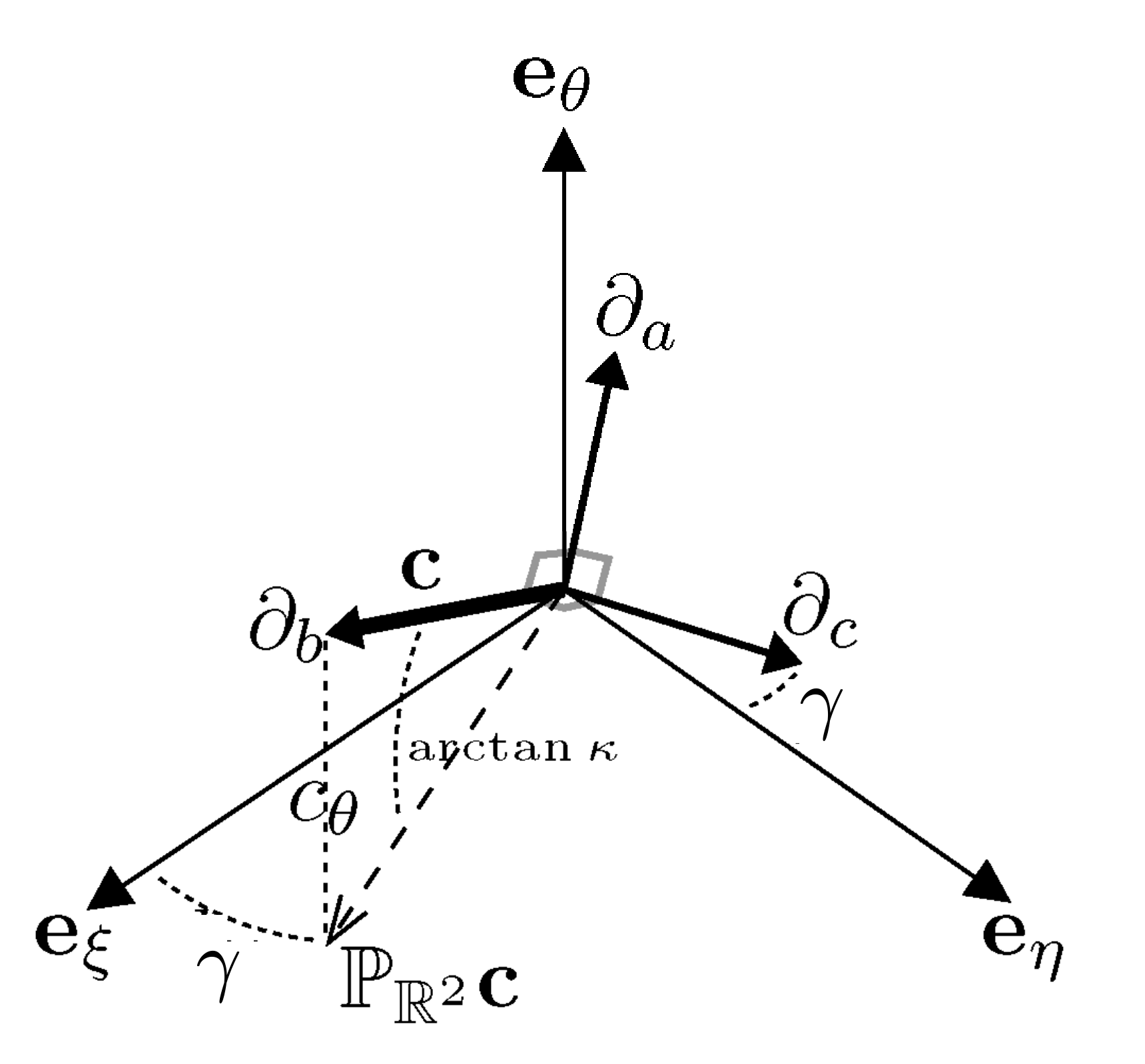}}
\caption{The gauge coordinate frame $\{\partial_{a},\partial_{b},\partial_{c}\}$ illustrated with respect to the basis of left-invariant vector fields $\{\partial_{\theta},\partial_{\xi},\partial_{\eta}\}$. Here we note that the curvature estimation $\kappa$ is given by (\ref{approach2}) and $\{\partial_{a},\partial_{b},\partial_{c}\}$ are given by (\ref{gauge}), where $\partial_{b}$ is the direction determined by the optimal exponential curve fit (\ref{minprob}) to the data. The angle $\gamma=\gamma(U)(g,t)=\arg (c_{*}^{\xi}(g,t) + i \, c_{*}^{\eta}(g,t)) $ intuitively tells us how ``horizontal'' the orientation score is at location $g \in SE(2)$ at time $t>0$. }\label{fig:gaugecoordinates}
\end{figure}
For geometric understanding it helps to consider the Gauge tangent-vectors in ball-coordinates with respect to the basis of left-invariant vector fields $\{\partial_{\theta},\partial_{\xi},\partial_{\eta}\}$ so that it becomes obvious which rotation in $SO(3)$ (or rather which class of rotations in $SO(3)/SO(2)\equiv S^{2}$, if we do not distinguish between directions in plane $(\ul{c}_{*})^{\bot}$) is required to map the standard left invariant basis $\{\partial_{\theta},\partial_{\xi},\partial_{\eta}\}$ into the basis Gauge-coordinates. See Figure \ref{fig:gaugecoordinates}. The Gauge-coordinates in ball-coordinates read
\begin{equation} \label{gauge}
\left\{
\begin{array}{l}
\partial_{a}=-\cos \alpha \, \cos \gamma \, \partial_{\xi} - \cos \alpha \, \sin \gamma \, \partial_{\eta} + \beta \sin \alpha \, \partial_{\theta} \ ,\\
\partial_{b}= \sin \alpha  \, \cos \gamma \, \partial_{\xi} + \sin \alpha \, \sin \gamma \, \partial_{\eta} + \beta \cos \alpha \partial_{\theta} \ ,\\
\partial_{c}= - \sin \gamma \, \partial_{\xi} + \cos \gamma \, \partial_{\eta}\ ,
\end{array}
\right.
\end{equation}
where the Euler-angles read
\[
\begin{array}{l}
\alpha = \arccos c_{*}^{\theta} = \arccos \frac{\kappa}{\sqrt{\kappa^2+\beta^2}}\ , \\
\gamma = \arg (c_{*}^{\xi} + i \, c_{*}^{\eta}).
\end{array}
\]
Here the function $\gamma$ which maps $U$ to $\gamma(U)(g,t)=\arg (c_{*}^{\xi}(g,t) + i \, c_{*}^{\eta}(g,t))$
See Figure \ref{fig:dh}.

Now the diffusion generator diagonal along the left-invariant Gauge vector fields is given by
\[
\begin{array}{l}
D_{aa} (\partial_{a})^2 + D_{bb} (\partial_{b})^2 + D_{cc} (\partial_{c})^2 =
(
\begin{array}{ccc}
\beta \partial_{\theta} & \partial_{\xi} & \partial_{\eta}
\end{array}
)
M_{\alpha,\gamma}^{T} \left( \begin{array}{ccc} D_{aa} & 0 & 0 \\ 0 & D_{bb} & 0 \\ 0 & 0 & D_{cc}
\end{array}
\right)
M_{\alpha,\gamma}
\left(
\begin{array}{l}
\beta \partial_{\theta} \\
\partial_{\xi} \\
\partial_{\eta}
\end{array}
\right) .
\end{array}
\]
where {\scriptsize $M_{\alpha,\gamma}^{T}=
\left(
\begin{array}{ccc}
\sin \alpha & -\cos \alpha \cos \gamma & -\cos \alpha \sin \gamma \\
\cos \alpha & \cos \gamma \sin \alpha & \sin \alpha \sin \gamma \\
0 & -\sin \gamma & \cos \gamma
\end{array}
\right)$ }is the rotation matrix in $SO(3)$ which maps the tangent vector $\beta \partial_{\theta}$ to the tangent vector $\partial_{b}$.

Now by straightforward computation this leads to the following non-linear evolution equations on orientation scores
{\scriptsize
\begin{equation}\label{coherence2}
\hspace{-0.3cm}
\left\{ \!
\begin{array}{l}
\partial_t U(g,t) = (\!\begin{array}{lll} \beta \partial_\theta  \! \!&\! \! \partial_\xi \! \!&\! \! \partial_\eta \end{array} \!)
\frac{1}{\kappa^2+\beta^{2}} \times  \\
{\tiny
\left(
\begin{array}{ccc}
D_{bb}\kappa^2+D_{aa}\beta^2  & \kappa \beta (D_{bb}\!-\!D_{aa})\cos \gamma  & \kappa \beta (D_{bb}\!-\!D_{aa}) \sin \gamma \\
\kappa \beta (D_{bb}\!-\!D_{aa})\cos \gamma &  D_{cc} (\kappa^{2} \!+\!\beta^{2}) +((D_{bb}\!-\! D_{cc})\beta^{2}+(D_{aa}\!-\!D_{cc})\kappa^{2}) \cos^{2}\gamma & \cos \gamma  \sin \gamma  (\kappa^{2}(D_{aa} \!- \!D_{cc})+\beta^{2}(D_{bb} \!- \!D_{cc}))\\
\kappa \beta (D_{bb}\!-\!D_{aa}) \sin \gamma  &  \cos \gamma  \sin \gamma  (\kappa^{2}(D_{aa} \!- \!D_{cc})+\beta^{2}(D_{bb} \!- \!D_{cc})) & D_{bb}\beta^{2}\!+\!D_{aa}\kappa^{2} +((D_{cc} \!- \!D_{bb})\beta^{2} \!+\!(D_{cc}\!-\!D_{aa})\kappa^{2})\cos^{2}\gamma
\end{array}
\right)
}
\\
\times \left(
\begin{array}{l}
 \beta \partial_\theta \\ \partial_\xi \\ \partial_\eta \end{array} \right) U(g,t), \qquad \textrm{ for all } g \in SE(2),t>0,  \\[8pt]
U(g,t=0)=\mathcal{W}_{\psi}[f](g) \textrm{ for all }g \in SE(2),
\end{array}
\right.
\end{equation}
}
where for the sake of clarity we used short notation $D_{ii}=(D_{ii}(U))(g,t)$, for $i=a,b,c$.
Now again we set
\[
D_{bb}=1 \textrm{ and }(D_{aa}(U))(g,t)=(D_{cc}(U))(g,t)= e^{-\frac{(s(|U|)(g,t))^2}{c}}, c>0.
\]
Here we take (\ref{orstrength2}) as a measure for orientation strength. In the Gauge coordinates this measure can be written
\[
s(g,t)=\max\{-\Delta_{\ul{c}_{*}^{\bot}} |U(\cdot,t)|(g),0 \} =\max\{-\left((\partial_{a})^{2}|U(\cdot,t)|+(\partial_{c})^{2}|U(\cdot,t)|\right)(g),0\}.
\]
and the conductivity matrix in (\ref{coherence2}) simplifies to
\[
{\small
\frac{1}{\beta^{2}+\kappa^{2}}
\left(
\begin{array}{ccc}
\kappa^2+D_{aa}\beta^2  & \kappa \beta (1-\!D_{aa})\cos \gamma  & \kappa \beta (1\!-\!D_{aa}) \sin \gamma \\
\kappa \beta (1\!-\!D_{aa})\cos \gamma &  D_{aa} (\kappa^{2} \!+\!\beta^{2}) +(1\!-\! D_{aa})\beta^{2} \cos^{2}\gamma & \cos \gamma  \sin \gamma  \beta^{2}(1 \!- \!D_{aa})\\
\kappa \beta (1\!-\!D_{aa}) \sin \gamma  &  \cos \gamma  \sin \gamma  \beta^{2}(1 \!- \!D_{aa}) & \beta^{2}\!+\!D_{aa}\kappa^{2} +(D_{aa} \!- \!1)\beta^{2}  \cos^{2}\gamma
\end{array}
\right)
}.
\]
See Figure \ref{fig:GFs} for an illustration of the special case $D_{bb}$ is constant, $Daa=D_{cc}=0$, $\gamma=0$, which despite the strong degree of degeneracy still leads to a smooth and useful Green's function since the H\"{o}rmander condition, recall subsection \ref{ch:hoermander}, is satisfied. \\ \\
\textbf{Remark:} \\
Although not discussed here it is worthwhile to consider the components $\{C^{ij}(U)(g,s)\}$ of the adaptive conductivity matrix with respect to the basis of left-invariant vector fields $\{\mathcal{A}_{1},\mathcal{A}_{2},\mathcal{A}_{3}\}:=\{\partial_{\theta},\partial_{\xi},\partial_{\eta}\}$ on $SE(2)$ as the inverse-components of a metric attached to the graph of $g \mapsto U(g,s)$. So rather than changing the constant conductivity by an orientation score adaptive conductivity within the diffusion equation (as we did in this section) we could change the constant left-invariant metric (\ref{metric}) by an orientation score adaptive metric $G(U)(g,t)= \sum \limits_{i=1}^{3}\sum \limits_{j=1}^{3}G_{ij}(U)(g,t){\rm d}\mathcal{A}^{i} \otimes {\rm d}\mathcal{A}^{j}$, where $G_{ij}(U)(g,t)$ are the components of the inverse matrix of $[G^{ij}(U)(g,t)]:=[C^{ij}(U)(g,t)]$. So in stead of (\ref{coherence2}) we could consider the Laplace-Beltrami flow
{\small
\[
\left\{ \!
\begin{array}{l}
\partial_t U(g,t) = \frac{1}{\sqrt{\det{G(U)(g,t)}}} \sum \limits_{i=1}^{3}\sum \limits_{j=1}^{3}\mathcal{A}_{i}\left\{\sqrt{\det{G}(U)(g,t)} \, G^{ij}(U)(g,t) \, \mathcal{A}_{j}  U(g,t)\right\}(g,t)  \   \textrm{ for all } g \in SE(2),t>0,  \\[8pt]
U(g,t=0)=\mathcal{W}_{\psi}[f](g) \textrm{ for all }g \in SE(2).
\end{array}
\right.
\]
}
This only means (by the product rule for differentiation) that we must add to the righthand side in the PDE in (\ref{coherence2}) which equals
\[
\sum \limits_{i=1}^{3}\sum \limits_{j=1}^{3} \mathcal{A}_{i}\left\{G^{ij}(U)(g,t) \mathcal{A}_{j}  U(g,t)\right\}(g,t),
\]
the following terms: $+\frac{1}{\sqrt{\det{G(U)(g,t)}}} \sum \limits_{i=1}^{3} \mathcal{A}_{i}\left\{\sqrt{\det{G}(U)(g,t)}\right\}(g,t) \sum \limits_{j=1}^{3} G^{ij}(U)(g,t) \mathcal{A}_{j}  U(g,t)$.

\section{Towards a graphical eraser: Morphology PDE's on $SE(2)$. \label{ch:erasor}}

\begin{figure}
\centerline{
\includegraphics[width=0.75\hsize]{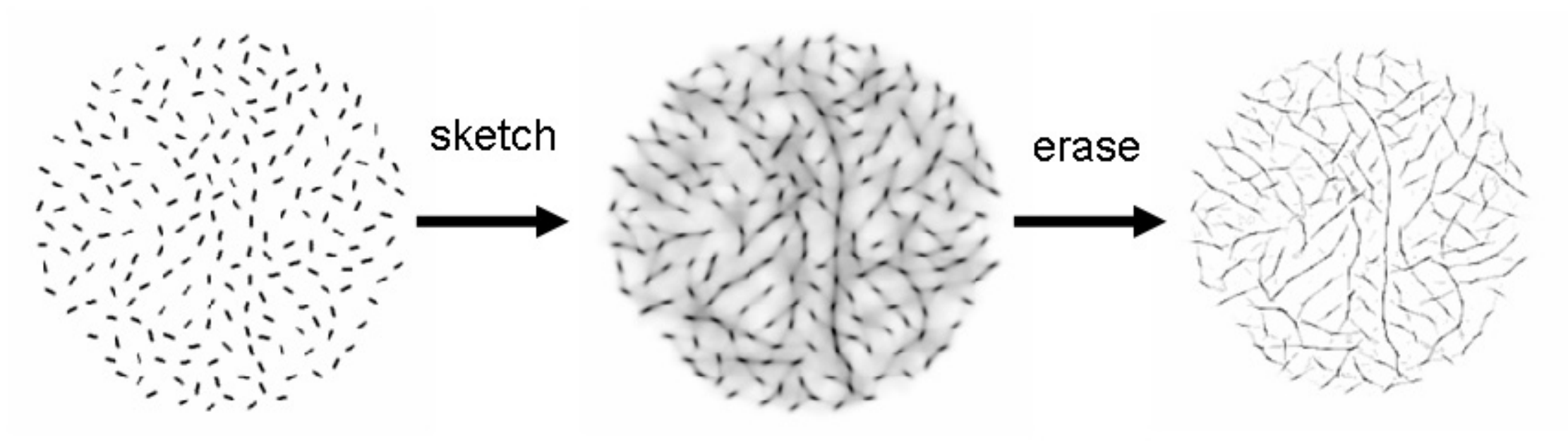}
}
\caption{The sections \ref{ch:OS}, \ref{ch:2}, \ref{ch:leftinvimageproc},  \ref{ch:diffSE2}, \ref{ch:Cartan}, \ref{ch:CED} provide the theory for automatic sketching process by means of left-invariant evolution equations on invertible orientation scores, where the general scheme is explained in Figure \ref{fig:leftinv}. The sections \ref{ch:elastica}, \ref{ch:elasticatrue}, \ref{ch:CED} and Appendix \ref{ch:app}, \ref{ch:app2}, \ref{ch:app3} serve the practical purpose of developing an automated erasor by means of left-invariant curve extraction. Section \ref{ch:erasor} serves the practical goal to create a continuous process from a sketched image to an image consisting only of these curves, where we again use the general scheme explained in Figure \ref{fig:leftinv}, but where we replace the (non-linear) left-invariant evolution equations (\ref{coherencesimple}) and (\ref{coherence2}) (in section \ref{ch:CED}) by Hamilton-Jakobi equations (also known as morphology equations in image processing) (\ref{er1alpha}).}\label{fig:erasor}
\end{figure}

Sofar we considered automatic line and contour enhancement via linear and non-linear left-invariant diffusion equations on (invertible) orientation scores of 2D-images. In this way we obtained an automated ``graphical sketcher''. What is missing is an automated ``graphical eraser'' which erases brush strokes which are to far away from the the zero crossings $\partial_{\theta}U$ and $\partial_{\eta}U$ and which make the completion fields to broad. See Figure \ref{fig:erasor}.

Moreover, our oriented wavelets are not ``perfect'' line detectors.
They often yield too broad distributions $\theta \mapsto U_{f}(\ul{x},e^{i\theta})$ at fixed positions $\ul{x}\in \R^2$ where a line/contour is present. Therefore it is often desirable
to erase (i.e. to apply erosion) in both $\theta$ and $\eta$ direction.

In the previous chapter we have seen that the modes of the direction process are for reasonable parameter settings ($\frac{D_{11}}{\alpha}$ small) closely approximated by the intersection of the smooth surfaces $\{g \in SE(2)\; |\; \partial_{\theta}U(g)=0 \}$ and $\{g \in SE(2) \; |\; \partial_{\eta}U(g)=0\}$. 
The gradient $\nabla_{2}U=(\partial_{\theta}U,\partial_{\eta}U)$ of a (diffused\footnote{Since thinning and erosion algorithms can be useful both before and after a diffusion step, we simply write $U \in \mathbb{L}_{2}(SE(2))$. Here $U$ could be $U=|\mathcal{W}_{\psi}f|$ or $U=\Phi_{t}(|\mathcal{W}_{\psi}f|)$, where $t$ is the stopping time of the non-linear left-invariant diffusion or $U$ could be the completion field $U=(A-\gamma I)^{-1}|\mathcal{W}_{\psi}f| (A^{*}-\gamma I)^{-1}|\mathcal{W}_{\psi}f|$ obtained by linear left-invariant diffusion. }) orientation score $U \in \mathbb{L}_{2}(SE(2))$
locally points to these zero-crossings.

Therefore we propose the following 2 PDE-systems on \emph{real-valued } (processed) orientation scores $U:SE(2) \to \R$ 
\begin{equation} \label{er1alpha}
\left\{
\begin{array}{l}
\partial_{t}W(g,t)= -\|(C \nabla_{2}) W(g,t)\|^{2 \zeta}= -((c_{1})^{2} (W_{\theta}(g,t))^2 + (c_{3})^2 (W_{\eta}(g,t))^{2})^{2\zeta}, \\
W(g,0)=|U(g)|,  \textrm{ with }c_{1}=D_{11}\textrm{ and }c_{3}=D_{33},
\end{array}
\right.
\end{equation}
$\zeta \in [\frac{1}{2},1]$, with $C=\textrm{diag}(c_{1},c_{3})$. For $\zeta=1$ the solution is given by
\begin{equation} \label{viscosity}
W(g,t)=(k_t \ominus_{SE(2)} U)(g) :=\inf \limits_{h \in \mathcal{S}_{g_{0}}}[ U(h)-k_{t}^{c_{1},c_{2}}(g^{-1} h)].
\end{equation}
where the morphology kernel $k_t$ is given by
\begin{equation} \label{morphkern}
k_t(g)=\frac{K_{t}^{D_{11},D_{33}}(g)}{K_{t}^{D_{11},D_{33}}(e)}
\end{equation}
where $K_{t}^{D_{11},D_{33}}(g)$ equals (after rotation by $\pi/2$ in each spatial plane, mapping $\xi$ to $\eta$) the diffusion kernel on $SE(2)$ studied in section \ref{ch:diffSE2}, Theorem \ref{th:exaxtdiff}, (with analytic approximations in section \ref{ch:heisapprox}).

Here we note that the corresponding dilation operator is defined by
\[
(K \oplus_{SE(2)} U)(g):= \sup \limits_{h \in SE(2)} \left[ K(h^{-1}g)+U(h) \right]
\]
which is the equivalent of a $SE(2)$-convolution, (\ref{Gconv}), where we replaced the $(+,\cdot)$-algebra by the $(\max,+)$ algebra.
\begin{remark} Due to the non-commutative nature of $SE(2)$ (and H\"{o}rmander's theorem, \cite{Hoermander}) both the diffusion kernel in Theorem \ref{th:exaxtdiff} and the corresponding kernel $k_t$ given by (\ref{morphkern}), in contrast to their well-known analogues on $\R^2$:
\begin{equation} \label{MorphonR2}
K_{t}(x,y)=\frac{1}{4\pi t}e^{-\frac{x^2+y^2}{4t}} \textrm{ and }k_{t}(x,y)=\log \frac{K_{t}(x,y)}{K_{t}(0,0)}=\frac{x^{2}+y^{2}}{4t},  (x,y) \in \R^2,
\end{equation}
are \emph{not separable} along the exponential curves, respectively, in the $(+,\cdot)$-algebra and $(\max,+)$ algebra. For if a differentiable kernel $K$ is separable $K(a,b)=k_{1}(a)k_{2}(b)$ (resp. $k(a,b)=k_{1}(a)+k_{2}(b)$) and isotropic $(a\partial_b -b \partial_{a})K(a,b)=0$ (resp. $(a\partial_b -b \partial_{a})K(a,b)=0$) then clearly it must be equal to $K(a,b)=\mu e^{\frac{a^{2}+b^2}{\lambda}}$ and $k(a,b)=\gamma (a^{2} +b^{2})$, for some separation constants $\mu, \lambda,\gamma >0$.
\end{remark}
Before we motivate our conjecture we first explain why the standard approach on finding viscosity solutions of morphology equations on $\R^2$ by means of the Cramer transform is not easily generalized to $SE(2)$. The homomorphism between dilation/erosion and diffusion/inverse diffusion is given by the Cramer transform $\mathcal{C}=\gothic{F} \circ \log \circ \mathcal{L}$, \cite{Akian}, \cite{Burgeth}, which is a concatenation of the multi-variate Laplace-transform, logarithm and the Fenchel transform (mapping a convex function $c:\R^n \to \overline{\R}$ onto the function $\ul{x} \mapsto [\gothic{F}c](\ul{x})=\sup_{\ul{y}}[\ul{y} \cdot \ul{x}-c(\ul{x})]$). This is due to the fact that
\[
\mathcal{C}(f*g)=\gothic{F} \log \mathcal{L}(f*g)=\gothic{F}(\log \mathcal{L}f+\log \mathcal{L}g)=\gothic{F} \log \mathcal{L} f \oplus \gothic{F} \log \mathcal{L} g
=\mathcal{C}f \oplus \mathcal{C}g.
\]
Now, since $SE(2)$ is non-commutative the irreducible representations
are no-longer 1 dimensional and the Fourier/Laplace transform on $SE(2)$ becomes relatively complicated (see \cite{DuitsR2006AMS}){App. B}). Moreover, it is not clear how the Fenchel transform should be generalized to $SE(2)$.

\subsection{Morphology and Hamilton-Jakobi theory \label{ch:MorphJakobi}}

In the calculus of variations and the corresponding Hamilton-Jakobi theory on some finite dimensional manifold $G$, with local coordinates $\{g_{i}\}_{i=1}^{n}$, one usually starts with a Lagrangian $L:\R^{+} \times G \times T(G) \to \R^{+}$ of class $C^2$,
giving rise to the following energy on curves
\begin{equation}\label{energy}
\mathcal{E}(\gamma)= \int \limits_{t_0}^{t_{1}} L(t,\gamma(t),\gamma'(t)) \; {\rm dt},
\end{equation}
where, for now, we shall assume non-degeneracy of $L$ with respect to its dependence on $T(G)$, i.e.
\begin{equation} \label{nondeg}
\det \left(\nabla^{2}_{\dot{g}} L(t,g,\dot{g}) \right)\neq 0 \textrm{ for all }t>0, g \in G .
\end{equation}
so that we can express $\dot{g}^{i}$ as functions of the canonical variables $(t,g^{j}, p_{i})$:
\begin{equation} \label{phigp}
\dot{g}^{i}=\phi(t,g^{j}, p_{i}),
\end{equation}
where the components momentum $p=\nabla_{\dot{g}}(L(t,g,\dot{g}))$ are given by $p_{i}= \frac{\partial L(t,g^j,\dot{g}^{j})}{\partial \dot{g}^{i}}$.

Here one should make a clear distinction between the parameter dependent non-homogeneous case (recall the elastics in section \ref{ch:elastica}, where traveling time coincides with spatial arc-length, $t=s$ and where $G=SE(2)$, $L(s,[\gamma],[\dot{\gamma}])=((\dot{\theta}(s))^{2}+ \epsilon )\|\dot{x}(s)\|=(\dot{\theta}(s))^{2}+ \epsilon $) and the parameter independent homogeneous case, recall the geodesics in section \ref{ch:geodesics} with
\begin{equation}
G=SE(2)/Y \textrm{ and }  L(s,[\gamma],[\dot{\gamma}])=\sqrt{(\dot{\theta}(s))^{2}+ \epsilon \|\dot{x}(s)\|^{2} }.
\end{equation}
Moreover, one should make a clear distinction between the cases where $t>0$ is an intrinsic coordinate on the manifold and the cases where time $t>0$ is an independent coordinate, see \cite{Rund}{p.44-p48}.

For now we will consider the non-homogeneous case with time as an independent coordinate. Soon we will consider $t>0$ as an intrinsic coordinate on the manifold, which will turn out to be important to relate the direction process to morphological equations on $SE(2)$, whereas the relation to contour enhancement processes and morphology on $SE(2)$ requires $t>0$ to play an independent role.
\begin{theorem}
A 1-parameter family of hypersurfaces $S(t,g)=\sigma$ is geodesically equidistant with respect to a \emph{non-degenerate } Lagrangian, that is $\det \left(\nabla^{2}_{\dot{g}} L(t,g,\dot{g}) \right)\neq 0$ and $L(t,\gamma(t),\gamma'(t))=\frac{d\sigma}{dt}(t)$ for all congruency curves $\gamma$ (satisfying $\nabla_{\dot{g}} L(t,\gamma(t),\gamma'(t))=\nabla_{g}S (t,\gamma(t))$) if and only if
\begin{equation} \label{Jakobieq}
\frac{\partial S}{\partial t} + H(t, g, \nabla S)=0,
\end{equation}
where the uniquely corresponding non-vanishing Hamiltonian function $H \in C^{2}(\R^{+}\times G \times T(G)', \R^{+})$ is given by
\begin{equation} \label{Hamiltonian}
H(t,g,p)=-L(t,g,\phi(t,g,p))+ \langle p,\phi(t,g,p) \rangle.
\end{equation}
In particular the characteristic function $W: \R^{+} \times G $ given by
\[
W(t,g ; t_0,g_0)= \inf \{\int_{t_{0}}^{t} L(t,\gamma(t),\dot{\gamma}(t)) \, {\rm d}t \; |\; \gamma \in C^{2}(\R^{+},G) \textrm{ with }\gamma(t)=g, \gamma(t_0)=g_0 \} , g_0 \in G, t>0,
\]
satisfies the Hamilton-Jakobi equation. Along congruency curves the following fundamental equations hold
\begin{equation} \label{Hamfundamental}
\dot{p}=-\nabla_{g}H(g,p) \textrm{ and }\dot{g}=\nabla_{p} H(g,p).
\end{equation}
\end{theorem}
For proof and more background see \cite{Rund}{p.12--25}.
\begin{theorem}
The Hyper-surface
\begin{equation}\label{omhullende}
\{g \in G \; |\; S(t,g)= \sigma + R\}, \qquad R>0,
\end{equation}
where $S(t,g)$ is a solution of the Hamilton-Jakobi equation, is the envelope of the set of geodesic spheres
\begin{equation}\label{geodesicsphere}
S_{(g_0,t_{0}),R}  \{(t,g) \in \R^{+} \times G \; |\; W(t,g ; t_0,g_0)=R\}
\end{equation}
of radius $R>0$ centered on the hyper surface $\{g \in G \; |\; S(t,g)=\sigma\}$. If moreover,
the Weierstrass excess function $E: \R^{+} \times G \times T(G) \times T(G) \to \R$ given by
\[
E(t,g^{j},\dot{g}^{j},\dot{\overline{g}}^{j})=L(t,g^j,\dot{\overline{g}}^{j})-L(t,g^j,\dot{g}^{j}) -(\dot{\overline{g}}^{j}-\dot{g}^{j}) \frac{\partial L(t,g^{j},\dot{g}^{j})}{\partial \dot{g}^{j}}>0
\]
is strictly positive, where $\dot{\overline{g}}^{j}$ denotes the tangent vector an arbitrary extremal curve $\overline{\Gamma}$ issuing from $g_0$ to an arbitrary point $\overline{g} \in G$ such that
\[
\int_{g_0, \Gamma}^{\overline{g}} L(t,g^j, \dot{\overline{g}}^{j}) {\rm d}t = R, \qquad R>0,
\]
we have that the hyper surfaces (\ref{omhullende}) are supporting hyper surfaces of the geodesic spheres centered on the hyper surface $\{g \in G \; |\; S(t,g)=\sigma\}$.
\end{theorem}
This results puts an analogy between morphological PDE's, which are of the type (\ref{Hamiltonian}), and Huygens' principle. Now consider $G=SE(2)$ and put the following Cartan connection on $SE(2)$ given by $\omega(X_g)={\rm d}\xi(X_g) \partial_{x}$, $X_g \in T_{g}(SE(2))$, this means that the horizontal part of the tangent space is given by the span $\{\partial_{\theta},\partial_{\eta}\}$.
Set
$\mathcal{L}(s,[\gamma],[\gamma'])=\frac{1}{4 D_{11}}(\frac{d \theta}{ds})^{2}+ \frac{1}{4 D_{33}} (\frac{d \eta}{ds})^{2} $ and we get
\[
H(p_{1},p_{3})=-L(2D_{11}p_{1},2D_{33}p_2)+ 2 D_{11} p_{1}^{2} + 2D_{33} p_{2}^{2}=D_{11} p_{1}^{2} + D_{33} p_{2}^{2}
\]
and thereby the morphological PDE-system (\ref{er1}) coincides with the Hamilton-Jakobi equation (\ref{Hamiltonian}), with $\zeta=1$, and by the above theorems (using the moving frame $\{\xi,\eta,\theta\}$ of reference to ensure left-invariance) we see that the viscosity solution is indeed of the type (\ref{viscosity}).  Similarly, if we consider
$G=SE(2)$, with Cartan connection $\omega(X_g)={\rm d}\eta(X_g) \partial_{y}$, $X_g \in T_{g}(SE(2))$, this means that the horizontal part of the tangent space is given by the span $\{\partial_{\theta},\partial_{\xi}\}$.
\[
H(p_{1},p_{3})=-L(2D_{11}p_{1},2D_{22}p_2)+ 2 D_{11} p_{1}^{2} + 2D_{22} p_{2}^{2}=D_{11} p_{1}^{2} + D_{22} p_{2}^{2}
\]
and Hamilton-Jakobi system
\begin{equation} \label{er1}
\left\{
\begin{array}{l}
\partial_{t}W(g,t)= -D_{11} (W_{\theta}(g,t))^2 - D_{22} (W_{\xi}(g,t))^{2},\\
W(g,0)=U(g),  \ , g \in G, t>0
\end{array}
\right.
\end{equation}
It is still the question though whether the structure element (i.e. the max-plus convolution kernel)
$k_t(g)=\frac{1}{t}k(g)$, with $k(g)=W(1,g,0,g_0)$, indeed satisfies $k_t(g)=-\log \frac{K^{D_{11},D_{22}}_{t}(g)}{K^{D_{11},D_{22}}_{t}(e)}$ ? \\
 \\
Considering the usual $\R^2$-case, i.e. morphological systems on images, this relation does hold (\ref{MorphonR2}) and it is the question whether the commutative nature of $\R^2$ plays a
crucial role. The next results indicate that the non-commutative nature of $SE(2)$ does not cause problems in this respect. In fact they subscribe our conjecture.

\subsubsection{Hamilton-Jakobi theory and (lifted) elastica curves}

Recall from section \ref{ch:elastica} that the elastica curves correspond to the modes of the direction process. As in a direction process a unit speed grey-value particle is always suppose to move in $\xi$-direction it makes sense to identify the arc-length $s>0$ with $\xi$. By doing this we use the temporal parameter $t\equiv s \equiv \xi$ as a spatial parameter in $SE(2)$. As pointed out in \cite{Rund}{p.44-48} this requires a different approach in Hamilton Jakobi-theory.

We rewrite the Lagrangian
\begin{equation} \label{Lagrangianxiands}
L(s,\theta(s),\xi(s),\eta(s),\dot{\theta}(s),\dot{\xi}(s),\dot{\eta}(s))= \mp (\frac{1}{4 D_{11}}(\dot{\theta}(s))^{2}+\alpha \cdot 1)  \; \equiv ((\dot{\theta}(s))^{2}+\epsilon) \langle {\rm d}\xi,\dot{x}(s)\rangle,
\end{equation}
with $\langle {\rm d}\xi,\dot{x}(s)\rangle=\|\dot{x}(s)\|=1$, in (\ref{energy}) as follows
\[
\begin{array}{l}
L^{*}([g(s)],[g'(s)])=L^{*}(\theta(s),s,\eta(s); \dot{\theta}(s),\dot{\xi}(s),\dot{\eta}(s) ) \cdot \dot{\xi}(s):= \\
L(s,\theta(s),\eta(s), \frac{\dot{\theta}(s)}{\dot{\xi}(s)}, \frac{\dot{\eta}(s)}{\dot{\xi}(s)} ) \cdot \dot{\xi}(s),
\end{array}
\]
and therefore the components of the canonical impuls $p^{*}$ of the new Lagrangian $L^{*}$ equal
\[
\left\{
\begin{array}{l}
p^{*}_{1}=\frac{\partial L^{*}}{\partial \dot{\theta}}=\frac{\partial L}{\partial \dot{\theta}}= \mp \frac{\dot{\theta}}{2D_{11}}=p_{1},   \\
p^{*}_{2}= L-\dot{\theta}(s) p_{1} \\
p^{*}_3=p_{3}=0
\end{array}
\right.
\]
and the Hamilton-Jacobi equation on $G=SE(2)$, where we again restrict ourselves to horizontal curves and horizontal subspaces spanned by $\{\partial_{\theta},\partial_{\xi}\}$, corresponding to the Lagrangian (\ref{Lagrangianxiands}) is:
\begin{equation} \label{JacobiEl}
H^{*}([g],p)=p_{2}+H([g],p)=p_{2} \mp (D_{11} p_{1}^{2} -\alpha)= 0
\end{equation}
So the corresponding Jakobi-Hamilton equation is given by
\begin{equation}\label{Jakobieqelastica}
\frac{\partial W}{\partial \xi}= \pm \left(D_{11} \left( \frac{\partial W}{\partial \theta} \right)^2 -\alpha \right)
\end{equation}
If we now drop our identification between $\xi$ and $s$ this yields the following non-linear morphology process on orientation scores:
\begin{equation}\label{Jakobieqelastica}
\left\{
\begin{array}{l}
\frac{\partial W}{\partial s}=\frac{\partial W}{\partial \xi} \pm  D_{11} \left( \frac{\partial W}{\partial \theta} \right)^2  \\
W(\cdot,s=0)=U_{f}
\end{array}
\right.
\end{equation}
which is the morphological equivalent of the forward Kolmogorov equation (\ref{ForwardKolmogorovMumford}) of Mumford's direction process (\ref{Mumford}).

For the corresponding Heisenberg approximation of Mumford's direction process, where the curve-length $s$ parameterized is replaced by $x$ rather than $\xi$ and where the Lagrangian is simply given by $L(x,y(x),\theta(x))=(\theta'(x))^{2} +\epsilon$, $\theta(x)=y'(x)$, one can follow the same scheme.
Following \cite{Rund}{p.44-48} we for the moment introduce an independent time variable $\tau$, with $x'(\tau)\neq 0$. Later on we set $\tau=x$
\[
\begin{array}{l}
L(x,y(x),\theta(x))= (\dot{\theta}(x))^{2} +\epsilon, \qquad \theta(x)=\dot{y}(x),  \\[7pt]
L^{*}(x(\tau),y(\tau),\theta(\tau))=x'(\tau) \, L(x(\tau),y(\tau),\theta(\tau), \frac{y'(\tau)}{x'(\tau)},\frac{\theta'(\tau)}{x'(\tau)}),
\end{array}
\]
\mbox{canonical impuls vectors $(i=1: \theta, i=2: x, i=3: y)$, recall $\epsilon=4 \alpha D_{11}$}
\[
\left\{
\begin{array}{l}
p_{1}^{*}= \frac{\partial L^{*}}{\partial \theta'}=\frac{\partial L}{\partial \dot{\theta}}= 2 \dot{\theta}(x) \\
p_{2}^{*}=\frac{\partial L^{*}}{\partial x'} =L -p_{1}\dot{y}(x) -p_2 \dot{\theta}(x) \\
p_{3}^{*}=p_3=0
\end{array}
\right.
\]
and consequently the Hamiltonian is given by
\[
H(x,y,\theta, p_{1},p_{2})=-L+\dot{y}(x) p_3 + \frac{p_{1}^{2}}{2 }
=
\frac{p_{1}^{2}}{4 } +\theta p_{3} -\epsilon
\]
\mbox{yielding Hamilton-Jakobi equation: $p_{2}^{*}+H(x,p_{1},p_{3})=0$}:
\begin{equation} \label{JakobiH}
\partial_{x} S+ \theta \partial_{y}S= + \frac{1}{4 } (\partial_{\theta}S)^2 +\epsilon \Leftrightarrow \hat{A}_{2}S= \frac{1}{4}(\hat{A}_{1}S)^2 +\epsilon
\end{equation}
which is again related to the exact case (\ref{JacobiEl}) by replacing $\cos \theta$ by $1$ and $\sin \theta$ by $\theta$, i.e. replacing $\mathcal{A}_{i}$ by $\hat{A}_{i}$, recall (\ref{approxVF}). The characteristic function (computed by the B-spline solutions, recall (\ref{modes}) now equals :
\[
\begin{array}{l}
\min \left\{ \int \limits_{0}^{x} ((y''(\tau))^2+ \epsilon){\rm d\tau}\; |\; y(0)=0,y'(0)=0, y'(x)=\theta, y(x)=y \right\} \\
= 4 \frac{3y^2 +3xy\theta +x^2 \theta^2}{x^3} +\epsilon x = D_{11} \log  \frac{K_{s}^{\alpha, D_{11}}(x,y,\theta)}{K_{s}^{\alpha,D_{11}}(x,0,0)}
\end{array}
\]
The canonical equations (\ref{Hamfundamental}) for the congruency curves through the geodesically equidistant surfaces $\{(x,y,e^{i\theta}) \in SE(2) \; |\; S(x,y,\theta)=\sigma\}$ in the Heisenberg-approximation are now given by
\[
\left\{
\begin{array}{l}
\frac{dp_i}{dt} = - \frac{\partial H^{*}(x,p)}{\partial x_i},  \\
\frac{dx_i}{dt} =  \frac{\partial H^{*}(x,p)}{\partial p^{*}_i} \\
\end{array}
\right. \Rightarrow
\left\{
\begin{array}{l}
\theta'(x)= \mp 2 D_{11} p_{1}(x) \\
x'(x)=1 \\
y'(x)= \theta(x) \\
\end{array}
\right.
\textrm{ and }
\left\{
\begin{array}{l}
p_{1}'(x)=p_3(x) \\
p_{2}'(x)= 0 \\
p_{3}'(x)= 0
\end{array}
\right.
\]
with $x_{1}=\theta, x_{2}=x, x_3=y$,
from which we directly deduce
$\rightarrow y''''(x)=0, \theta(x)=y'(x)$ yielding indeed the B-spline solutions (\ref{modes}). The impuls vector $p$ along a $B$-spline mode starting at $(x_0,y_0,\theta_0)$ ending at
$(x_{1},y_{1},-\theta_{1})$, $x_{0}<x_{1}$ is given by
\[
p(x)=\left(
\begin{array}{l}
p_{1}(x)= \frac{1}{2D_{11}}\left(\frac{6 y_1 +2 x_{1}(\theta_{1}-2 \theta_0)}{x_1^{2}}-\frac{6 x(2y_{1}+x_{1}(\theta_{1}-\theta_{0}))}{x_1^3} \right) \\
p_{2}(x)= 1 \\
p_{3}(x)= - \frac{6 (2y_{1}+x_{1}(\theta_{1}-\theta_{0}))}{2D_{11}x_1^3}
\end{array}
\right)=
\]
In the Heisenberg approximation case we can compute the characteristic function by means of the $B$-spline modes (\ref{modes}) yielding
\[
\begin{array}{ll}
W(x,y,\theta; 0,0,0) &= \min \left\{
\int \limits_{0}^{x} \frac{1}{4 D_{11}}( y''(x))^2 {\rm d}x \; |\;
y(0)=y, y(0)=0, y'(0)=0, y'(0)=\theta
\right\} \\
 &=  \frac{3y^{2}+3x y \theta +x^2 \theta^2}{D_{11} x^3}
 \end{array}
\]
which can easily be checked to satisfy the Hamiltonian equation (\ref{JakobiH}).

Finally, we note that the structure element of the morphology related to the direction process is
again related to the direction process \emph{resolvent} Green's function by means of
\begin{equation}
W(x,y,\theta; 0,0,0)= -4 \, \cdot  \log \left(\frac{R^{\alpha,D_{11}}(x,y,\theta)}{R^{\alpha, D_{11}}(x,0,0)}\right).
\end{equation}
\begin{remark} Consider the corresponding PDE-systems for a direction process and morphology on $\R^2$, then we again have
\[
\begin{array}{l}
(\partial_{x}-D \partial_{y}^{2} -\alpha I)R_{\alpha, D}(x,y)=\alpha \delta_{0,0} \rightarrow R_{\alpha, D}(x,y)=\frac{\alpha}{\sqrt{4\pi s}} e^{-\frac{y^2}{4s}-\alpha x}
k_x(y)= -\log \frac{R_{\alpha, D }(x,y)}{R_{\alpha, D_{22}}(x,0)}= \frac{y^2}{4x},
\end{array}
\]
and $k_x(y)=\frac{y^2}{4x}$ is the viscosity solution of the Hamilton-Jakobi equation $\partial_{x} k(x,y) = \partial_{y}^{2} k(x,y)$, with $k(0,y)=\delta_{0}^{c}$, with $\delta_{0}^{c}$ the denotes the convex Dirac function given by $\delta_{0}^{c}(y)=\infty$ if $y\neq 0$ and $\delta_{0}^{c}(y)=0$ if $y=x$.
\end{remark}
Before we proceed with the homogeneous case we put a relevant observation on the exact case.
\begin{remark}
Despite the fact that $\kappa(s)=\dot{\theta}(s)$ for horizontal curves, the exact elastica equation (\ref{elastics}) is \emph{not} of the standard Lagrangian type
\mbox{$
\partial_{g_i} \mathcal{L} + \partial_{t} \{\partial_{g'_i} \mathcal{L} - \partial_{t} \partial_{g_{i}''}\mathcal{L}\}=0.
$}. The reason for this is that the arclength parameter $s>0$ is \emph{curve dependent}. Recall that in the only relevant horizontal pertubations are:
\[
\begin{array}{l}
\ul{x}_{NEW}(s)=\ul{x}(s)+ h \delta \ul{n}(s) \\
\theta_{NEW}(s)=\theta(s) + \arctan \frac{ \epsilon \delta(s) \dot{\theta}(s)}{1+\epsilon \dot{\delta}(s)}
\end{array}
\]
we have $\frac{ds_{NEW}}{ds}(s)=1- h \delta(s) \dot{\theta}(s) +O(h^2)$. Therefore in order to embed elastica curves in the standard Euler-Lagrange/Hamiltonian theory
we must follow \cite{MarkusThesis}{App. C} and rewrite the Lagrangian in a curve dependent way
\[
\mathcal{L}(\eta(t),\xi(t),\eta(t), \xi'(t),\eta'(t),\theta'(t))= \frac{(\theta'(t))^2}{s'(t)} + \epsilon s'(t),
\]
where $s'(t)=\|\ul{x}'(t)\|=\sqrt{(\xi'(t))^2+(\eta'(t))^2}$ and in our notation we distinguish between $\theta'(t)= \frac{d}{dt}\theta(s(t))$ and $\dot{\theta}(s)=\frac{d}{ds} \theta(s)$.
%
The curvature along the curve $t\mapsto g(t)=(\xi(t)\cos(\theta(t))- \eta(t)\sin(\theta(t)), \xi(t)\sin(\theta(t))+ \eta(t)\cos(\theta(t)),\theta(t))$ at time $t>0$ equals
\[
\kappa(s(t))= \frac{\xi'(t)\eta''(t)-\xi''(t)\eta'(t)}{(s'(t))^3}
\]
and consequently
\[
\begin{array}{ll}
\partial_{\xi'}s'= \frac{\xi'}{s'} \textrm{ and } \partial_{\eta'}s'= \frac{\eta'}{s'}, &
\partial_{\xi''}s'= \partial_{\eta''}s'=0 \\
\partial_{\xi'}\kappa = -2 \frac{\kappa}{(s')^2} \xi' + \frac{s''}{(s')^4} \eta', \partial_{\eta'}\kappa = -2 \frac{\kappa}{(s')^2} \eta' - \frac{s''}{(s')^4} \xi', &
\partial_{\xi''}\kappa= -\frac{\eta'}{(s')^3}, \  \partial_{\eta''}\kappa= -\frac{\xi'}{(s')^3}.
\end{array}
\]
and thereby as pointed out in \cite{MarkusThesis}{App. C} the standard Euler Lagrange equations yield
\[
\left\{
\begin{array}{l}
\partial_{t} (\partial_{\xi'}\mathcal{L} - \partial_{t} \partial_{\xi'}\mathcal{L})=0  \\
\partial_{t} (\partial_{\eta'}\mathcal{L} - \partial_{t} \partial_{\eta''}\mathcal{L})=0
\end{array}
\right. \desda \left(\epsilon \kappa +\kappa^{3} - 2 \cdot \frac{k' s''-k'' s'}{s'}\right)
\left(
\begin{array}{l}
-\eta'(t) \\
\xi'(t)
\end{array}
\right) \desda 2\ddot{\kappa}(s)+\kappa^{3}(s) +\epsilon \kappa(s)=0.
\]
Now for horizontal curves $\angle \dot{x}(s(t))=\angle \ul{x}'(t)=\theta(s(t))$ one can also express the curvature along the curve at time $t>0$
\[
\kappa(s(t))= \dot{\theta}(s(t))=\frac{\theta'(t)}{s'(t)}
\]
yielding the equivalent angular Euler-Lagrange equation
\[
\epsilon \partial_{\dot{\theta}} L -\dot{\theta}L - \frac{d^2}{ds^2}\{ \partial_{\dot{\theta}}L\}=2\ddot{\kappa}(s)+\kappa^{3}(s) +\epsilon \kappa(s)=0, \textrm{ with }\kappa(s)=\dot{\theta}(s).
\]
\end{remark}

\subsubsection{The homogenous case: Hamilton-Jakobi theory and geodesics in $SE(2)$}

Sofar we restricted ourselves to variational problems based on a non-degenerate Lagrangian, i.e. a Lagrangian satisfying (\ref{nondeg}) on $SE(2)$. For the elastica curves studied in section \ref{ch:elastica} this was fine, however for the geodesics studied in section \ref{ch:geodesics} it is not. With this respect we note that if a Lagrangian $L(g,\dot{g})$ is homogeneous, that is if $L(g,\lambda \dot{g})=\lambda L(g,\dot{g})$, then $\nabla_{\dot{g}} L(g,\dot{g}) \cdot \dot{g}=L(g,\dot{g})$ and thereby $\det ( \nabla_{\dot{g}} \nabla_{\dot{g}} L)=0$. So we can not express $\dot{g}$ in the canonical variables $g,p$ like we did in (\ref{phigp})!

Therefore the non-degeneracy condition (\ref{nondeg}) is replaced by another non-degeneracy condition: 
\[
\det (g_{ij}(g,\dot{g}))\neq 0,
\]
where the fundamental tensor $g_{ij} {\rm d}g^{i} \otimes {\rm d}g^{j}$ is given by
\[
g_{ij}(g, \dot{g}) =\frac{1}{2} \frac{\partial^{2} L(g,\dot{g})}{\partial \dot{g}^{i} \partial \dot{g}^{j}}.
\]
and the canonical variables $p=p_{i}{\rm d g}^{i}$ in the non-homogeneous case are now replaced by the following variables
\[
y= y_{i} {\rm d y}^{i}, \textrm{ where }y_{i}=g_{ij}(g,\dot{g}) \, \dot{g}^{j} = L(g, \dot{g}) \frac{\partial L(g, \dot{g})}{ \partial \dot{g}^{i}}.
\]
allowing us to rewrite the energy in
\[
\mathcal{E}(\gamma)= \int \limits_{s_0}^{s_{1}} \sqrt{g_{ij}(\gamma(s),\dot{\gamma}(s))\; \dot{\gamma}^i(s) \dot{\gamma}^{j}(s) } \; {\rm ds},
\]
where $g_{ij}(g,\dot{g})g_{jk}(g,\dot{g})=\delta_{ik}$, which is the length integral in the manifold $G$ in Riemannian geometry. The Hamiltonian now reads
\begin{equation}
H^{2}(g^{k},y_{h})=g^{ij}(g^k,y_h) y_i y_j
\end{equation}
and consequently we get
\[
\begin{array}{l}
g^{ij}(g^k,y_h)= \frac{1}{2} \frac{\partial^{2} H^{2}(g^k,y_h)}{\partial y_i \partial y_j}, \\[8pt]
\dot{g}^{i}= H(g^{i},g_{ij}(g,\dot{g})\, \dot{g}^{j}) \frac{\partial H(g^{i},g_{ij}(g,\dot{g}) \, \dot{g}^{j})}{\partial y_i}= \left. H(g^{i},y_i) \frac{\partial H(g^{i},y_i)}{\partial y_i} \right|_{y_i=g_{ij}(g,\dot{g}) \, \dot{g}^{j}}
\end{array}
\]
from which we deduce the following relation between Hamiltonian and Lagrangian:
\[
H(g^{h},y_{h})=L(g^h, \dot{g}^{h}), \qquad \textrm{ with } y_h=g_{hj}(g,\dot{g}) \, \dot{g}^{j}.
\]
As a result, for details see \cite{Rund}{p.166-170} the corresponding Hamilton Jabobi equations now become
\[
H(g,\nabla S)= \pm 1.
\]
Now we return to our case of interest, where we equip the principal fiber bundle $P_X=(SE(2),SE(2)/X,\pi,R)$, $X=\{(x,0,0)\; |\; x \in \R\}$ and where we recall $\pi(g)=[g]$, $R_g h=hg$ set connection $\omega(X_g)=\langle{\rm d}\xi,X_g \rangle \partial_{x}$ so that with this convention the horizontal part of the tangent spaces is given by
\[
\textrm{span}\{\partial_{\theta},\partial_{\eta}\} \subset T(SE(2))
\]
on which we apply the homogeneous Lagrangian:
\begin{equation}
\label{Lagrangianxiands}
L(s,\theta(s),\eta(s),\dot{\theta}(s),\dot{\eta}(s))= \sqrt{ \frac{1}{4D_{11}}(\dot{\theta}(s))^{2}+ \frac{1}{4D_{22}} (\dot{\eta}(s))^{2}},
\end{equation}
whose corresponding geodesics were studied in section \ref{ch:geodesics}. The corresponding Hamiltonian is given by $H([g],y)=\sqrt{D_{11}(y_{1})^{2}+D_{22}(y_{2})^{2}}$ and the Hamilton-Jakobi equation now reads
\[
1= \pm \sqrt{ \left( \frac{\partial S}{\partial \theta} \right)^2 +\left( \frac{\partial S}{\partial \eta} \right)^2},
\]
now by setting $W(g,t)=S(g)P(T=t)=\alpha e^{-\alpha t} S(g)$ we get the following morphological system
\[
\left\{
\begin{array}{l}
\partial_{t}W= \sqrt{D_{11}\left(\frac{\partial W}{\partial \theta}\right)^2+D_{22}\left(\frac{\partial W}{\partial \eta}\right)^2 } \\
W(g,0)=U(g)
\end{array}
\right.
\]
which is exactly (\ref{er1alpha}) for $\zeta=\frac{1}{2}$.

\subsection{Graphical thinning}

Another approach for narrowing down the orientation scores (and/or the corresponding completion fields) around the zero crossings of $\partial_{\theta}U$ and $\partial_{\eta}U$ of some real-valued function $U:SE(2) \to \R$ is what we call ``graphical thinning''. Here oriented grey-value particles are transported to the modes, rather than being erased. We propose the following two models for graphical thinning:
{\small
\begin{equation} \label{th1}
\left\{
\begin{array}{l}
\partial_{t}W(g,t)= -C\nabla_{2} \cdot (W(g,t) \, C \nabla_{2}W(g,t)) = -c^{2} (W(g,t) \, W_{\theta \theta}(g,t) + W_{\theta}^{2}(g,t))-(W(g,t) W_{\eta \eta}(g,t)+ W_{\eta}^{2}(g,t)) , \\  \ \  
W(g,0)=U(g)
\end{array}
\right.
\end{equation}
}
and its linear counterpart
{\small
\begin{equation} \label{th2}
\left\{
\begin{array}{l}
\partial_{t}W(g,t)=  -C\nabla_{2} \cdot (W(g,t) \, C \nabla_{2}U|)(g) \\= -c^{2} (W(g,t) \, U_{\theta \theta}(g) + W_{\theta}(g,t) U_{\theta}(g)-(W(g,t) \, U_{\eta \eta}(g) + W_{\eta}(g,t)U_{\eta}(g)) ,  \\
W(g,0)=U(g).
\end{array}
\right.
\end{equation}
}
These PDE's have the advantage that the total $\mathbb{L}_{1}$-norm is preserved. However, these equations are unstable near the zero crossings of $\partial_{\theta}W$ and $\partial_{\eta}W$, which causes serious problems in practice.

\section{Acknowledgements}

The Dutch organization for Scientific research is gratefully acknowledged for financial support.

The authors wish to thank the following persons from the biomedical engineering department Eindhoven University of Technology  dr. Markus van Almsick (Chapter \ref{ch:OS}, Chapter \ref{ch:2}, Chapter \ref{ch:erasor}, Chapter \ref{ch:elastica}), prof. Luc Florack (Chapter \ref{ch:OS}) and ir. Gijs Huisman (Chapter \ref{ch:elastica} and Appendix \ref{ch:app2}) for their contributions and discussions on this manuscript.

Furthermore the authors wish to thank the following persons from the Department of mathematics and computing science at Eindhoven University of technology: prof. Mark Peletier (Chapter \ref{ch:elastica}, Chapter \ref{ch:elasticatrue} and Appendix \ref{ch:app3}), drs. Yves van Gennip (Chapter \ref{ch:elastica}, Chapter \ref{ch:elasticatrue} and Appendix \ref{ch:app3}) and dr. Olaf Wittich (Chapter \ref{ch:diffSE2} and Appendix \ref{ch:app2}), dr. Tycho van Noorden (Chapter \ref{ch:geodesics})
prof. Jan de Graaf (Chapter \ref{ch:OS} and Chapter \ref{ch:diffSE2}) and dr. Tom ter Elst (currently in University of Auckland, New Zealand, Chapter \ref{ch:2}, subsection \ref{ch:hoermander}, Appendix \ref{ch:terElst}) and ir. Maurice Duits (currently in University of Leuven, Belgium, Chapter \ref{ch:OS}) for several important corrections and suggestions included in this manuscript.

\appendix

\section{Derivation of the geodesics by means of reduction of Pfaffian systems using Noether's Theorem \label{ch:app}}

Next we apply Bryant \and Griffiths approach \cite{Bryant} on the Marsden-Weinstein reduction for Hamiltonian systems \cite{Marsden} admitting a Lie group of symmetries on Euler-Lagrange equations associated to the functional $\int \sqrt{\kappa^{2}(s)+ \epsilon} {\rm d}s$, to explicitly derive the solution curves $s \mapsto \gamma(s)$ in $SE(2)$. Recall that in section \ref{ch:geodesics} we derived the curvature of the minimizer of $\int \sqrt{\kappa^{2}(s)+ \epsilon} {\rm d}s$ by solving an ODE for $\kappa$ that we derived from Euler-Lagrange minimization. Here we will derive the same equation, taking into account the restriction to horizontal curves, in a much more structured way, avoiding extensive computations, by means of symplectic geometry. Moreover we will derive an important underlying conservation law and by the Marsden-Weinstein reduction we will derive the curves themselves. Similar to section \ref{ch:elasticatrue} we will not restrict ourselves to curves of fixed length.

Consider the manifold  $Q= SE(2) \times \R^{+} \times \R \times \R$ with coordinates $(x,y,e^{i\theta},\sigma, \kappa,t)$, where $\sigma= \|\ul{x}'(t)\|$ so that ${\rm d}s= \sigma {\rm d}t$.
On $Q$ we consider the Pfaffian equations
\begin{equation} \label{Pfaffhor}
\begin{array}{l}
\theta^{1}:={\rm d}\xi - \sigma {\rm d}t = 0 , \qquad \sigma >0, \xi=x\cos \theta +y \sin \theta, \eta=-x\sin \theta + y \cos \theta, \\
\theta^{2}:={\rm d}\eta =0, \\
\theta^{3}:={\rm d}\theta -\kappa \sigma {\rm d}t =0,
\end{array}
\end{equation}
note that these Pfaffian equations uniquely determine the horizontal part $I(Q)$ of the dual tangent space $T^{*}(Q)$, where we recall that along horizontal curves we have
$\frac{{\rm d}\theta}{{\rm ds}}=\sigma^{-1}\frac{{\rm d}\theta}{{\rm dt}}= \kappa$, $\langle {\rm d}\eta,\ul{x}'(t)\rangle=0$, $\langle {\rm d}\xi,\ul{x}'(t)\rangle=\sigma$.

We would like to minimize the energy $\int \sqrt{\kappa^{2}+\epsilon} \sigma {\rm d}t$ under the side conditions (\ref{Pfaffhor}), then the gradient of the energy should be linearly dependent on the gradient of the side condition and therefor we set
\[
\psi=\sqrt{\kappa^{2}+\epsilon} \sigma {\rm d}t
 + \lambda_{1}({\rm d\theta -\kappa \sigma {\rm d}t})+\lambda_{2}({\rm d}\xi -\sigma {\rm d}t) +\lambda_{3}{\rm d}\eta
\]
where $\lambda_{1}$, $\lambda_{2}$, $\lambda_{3}$ are Lagrange multipliers.
Formally speaking, we consider the affine sub-bundle $Z=\{\;Z_q |\; q \in Q\}\equiv Q \times T(SE(2))^{*}$ of $T^{*}(Q)$ determined by
\[
\begin{array}{l}
Z_q = \{\left. \sqrt{\kappa^{2}+\epsilon} \sigma {\rm d}t \right|_{q} \in I_{q} \subset T^{*}_{q}(Q) \}, \\
Z\equiv Q \times T(SE(2))^{*} \ \textrm{ by the isomorphism } \ (q,\LLL) \leftrightarrow \left. \sqrt{\kappa^{2}+\epsilon} \sigma {\rm d}t \right|_{q} +  \sum \limits_{k=1}^3\left. \lambda_{k}\theta^{k} \right|_{q} 
\end{array}
\]
Next we compute the exterior derivative of $\psi$ :
\[
\begin{array}{ll}
{\rm d}\psi &= \sqrt{\kappa^{2}+\epsilon} {\rm d}\theta \wedge {\rm d}t + \frac{\kappa \sigma}{ \sqrt{\kappa^{2}+\epsilon}} {\rm d}\kappa \wedge {\rm d}t + \lambda_{2} {\rm d}\theta \wedge {\rm d}\eta + {\rm d}\lambda_{2} \wedge {\rm d}\xi - {\rm d}\lambda_{2} \wedge \sigma {\rm d}t - \lambda_{3} {\rm d}\theta \wedge {\rm d}\xi \\
 & -\lambda_{2} {\rm d}\sigma \wedge {\rm d}t + {\rm d}\lambda_{3} \wedge {\rm d}\eta + {\rm d}\lambda_{1} \wedge {\rm d}\theta - \kappa \sigma {\rm d}\lambda_{1}\wedge {\rm d}t - \sigma \lambda_{1} {\rm d}\kappa \wedge {\rm d}t - \kappa \lambda_{1} {\rm d}\sigma \wedge {\rm d}t
\end{array}
\]
where we used the following two equalities
\[
\begin{array}{l}
{\rm d}{\rm d}\xi ={\rm d}(\cos \theta {\rm d}x +\sin \theta {\rm d}y)=-\sin \theta {\rm d}\theta \wedge {\rm d}x +\cos \theta {\rm d}\theta \wedge {\rm d}y={\rm d}\theta \wedge (-\sin \theta {\rm d}x +\cos \theta {\rm d}y)= {\rm d\theta}\wedge {\rm d\eta}, \\
{\rm d}{\rm d}\eta = -{\rm d\theta}\wedge {\rm d\xi}.
\end{array}
\]
The exterior derivative ${\rm d}\psi$ determines the characteristic curves (in our case the geodesics) by means of
\[
\begin{array}{l}
\gamma'(t) \rfloor {\rm d}\psi_{\gamma(t)} = 0, \qquad \textrm{and }
\gamma^{*}{\rm dt} \neq 0.
\end{array}
\]
So the Pfaffian equations for decent parameterizations satisfying $\gamma^{*}{\rm dt} \neq 0$ are given by
\begin{equation} \label{Pfaff}
\left\{
\begin{array}{l}
\partial_{\lambda_{1}} \rfloor {\rm d}\psi ={\rm d}\theta- \kappa \sigma {\rm d}t=0 \\
\partial_{\lambda_{2}} \rfloor {\rm d}\psi ={\rm d}\xi- \sigma {\rm d}t=0  \\
\partial_{\lambda_3} \rfloor {\rm d}\eta =0 \\[7pt]
\partial_{\sigma} \rfloor {\rm d}\psi = (\sqrt{\kappa^2 +\epsilon}-\lambda_{1}\kappa-\lambda_2) {\rm d}t =0 \\
\partial_{\kappa} \rfloor {\rm d}\psi = \sigma (\kappa (\kappa^2 + \epsilon)^{-1/2} -\lambda_1 ) {\rm d}t =0 \\[7pt]
-\partial_{\theta} \rfloor {\rm d}\psi= {\rm d}\lambda_{1} -\lambda_{2}{\rm d}\eta +\lambda_{3}{\rm d}\xi =0 \\
-\partial_{\xi} \rfloor {\rm d}\psi= {\rm d}\lambda_{2} -\lambda_{3}{\rm d}\theta=0 \\
-\partial_{\eta} \rfloor {\rm d}\psi= {\rm d}\lambda_{3} +\lambda_{2}{\rm d}\theta=0 \\
\end{array}
\right. \ .
\end{equation}
The first three equations represent the horizontality restriction, the 4th en 5th equation represent the Euler-Lagrange optimization of the energy and the last three equations provide the Lagrange multipliers (recall (\ref{zform}))
\begin{equation} \label{lambdas}
\left\{
\begin{array}{l}
\lambda_{1}= \frac{\kappa}{\sqrt{\kappa^{2}+\epsilon}}=z \\
\lambda_{2}= -\sqrt{\epsilon}\sqrt{1-z^2} \\
{\rm d}z +\lambda_{3} \sigma{\rm d}t={\rm d}z +\lambda_{3}{\rm d}s \Rightarrow \lambda_{3}=-\dot{z},
\end{array}
\right.
\end{equation}
by employing Noether's theorem and an invariance group of symmetries of $\psi$ (which is defined as a group acting on $Q$ such that the induced action $\eta$ on $T^{*}(Q)$ satisfies $\eta_{g}(Z)=Z$ for all $g \in G$) as will briefly explain next. In our case the action $\eta$ on $T^{*}(Q)$ is induced by the action of $SE(2)$ acting on itself) of the minimization problem.

Noether's theorem (which provides a conservation law on momentum) says that the momentum mapping $m:Z \to SE(2)^{*}$ given by
\[
\langle m(p), \xi \rangle = (\xi \rfloor \psi)(p), \qquad p \in Z,
\]
is constant along the characteristic curves. Note that computation of the Lie-derivative of $\psi$ along a characteristic curve gives
\[
0=\mathcal{L}_{\xi}\mathcal{A}=\mathcal{A} \rfloor {\rm d}\psi +d(\mathcal{A} \rfloor \psi)= \mathcal{A} \rfloor {\rm d}\psi, \qquad \textrm{ for all left-invariant vector fields }\mathcal{A} \in \mathcal{L}(SE(2)),
\]
which explains the last three equalities in (\ref{Pfaff}).

The momentum mapping is invariant under the co-adjoint representation $\textrm{Ad}^{*}$ (this is the representation dual to the adjoint representation (\ref{adjointrep}))
\begin{equation} \label{Adrepinmomentmap}
m(\eta_{g}(p))= (\textrm{Ad}_{g^{-1}})^{*} m(p)
\end{equation}
which follows from the fact that $\eta_{g}^{*}\psi=\psi$ and $(\eta_{g})_{*}\xi=(Ad_{g^{-1}})_{*} \xi$.
Consequently, the characteristic curves are contained in the co-adjoint orbits.
It can be verified that the co-adjoint orbits of $SE(2)$ are given by
\[
\lambda_{2}^{2}+\lambda_{3}^{2}=c^2 \epsilon \geq 0, \qquad c>0,
\]
so we get the following \emph{preservation law} that holds along the characteristic curves
\begin{equation}\label{law}
(\dot{z}(s))^2 +\epsilon -c^2 \epsilon = \epsilon (z(s))^2, \qquad s>0,
\end{equation}
where the normalized curvature $z(s)=\frac{\kappa(s)}{\sqrt{\kappa^{2}(s)+\epsilon}}$ satisfies $|z|<1$ and indeed this formula follows by integration of (\ref{zdiv}), since
\[
\ddot{z}=\epsilon z \Rightarrow \dot{z} \ddot{z} = \epsilon \dot{z} z \rightarrow (\dot{z}(s))^2 = \epsilon (z(s))^2 +C, C \in \R.
\]
As observed by Bryant \and Griffiths \cite{Bryant}{p.543-544} (with slightly different conventions) the last three equations of (\ref{Pfaff}) can be written
\[
{\rm d} \hat{\lambda} = \hat{\lambda} g^{-1}{\rm d}g \desda {\rm d}(\hat{\lambda} \cdot g^{-1})={\rm d} \; \widehat{\left((\textrm{Ad}_{g^{-1}})^* \cdot \lambda\right)} =0
\]
where $\hat{\lambda}=(-\lambda_{3}, \lambda_{2}, \lambda_{1})$ and where the matrix form of the Cartan-connection equals,
\[
g^{-1}{\rm d}g =\left(
\begin{array}{ccc}
\cos \theta & -\sin \theta & x \\
\sin \theta & \cos \theta & y \\
0 & 0 & 1
\end{array}
 \right)^{-1} {\rm d} \left(
\begin{array}{ccc}
\cos \theta & -\sin \theta & x \\
\sin \theta & \cos \theta & y \\
0 & 0 & 1
\end{array}
 \right)= \left(
\begin{array}{ccc}
0 & -{\rm d} \theta & {\rm d}\xi \\
{\rm d}\theta & 0 & {\rm d}\eta \\
0 & 0 & 0
\end{array}
 \right)
\]
where both Lie-algebra and Lie-group are embedded in the group of invertible $3\times 3$ matrices. Consequently, by Noether's theorem we have $\hat{\lambda}=\hat{\mu} \cdot g$, for some constant $\hat{\mu}=(-\mu_{3},\mu_{2},
\mu_{1})$, or more explicitly we have
\[
\left\{
\begin{array}{l}
z=\mu_{1} -\mu_{3} x + \mu_{2} y \\
\dot{z}= -\mu_{3} \cos \theta + \mu_{2} \sin \theta \\
\sqrt{\epsilon(1-z^2)}= \mu_{3} \sin \theta + \mu_{2} \cos \theta, \qquad \textrm{with } \mu_{2}^2+\mu_{3}^{2}=c^2 \epsilon.
\end{array}
\right.
\]
Next we choose
\[
h_{0}=\left(
\begin{array}{ccc}
-\frac{\mu_3}{c \sqrt{\epsilon}} & -\frac{\mu_{2}}{c \sqrt{\epsilon}} & \frac{\mu_1 \mu_3}{c^2 \epsilon} \\
\frac{\mu_2}{c \sqrt{\epsilon}} & -\frac{\mu_3}{c \sqrt{\epsilon}} & -\frac{\mu_1 \mu_2}{c^2 \epsilon} \\
0 & 0 & 1
\end{array}
\right)^{-1}\in SE(2)
\]
so that $\hat{\mu}\cdot h_{0}^{-1}=(\sqrt{\epsilon}c,0,0)$ and use left-invariance $g=h_{0}^{-1}\tilde{g}$, $\tilde{g}\equiv (\tilde{x},\tilde{y},e^{i\tilde{\theta}})$ then we get $\hat{\lambda}=\hat{\mu}\cdot g = (c\sqrt{\epsilon},0,0)\cdot \tilde{g}$, i.e.
\[
\begin{array}{l}
\tilde{x}= \frac{z}{c\sqrt{\epsilon}}, \ \
c \sqrt{\epsilon} \cos \tilde{\theta} =c \sqrt{\epsilon}\dot{\tilde{x}} =\dot{z}, \textrm{ and }-\sqrt{\epsilon(1-z^2)} = -\sqrt{\epsilon}c \sin \tilde{\theta}= c \sqrt{\epsilon}\dot{\tilde{y}}
\end{array}
\]
and consequently we have
\[
\begin{array}{l}
\tilde{x}(s)= (\sqrt{\epsilon}c)^{-1} z(s) \\
\tilde{y}(s)=\tilde{y}(0) + \frac{1}{c} \int_{0}^{s} \sqrt{1-(z^{2}(\tau))}\, {\rm d}\tau \\
\tilde{\theta}(s)= \angle (\dot{\tilde{\ul{x}}}(s),\ul{e}_{x})= \tilde{\theta}(0)+ \int \limits_{0}^{s} \kappa(\tau)\, {\rm d}\tau.
\end{array}
\]
So by means of (\ref{z}) we get the solution $g(s)=(x(s),y(s),\theta(s))=h_{0}^{-1}(\tilde{x}(s),\tilde{y}(s),\tilde{\theta}(s))$, i.e.
\begin{equation} \label{solutiongeodesic}
\left\{
\begin{array}{l}
x(s)= \frac{\mu_{1}\mu_3}{c^2 \epsilon} -\frac{\mu_3}{c\sqrt{\epsilon}} \tilde{x}(s)-\frac{\mu_2}{c\sqrt{\epsilon}} \tilde{y}(s)  \\
y(s)= \frac{-\mu_{1}\mu_2}{c^2 \epsilon} +\frac{\mu_2}{c\sqrt{\epsilon}} \tilde{x}(s)-\frac{\mu_3}{c\sqrt{\epsilon}} \tilde{y}(s)  \\
\theta(s)=\tilde{\theta}(s)+\arccos \left( -\frac{\mu_3}{c\sqrt{\epsilon}}\right)
\end{array}
\right.
\textrm{ with }
\left\{
\begin{array}{l}
\tilde{x}(s)= \frac{z_0}{\sqrt{\epsilon}c} \cosh (\sqrt{\epsilon}s)+ \frac{z_0'}{c \epsilon} \sinh(\sqrt{\epsilon s}) \\
\tilde{y}(s)=\tilde{y}_{0} +\frac{1}{c} \int_{0}^{s} \sqrt{1-c^2 (\tilde{x}(\tau))^2 \epsilon}\; {\rm d}\tau \\
\tilde{\theta}(s)= \arccos \left(\frac{z_0}{c}\sinh(\sqrt{\epsilon \,  s})+ \frac{z_0'}{c\sqrt{\epsilon}}\cosh (\sqrt{\epsilon}s) \right), \qquad
\end{array}
\right.
\end{equation}
where $c=\sqrt{1+ \frac{(z_0')^2}{\epsilon}-z_0^2}$.

Now we have $6$ unknown parameters $\mu_{1},\mu_3, z_0, z_{0}', \tilde{y}(0),L$, (note that $\mu_{2}$ is not unknown since $\mu^{2}_{2}+\mu_{3}^{2}=c^2 \epsilon$ and $c=\sqrt{1+\frac{(z_0')^2}{\epsilon}-(z_0)^{2}}$) to ensure the given boundary conditions
\[
\left\{
\begin{array}{l}
g(0)=(x(0),y(0),e^{i\theta(0)})=g_0:=(x_0,y_0,e^{i\theta_{0}}), \\
g(L)=(x(L),y(L),e^{i\theta(L)})=g_1=(x_1,y_1,e^{i\theta_{1}})
\end{array}
\right.
\]
By means of left-invariance we can always make sure (by multiplying from the left with $g_{1}^{-1}$) that $g_{1}=e$, so $\theta_{1}=0,x_{1}=0,y_1=0$.

In this case straightforward and intense computations yield
\begin{equation}\label{pars}
\begin{array}{ll}
\begin{array}{l}
\mu_{1}=z_0+ \mu_3 x_0 -\mu_{2} y_0, \\
\mu_{2}=c \sqrt{\epsilon} \sin (\arccos \left(\frac{-\mu_3}{c\sqrt{\epsilon}}\right)), \\
\mu_3= -z_0' \cos \theta_0 + \sqrt{\epsilon} \sin \theta_{0} \sqrt{1-z_0^{2}}, \\
c=\sqrt{\frac{\mu_{2}^{2}+\mu_{3}^{2}}{\epsilon}}=\sqrt{1+\frac{(z_0')^2}{\epsilon}-(z_0)^{2}} \\
\end{array}
 &
\begin{array}{l}
\tilde{y}(0)= \frac{-\mu_3 y_0 -\mu_2 x_0}{c \sqrt{\epsilon}}, \\
{\tiny L= \left\{
\begin{array}{l} \frac{-1}{\sqrt{\epsilon}} \log \left( \frac{\mu_3}{\sqrt{\epsilon z_0}}\right) \textrm{if }c=1  \\
\\ \frac{1}{\sqrt{\epsilon}} \log \left( \frac{-\mu_3 + \sqrt{-(z_0')^{2}+ (z_0)^2 \epsilon +\mu_3^2}}{z_0'+z_0 \sqrt{\epsilon}}\right) \textrm{ if }c > 1 , \mu_3 <0, z_0'+z_0 \sqrt{\epsilon}>0 \\
 \frac{1}{\sqrt{\epsilon}} \log \left( \frac{-\mu_3 - \sqrt{-(z_0')^{2}+ (z_0)^2 \epsilon +\mu_3^2}}{z_0'+z_0 \sqrt{\epsilon}}\right) \textrm{ if }c < 1, \mu_3 <0, z_0'+z_0 \sqrt{\epsilon}<0
\end{array} \right. } 
\end{array}
\end{array}
\end{equation}
So all parameters are now expressed in the two unknown $z_0$ and $z_{0}'$ which are determined by the two remaining boundary conditions:
\begin{equation} \label{twoeq}
\left\{
\begin{array}{l}
\frac{\mu_{1} \mu_{3}}{c^2 \epsilon} - \frac{\mu_3}{c\sqrt{\epsilon}} \tilde{x}(L)- \frac{\mu_2}{c \sqrt{\epsilon}} \tilde{y}(L)=x_1, \\
-\frac{\mu_1 \mu_2}{c^2 \epsilon} + \frac{\mu_{2}}{c\sqrt{\epsilon}} \tilde{x}(L) - \frac{\mu_3}{c\sqrt{\epsilon}} \tilde{y}(L)= y_1.
\end{array}
\right.
\end{equation}
Now since $SE(2)$ is a symmetric space \cite{Jost} all points can be connected by a geodesic and we may expect that there indeed exist $z_0$ and $z_0'$ such that (\ref{twoeq}) holds.
Consequently, the singularities (which cause extreme problems in the numerical shooting algorithm (\ref{NR}) of section \ref{ch:geodesics}) where $z(s_{max})=1$ occur always at $s_{max} \geq L$ (and if $\mu_3 \neq c \sqrt{\epsilon}$ then $s_{max}>L$). Next we explicitly verify that $s_{max} \geq L$ in 2 cases.

In case $c>1, \mu_{3}<0$ and $z_0'+\sqrt{\epsilon} z_0>0$ we have
\[
\begin{array}{l}
e^{\sqrt{\epsilon} L}= \frac{-\mu_3 +\sqrt{-(z_0')^2 + (z_0)^{2}\epsilon +\mu_3^2}}{z_0' + z_0 \sqrt{\epsilon}}\ , \\
e^{\sqrt{\epsilon} s_{max}}= \frac{-\sqrt{\epsilon} +\sqrt{(z_0')^2 - (z_0)^{2}\epsilon +\epsilon}}{z_0' + z_0 \sqrt{\epsilon}} = \frac{\sqrt{\epsilon}(1+c)}{z_0'+z_0 \sqrt{\epsilon}}
\end{array}
\]
and indeed $-\mu_3 + \sqrt{-(z_0')^2 + (z_0)^{2}\epsilon +\mu_3^2}< 2 \sqrt{\epsilon}< (1+c) \sqrt{\epsilon}$ so $L<s_{max}$.

In case $c<1, \mu_{3}>0$ and $z_0'+\sqrt{\epsilon} z_0<0$ we have
\[
\begin{array}{l}
e^{\sqrt{\epsilon} L}= \frac{-\mu_3 -\sqrt{-(z_0')^2 + (z_0)^{2}\epsilon +\mu_3^2}}{z_0' + z_0 \sqrt{\epsilon}}\ , \\
e^{\sqrt{\epsilon} s_{max}}= \frac{-\sqrt{\epsilon} +\sqrt{(z_0')^2 - (z_0)^{2}\epsilon +\epsilon}}{z_0' + z_0 \sqrt{\epsilon}} = \frac{\sqrt{\epsilon}(1+c)}{|z_0'+z_0 \sqrt{\epsilon}|}=\frac{\sqrt{\epsilon}(1+c)}{-(z_0'+z_0 \sqrt{\epsilon})}
\end{array}
\]
and indeed we have $e^{\sqrt{\epsilon} s_{max}} \geq e^{\sqrt{\epsilon} L}$, since $\mu_3 + \sqrt{-(z_0')^2 + (z_0)^{2}\epsilon +\mu_3^2}\leq c \sqrt{\epsilon} + \sqrt{\epsilon (1-c^2)+c^2\epsilon}=\sqrt{\epsilon}(1+c)$. Equality is obtained if $\mu_3=c \sqrt{\epsilon}$.

\begin{figure}
\centerline{
\includegraphics[width=0.45\hsize]{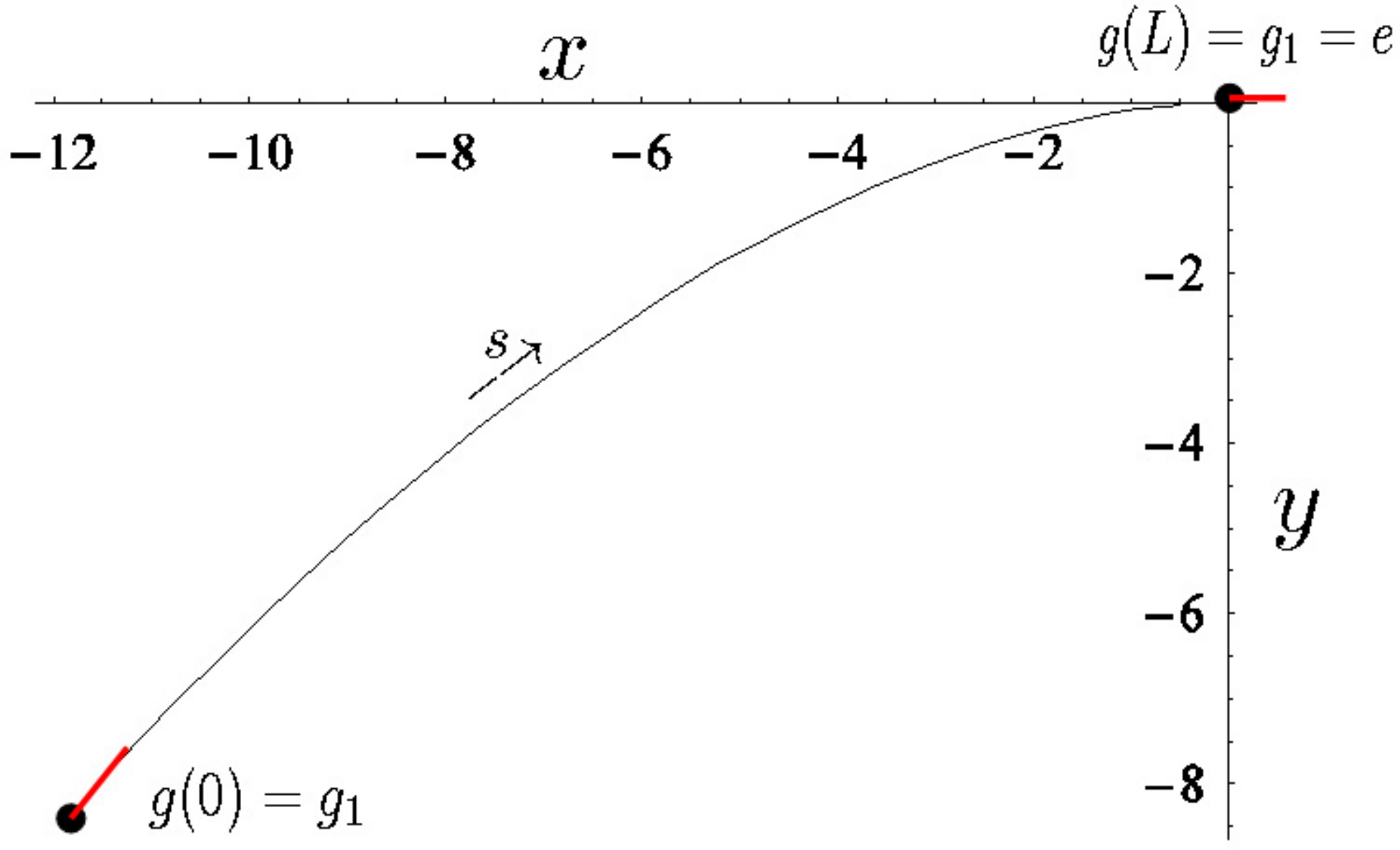} \hfill
\includegraphics[width=0.45\hsize]{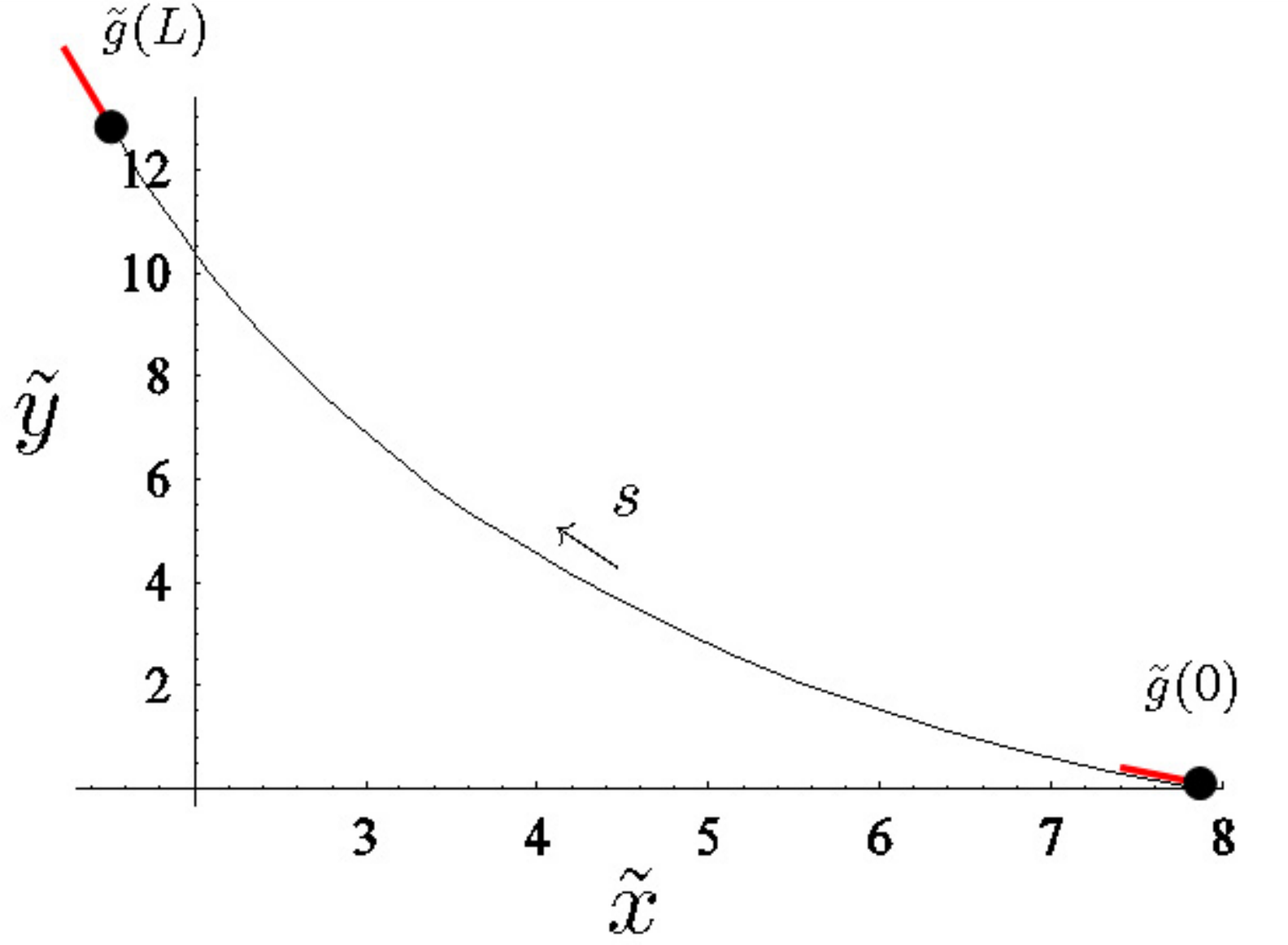}
}
\caption{Left figure: Illustration of a geodesic $s \mapsto g(s)$ computed by (\ref{solutiongeodesic}) and its affine relative $s \mapsto \tilde{g}(s)=h_{0}^{-1} g(s)$. Parameter settings
$x_0=-11.868$, $y_0=-8.44337$, $\theta_0= 51.95^{\circ}$, $x_1=y_1=\theta_1=0$, $L=15$, $\epsilon=0.0125$, $z_0=-0.1641$, $z_0'=0.0183$, $c=1$. }\label{fig:geodesics}
\end{figure}

\begin{figure}
\centerline{
\includegraphics[width=0.35\hsize]{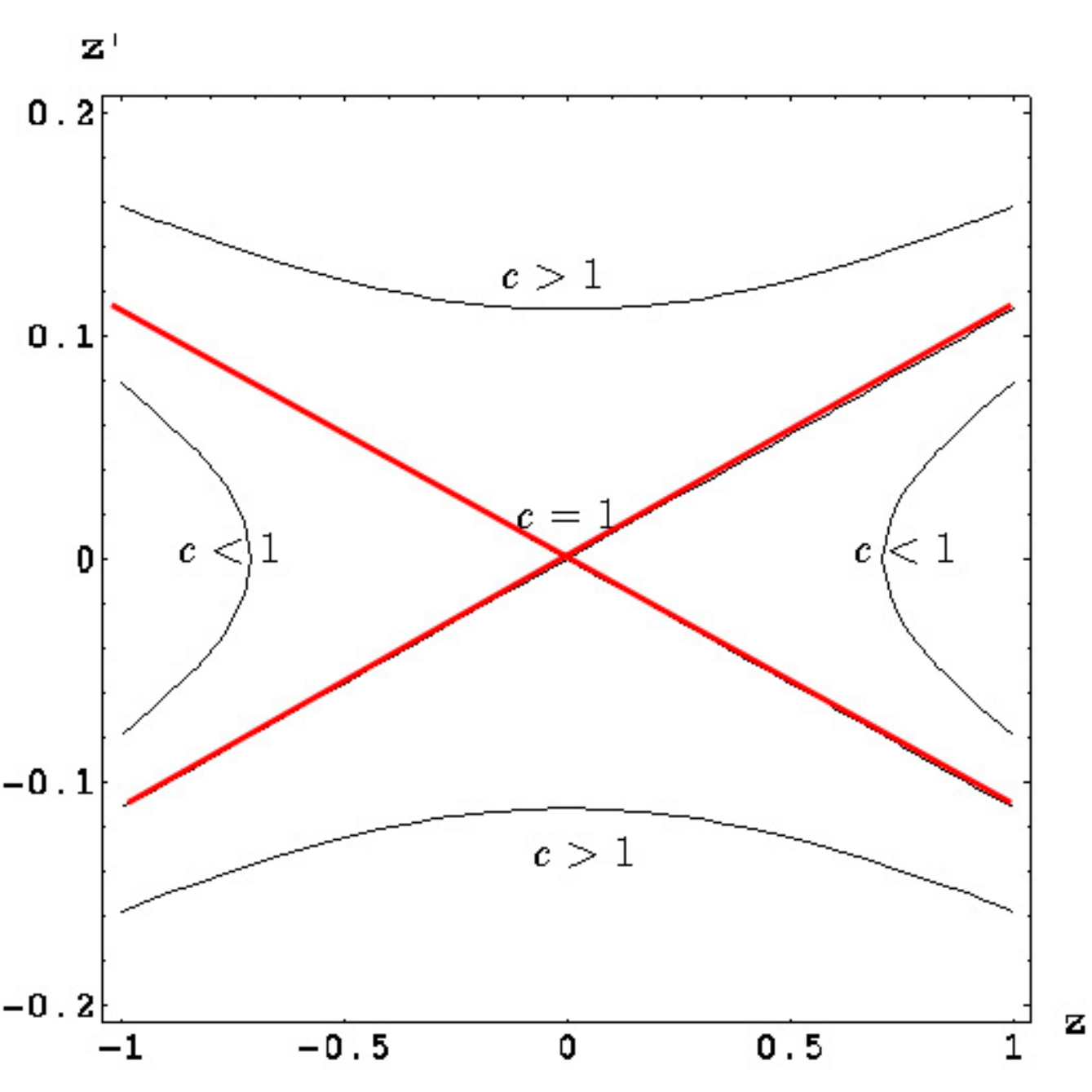} \hfill
\includegraphics[width=0.55\hsize]{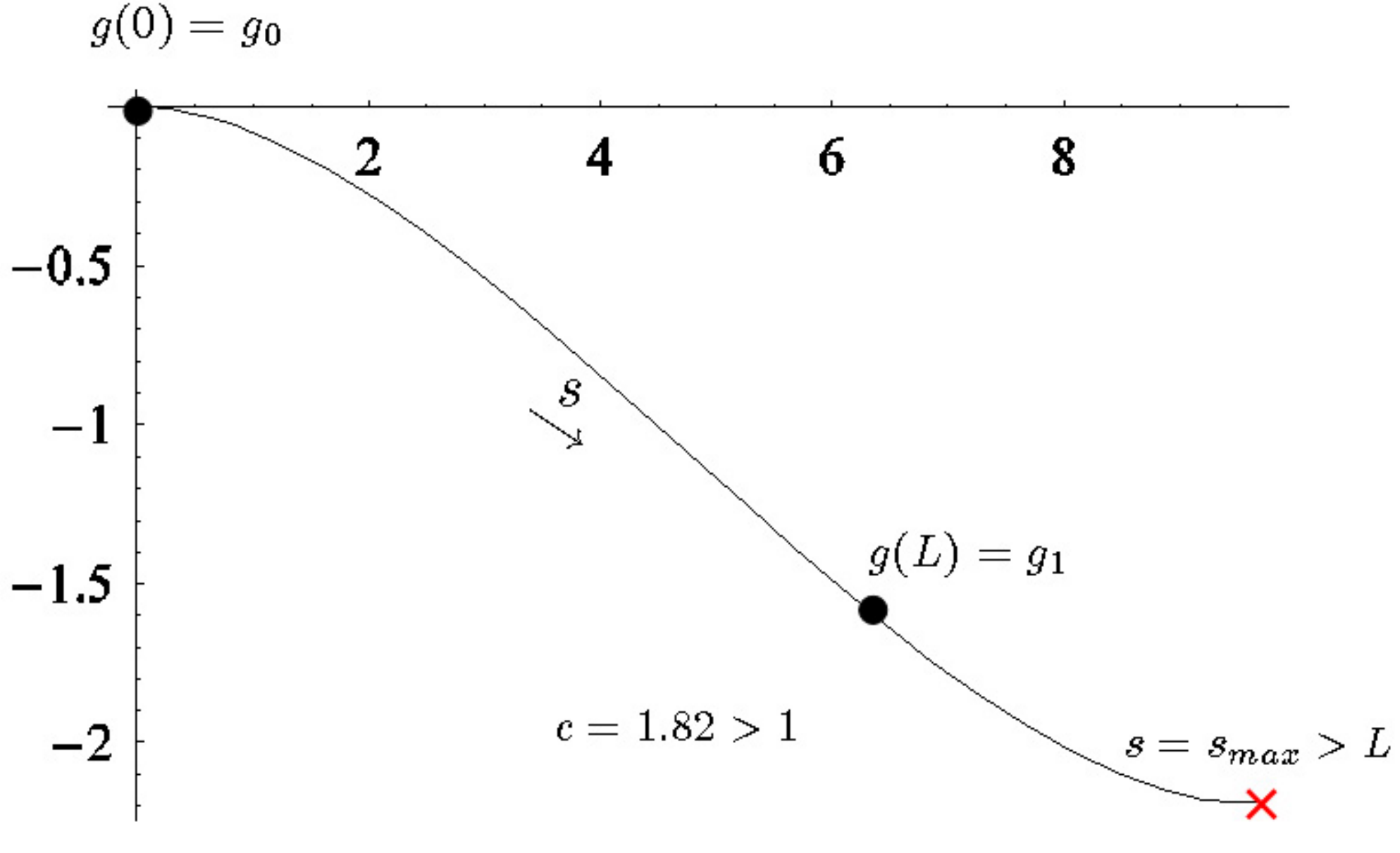}
}
\caption{Left: A phase plot of $z(s)=\frac{\kappa(s)}{\kappa^{2}(s) +\epsilon}$ and $\dot{z}(s)$ clearly indeed reveals that all paths will eventually end up at $z=1$ where solutions brake down because of infinite curvature. Except for the cases where the initial condition is such that $c=\sqrt{1-\epsilon^{-1}(\dot{z}(0))^2-(z(0))^2}=1$ then solutions stay at $c=1$ but reach the line $z=1$ only for $s\to \infty$.
Right: These infinite curvature singularities always take place at $s_{max}>L$, so this causes no problems in our exact analytic solutions (in contrast to the numerical shooting algorithm (\ref{NR}).).}
\label{fig:geodesics2}
\end{figure}

See Figure \ref{fig:geodesics} and see Figure \ref{fig:geodesics2}.
\pagebreak
\section{Completion measures, Brownian bridges and geodesics on $P_Y=(SE(2),SE(2)/Y, \pi ,R )$ \label{ch:app2}}

In this chapter we shall again work with the principal fiber bundle $P_Y=(SE(2),SE(2)/Y, \pi ,R)$, $\pi(g)=[g]=gY$, $Y=\{(0,y,0)\; |\; y \in \R\}$, $R_g h=hg$, equipped with Cartan connection $\omega:T(SE(2) \to T_{e}(SE(2))$ given by $\omega(X_g)=\langle {\rm d}\eta,X_{g}\rangle \partial_{x}$ so that the horizontal part $\mathcal{H}$ of the tangent space $T(SE(2)$ of $SE(2)$ is (by definition) the kernel of the connection $\omega$ which is spanned by
\[
\mathcal{H}=\textrm{span}\{\partial_{\theta},\partial_{\xi}\}.
\]
Recall from section \ref{ch:fiberbundles} that the tangent vectors along horizontal curves, recall Definition \ref{def:horizontal}, are always in the horizontal part of the tangent space.

Sofar we considered completion fields as collision probability of a forward direction process and a backward direction process. However, the direction process is a stochastic process for contour completion. In this section, however, we shall replace the (exact) resolvent Green's functions $R_{\alpha}^{D_{11}}$ for contour completion by the Green's function of contour enhancement $R_{\alpha}^{D_{11},D_{22}}$. The corresponding
completion distribution on $SE(2)$ between $e=(0,0,e^{i0})\in SE(2)$ and $h \in SE(2)$ is now given by
\begin{equation} \label{completiondistr}
g \mapsto R_{\alpha}^{D_{11},D_{22}}(g)R_{\alpha}^{D_{11},D_{22}}(g^{-1}h).
\end{equation}
Next, we show it is related to the Brownian bridge measure by condition time integration, i.e. Laplace transform.

The Brownian bridge measure on $SE(2)$ is given by
\[
\mathbb{Q}^{e,h}_{0,t}(A) = \int \limits_{0}^{t} \int \limits_{A} \mathcal{K}_{t,s}^{h}(g) {\rm d}\mu_{SE(2)}(g) {\rm d}s,
\]
where $A$ is a measurable set within $SE(2)$ supported by the set of paths starting at time $0$ at $e$ and ending up at $h$ at time $t>0$ and where the density $\mathcal{K}_{s,t}^{h}:SE(2) \to \R^{+}$ is given by
\begin{equation} \label{brownianbridge}
\mathcal{K}_{s,t}^{h}(g)= \frac{K_{s}^{D_{11},D_{22}}(g)K_{t-s}^{D_{11},D_{22}}(g^{-1}h)}{K_{t}^{D_{11},D_{22}}(h)},
\end{equation}
where $K_{s}^{D_{11},D_{22}}$ denotes the heat-kernel determined in subsection \ref{ch:exactheatkernels}. The kernel $\mathcal{K}_{s,t}^{h}$ represents the probability density that a random walker which started at $e$ and ends up in $h$ passes $g$. Note that it is a conditional measure, since it only considers paths with start at $e$ and end at time $t$ at $h$, which explains the denominator. We associate to such Brownian bridge measure the following unconditional measure $\mathcal{K}^{h, UC}_{s,t}$ by removing the denominator in (\ref{brownianbridge}):
\[
\mathcal{K}_{s,t}^{h,UC}(g)= K_{s}^{D_{11},D_{22}}(g)K_{t-s}^{D_{11},D_{22}}(g^{-1}h).
\]
Then we have the following result by Wittich, see \cite{wittich} for more details and proof:
\begin{theorem}
Let $M$ be a Riemannian manifold with distance function $d_{M}:M \to \R^{+}$ with positive curvature. Let $p \in M$ and choose $r>0,\epsilon_0>0$ such that $R +\epsilon_0< r(p)$, where $r(p)>0$ is some number such that $B(p,r):=\{x \in M \; :|; d_{M}(p,x)<r\}$ is strongly convex, that is for any two points within such a ball there is a unique minimizing geodesic whose interior is again contained within the ball. Let $C \subset B(p,r)$ be an arbitrary closed subset. Write
\[
\kappa(C):= \max \left\{ 1, \sup \limits_{q \in C, \sigma \in \bigwedge_{q}^{2}(M) K(\sigma)} \right\},
\]
where $\bigwedge_{q}^{2}(M)$ denotes the set of anti-symmetric bilinear forms on the tangent space $T_{q}(M)$, which is isomorphic to the set of oriented 2D-planes in $T_{q}(M)$. Let $\mathbb{Q}^{p,q}_{0,t}$ denote the Brownian Bridge measure on $M$ supported by the set $\Omega(p,q,t)$ of continuous paths starting at time $0$ at $p$ and ending up at $q \in M$ at time $t>0$ and let $\mathbb{Q}^{UC, p,q}_{0,t}$ the associated unconditional measure, i.e. the conditional measure $\mathbb{Q}^{p,q}_{0,t}$ and unconditional measure $\mathbb{Q}^{UC, p,q}_{0,t}$ are given by
\[
\begin{array}{l}
\mathbb{Q}^{p,q}_{0,t}(A)= \lim \limits_{\eta \to 0} \frac{\mathbb{W}_{M}^{p}\left(A \cap \{w(t) \in B(q,\eta)\}\right)}{\mathbb{W}_{M}^{p}\left( \{w(t) \in B(q,\eta)\}\right)}, \\
\mathbb{Q}^{UC, p,q}_{0,t}(A)= \lim \limits_{\eta \to 0} \mathbb{W}_{M}^{p}\left(A \cap \{w(t) \in B(q,\eta)\}\right),
\end{array}
\]
where $\mathbb{W}^{p}_{M}$ denotes the well-known Wiener measure on $M$ centered at $p$.
Then for all $q \in C$ there is a unique geodesic $\gamma^{p,q,t}$, parameterized by with constant velocity $v(s)=\frac{d_{M}(p,q)}{t}$, joining $p$ and $q$. Furthermore, there is some $\delta>0$ such that for all $\epsilon >0 $ with $\epsilon <\epsilon_0$ and for all $q \in C$,
\[
\begin{array}{l}
\mathbb{Q}^{p,q}_{0,t}(B(\epsilon,p,q,t))\leq 2 e^{-\frac{2(R(\kappa(C))-\delta)\epsilon^2}{t}}, \\
\mathbb{Q}^{UC; p,q}_{0,t}(B(\epsilon,p,q,t))\leq 2 e^{-\frac{2(R(\kappa(C))-\delta)\epsilon^2}{t}} e^{-\frac{d_M(p,q)}{2t}}
\end{array}
\]
where \mbox{$B(\epsilon,p,q,t):= \{\omega \in \Omega(p,q,t)\; :\; \sup \limits_{s \in [0,t)} d(\omega(s),\gamma^{p,q,t}(s)) \geq \epsilon \}$} and $R(\kappa(C)):= \frac{\sqrt{\kappa(C)}\, d_{M}(p,q)}{2} \cot \left( \frac{\sqrt{\kappa(C)}\, d_{M}(p,q)}{2}\right)$.
\end{theorem}
Consequently,  the Brownian Bridge measure tends, as $t \downarrow 0$  to the point measure supported by the geodesic $\gamma^{p,q,t}$, parameterized proportional to arc-length with constant velocity $v(s)=\frac{d_{M}(p,q)}{t}$.

Set $M=SE(2)$ or rather $M=P_Y$, since we again restrict ourselves to horizontal curves.
Then we set
\begin{equation} \label{epsilonCitti}
\begin{array}{ll}
d_{SE(2)}(g,h) &= d_{SE(2)}(e, g^{-1}h)=
\inf \left\{\int_{0}^{1} \sqrt{\kappa^{2}(s) + \epsilon} \, {\rm d}s =\int_{0}^{1} \sqrt{ (\theta'(t))^{2} + \epsilon \|\ul{x}'(t)\|^2 }\; {\rm d}t \; |\;  \right. \\
 &\left. \gamma \textrm{ is a smooth horizontal curve connecting }g \textrm{ and }g_{0}, \gamma(0)=a, \gamma(1)=g^{-1}h\right\}, \qquad \epsilon=\frac{D_{11}}{D_{22}},
\end{array}
\end{equation}
where we recall $\frac{ds}{dt}=\|\ul{x}'(t)\|$ and we recall that for horizontal curves we have $\kappa(s)=\frac{d\theta}{ds}$.

Now we note that the explicit relation between the completion measure $\mu^{e,h}(A)$ induced by the completion distribution (\ref{completiondistr}) and the unconditional Brownian-Bridge measure $\mathbb{Q}^{UC; e,h}_{0,t}$ is given by
\[
\begin{array}{ll}
\mu^{e,h}(A) &= \int \limits_{A} R_{\alpha}^{D_{11},D_{22}}(g) R_{\alpha}^{D_{11},D_{22}}(g^{-1}h) {\rm d}\mu_{SE(2)}(g) \\
             &= \alpha^{2} \int \limits_{A} \mathcal{L} (t \mapsto K_{t}^{D_{11},D_{22}}(g))(\alpha) \mathcal{L} (t \mapsto K_{t}^{D_{11},D_{22}}(g^{-1}h))(\alpha) \; {\rm d}\mu_{SE(2)}(g)\\
             &= \alpha^{2} \int \limits_{A} \mathcal{L} (K_{\cdot}^{D_{11},D_{22}}(g)\star K_{\cdot}^{D_{11},D_{22}}(g^{-1}h))(\alpha) \; {\rm d}\mu_{SE(2)}(g) \\
             &= \alpha^{2} \mathcal{L}\left(t \mapsto \int \limits_{A} \left\{ \int \limits_{0}^{t} K_{s}^{D_{11},D_{22}}(g) K_{t-s}^{D_{11},D_{22}}(g) \, {\rm d}s \right\} \; {\rm d}\mu_{SE(2)}(g)\right)(\alpha) \\
             &= \alpha^{2} \mathcal{L} \left(t \mapsto \int \limits_{0}^{t} \left\{\int \limits_{A} K_{s}^{D_{11},D_{22}}(g) K_{t-s}^{D_{11},D_{22}}(g) {\rm d}\mu_{SE(2)}(g)\right\}\, {\rm d}s\right)(\alpha)  \; \\
             &= \alpha^2 \mathcal{L}(t \mapsto \mathbb{Q}^{UC;e,q}_{0,t}(A))(\alpha).
\end{array}
\]
Now we recall that the diffusion generator $\partial_{\xi}^{2} +\partial_{\theta}^{2}$ is hypo-elliptic (it satisfies the H\"{o}rmander condition, recall subsection \ref{ch:hoermander}) as a result the diffusion kernel satisfies similar estimates, that hold for Green's functions of elliptic operators:
\[
|K_{t}(g)| \leq \frac{c}{t^2} e^{-\frac{d_{SE(2)}(g,e)}{4t}},
\]
for some $c>0$, \cite{TerElst},\cite{Hebisch},\cite{Citti}, for more details on this particular case see Appendix \ref{ch:terElst}, and consequently we have for $A=B(\epsilon,e,h,t)$
\[
\begin{array}{ll}
\mu^{e,h}(A=B(\epsilon,e,h,t)) &= \alpha^2 \mathcal{L}(t \mapsto \mathbb{Q}^{UC;e,q}_{0,t}(A))(\alpha) \leq \alpha \int_{0}^{\infty} c e^{-\frac{(d_{SE(2)}(e,h))^2 + 2(R-\delta)\epsilon^2}{4t}} t^{-2}\, e^{-\alpha t} {\rm dt}, \\
 &= \frac{\alpha^{2}\sqrt{\alpha}K_{1}(\sqrt{(d_{SE(2)}(e,h))^2 + 2(R-\delta)\epsilon^2}\sqrt{\alpha})}{\sqrt{(d_{SE(2)}(e,h))^2 + 2(R-\delta)\epsilon^2}},
 \end{array}
\]
with $R=(1/2)\sqrt{\kappa(C)}d_{SE(2)}(e,h) \cot ((1/2) \sqrt{\kappa(C)} d_{SE(2)}(e,h) )$ and $\kappa(C)= \sup \limits_{q \in C} K_{q}(\sigma)$
where $K(\sigma)$ is the horizontal curvature (i.e. the sectional curvature in the plane spanned by $\{\partial_{\theta},\partial_{\xi}\}$.

Consequently, we see that if the expected life time $E(T)=\frac{1}{\alpha}$ tends to $0$ the completion measure between two delta distributions $\delta_{e}$ and $\delta_{h}$ tends to the point measure $\delta_{\gamma^{e,h}}$ supported by the unique geodesic minimizing $\int \sqrt{\frac{1}{4D_{11}}\kappa^{2}+ \frac{1}{4D_{22}}}\; {\rm d}s$ (which we explicitly derived in section \ref{ch:geodesics}) connecting $e$ and $h$.

\section{``Snakes'' in SE(2) based on completion fields of orientation scores \label{ch:app3}}

In this section we will formulate a variational problem, where the energy consists of two parts, with the goal of finding a sufficiently horizontal curve with given beginning and ending that fits the data $C:SE(2) \to \R^{+}$, where $C$ denotes for example the completion distribution of an orientation score (\ref{completionfield}). Here the internal energy of such a curve is the elastica functional $\int \kappa^{2}(s) {\rm ds}$ and the external energy of the curve, which takes care that the curve fits the data, is the total integral $\int e^{- C(\gamma(s))} {\rm d}s$, where $s>0$ denotes the arclength in $\R^2$ of the projected curve $\ul{x}= P_{\R^2}\gamma$. So this is just a direct generalization of the variational methods for so-called snakes in image analysis where a curve in $\R^2$ is supposed to fit the image data $(x,y)\mapsto f(x,y)$ with given restrictions on the internal energy of the curve which usually consists of a length and curvature penalization. However, the ``snakes'' on $SE(2)$ have in principle the advantage that they can deal with crossing contours.

Next we will derive the corresponding Euler-Lagrange equation, but before we can continue we need a small result on differentiating a function $SE(2)$ along a horizontal curve.
\begin{lemma} \label{horprop}
A smooth curve in $SE(2)$ given by $s \mapsto \gamma(s)= \xi(s)\ul{e}_{\xi}(s)+\eta(s) \ul{e}_{\eta}(s) + \theta(s) \ul{e}_{\theta}(s)$, with $\gamma(s)=(x(s),e^{i\theta(s)}) \in SE(2)$, $\ul{e}_{\xi}(s)=\cos \theta(s) \ul{e}_{x}+
\sin \theta(s) \ul{e}_{y}$, $\ul{e}_{\eta}(s)=-\sin \theta(s) \ul{e}_{x}+
\cos \theta(s) \ul{e}_{y}$ and $s>0$ the arclength parameter of the projected curve $x=\mathbb{P}_{\R^2}\gamma$ on the spatial plane, is a horizontal curve in $SE(2)$ iff
$\frac{d \eta}{ds}=-\xi \kappa$. Moreover, for such curves we have $\frac{d\xi}{ds}-\kappa \eta=1$. Now if we differentiate a smooth function on $C:SE(2) \to \R$ along a horizontal curve we get
\begin{equation} \label{alonghorcurve}
\frac{d}{ds} C(\gamma(s))= \partial_{\xi} C(\gamma(s)) + \kappa(s) \partial_{\theta} C(\gamma(s))
\end{equation}
\end{lemma}
\textbf{Proof }
By straightforward differentiation we get
\[
\dot{\gamma}(s)= \frac{d}{ds}( \xi(s)\ul{e}_{\xi}(s)+\eta(s) \ul{e}_{\eta}(s) + \theta(s) \ul{e}_{\theta})= (\dot{\xi}(s) -\kappa(s) \eta(s))\ul{e}_{\xi}(s) +
(\dot{\eta}(s) +\kappa(s) \xi(s))\ul{e}_{\eta}(s) +\dot{\theta}(s)\ul{e}_{\theta},
\]
which is horizontal iff $\frac{d \eta}{ds}=-\xi \kappa$ in which case we have
$1=\| \frac{d\ul{x}}{ds} \|=|\frac{d\xi}{ds}-\kappa \eta|$. The sign follows from the restriction $\angle (\dot{x}(s),\ul{e}_{x})=\theta(s)$. Finally we note that by the chain-law we have
\[
\begin{array}{ll}
\frac{d}{ds} C(\gamma(s)) &= \langle C_{\xi}(\gamma(s)) {\rm d}\xi + C_{\eta}(\gamma^{*}(s)) {\rm d}\eta + C_{\theta}(\gamma^{*}(s)){\rm d}\theta, \dot{\gamma}(s) \rangle \\
 &=  \left(C_{\xi}(\gamma(s)) \left(\frac{d\xi}{ds}-\kappa(s)\eta(s) \right)+ C_{\eta}(\gamma(s)) \left(\frac{d\eta}{ds}+\kappa(s)\xi(s) \right) +C_{\theta}(\gamma(s)) \kappa(s)\right) \\
  & =  \left(C_{\xi}(\gamma(s)) \left(\frac{d\xi}{ds}-\kappa(s)\eta(s) \right)+C_{\theta}(\gamma(s)) \kappa(s)\right) \hfill \Box
\end{array}
\]
\begin{theorem} \label{lemmabelangrijk2}
Let $g_{0}, g_{1} \in SE(2)$. Let $L>0$ be fixed. Consider a smooth positive function $C: SE(2) \to \R^{+}$ and define $\tilde{C}:SE(2) \to \R^{+}$ by $\tilde{C}(g)=e^{-C(g)}$. Then a local maximum (or mode) of the following variational problem
\[
\begin{array}{ll}
\arg \min \{ \mathcal{E}_{L}(\gamma) &=\int \limits_{0}^{L} (e^{-C(\gamma(s))}+\kappa^{2}(s)) \|\dot{\ul{x}}(s)\|\; {\rm d}s \;|\;s \mapsto \gamma(s) \textrm{ is a smooth horizontal curve with total length }L \\
 &\textrm{ connecting } \gamma(0)=g_0 \in SE(2) \textrm{ and }\gamma(L)=g_0 \in SE(2),\ul{x}=\mathbb{P}_{\R^2} \gamma  \}
\end{array}
\]
is a horizontal curve $s \mapsto \gamma^{*}(s)=(\ul{x}^{*}(s),e^{i \theta^{*}(s)})$ which satisfies the following Euler-Lagrange equation:
\begin{equation} \label{extremalcurves}
\kappa^{2} + 3\kappa'(s) +\tilde{C}_{\eta}(\gamma^{*}(s))-\frac{d}{ds}\tilde{C}_{\theta}(\gamma^{*}(s))= \lambda \kappa(s),
\end{equation}
where $\lambda$ is a Lagrange multiplier in $\R$, due to the restriction to curves with length $L$. 
\end{theorem}
\textbf{Proof }
Consider the following (normal) deviation on the horizontal curve
\begin{equation} \label{dev1}
\gamma^{*}(s)=(\ul{x}^{*}(s),\theta^{*}(s)) \mapsto  \gamma_{\epsilon \delta}(s)= (\ul{x}^{*}(s)+\epsilon \delta(s) \ul{e}_{\eta}(s), \theta^{*}(s)+ \arctan \frac{\epsilon\delta'(s)}{1-\kappa(s)\epsilon \delta(s)}), \epsilon >0
\end{equation}
with $\delta$ smooth, compactly supported within $(0,L)$, so $\delta(0)=\delta(L)=\delta'(0)=\delta'(L)=0$ and
\begin{equation}\label{lp}
\begin{array}{l}
\int \limits_{0}^{L} \kappa(s) \delta(s) {\rm d}s=0 \textrm{ and }
\int \limits_{0}^{L} \kappa(s) \tilde{C}(\gamma^{*}(s))\, \delta(s) {\rm d}s=0
\end{array}
\end{equation}
where $s>0$ represents the arc-length of the curve $\ul{x}^{*}=\mathbb{P}_{\R^2} \gamma^{*}$. Note that the total length of the deviated curve equals $L + O(\epsilon^{2})$ since by (\ref{lp}) we have
\[
\left. \partial_{\epsilon} \int \limits_{0}^{L} \|\frac{d}{ds}(x^{*}(s) +\epsilon \delta(s) \ul{e}_{\eta}(s)) \| {\rm d}s \right|_{\epsilon=0}= \int \limits_{0}^{L} \kappa(s) \delta(s) {\rm d}s=0
\]
Moreover, this deviation yields a new horizontal curve $\gamma_{\epsilon\delta}$ with 
length $L$ connecting $g_{0}$ and $g_{1}$. To this end we note that
$\angle ((1-\epsilon\delta(s) \kappa(s))\ul{e}_{\xi}(s) +\delta '(s) \ul{e}_{\eta}(s) ,\ul{e}_x)= \theta(s) + \arctan \left( \frac{\epsilon \delta'(s)}{1- \epsilon \delta(s) \kappa(s) }\right)$, with $\kappa(s)=\theta'(s)$.
Note that by first order Taylor-expansion of the integrand around $g(s)$ we get
{\small
\[
\begin{array}{ll}
\lim \limits_{\epsilon \to 0} \frac{\mathcal{E}_{L}(\gamma_{\epsilon\delta})-\mathcal{E}_{L}(\gamma^{*})}{\epsilon} &
= \lim \limits_{\epsilon \to 0} \frac{1}{\epsilon}\int_{0}^{L} \left\{\left(\tilde{C}(\gamma^{*}(s))+ \epsilon \delta(s) \left. \frac{\partial \tilde{C}}{\partial \eta} \right|_{\gamma^{*}(s)} +\arctan \{ \frac{\epsilon \delta'(s)}{1-\epsilon \delta(s)\theta'(s)}\}\left. \frac{\partial \tilde{C}}{\partial \theta} \right|_{\gamma^{*}(s)} \right)(1-\epsilon \delta(s) \kappa(s) )\right.
\\ & \left. \qquad -\tilde{C}(\gamma^{*}(s))   (\kappa(s)+\epsilon \delta''(s) +\epsilon \delta(s) \kappa^{2}(s))^2(1-\epsilon \delta(s) \kappa(s))-\kappa^{2}(s) +O(\epsilon^2)\right\}
\; {\rm d}s
\\ &=  \int_{0}^{L} \delta(s)\left( \left. \frac{\partial \tilde{C}}{\partial \eta} \right|_{
\gamma^{*}(s)} -\frac{d}{ds}\left. \frac{\partial \tilde{C}}{\partial \theta} \right|_{\gamma^{*}(s)}- \kappa(s)\tilde{C}(\gamma^{*}(s)) + 2\kappa(s)+\kappa^{3}(s) \right){\rm d}s  \\
&=  \int_{0}^{L} \delta(s)\left( \left. \frac{\partial \tilde{C}}{\partial \eta} \right|_{
\gamma^{*}(s)} -\frac{d}{ds}\left. \frac{\partial \tilde{C}}{\partial \theta} \right|_{\gamma^{*}(s)} + 2\kappa''(s) +\kappa^{3}(s)\right){\rm d}s
\end{array}
\]
}
Then since $\gamma^{*}$ is the minimizer and since $\delta$ satisfies the side conditions (\ref{lp}) the gradient of the Energy should be linear dependent on the gradients of the side conditions (which equal $\kappa$ and $\tilde{C}(\gamma^{*}) \kappa$) and the Euler-Lagrange equations read:
\begin{equation}\label{one}
2 \kappa''(s)+\kappa^{3}(s)+\left.
\frac{\partial \tilde{C}}{\partial \eta} \right|_{\gamma^{*}(s)}-\frac{d}{ds}\left\{ \left. \frac{\partial \tilde{C}}{\partial \theta}\right|_{\gamma^{*}(s)} \right\}= q_1(s),
\end{equation}
\begin{equation} \label{qs}
\textrm{with }q_{1}(s)=\lambda_{1}\kappa(s) + \lambda_{2} \tilde{C}(\gamma^{*}(s)) \kappa(s),
\end{equation}
where $\lambda_{1}$ and $\lambda_{2}$ are some Lagrange multipliers. We will show the final step $\lambda_{2}=0$, by systematically checking for all horizontal curve pertubations.

The same technique can be applied to horizontal curve deviations of the type
\begin{equation} \label{dev2}
(\ul{x}^{*}(s),\theta^{*}(s)) \mapsto (\ul{x}^{*}(s) + \epsilon\delta(s) \ul{e}_{\xi}(s), \theta^{*}(s)+\arctan \frac{\epsilon \delta(s) \kappa(s)}{1+\epsilon \delta'(s)})
\end{equation}
where $\delta$ is an arbitrary compactly supported smooth function within $(0,L)$, but this is just a re-parametrization of the same curve and gives (\ref{alonghorcurve}).

Finally, we consider horizontal curve deviations of the type
\begin{equation} \label{dev3}
(\ul{x}(s),\theta(s)) \mapsto (\ul{x}(s)+\epsilon \tilde{\delta}(s) \ul{e}_{\eta}(s) -\frac{\epsilon\tilde{\delta}'(s)}{\kappa(s)}\ul{e}_{\xi}(s), \theta(s))
\end{equation}
where $\tilde{\delta}$ is an arbitrary compactly supported smooth function within $(0,L)$ with  $\int \limits_{0}^{L} \tilde{\delta}(s) \, {\rm d}s=0$
which are again length preserving up to $O(\epsilon^{2})$:
\[
\left. \partial_{\epsilon} \int \limits_{0}^{L}  \| \frac{d}{ds}(\ul{x}(s)+\epsilon \delta(s) \ul{e}_{\eta}(s) -\frac{\epsilon \delta'(s)}{\kappa(s)}\ul{e}_{\xi}(s)) \| {\rm d}s \right|_{\epsilon=0}=
\int \limits_{0}^{L} \delta(s)\kappa(s)+ \left( \frac{\delta'(s)}{\kappa(s)} \right)'(s){\rm d}s=0+\frac{\delta'(L)}{\kappa(L)}-\frac{\delta'(0)}{\kappa(0)}=0,
\]
and which yield
\begin{equation}\label{three}
\frac{d}{ds}\left\{ \kappa^{-1}(s) \left. \frac{\partial \tilde{C}}{\partial \xi} \right|_{\gamma^{*}(s)}\right\}+ \left.\frac{\partial \tilde{C}}{\partial \eta} \right|_{\gamma^{*}(s)}-\frac{d}{ds}\left(
\frac{ \frac{d}{ds}\tilde{C}(\gamma^{*}(s))}{\kappa(s)}\right)=q_3(s)
\end{equation}
with $q_{3}=\lambda_{3} \kappa(s)$
Now divide equation (\ref{alonghorcurve}) by $\kappa$ and differentiate with respect to $s>0$ and by (\ref{one}) we may substitute $(\tilde{C}_{\theta})'=\tilde{C}_{\eta}-q_{1} + 2\kappa'''+\kappa^3$ this yields the following equation
\begin{equation} \label{threetilde}
\frac{d}{ds}\left\{ \kappa^{-1}(s) \left. \frac{\partial C}{\partial \xi} \right|_{\gamma^{*}(s)}\right\}+ \left.\frac{\partial C}{\partial \eta} \right|_{\gamma^{*}(s)}-\frac{d}{ds}\left(
\frac{ \frac{d}{ds}\tilde{C}(\gamma^{*}(s))}{\kappa(s)}\right)=q_{1}(s)- 2\kappa'''(s) +\kappa^{3}(s)
\end{equation}
from which we deduce $q_{1}(s)=q_{3}(s)$ and thereby $\lambda:=\lambda_{1}=\lambda_{3}$ and $\lambda_{2}=0$ from which the result follows.
$\hfill \Box$

\section{Semigroups generated by subcoercive operators on SE(2) \label{ch:terElst} }

In this chapter we shall apply the general theory in \cite{TerElst3} on weighted subcoercive operators on Lie groups, to our case of interest: The Forward Kolmogorov equation of the contour enhancement processes (\ref{Cittiproc}). From this general theory we will deduce
that the closure of the generator of a contour enhancement process indeed generates a holomorphic semi-group with a smooth and fast decaying Green's function on $SE(2)$ (which we explicitly derived in section \ref{ch:diffSE2}). Moreover, we will derive Gaussian estimates for both exact and approximate kernel (derived in subsection \ref{ch:heisapprox}) and we will put our approach of approximation in a more general context yielding a continuous family of holomorphic semi-groups connecting the exact semigroup in subsection \ref{ch:exactheatkernels} and its ``Heisenberg''-approximation in subsection \ref{ch:heisapprox}. We will show that the Gaussian estimates for the kernel $\alpha^{-1}\hat{R}_{\alpha}(x,y,\theta)$ are surprisingly sharp for the approximate resolvent case if $\alpha \downarrow 0$ (i.e. infinite lifetime) which is exactly given by (\ref{Happrox}). This indicates that the Gaussian estimates of the exact kernel (nice for computations in the spatial domain) can be used as reasonable, somewhat rough, approximations of the exact convolution kernels if $D_{11} << D_{22}$. Furthermore, these Gaussian estimates can be used for taking regularized/Gaussian derivatives on orientation scores (similar to (\ref{GDOS1}) and (\ref{GDOS2}) where $D_{22}=D_{33}$) for the (horizontal) case $D_{33}=0$.

Let $G$ be a Lie group, with Lie Algebra $T_{e}(G)$  of dimension $d$ with basis $\{A_{1},\ldots, A_{d}\}$. Let \[\{A_{1},\ldots, A_{d'}\}\subset \{A_{1},\ldots, A_{d}\}, \qquad d' \leq d,
\]
be an algebraic basis of the same Lie algebra, that is there exist an integer $r$ (called the rank of the algebraic basis) such that
\begin{equation} \label{filtration}
\gothic{g}_{1}:= \textrm{span}\{A_{1},\ldots, A_{d'}\},
\gothic{g}_{2}:=\textrm{span}\{[\gothic{g}_{1},\gothic{g}_{1}]\}, \  \ldots \ ,
\gothic{g}_{r}= \textrm{span}\{[\gothic{g}_{r-1},\gothic{g}_{r-1}]\}=T_{e}(G).
\end{equation}
Now for each element in $A \in T_{e}(G)$ there exists a \emph{minimum} integer $k \in \mathbb{N}$ such that $A \in \gothic{g}_{k}$.
We shall refer to this integer as the weight of $A$. In particular the weights of the basis elements $A_{i}$ will be denoted by $w_i$ for $i=1,\ldots, d$. Note that this particular convention of assigning weights, implies $w_{i}=1$ iff $i=1,\ldots,d'$. We stress that this particular convention coincides with a special case of a reduced weighted algebraic basis  $\{A_{1},\ldots A_{d'}\}$ in \cite{TerElst3}, where the filtration $\{\gothic{g}_{\lambda}\}_{\lambda \geq 0}$ should satisfy $\gothic{g}_{\lambda}=\{0\}$ if $\lambda<1$, $\gothic{g}_{\lambda} \subset \gothic{g}_{\mu}$ for all $\lambda \leq \mu$ and $[\gothic{g}_{\lambda},\gothic{g}_{\mu}]\subset \gothic{g}_{\lambda+\mu}$ for all $\lambda,\mu \geq 0$ and $\gothic{g}_{r}=0$ for large $r$, need not be given by (\ref{filtration}) and where the weights $\{w_{i}\}_{i=1,\ldots,d'}$ need not be equal to one, but should satisfy $A_{i} \notin \bigcup \limits_{\lambda<w_i} \gothic{g}_{\lambda}$ for all $i=1,\ldots, d'$.

Let $J(d')$ denote the space of multi-indices associated to the algebraic basis, {\small $J(d')= \bigcup \limits_{n=0}^{\infty} \bigoplus \limits_{k=0}^{n}\{1,\ldots,d'\}^{k}$} and for all
$\alpha=(i_{1},\ldots,i_{n}) \in J(d')$ we associate the the Lie algebra element $A_{\alpha}= A_{i_{1}}\cdots A_{i_{n}}$ and the weighted length
\[
\|\alpha\|=\sum \limits_{k=1}^{n} w_{i_k},
\]
where $n$ is the Euclidean length of $\alpha$ which we shall denote by $|\alpha|=n$.
If $C:J(d')\to \mathbb{C}$ is such that $C(\alpha)=0$ if $\|\alpha\|>m$, for some integer $m \in \mathbb{N}$ and if there exists an multi-index $\alpha \in J(d')$, with $\|\alpha\|= m$, such that $C(\alpha)\neq 0$, then $C$ is called an $m$-th order form. To each $m$-th order form we associate the $m$-th order left invariant operator
\[
\begin{array}{l}
\mathcal{A}_{C}:= \sum \limits_{\alpha \in J(d'), \alpha=(i_{1},\ldots, i_{n})}  C(\alpha) \mathcal{A}^{\alpha} \\
\textrm{ with }\mathcal{A}^{\alpha}=\mathcal{A}_{i_{1}}\ldots \mathcal{A}_{i_{n}}={\rm d}\mathcal{R}(A_{i_{1}})\ldots {\rm d}\mathcal{R}(A_{i_n}).
\end{array}
\]
\begin{definition}
Then $C$ is said to be a $G$-weighted subcoercive form if $\frac{m}{w_i} \in 2 \mathbb{N}$ for all $i=1,\ldots, d'$ and there exist $\mu >0$, $\nu \in \R$ such that the G{\aa}rding inequality holds
\begin{equation}
\textrm{Re} \left\{ (\phi, \mathcal{A}_{C}\phi)_{\mathbb{L}_{2}(SE(2))}\right\} \geq \mu \left(\max \limits_{|\alpha| \leq \frac{m}{2}} \|\mathcal{A}_{\alpha}\|_{\mathbb{L}_{2}(SE(2))} \right)^{2} - \nu \|\phi\|^{2}_{\mathbb{L}_{2}(SE(2))}.
\end{equation}
\end{definition}
Associated to group $G$ and reduced weighted algebraic basis $\{A_{1},\ldots, A_{d'}\}$  one can construct a homogeneous Lie-Algebra $G_{0}$, \cite{Nagel}, with dilations $(\gamma_{t})_{t \geq 0}$
by means of
\begin{equation} \label{homogen}
\begin{array}{l}
[A,B]_{t}= \gamma_{t}^{-1} ([\gamma_{t}(A), \gamma_{t}(B)]) \desda \gamma_{t}([A,B]_t)=[\gamma_{t}(A),\gamma_{t}(B)], \\
\textrm{ where } \gamma_{t}(A_{i})=t^{w_i} A_{i} \textrm{ for }i=1,\ldots,d.
\end{array}
\end{equation}
Now $(T_{e}(G),[\cdot,\cdot]_{t=1})$ is the original Lie Algebra $(T_{e}(G),[\cdot,\cdot])$ and $(T_{e}(G),[\cdot,\cdot]_{0}:=[\cdot,\cdot]_{t \downarrow 0})$ is a homogeneous Lie Algebra
with dilations $(\gamma_{t})_{t>0}$ which is uniquely determined by the filtration corresponding to the reduced algebraic basis. It can be shown that the reduced weighted algebraic basis $\{A_{1},\ldots, A_{d'}\}$ is a reduced weighted algebraic basis for the Lie algebra $(T_{e}(G),[\cdot,\cdot]_{t})$ for all $t>0$. The group simply connected group $G_{t}$ is generated by the Lie algebra $(T_{e}(G),[\cdot,\cdot]_{t})$ via the exponential mapping. The left-invariant vector fields $\mathcal{A}_{i}^{t}$ on $G_{t}$ are given by
\[
(\mathcal{A}_{i}^{t}\phi)(g):= ({\rm d}\mathcal{R}(A_{i})\phi)(g)= \left. \frac{d}{ds} \phi(g \exp_{t}(s A_{i})) \right|_{s=0}.
\]
Now the Lie algebra $(T_{e}(G),[\cdot,\cdot]_{t})$ can be equipped with the following modulus
\begin{equation} \label{modulus}
|g|_{t}'=d_{t}'(g,e)= \inf \left\{\delta>0 \; : \; \exists_{\gamma \in C_{t}(\delta)} \; :\; \gamma(0)=e, \gamma(1)=g  \right\},
\end{equation}
where $C_{t}(\delta)$ equals the space of all absolutely continuous curves $\gamma$ with tangent vectors in the plane spanned by $\{\mathcal{A}_{1}^{t},\ldots, \mathcal{A}_{d'}^{t}\}$ such that
\[
\dot{\gamma}(s)= \sum \limits_{i=1}^{d'} \gamma_{i}(s) \left. \mathcal{A}_{i}^{t} \right|_{\gamma(s)}, \textrm{ with } |\gamma_{i}(s)|=|\langle \left. {\rm d}\mathcal{A}^{i}\right|_{\gamma(s)},\dot{\gamma}(s)\rangle| <\delta^{\omega_{i}}, \textrm{ for all }i=1,\ldots,d' \textrm{ and }s>0.
\]
If we return to our special choice of filtration (\ref{filtration}) all weights of the algebraic basis elements are equal 1 in which case $(T_{e}(G),[\cdot,\cdot]_{0})$ is a nilpotent Lie algebra of the same rank $r$, \cite{TerElst3}{ Lemma 3.10, p.106}.  Using
$(T_{e}(G),[\cdot,\cdot]_{0})$ as a ``local approximation'' of $(T_{e}(G),[\cdot,\cdot])$ the authors in \cite{TerElst3} obtained their general result \cite{TerElst3}{Thm 1.1. p.93}. Next we give a brief summary for the special case of our interest $(H,G,\mathcal{U})=(\mathbb{L}_{2}(G),G,\mathcal{R})$, $G=SE(2)$, with $\mathcal{R}$ the right regular representation whose derivative is the isomorphism between $T_{e}(G)$ and the Lie algebra of left invariant vector fields $\mathcal{L}(G)$, recall (\ref{leftinv3}).
\begin{theorem}\label{th:terElst}
Let $C$ be a $G$-weighted subcoercive form defined on a Lie group $G$, with Haarmeasure $\mu_{G}$. Then the closure of $-\mathcal{A}_{C}$ generates a holomorphic semigroup $s \mapsto S_s$ on $\mathbb{L}_{2}(G)$ which has a fast decreasing kernel in $K_{s} \in \mathbb{L}_{1}(SE(2)) \cap C^{\infty}(G)$ such that
\begin{equation} \label{star}
\mathcal{A}^{\alpha} S_{s}U= \int \limits_{G} (A^{\alpha} K_{s})(h^{-1}g)U(h) {\rm d}\mu_{G}(h), \textrm{ for all }\alpha \in J(d'), \textrm{ and all }U \in \mathbb{L}_{2}(G)
\end{equation}
and for all $\alpha$ there exist $b,c>0$ such that
\begin{equation}\label{ineq}
|\mathcal{A}^{\alpha} K_{s} (g)| \leq c s^{- m^{-1}(\|\alpha\|+\sum \limits_{i=1}^{d}w_{i})} e^{- b \left(\frac{|g|'_{1}}{s}\right)^{\frac{1}{m-1}}},
\end{equation}
for all $g\in G$ and all $s>0$, where $|g|'_{1}$ is given by (\ref{modulus}) with $t=1$.
\end{theorem}
\subsection{Application to Forward Kolmogorov Equation of Contour Enhancement Process }
Now to apply this result to the case of the Forward Kolmogorov equation for the contour enhancement process we set
\[
\begin{array}{l}
G=SE(2), T_{e}(SE(2))=\textrm{span}\{A_{1},A_{2},A_{3}\}=\textrm{span}\{\partial_{\theta},\partial_{x},\partial_{y}\}, \\ \mathcal{L}(SE(2))=\textrm{span}\{\mathcal{A}_{1},\mathcal{A}_{2},\mathcal{A}_{3}\}=\textrm{span}\{\partial_{\theta},\partial_{\xi},\partial_{\eta}\}.
\end{array}
\]
As algebraic basis we use $\{A_{1},A_{2}\}$ since the corresponding left-invariant vector fields $\{\mathcal{A}_{1},\mathcal{A}_{2}\}$ span the horizontal part of the tangentspaces with respect to the Cartan-Maurer form (\ref{Maurer}), which is the only natural choice if it comes to the restriction of curves in $SE(2)$ which arise as ``lifts'' from their projections on the spatial plane. The rank of this algebraic basis is equal to 2 and we have
\[
\gothic{g}_{1}= \textrm{span}\{\partial_{\theta},\partial_{\xi}\} \textrm{ and }\gothic{g}_{2}=[\gothic{g}_{1},\gothic{g}_{1}]=\textrm{span}\{\partial_{\theta},\partial_{\xi},\partial_{\eta}\}=\mathcal{L}(SE(2)),
\]
so $w_{1}=w_{2}=1$ and $\partial_{\eta}\in \gothic{g}_2 \Rightarrow w_{3}=2$.
Now $C(\alpha)=-\alpha_{1}^{2}-\alpha_{2}^{2}$ is a $SE(2)$ weighted subcoercive form, since $m=2$ and $w_{1}=w_{2}=1$ and moreover the G{\aa}rding inequality holds with $\mu=1$, $\nu=0$ since by partial integration
\[
-(\phi,(\partial_{\theta}^{2}+\partial_{\xi}^{2})\phi)_{\mathbb{L}_{2}(SE(2))} \geq \max \{ \|\partial_{\theta}\phi\|^{2}_{\mathbb{L}_{2}(SE(2))},\|\partial_{\xi}\phi\|^{2}_{\mathbb{L}_{2}(SE(2))} \}
\]
for all test functions $\phi \in \mathcal{D}(\Omega_{e})$, where $\Omega_{e}$ is some open environment around the unity in $SE(2)$. Consequently, by Theorem \ref{th:terElst}, the closure of $\partial_{\theta}^{2}+\partial_{\xi}^{2}$ generates a holomorphic semigroup $t \mapsto S_{t}$ on $\mathbb{L}_{2}(SE(2))$ with a fast decreasing kernel $K_{t} \in \mathbb{L}_{2}(SE(2)) \cap
\mathbb{L}_{1}(SE(2))$ such that (\ref{star}) holds. Moreover, the kernel satisfies the following estimate:
\[
|K_{t}(g)| \leq c t^{-2} e^{-b \frac{(|g|_{1}')^2}{t}},
\]
with locally $(|g|')^2\equiv \xi^{2}+\theta^{2}+|\eta|$, since $(\theta,\xi,\eta)$ are the coordinates of the second kind in $SE(2)$. The kernel satisfies $(\partial_{t} +\partial_{\theta}^{2}+\partial_{\xi})K_{t}=\delta_{0}^{t} \otimes \delta^{g}_{e}$ as distributions.

This can be generalized to $D_{11}>0$, $D_{22}>0$ and $C(\alpha)=-D_{11}\alpha_{1}^{2}-D_{22}\alpha_{2}^{2}$ yeilding
\[
K_{t}(x,y,\theta) \leq \frac{1}{4\pi t^2 D_{11} D_{22}} e^{-\left\{ \frac{(x\cos \theta+ y \sin \theta)^2}{D_{22}}+  \frac{\theta^{2}}{D_{22}}+ \frac{|-x \, \sin \theta+ y \cos \theta|}{\sqrt{D_{11}D_{22}}}\right\}\frac{1}{4t}},
\]
for all $x,y \in \R$ and all $\theta \in [0,2\pi)$, where we note that the constant $b$ in the equality (\ref{ineq}) is independent on $D_{11}, D_{22}$ so $b=\frac{1}{4}$.
This estimate also coincides with the estimate by Citti and Sarti \cite{Citti}{ Thm 5.1}.

Next we estimate will investigate the sharpness of the Gaussian estimates. Now since the moduli $|\cdot|_{t}$ are locally equivalent \cite{TerElst3} and since we have a simple exact formula for the  resolvent Heisenberg approximation kernel (\ref{Happrox}) in the spatial domain we choose to study the sharpness of the estimate of the kernels on the group $G_{0}$ which will turn out to be isomorphic to the Heisenberg group $H_3$. Here we shall devide the analysis in two steps. First we shall show that the approach in subsection \ref{ch:heisapprox} is a special case of the above homogenization (\ref{homogen}) of the Lie algebra, yielding a \emph{continuum of semigroups} between the exact case studied in subsection \ref{ch:exactheatkernels} and the Heisenberg approximation case studied in subsection \ref{ch:heisapprox}. Then we shall consider the logarithmic weighted modulus which is locally equivalent to the modulus (\ref{modulus}) and derive surprisingly sharp Gaussian estimates both from above and below of the resolvent Heisenberg kernel (\ref{Happrox}).

Following the general scheme we define the dilation on the algebra by $\gamma_{t}:T_{e}(SE(2)) \to T_{e}(SE(2))$ by $\gamma_{t}(c^{1}A_{1}+c^{2}A_{2}+ c^{3}A_{3})=t \, c^{1}A_{1}+t \, c^{2}A_{2}+ t^{2} c^{3}A_{3}$. Furthermore we define the corresponding dilation on the group by $\tilde{\gamma}_{t}(x,y,e^{i\theta})=(\frac{x}{t},\frac{y}{t^2},e^{i\frac{\theta}{t}})$. Now $(T_{e}(SE(2)),[\cdot,\cdot]_{t})$ with $[\cdot,\cdot]_{t}$ given by (\ref{homogen}) is a Lie-algebra with corresponding simply connected group $(SE(2))_t =\exp_{t}(T_{e}(SE(2)))$. Note that the dilation on the Lie-algebra coincides with the pushforward of the dilation on the group $\gamma_{t}=(\tilde{\gamma})_{*}$ and the left invariant vector fields on $(SE(2))_t$ are given by
\[
\left. \mathcal{A}_{i}^{t} \right|_{g}= (\tilde{\gamma}_{t}^{-1} \circ L_{g} \circ \tilde{\gamma}_{t})_{*}A_{i},
\]
for all $t \in (0,1]$ and a brief computation yields
\[
\begin{array}{ll}
\left. \mathcal{A}_{i}^{t} \right|_{g} \phi &=(\tilde{\gamma}_{t}^{-1} \circ L_{g} )_{*} (\tilde{\gamma}_{t})_{*} A_{i} \phi =(\tilde{\gamma}_{t}^{-1} \circ L_{g} )_{*} \gamma_{t} (A_{i}) \phi  =t^{w_{i}}(\tilde{\gamma}_{t}^{-1} \circ L_{g} )_{*} (A_{i}) \phi  =t^{w_{i}} (\tilde{\gamma}_{t}^{-1})_{*} \left. \mathcal{A}_{i} \right|_{g} \phi \\
 &=t^{w_{i}} \left. \mathcal{A}_{i} \right|_{\tilde{\gamma}^{-1}_{t}g}(\phi \circ \tilde{\gamma}_{t})
\end{array}
\]
for all smooth complex-valued functions $\phi$ defined on a small open environment around $g \in SE(2)$.
So we see that for all $g=(x,y,e^{i\theta}) \in SE(2)$ we have
\[
\begin{array}{ll}
\left. \mathcal{A}_{1}^{t} \right|_{g} &=\frac{1}{t} ( t \partial_{\theta})=\partial_{\theta} \\
\left. \mathcal{A}_{2}^{t} \right|_{g} &= t \left( \frac{\cos(\theta t)}{t} \partial_{x}+ \frac{\sin \theta t}{t^2}\partial_{y} \right) = \; \cos(\theta t) \partial_{x} + \frac{\sin (\theta t)}{t} \partial_{y} \\
\left. \mathcal{A}_{3}^{t} \right|_{g} &= t^{2}\left( -\frac{\sin(\theta t)}{t} \partial_{x} + \frac{\cos (\theta t)}{t^{2}}\partial_{y}\right) = \; - t \sin(\theta t) \partial_{x} + \cos(\theta t)\partial_{y}
\end{array}
\]
and indeed $\left.\mathcal{A}_{i}^{t}\right|_{e}=\left. \mathcal{A}_{i}\right|_{e} =A_{i}$ and
\[
\begin{array}{lll}
\
[\mathcal{A}_{2}^{t},\mathcal{A}_{3}^{t}]=0,
& \ [\mathcal{A}_{2}^{t},\mathcal{A}_{3}^{t}]=t^2 \mathcal{A}_{2}^{t}, &
\ [\mathcal{A}_{1}^{t},\mathcal{A}_{3}^{t}]=t^2 \mathcal{A}_{3}^{t},
\end{array}
\]
and by taking the limit $t \downarrow 0$ we see that the homogeneous contraction $(SE(2))_{0}=\lim \limits_{t \downarrow 0} (SE(2))_t$ is isomorphic to $H_{3}$ and the space of corresponding left-invariant vector fields equals $\mathcal{L}(H_{3})=\textrm{span}\{\partial_{\theta}, \partial_{x}+\theta \partial_{y},\partial_{y}\}$ which is indeed a nilpotent Lie group of rank $2$.
This nilpotent Lie group isomorphic to $H_{3}$ is a subgroup of the five dimensional group $H_{5}$ of Heisenberg type that arises by approximating $\cos \theta \approx 1$ and $\sin \theta \approx \theta$ whose left-invariant vector fields are given by
\[
\begin{array}{lll}
\hat{A}_{1}=\mathcal{A}_{1}^{0}=\partial_{\theta} , & \hat{A}_{2}=\mathcal{A}_{2}^{0}=\partial_{x}+\theta \partial_{y}, & \hat{A}_{4}= \partial_{y}, \\ \hat{A}_{3}=-\theta \partial_{x}+ \partial_{y}, & \hat{A}_{5}= \partial_{x}. & \
\end{array}
\]
This Lie-algebra $\mathcal{L}(H_{5})=\textrm{span}\{\hat{A}_{1},\hat{A}_{2},\hat{A}_{3},\hat{A}_{4},\hat{A}_{5}\}$ is isomorphic to the matrix-algebra
\[
\sum \limits_{i=1}^{5} a^{i} \hat{A}_{i} \leftrightarrow
\left(
\begin{array}{cccc}
0 & a^{1} & a^{4} & a^{5} \\
0 & 0 & a^{2} & a^{3} \\
0 & 0 & 0 & 0 \\
0 & 0 & 0 & 0
\end{array}
\right)=: \sum \limits_{i=1}^{5} a^{i} E_{i} =:B
\]
whose exponent is given by
{\small
\[
\exp(t B)=1+tB+ \frac{t^{2}}{2} B^2 = \left(
\begin{array}{cccc}
1 & t \, a^{1} & t \, a^{4} + \frac{1}{2} t^{2} a^{1} a^{2}& t \, a^{5} +\frac{1}{2} t^2 a^{1} a^{3} \\
0 & 1 & t \, a^{2} & t \, a^{3} \\
0 & 0 & 1 & 0 \\
0 & 0 & 0 & 1
\end{array}
\right).
\]
}
This isomorphism enables us to quickly relate the coordinates of first kind to the coordinates of the second kind in $H_{3}=(SE(2))_{0}$ without explicitly using the CBH-formula :
\[
\begin{array}{l}
(x,y,\theta)= \exp_{0}(\alpha^{3} A_{3}) \exp_{0}(\alpha^{2} A_{2})\exp(\alpha^{1} A_{1})= \exp_{0}(\beta^{1}A_{1}+ \beta^{2}A_{2} + \beta^{3}A_{3}) \desda \\[8pt]
\beta^{1}=\alpha^{1}=\theta,  \
\beta^{3}+\frac{1}{2}\beta^{1}\beta^{2}= \alpha^{3}=y, \
\beta^{2}=\alpha^{2}= x \\
\end{array}
\]
so we see that the coordinates of the first kind on $(SE(2))_0$ read
\[
\beta^{1}=\theta, \beta^{2}=x \textrm{ and }\beta^{3}=y- \frac{1}{2}x \theta
\]
and as a result the weighted modulus on $(SE(2))_{0}$ associated to the filtration (\ref{filtration}) is given by
\[
|g|_{0}=\sqrt{\theta^{2}+x^2+|y-\frac{1}{2}x\theta|},
\]
now by \cite{TerElst3} Prop.6.1 there exists a $c\geq 1$ and an $\epsilon \in (0,1]$ such that for all $a \in T_{e}(SE(2))$ with $\|a\|\leq \epsilon$ such that $c^{-1}|a|_{t}\leq |\exp_{t}(a)|_{t}'\leq c |a|_{t}$, where the weighted modulus is given by $|a|_t= |\sum \limits \beta^{i}_{t} A_{i}^{t}|= \sqrt{(\beta^{1}_{t})^{2/w_1} +(\beta^{2}_{t})^{2/w_{2}} + |\beta^{3}_{t}|^{2/w_{2}}}$.
As a result we have for $t=0$, $\beta^{k}_{t=0}=\beta^{k}$, $k=1,2,3$ that
\[
\begin{array}{l}
|(x,y,\theta)|' \geq \frac{1}{c} (x^{2}+\theta^{2}+|y -\frac{1}{2}x\theta|) \Rightarrow \\ \hat{K}_{s}^{D_{11},D_{22}}(x,y,e^{i\theta}) \leq \frac{1}{4\pi s^2} e^{-\frac{(|(x,y,\theta)|')^2}{4s}} \leq \frac{1}{4\pi s^2} e^{-\frac{1}{c^2 4s}(x^2 +\theta^{2}+ |y -\frac{1}{2}x\theta|) },
\end{array}
\]
so that by integration over traveling time $s>0$ we find
\[
\lim \limits_{\alpha \downarrow 0} \hat{R}^{D_{11},D_{22}}_{\alpha}(x,y,\theta)= \int \limits_{0}^{\infty} \hat{K}_{s}^{D_{11},D_{22}}(x,y,\theta) \, {\rm d}s \leq
\frac{1}{\pi D_{11}D_{22}} \frac{c^{2}}{\frac{x^{2}}{D_{22}}+ \frac{\theta^{2}}{D_{11}}+ \frac{|y-\frac{1}{2}x\theta|}{\sqrt{D_{11}D_{22}}}}.
\]
Now if we consider the exact solution
{\small
\[
\begin{array}{ll}
\frac{1}{4\pi D_{11}D_{22}}\frac{1}{\frac{x^{2}}{D_{22}}+ \frac{\theta^{2}}{D_{11}}+ \frac{|y-\frac{1}{2}x\theta|}{\sqrt{D_{11}D_{22}}}}
&\leq
\int \limits_{0}^{\infty} \hat{K}_{s}^{D_{11},D_{22}}(x,y,\theta) \, {\rm d}s = \frac{1}{\pi D_{11}D_{22}} \frac{1}{\sqrt{\left(\frac{x^{2}}{D_{22}}+\frac{\theta^{2}}{D_{11}}\right)^2+ 16 \frac{|y-\frac{1}{2}x\theta|^{2}}{D_{11}D_{22}}}} \\
 &\leq \frac{1}{\pi D_{11}D_{22}} \frac{\sqrt{2}}{\frac{x^{2}}{D_{22}}+ \frac{\theta^{2}}{D_{11}}+ \frac{|y-\frac{1}{2}x\theta|}{\sqrt{D_{11}D_{22}}}}
\end{array}
\]
}
where we note that for all $a,b>0$ one has $a+b \geq \sqrt{a^{2}+b^{2}} \geq \frac{1}{\sqrt{2}}(a+b)$, then we see that $c=\sqrt[4]{2}\approx 1.19 >1$ indeed yields a Gaussian upper bound for the exact Heisenberg kernel for $\alpha \downarrow 0$, whereas $c=0.5$ yields a Gaussian lower-bound for the same kernel.
\end{document}